\documentclass[titlepage,11pt]{report} % use larger type; default would be 10pt
\usepackage[utf8]{inputenc} % set input encoding (not needed with XeLaTeX)
\usepackage{amsmath,amsfonts,amssymb,amsthm}
\usepackage{geometry} % to change the page dimensions
\usepackage{graphicx} % support the \includegraphics command and options
\usepackage{booktabs} % for much better looking tables
\usepackage{array} % for better arrays (eg matrices) in maths
\usepackage{paralist} % very flexible & customisable lists (eg. enumerate/itemize, etc.)
\usepackage{verbatim} % adds environment for commenting out blocks of text & for better verbatim
\usepackage{subfig} % make it possible to include more than one captioned figure/table in a single float
\usepackage{fancyhdr} % This should be set AFTER setting up the page geometry
\usepackage{hyperref}
\usepackage{sectsty}
\usepackage{fullpage}
\usepackage{minitoc}
\allsectionsfont{\sffamily\mdseries\upshape} % (See the fntguide.pdf for font help)
\usepackage[nottoc,notlof,notlot]{tocbibind} % Put the bibliography in the ToC
\usepackage[titles,subfigure]{tocloft} % Alter the style of the Table of Contents

\newtheorem{proposition}{Proposition}
\newtheorem{theorem}{Theorem}

\newtheorem{lemma}{Lemma}

\usepackage{chngcntr}
\counterwithin{proposition}{chapter}
\counterwithin{theorem}{chapter}
\counterwithin{definition}{chapter}
\counterwithin{lemma}{chapter}

\begin{document}

%- une synthèse en une quinzaine de pages des recherches post‐thèse et des perspectives
%qui seront développées dans le mémoire HDR.

\title{Statistical models and \\probabilistic methods \\on Riemannian manifolds}
\author{Salem Said -- CNRS, Université de Bordeaux}
\date{}

\maketitle

\dominitoc
\tableofcontents

%\chapter{Introduction\,: why read this thesis\,?}

\chapter*{A guide to this thesis} 
This thesis reflects the major themes of my work, which I have carried out in the past four years.\hfill\linebreak It does leave out some of this work, especially  on the subject of warped information metrics. I~hope readers may find time for at least a glance at this ``missing" part, (for example, in~\cite{warped}). However, the thesis is rather self-contained, and I feel that the best way of reading it is just from beginning to end, uninterrupted.

At any rate, I would like to ask readers to begin with Chapter \ref{introduction}. Then, once this is done, they can skip to Chapters \ref{gaussian} and \ref{bayesian}, which I would like to ask them to read together, or go on to\hfill\linebreak Chapter \ref{barycentre}, which is quite independent from following. The same goes for Chapter \ref{stocha}, which can be read right after Chapter \ref{introduction}, provided just a little bit of familiarity with Chapter \ref{gaussian}. 

Each chapter begins with a table of contents, followed by a sort of ``abstract", which provides some additional details, on the table of contents, and points to some of the more interesting results. I have done my best to avoid the thesis being a copy-paste of published research papers. Chapter \ref{gaussian} uncovers several new connections between Riemannian Gaussian distributions and random matrix theory, while Chapter \ref{bayesian} is entirely made up of previously unpublished material. The other chapters stick more closely to my existing papers (published, or under review), although I have made a consistent effort to improve the presentation, and to include useful background and historical discussion.

I hope that readers will find it stimulating, to read an ``original thesis". On the other hand, exploring new ideas exposes one to the risk of making mistakes (of various magnitude), and I also hope these are duly pointed out, and the appropriate criticism is served up, without restraint. On the whole, writing this thesis has been a humbling experience for me. I have found out, once and again, that I was unable to answer questions or to prove statements,\hfill\linebreak even when they seemed very natural. Chapter \ref{chapopen}, a short final chapter, contains a list of such ``open problems" (they are open to me, but others may find them easy).   
%humbling experience
%open problems

I should acknowledge the input of many colleagues, who have shaped the ideas layed out in the following. Chapter \ref{barycentre} was born out of discussions with Marc Arnaudon, and Chapter \ref{stocha} relies heavily on joint work with Alain Durmus, Pablo Jimenez, and Eric Moulines~\cite{colt}\cite{aistats}. 

For Chapter \ref{gaussian}, the idea of a useful connection between Riemannian Gaussian distributions and random matrix theory was first suggested to me by my colleague Yannick Berthoumieu. During the summer of 2020, I worked on this idea with Cyrus Mostajeran and Simon Heuveline. Later on, when I was nearly finished writing this thesis, I was very excited to discover the work of Leonardo Santilli and Miguel Tierz~\cite{tierz}, who were simultaneously developing the same idea. It is really a great satisfaction to see a whole project unfold out of an ``innocent" discussion. \hfill\linebreak
For this, I want to thank all of the colleagues I just mentioned. 

Perhaps nobody will ever write a better preface than Cervantes, whose following famous words certainly apply here are now. \textit{Idle reader\,: Without my swearing to it, you can believe that I would like this book, the child of my understanding, to be the most beautiful, the most brilliant, and the most discreet that anyone could imagine. But I have not been able to contravene the natural order; in it, like begets like}.

\chapter{Notation and background} \label{introduction}

\minitoc
\vspace{0.1cm}

{\small Certainly, this thesis is intended for specialised readers, who are already familiar with the basics of Riemannian geometry. This first chapter is not a stand-alone introduction to Riemannian geometry, but merely hopes to help the readers ease into the material in subsequent chapters\,: by recalling some elementary notions in Riemannian geometry, I hope to find a shared language with my benevolent readers. Some original, or even unpublished, material is still included. As discussed in the following\,:
\begin{itemize}
  \item \ref{sec:levi} -- \ref{sec:taylor} lead up to the second-order Taylor expansion of a function defined on a Riemannian manifold. Proposition \ref{prop:geodesic} of \ref{sec:taylor} states that geodesic curves are exactly those curves which admit a Taylor expansion of any $C^2$ function, in terms of its gradient and Hessian.
\item \ref{sec:pca} recalls real Grassmann manifolds, and the PCA objective function. It presents a calculation of the gradient and Hessian of this function, based on the symmetric space structure of real Grassmann manifolds.
\item \ref{sec:retractions} and \ref{sec:grassret} contain some original material, from~\cite{colt}. In \ref{sec:retractions}, the original concept of regular retraction (a retraction which avoids the cut locus) is introduced. In \ref{sec:grassret}, the usual ``projection retraction" on the real Grassman manifold is shown to be a regular retraction.
\item \ref{sec:squareddistance} recalls the usual metric and Hessian comparison theorems of Riemannian geometry.
\item \ref{sec:robust} is an application of \ref{sec:squareddistance}, studied in~\cite{colt}. It introduces the robust Riemanian barycentre of a probability measure on a Hadamard manifold, proving its existence and uniqueness.
\item \ref{sec:volume} is concerned with Riemannian volume\,: integration in geodesic spherical coordinates, volume comparison, and integral formulae for Riemannian symmetric spaces. All of these will be very important, in Chapters \ref{gaussian} and \ref{bayesian}.
\item \ref{sec:geolemma} provides the previously unpublished Propositions \ref{prop:geometriclemma} and \ref{prop:geodesiclemma}, which may be used to compute geodesics in symmetric spaces.
\end{itemize}
}
\vfill
\pagebreak

\section{The Levi-Civita connection} \label{sec:levi}
A smooth (0,2)-tensor field $g$, on a finite-dimensional smooth manifold $M$, is called a Riemannian metric, if the bilinear form
\begin{equation} \label{eq:metric}
\langle u,\!v\rangle_{\scriptscriptstyle x} =   g_{\scriptscriptstyle x}(u,v)  \hspace{1cm} u\,,v \in T_xM
\end{equation}
is a true scalar product, for each $x \in M$. In this case, for $u \in T^{\phantom{*}}_xM$ and $t \in T^*_xM$, the identities
$$
  (u^{\flat},v) = \langle u,\!v\rangle_{\scriptscriptstyle x} \hspace{0.3cm}\mbox{and}\hspace{0.3cm}
  \langle t^{\,\scriptscriptstyle \#},\!v\rangle_{\scriptscriptstyle x}=  \left(t,v\right) \hspace{1cm} v \in T^{\phantom{*}}_xM
$$
uniquely define $u^{\scriptscriptstyle \flat} \in T^*_xM$ and $t^{\,\#} \in T^{\phantom{*}}_xM$. By a useful abuse of notation,
\begin{equation} \label{eq:pairing}
u^{\flat} = g(u) \hspace{1cm}  t^{\,\scriptscriptstyle \#} = g^{\scriptscriptstyle -1}(t)
\end{equation}
The Levi-Civita connection of $g$ is the unique affine connection $\nabla$ which is metric, so that
\begin{equation} \label{eq:metricconnection}
  \nabla g = 0
\end{equation}
and tortionless, so that the exterior derivative $d\theta$, of any $1$-form $\theta$, reads
\begin{equation} \label{eq:tortionless}
d\theta(X,Y) =  \nabla_X\theta(Y) - \nabla_Y\theta(X) 
\end{equation}
for vector fields $X$ and $Y$. In effect, by (\ref{eq:metricconnection}) and (\ref{eq:tortionless}), the $(1,1)$-tensor field $\nabla X$ (the covariant derivative of the vector field $X$), decomposes into self-adjoint and skew parts,
\begin{equation} \label{eq:koszul}
  2\,\langle\nabla_YX,Z\rangle = \mathcal{L}_Xg(Y,Z) + dX^{\flat}(Y,Z)
\end{equation}
where $\mathcal{L}_Xg$ denotes the Lie derivative of the metric $g$ along $X$, the ``the linear elasticity tensor" (the equivalence between (\ref{eq:metricconnection})--(\ref{eq:tortionless}) and (\ref{eq:koszul}) is the content of Koszul's theorem).

Given local coordinates $(x^i\,;i=1,\ldots,n)$ on an open $U \subset M$, there is a coordinate frame $(\partial_i)$, along with a coframe $(dx^i)$ --- of course,  $\partial_i$ stands for $\left.\partial\middle/\partial x^i\right.$. In terms of these coordinates, the metric $g$ takes on the form of a length element
\begin{equation} \label{eq:lengthelement}
 g = g_{ij}\,dx^i\otimes dx^j \hspace{1cm} g_{ij} = \langle \partial_i\,,\partial_j\rangle
\end{equation}
and covariant derivatives may be expressed, in coordinate form,
\begin{equation} \label{eq:christoffel}
  \nabla X = \left\lbrace\partial_j X^i + \Gamma^i_{jk} X^k\right\rbrace \partial_i \otimes dx^j \hspace{1cm} \nabla_{\partial j}\,\partial_k = \Gamma^i_{jk}\,\partial_i
\end{equation}
using the Christoffel symbols $(\Gamma^i_{jk})$. 
\section{Parallel transport and geodesics} \label{sec:geod}
A vector field $X$, along a smooth curve $c:I\rightarrow M$, defined on some interval $I \subset \mathbb{R}$, is a map $X:I \rightarrow TM$ such that 
$\pi \circ X = c$ ---  of course, $\pi:TM \rightarrow M$ denotes the canonical projection. The Levi-Civita connection $\nabla$ can be used to compute the covariant derivative of $X$ along $c$, itself a vector field along $c$, here denoted $\nabla_{\dot{c}\,}X$. In local coordinates,
\begin{equation} \label{eq:dotX}
\nabla_{\dot{c}\,}X(t) = \left \lbrace \frac{d}{dt}X^i(t) + (\Gamma^i_{jk}\circ c(t))\, \dot{c}^{\,j}(t)X^k(t)\right\rbrace (\partial_i\circ c)(t)
\end{equation}
and this suggests writing $\nabla_{\dot{c}\,}X = \nabla_{t\,}X$ or even $\dot{X}$, when $c$ is understood from the context. Now, $X$ is called parallel along $c$ if $\nabla_{\dot{c}\,}X = 0$. From (\ref{eq:dotX}), this means that the components $X^i(t)$ satisfy a linear differential equation with smooth coefficients.

Thus, if $X$ is parallel along $c$, then $X$ is completely determined by its value at any instant, say $t_{\scriptscriptstyle o} \in I$. Equivalently, if $v$ is tangent to $M$ at $c(t_{\scriptscriptstyle o})$, then there exists a unique parallel vector field $X$ along $c$, with $X(t_{\scriptscriptstyle o}) = v$. It follows that, for $t \in I$, there exists a linear operator $\Pi^t_{t_{\scriptscriptstyle o}}$ which maps $T_{c(t_{\scriptscriptstyle o})}M$ onto $T_{c(t)}M$, by $\Pi^t_{t_{\scriptscriptstyle o}}(v) = X(t)$. This linear operator $\Pi^t_{t_{\scriptscriptstyle o}}$ is called parallel transport along $c$, from  $c(t_{\scriptscriptstyle o})$ to $c(t)$, and has the following properties,
\begin{equation} \label{eq:hemigroup}
  \text{hemigroup property} \hspace{1cm} \Pi^t_{t_{\scriptscriptstyle o}} = \Pi^t_{t_{\scriptscriptstyle 1}} \circ \Pi^{t_{\scriptscriptstyle 1}}_{t_{\scriptscriptstyle o}}
\end{equation}
\begin{equation} \label{eq:isometry}
\phantom{abcd}\text{isometry property} \hspace{1cm} \Vert\Pi^t_{t_{\scriptscriptstyle o}}(v)\Vert_{\scriptscriptstyle c(t)} = \Vert v\Vert_{\scriptscriptstyle c(t_{\scriptscriptstyle \scriptscriptstyle o})}
\end{equation}
where $\Vert \cdot \Vert_x$ is the norm associated with the scalar product in (\ref{eq:metric}), for any $x \in M$. Clearly, if one knows how to compute parallel transports, then one is able to recover covariant derivatives,
\begin{equation} \label{eq:pinabla}
\nabla_{\dot{c}\,}X(t_{\scriptscriptstyle o}) = \left.\frac{d}{dt}\right|_{t=t_{\scriptscriptstyle o}}\Pi^{t_{\scriptscriptstyle o}}_{t}(X(t))
\end{equation}
A smooth curve $c:I\rightarrow M$ is called a geodesic curve, if its velocity vector field $\dot{c}$ is parallel. This means that $c$ satisfies the geodesic equation,
\begin{equation} \label{eq:geodesicequation}
\nabla_{\dot{c}\,}\dot{c} = 0 \text{ or } \ddot{c} = 0
\end{equation}
Written out in local coordinates, this is a non-linear ordinary differential equation,
\begin{equation} \label{eq:accelerationcoordinates}
\frac{d^{\scriptscriptstyle\hspace{0.03cm}2}}{dt^{\scriptscriptstyle 2}}c^i(t) + \Gamma^i_{jk}(c(t))\frac{d}{dt}c^{\,j}(t)\frac{d}{dt}c^k(t) = 0
\end{equation}
If its solutions $c(t)$ exists at all finite $t \in \mathbb{R}$, for any initial conditions $c(t_{\scriptscriptstyle o}) = x$ and $\dot{c}(t_{\scriptscriptstyle o})=v$, then the metric $g$ on the manifold $M$ is called geodesically complete. In this case, the Riemannian exponential map $\mathrm{Exp}:TM\rightarrow M$,  given by $\mathrm{Exp}_x(v) = c(1)$ is well-defined. 

The geodesic equation (\ref{eq:geodesicequation}) states that the curve $c$ has zero acceleration (just like a particle in free motion). This means that geodesic curves are extremals of the energy functional
\begin{equation} \label{eq:energy}
  E(c) = \int_I\,\Vert \dot{c}(t)\Vert^2\,dt
\end{equation}
and that re-parameterised geodesic curves (of the form $c\circ \varphi$ where $c$ is a geodesic, and $t^\prime = \varphi(t)$ a new parameterisation) are extremals of the length functional
\begin{equation} \label{eq:length}
  L(c) = \int_I\,\Vert \dot{c}(t)\Vert\,dt
\end{equation}
This leads to the notion of Riemannian distance, which will be discussed in \ref{sec:squareddistance} below.

\section{Taylor expansion of a function} \label{sec:taylor}
Let $f:M \rightarrow \mathbb{R}$ be a $C^2$ function, denote $df$ its differential. The gradient of $f$ is the vector field
\begin{equation} \label{eq:gradient}
\mathrm{grad}\,f = g^{\scriptscriptstyle -1}(df)
\end{equation}
In the notation of (\ref{eq:pairing}). The Hessian of $f$ is the $(1,1)$-tensor field
\begin{equation} \label{eq:hessian}
 \mathrm{Hess}\,f = \nabla\,\mathrm{grad}\,f 
\end{equation}
The following proposition says that geodesic curves are exactly the curves which admit a Taylor expansion of any $C^2$ function, in terms of its gradient and Hessian.
\begin{proposition} \label{prop:geodesic}
 A smooth curve $c:I\rightarrow M$ is a geodesic curve, if and only if, for any $s,t \in I$, and any $C^2$ function $f:M \rightarrow \mathbb{R}$, 
\begin{equation} \label{eq:taylor1}
  (f\circ c)(t) = (f\circ c)(s) + \langle \mathrm{grad}\,f,\dot{c}\rangle_{\scriptscriptstyle c(s)}(t-s) + \frac{1}{2}\, \langle\mathrm{Hess}\,f\cdot\dot{c}\,,\dot{c}\rangle_{\scriptscriptstyle c(s)}(t-s)^{\scriptscriptstyle 2} + o(|t-s|^{\scriptscriptstyle 3})
\end{equation}
\end{proposition}
The proof of this proposition follows from the identity,
$$
\frac{d^{\scriptscriptstyle\hspace{0.03cm}2}}{dt^{\scriptscriptstyle 2}}(f\circ c)(t) = \frac{d}{dt}\langle \mathrm{grad}\,f,\dot{c}\rangle_{\scriptscriptstyle c(t)} = 
\langle\nabla_{\dot{c}\,} \mathrm{grad}\,f,\dot{c}\rangle_{\scriptscriptstyle c(t)} + \langle \mathrm{grad}\,f,\ddot{c}\rangle_{\scriptscriptstyle c(t)}
$$
Indeed, the last term is identically zero (\textit{i.e.}, for any $C^2$ function $f$), if and only if $\ddot{c}$ is identically zero, as in the geodesic equation (\ref{eq:geodesicequation}). 

In (\ref{eq:hessian}), $\mathrm{Hess}\,f$ is a $(1,1)$-tensor field. This tensor field is self-adjoint, and it is common practice to identify it with a $(0,2)$-tensor field, 
\begin{equation} \label{eq:hessianbis}
 \mathrm{Hess}\,f = \nabla\,df = \frac{1}{2}\mathcal{L}_{\mathrm{grad}\,f}\,g
\end{equation}
where the second equality follows from (\ref{eq:koszul}) and (\ref{eq:gradient}). This yields a lighter notation,
$$
\langle\mathrm{Hess}\,f\cdot u\,,v\rangle = \mathrm{Hess}\,f(u,v)
$$
Recall the Riemannian exponential map $\mathrm{Exp}$, from \ref{sec:geod} (always assume geodesic completeness). Proposition \ref{prop:geodesic} can be used to write down a Taylor expansion with Lagrange remainder,
\begin{equation} \label{eq:taylor2}
  f\left(\mathrm{Exp}_x(v)\right) = f(x) + \langle \mathrm{grad}\,f,v\rangle_{\scriptscriptstyle x} + \frac{1}{2}\,\mathrm{Hess}\,f_{\scriptscriptstyle c(t^*)}(\dot{c},\dot{c})
\end{equation}
where $c(t^*)$ is a point along the geodesic $c(t) = \mathrm{Exp}_x(t\,v)$, corresponding to an instant $t^* \in (0,1)$. \\[0.15cm]
\textbf{Remark\,:} writing (\ref{eq:hessianbis}) in local coordinates,
\begin{equation} \label{eq:hesscoordinates}
\mathrm{Hess}\,f = \left\lbrace \partial_{ ij} f - \Gamma^k_{ij\,} \partial_k f\right\rbrace dx^i\otimes dx^{\,j}
\end{equation}
The second derivatives $\partial_{ ij} f$ do not transform like a covariant tensor, but the Christoffel symbols correct for this problem, yeilding a true covariant tensor, $\mathrm{Hess}\,f$. A very nice way of saying this is that the Levi-Civita connection transforms second-order differentials, into covariant tensors. 

The concepts of second-order vectors and of second-order differentials are reviewed in~\cite{sanslarmes}, where they are used as a starting point for stochastic analysis in manifolds. 

\section{Example\,: The PCA objective function} \label{sec:pca}
The problem of principal component analysis consists in maximising the objective function
\begin{equation} \label{eq:pcaf}
  f(x) = \mathrm{tr}\left(x\Delta\right) \hspace{1cm} x \in \mathrm{Gr}_{\scriptscriptstyle \mathbb{R}}(p\,,q)
\end{equation}
where $\Delta$ is a symmetric positive-definite matrix, of size $(p+q)\times(p+q)$. The maximisation is over $x$ in the real Grassmann manifold $\mathrm{Gr}_{\scriptscriptstyle \mathbb{R}}(p\,,q)$, identified with a space of orthogonal projectors
\begin{equation} \label{eq:grassconst}
\mathrm{Gr}_{\scriptscriptstyle \mathbb{R}}(p\,,q)= \left \lbrace x \in \mathbb{R}^{\scriptscriptstyle (p+q)\times(p+q)}\,:x^{\dagger} - x = 0 \hspace{0.1cm},\hspace{0.1cm} x^2 - x = 0 \hspace{0.1cm},\hspace{0.1cm} \mathrm{tr}(x) = p\right\rbrace
\end{equation}
where $^\dagger$ denotes the transpose. Remarkably, it is possible to show that $\mathrm{Gr}_{\scriptscriptstyle \mathbb{R}}(p\,,q)$ is a submanifold of $\mathrm{S}(p+q)$, the affine space of symmetric matrices of size $(p+q)\times(p+q)$, with tangent spaces (the proof of this statement may be found in~\cite{eschenburg}), 
\begin{equation} \label{eq:tangentgrassconst}
T_x \mathrm{Gr}_{\scriptscriptstyle \mathbb{R}}(p\,,q)= \left\lbrace v \in \mathrm{S}(p+q)\,: xv + vx = v \right\rbrace
\end{equation}
It then follows that $\mathrm{Gr}_{\scriptscriptstyle \mathbb{R}}(p\,,q)$ is of dimension $pq$. Clearly, $\mathrm{Gr}_{\scriptscriptstyle \mathbb{R}}(p\,,q)$ admits of a Riemannian metric, which is the restriction of the trace scalar product of $\mathrm{S}(p+q)$,
\begin{equation} \label{eq:grasssp}
  \langle u,\!v\rangle_{\scriptscriptstyle x}  = \mathrm{tr}(uv) \hspace{1cm} u\,,v \in T_x\mathrm{Gr}_{\scriptscriptstyle \mathbb{R}}(p\,,q)
\end{equation}
By (\ref{eq:pairing}) and (\ref{eq:gradient}), it follows from (\ref{eq:pcaf}) that the gradient of $f(x)$ is given by
\begin{equation} \label{eq:pcapx}
\mathrm{grad}\,f(x) = \mathrm{P}_x(\Delta)
\end{equation}
where $\mathrm{P}_x:\mathrm{S}(p+q)\rightarrow \mathrm{S}(p+q)$ is the orthogonal projection onto $T_x \mathrm{Gr}_{\scriptscriptstyle \mathbb{R}}(p\,,q)$. Now, let $x = o$, the projector onto the span of the first $p$ vectors in the canonical basis of $\mathbb{R}^{\scriptscriptstyle p+q}$. One readily checks from (\ref{eq:tangentgrassconst}) that
\begin{equation} \label{eq:grassto}
T_o \mathrm{Gr}_{\scriptscriptstyle \mathbb{R}}(p\,,q)= \left\lbrace \tilde{\omega} = 
\left(\begin{array}{cc} 0_{\scriptscriptstyle p\times p} & \omega^\dagger \\[0.15cm] \omega & 0_{\scriptscriptstyle q\times q}\end{array}\right)\,;\, \omega \in \mathbb{R}^{\scriptscriptstyle q\times p}
 \right\rbrace
\end{equation}
Therefore, $\mathrm{P}_o(\Delta)$ is just $\Delta$ with its main diagonal blocks of size $p \times p$ and $q\times q$ set to zero.\hfill\linebreak Then, note that the orthogonal group $O(p+q)$ acts transitively on $\mathrm{Gr}_{\scriptscriptstyle \mathbb{R}}(p\,,q)$, by $g\cdot x = gxg^\dagger$\hfill\linebreak for $g \in O(p+q)$ and $x \in \mathrm{Gr}_{\scriptscriptstyle \mathbb{R}}(p\,,q)$, and that this action preserves the  Riemannian metric (\ref{eq:grasssp}). Therefore, one has the following alternative to (\ref{eq:tangentgrassconst})\footnote{By an abuse of notation, $g\cdot a= g\,a\,g^\dagger$, for any matrix $a$ of size $(p+q)\times(p+q)$.},% (By an abuse of notation, $g\cdot \tilde{\omega}= g\,\tilde{\omega}\,g^\dagger$),
\begin{equation} \label{eq:tangentgrassgroup}
T_x \mathrm{Gr}_{\scriptscriptstyle \mathbb{R}}(p\,,q)= \left \lbrace v = g\cdot \tilde{\omega}\,;\, g\cdot o = x\hspace{0.1cm},\hspace{0.1cm} \tilde{\omega} \in T_o \mathrm{Gr}_{\scriptscriptstyle \mathbb{R}}(p\,,q) \right\rbrace
\end{equation} 
where $g\cdot o = x$ simply means the first $p$ columns of $g$ span the image space of $x$. Since the action of $O(p+q)$ preserves the  Riemannian metric (\ref{eq:grasssp}), it easily follows
\begin{equation} \label{eq:grasspx}
\mathrm{P}_x(\Delta) = g \cdot \mathrm{P}_o(g^\dagger\cdot \Delta) \text{ for any $g$ such that $g\cdot o = x$}
\end{equation}
which can be used to evaluate the gradient of $f(x)$, in (\ref{eq:pcapx}). For the Hessian of $f(x)$, note that, according to Propositon \ref{prop:geodesic}, 
\begin{equation} \label{eq:pcahess}
\mathrm{Hess}\,f_x(v,v) = \left.\frac{d^{\scriptscriptstyle\hspace{0.03cm}2}}{dt^{\scriptscriptstyle 2}}\right|_{t=0} f\left(\mathrm{Exp}_x(tv)\right)
\end{equation}
Here, the Riemannian exponential can be transformed into a matrix exponential (see Proposition \ref{prop:geodesiclemma}, in \ref{sec:geolemma}). For $g \in O(p+q)$, note that $g \cdot o = o$ if and only if $g \in O(p) \times O(q) \subset O(p+q)$. Denote $\mathfrak{g}$ and $\mathfrak{k}$ the Lie algebras of $O(p+q)$ and $O(p) \times O(q) \subset O(p+q)$. Let $\mathfrak{p}$ denote the orthogonal complement of $\mathfrak{k}$ (with respect to the bilinear form $Q(\xi\hspace{0.02cm},\eta) = \mathrm{tr}(\xi\eta)$, for $\xi\,,\eta \in \mathfrak{o}(p+q)$). Then, 
\begin{equation} \label{eq:grassp}
\mathfrak{p} = \left\lbrace \hat{\omega} = 
\left(\begin{array}{cc} 0_{\scriptscriptstyle p\times p} & -\omega^\dagger \\[0.15cm] \omega & 0_{\scriptscriptstyle q\times q}\end{array}\right)\,;\, \omega \in \mathbb{R}^{\scriptscriptstyle q\times p}
 \right\rbrace
\end{equation}
From (\ref{eq:grassto}) and (\ref{eq:grassp}), it is clear there exists a canonical isomorphism $\pi_o: T_o \mathrm{Gr}_{\scriptscriptstyle \mathbb{R}}(p\,,q) \rightarrow \mathfrak{p}$ (just add a minus sign in front of $\omega^\dagger$ in (\ref{eq:grassto})). In terms of this isomorphism,
\begin{equation} \label{eq:grasslift}
 \mathrm{Exp}_x(tv) = \exp(t\, \hat{\omega}_{\scriptscriptstyle v})\cdot x \hspace{1cm} \hat{\omega}_{\scriptscriptstyle v} = g\cdot \pi_o(g^\dagger \cdot v)
\end{equation}
Replacing (\ref{eq:grasslift}) into (\ref{eq:pcahess}), the second derivative is easily computed,
\begin{equation} \label{eq:pcahessbis}
\mathrm{Hess}\,f_x(v,v) = \mathrm{tr}\!\left( \Delta\,\hat{\omega}^2_{\scriptscriptstyle v}\,x \right) +   
\mathrm{tr}\!\left( \Delta\,x\,\hat{\omega}^2_{\scriptscriptstyle v}  \right) -2\,\mathrm{tr}\!\left(\Delta\,\hat{\omega}_{\scriptscriptstyle v}\,x\, \hat{\omega}_{\scriptscriptstyle v} \right)
\end{equation}
\textbf{Remark\,:} a nice property of the linear map $v \mapsto \hat{\omega}_{\scriptscriptstyle v}$ is that $\mathrm{tr}(\hat{\omega}^2_{\scriptscriptstyle v}) = \langle v,v\rangle_{\scriptscriptstyle x\,}$.

\section{Regular retractions} \label{sec:retractions}
A retraction is a map $\mathrm{Ret}:TM \rightarrow M$, taking $v \in T_xM$ to $\mathrm{Ret}_x(v)$, and which verifies~\cite{absil}\cite{absilmalick},
\begin{equation} \label{eq:retdefinition}
  \mathrm{Ret}_x(0_x) = x \hspace{1cm} d\,\mathrm{Ret}_x(0_x) = \mathrm{Id}_x
\end{equation}
where $0_x \in T_xM$ is the zero vector in $T_xM$, and $\mathrm{Id}_x$ is the identity map of $T_xM$. While the Riemannian exponential $\mathrm{Exp}:TM \rightarrow M$ is itself a retraction\footnote{Recall that it is always assumed $M$ is geodesically complete.}, other retractions are often used as computationally cheap (or numerically stable) substitutes for the Riemannian exponential.

From (\ref{eq:retdefinition}), for any retraction $\mathrm{Ret}$, $\mathrm{Ret}_x$ agrees with $\mathrm{Exp}_x$ up to first-order derivatives. Further, $\mathrm{Ret}$ will be called geodesic, if $\mathrm{Ret}_x$ agrees with $\mathrm{Exp}_x$ up to second-order derivatives. This means the curve $c(t) = \mathrm{Ret}_x(tv)$ has zero initial acceleration\,: $\ddot{c}(0) = 0_{x\,}$, for any $v \in T_xM$ (in the notation of (\ref{eq:geodesicequation}). %In this case, $c(t)$ behaves like a geodesic, 

To compare a retraction $\mathrm{Ret}$ with the exponential $\mathrm{Exp}$, it is useful to introduce the maps
\begin{equation} \label{eq:PHI}
  \Phi_x : T_xM \rightarrow T_xM \hspace{1cm} \Phi_x(v) = \left(\mathrm{Exp}^{-1}_x \circ \mathrm{Ret}_x\right)(v)
\end{equation}
These maps are well-defined if $\mathrm{Ret}$ is regular. That is, if $\mathrm{Ret}_x(v) \notin \mathrm{Cut}(x)$ for any $v \in T_xM$ ($\mathrm{Cut}(x)$ denotes the cut locus of $x$, whose definition is recalled in \ref{sec:squareddistance}). In addition, they satisfy the following propositions.
\begin{proposition} \label{prop:retractions1}
  Let $\mathrm{Ret}:TM\rightarrow M$ be a regular retraction. Then, $\Phi_x : T_xM \rightarrow T_xM$ verify \\[0.1cm]
\emph{(a)} $\Phi_x(0_x) = 0_x$ and $\Phi^\prime_x(0_x) = \mathrm{Id}_x$ (the prime denotes the Fr\'echet derivative). \\[0.1cm]
\emph{(b)} $\Phi^{\prime\prime}_x(0_x)(v,v) = \ddot{c}(0)$, where the curve $c(t)$ is given by $c(t) = \mathrm{Ret}_x(tv)$.
\end{proposition}
\begin{proposition} \label{prop:retractions2}
  Let $\mathrm{Ret}:TM\rightarrow M$ be a regular retraction and $f:M \rightarrow \mathbb{R}$ be a $C^2$ function.  
\begin{equation} \label{eq:taylor2ret}
  f\left(\mathrm{Ret}_x(v)\right) = f(x) + \langle \mathrm{grad}\,f,\Phi_x(v)\rangle_{\scriptscriptstyle x} + \frac{1}{2}\,\mathrm{Hess}\,f_{\scriptscriptstyle \gamma(t^*)}(\dot{\gamma},\dot{\gamma})
\end{equation}
where $\gamma(t^*)$ is a point along the geodesic $\gamma(t) = \mathrm{Exp}_x(t\Phi_x(v))$, corresponding to some $t^* \in (0,1)$. 
\end{proposition}
As an application of Proposition \ref{prop:retractions1}, consider the following examples. \\[0.1cm]
\textbf{Example 1\,:} let $M = S^n \subset \mathbb{R}^{n+1}$, the unit sphere of dimension $n$, with its usual (round) Riemannian metric. The retraction $\mathrm{Ret}_x(v) = \left(x+v)\middle/\Vert x + v \Vert\right.$ ($\Vert \cdot \Vert$ is the Euclidean norm) is regular, and the maps $\Phi_x$ are given by
\begin{equation} \label{eq:phisphere}
 \Phi_x(v) = \arctan(\Vert v \Vert)\,\frac{v}{\Vert v \Vert} 
\end{equation} 
\textbf{Example 2\,:} let $M = U(d)$, the Lie group of $d \times d$ unitary matrices, with its bi-invariant metric $\langle u,\!v\rangle_{\scriptscriptstyle x} = -(1/2)\mathrm{tr}(uv)$. The retraction $\mathrm{Ret}_x(v) = \mathrm{Pol}(x+v)$ ($\mathrm{Pol}$ denotes the left polar factor) is regular, and the maps $\Phi_x$ are given by
\begin{equation} \label{eq:phiun}
 \Phi_x(v) = x\left(u \exp(i\arctan(\theta))\, u^\dagger\right)
\end{equation} 
where $^\dagger$ denotes the conjugate-tranpose, and $\omega = x^\dagger v$ has spectral decomposition $\omega = u(i\theta)u^\dagger$, where $u$ is unitary and $\theta$ is real and diagonal --- as one may expect, $\arctan(\theta) = \mathrm{diag}(\arctan(\theta_{ii}))$. \\[0.1cm]
\indent Now, (b) of Proposition \ref{prop:retractions1} implies the retractions in question are geodesic, since the Taylor expansion at zero of the arctangent only contains odd powers. Both of these retractions are based on orthogonal projection onto the manifold $M$, which is embedded in a Euclidean space. \\[0.1cm]
\textbf{Proof of Proposition \ref{prop:retractions1}\,:} note that (a) is immediate, by (\ref{eq:retdefinition}), and the fact that $\mathrm{Exp}$ is a retraction. To prove (b), note that 
\begin{equation} \label{eq:PHInorm}
  \Phi_x(v) = \tau^i(\mathrm{Ret}_x(v))\,\partial_i(x)
\end{equation}
where $(\tau^i\,;i=1,\ldots,n)$ are normal coordinates with origin at $x$, and where $\partial_i = \left.\partial\middle/\partial \tau^i\right.$. Since $\Phi_x$ is smooth (precisely, $C^2$), 
$$
\Phi^{\prime\prime}_x(0_x)(v,v) = \left. \frac{d^{\scriptscriptstyle\hspace{0.03cm}2}}{dt^{\scriptscriptstyle 2}}\right|_{t=0}\! \Phi_x(tv) 
$$
Thus, if $c(t) = \mathrm{Ret}_x(tv)$ and $c^i(t) = (\tau^i \circ c)(t)$, then
$$
\Phi^{\prime\prime}_x(0_x)(v,v) = \frac{d^{\scriptscriptstyle\hspace{0.03cm}2}}{dt^{\scriptscriptstyle 2}}c^i(0)\,\partial_i(x) = 
\left\lbrace\frac{d^{\scriptscriptstyle\hspace{0.03cm}2}}{dt^{\scriptscriptstyle 2}}c^i(0) + \Gamma^i_{jk}(c(0))\frac{d}{dt}c^{\,j}(0)\frac{d}{dt}c^k(0)\right\rbrace\partial_i(x)
$$
where the second equality holds since $\Gamma^i_{jk}(c(0)) = \Gamma^i_{jk}(x) = 0$, by the definition of normal coordinates. Comparing to (\ref{eq:geodesicequation}) and (\ref{eq:accelerationcoordinates}), it is clear $\Phi^{\prime\prime}_x(0_x)(v,v) = \ddot{c}(0)$. \\[0.1cm]
\textbf{Proof of Proposition \ref{prop:retractions2}\,:} this is a direct application of (\ref{eq:taylor2}), using $\mathrm{Ret}_x(v)  =\mathrm{Exp}_x(\Phi_x(v))$. \\[0.1cm]
\textbf{Remark\,:} the claims in Examples 1 and 2 above will not be proved in detail. Example 1 is quite elementary, and only requires one to recall that $\mathrm{Cut}(x) = \lbrace - x\rbrace$. For Example 2, the cut locus on a point $x$ in $U(d)$ is described in~\cite{sakai}, and (\ref{eq:phiun}) follows by a straightforward matrix calculation. 

\section{Example\,: a retraction for $\mathrm{Gr}_{\scriptscriptstyle \mathbb{R}}(p\,,q)$} \label{sec:grassret}
Let $\mathrm{St}_{\scriptscriptstyle \mathbb{R}}(p\,,q)$ denote the Stiefel manifold, whose elements are the $d \times p$ matrices $b$ with 
$b^\dagger b = \mathrm{I}_p$\hfill\linebreak ($\mathrm{I}_p$ is the $p \times p$ identity matrix, and $d = p+q$). Note that $T_{b\,}\mathrm{St}_{\scriptscriptstyle \mathbb{R}}(p\,,q) = \lbrace w\,: w^\dagger b + b^\dagger w = 0\rbrace$. For $w \in T_{b\,}\mathrm{St}_{\scriptscriptstyle \mathbb{R}}(p\,,q)$, let $[b] = bb^\dagger$ and $[w] = wb^\dagger + bw^\dagger$. If $v \in T_x \mathrm{Gr}_{\scriptscriptstyle \mathbb{R}}(p\,,q)$, one says that $(b,w)$ is representative of $(x,v)$, whenever $x = [b]$ and $v = [w]$. 

Recall that $x$ and $v$ may always be expressed $x = g\cdot o$ and $v = g\cdot \tilde{\omega}$, using (\ref{eq:tangentgrassgroup}). If $x = [b]$, then $g$ may be chosen 
$g = (b,b^{\scriptscriptstyle \perp})$ (the columns of $b^{\scriptscriptstyle \perp}$ span the orthogonal complement of the image space of $x$). Then, a direct calculation shows $v = [w_{\scriptscriptstyle v}]$, where $w_{\scriptscriptstyle v} = b^{\scriptscriptstyle \perp}\omega$. Now, define
\begin{equation} \label{eq:retgrassman1}
  \mathrm{Ret}_x(v) = \mathrm{Proj}\left( b + w_{\scriptscriptstyle v}\right) \text{ for some $b$ such that } x = [b]
\end{equation}
where $\mathrm{Proj}(h)$ denotes the orthogonal projector onto the span of the columns of $h$, for $h \in \mathbb{R}^{\scriptscriptstyle d\times p}$. This is well-defined, since it does not depend on the choice of $b$ and $b^{\scriptscriptstyle \perp}$, and is indeed a retraction, since it verifies (\ref{eq:retdefinition}).

For a nicer expression of (\ref{eq:retgrassman1}), identify each $x \in \mathrm{Gr}_{\scriptscriptstyle \mathbb{R}}(p\,,q)$ with its image space $\mathrm{Im}(x)$.\hfill\linebreak In other word, consider $\mathrm{Gr}_{\scriptscriptstyle \mathbb{R}}(p\,,q)$ as the space of all $p$-dimensional subspaces of $\mathbb{R}^d$. Then, 
\begin{equation} \label{eq:retgrassman}
\mathrm{Ret}_x(v) = \mathrm{Span}\left( b + w_{\scriptscriptstyle v}\right)
\end{equation}
where $\mathrm{Span}(h)$ denotes the span of the column space of $h \in \mathbb{R}^{\scriptscriptstyle d\times p}$. 
\begin{proposition} \label{prop:phigrass}\footnote{Whenever $a = (\alpha\,,0_{\scriptscriptstyle p\times q})^\dagger$, where $\alpha$ is $p\times p$ and diagonal, let $\arctan(a) = (\arctan(\alpha)\,,0_{\scriptscriptstyle p\times q})^\dagger$ where $\arctan(\alpha) = \mathrm{diag}(\arctan(\alpha_{ii}))$. For the proof on the following page, define $\cos(a)$ and $\sin(a)$ in the same way.}  The retraction $\mathrm{Ret}$ in (\ref{eq:retgrassman}) is regular, and the corresponding maps $\Phi_x$ (defined as in (\ref{eq:PHI}) are given by
\begin{equation} \label{eq:phigrass}
  \Phi_x(v) = \left[ b^{\scriptscriptstyle \perp}(r\arctan(a) s^\dagger) \right]
\end{equation}
for $x = [b]$ and $v = [b^{\scriptscriptstyle \perp}\omega]$, where $\omega$ has s.v.d. $\omega = ras^\dagger$ with $r \in O(q)$ and $s \in O(p)$.
\end{proposition}
As for Examples 1 and 2 in \ref{sec:retractions}, (b) of Proposition \ref{prop:retractions1} now implies $\mathrm{Ret}$ is a  geodesic retraction. 

\vfill
\pagebreak
\noindent \textbf{Proof of Proposition \ref{prop:phigrass}\,:} here, $x \in \mathrm{Gr}_{\scriptscriptstyle \mathbb{R}}(p\,,q)$ is identified with its image space, $\mathrm{Im}(x)$. Without loss of generality, it is assumed $p \leq q$. 

With $\Phi_x$ given by (\ref{eq:phigrass}), the aim will be to show that, for $x \in \mathrm{Gr}_{\scriptscriptstyle \mathbb{R}}(p\,,q)$ and $v \in T_x \mathrm{Gr}_{\scriptscriptstyle \mathbb{R}}(p\,,q)$,
\begin{equation} \label{eq:phigrassproof1}
  \mathrm{Exp}_x(\Phi_x(v)) = \mathrm{Ret}_x(v)
\end{equation}
In~\cite{sakai}, the cut locus of $x$ is obtained under the form
\begin{equation} \label{eq:cutgrass}
  \mathrm{Cut}(x) = \left\lbrace\mathrm{Exp}_x\left([b^{\scriptscriptstyle \perp}\omega]\right)\,;\,\omega = ras^\dagger\,,\,\Vert a\Vert_{\scriptscriptstyle \infty} = \frac{\pi}{2}\right\rbrace
\end{equation}
where $\omega= ras^\dagger$ is the s.v.d. of $\omega \in \mathbb{R}^{\scriptscriptstyle\,q\times p}$, with $r \in O(q)$ and $s \in O(p)$, and $\Vert a\Vert_{\scriptscriptstyle \infty} = \max_{ij}|a_{ij}|$. Since $\Vert \arctan(a)\Vert_{\scriptscriptstyle \infty} < \pi/2$, it follows from(\ref{eq:phigrass}) and (\ref{eq:phigrassproof1}) that $\mathrm{Ret}_x(v) \notin \mathrm{Cut}(x)$, so $\mathrm{Ret}$ is a regular retraction. Thus, to prove the proposition, one only has to prove (\ref{eq:phigrassproof1}).

Starting with the left-hand side of (\ref{eq:phigrassproof1}), let $\varphi = r \arctan(a) s^\dagger$, so $\Phi_x(v) = [ b^{\scriptscriptstyle \perp}\varphi]$. By the discussion before (\ref{eq:retgrassman1}), it follows that 
\begin{equation} \label{eq:proofgrass1}
 \Phi_x(v) = g \cdot \tilde{\varphi} \hspace{1cm} \text{(where $g = (b,b^{\scriptscriptstyle \perp})$)}
\end{equation}
However, then, by (\ref{eq:grasslift}),
\begin{equation} \label{eq:proofgrass2}
  \mathrm{Exp}_x(\Phi_x(v)) =   \exp(g\cdot \hat{\varphi})\cdot x = \left(g\exp(\hat{\varphi})\right)\cdot o
\end{equation}
where the second equality follows from $g^\dagger x = o$, using $g\cdot \hat{\varphi} = g\,\hat{\varphi}\,g^\dagger$. Using the s.v.d. of $\varphi$ ($\varphi = r\arctan(a) s^\dagger$), a straightforward matrix multiplication yields
\begin{equation} \label{eq:proofgrass3}
  \hat{\varphi} = k\cdot\hat{q} \hspace{0.5cm} \text{where } k = \left(\begin{array}{cc} s & 0 \\[0.1cm] 0 & r \end{array}\right) \;,\;
q = \arctan(a)
\end{equation}
Thus, from (\ref{eq:proofgrass2}) and (\ref{eq:proofgrass3}), using the fact that $k \in O(p) \times O(q)$, so $k\cdot o = o$ (or $k^\dagger \cdot o = o$),
$$
\mathrm{Exp}_x(\Phi_x(v)) = \left(g\,k \exp(\hat{q})\right)\cdot o  
$$
That is, by the group action property,
\begin{equation} \label{eq:proofgrass33}
\mathrm{Exp}_x(\Phi_x(v)) = g\,k\cdot\left(\exp(\hat{q})\cdot o\right)
\end{equation}
Now, let $b_o = (\mathrm{I}_p\,,0_{\scriptscriptstyle p\times q})^\dagger$, so $o = \mathrm{Span}(b_o)$ and
\begin{equation} \label{eq:proofgrass4}
 \exp(\hat{q})\cdot o = \mathrm{Span}\left( \exp(\hat{q})\,b_o\right)
\end{equation}
Then, let $a = (\alpha\,,0_{\scriptscriptstyle p\times q})^{\dagger\,}$, where $\alpha$ is $p \times p$ and diagonal. It will be shown below that
\begin{equation} \label{eq:proofgrass5}
  \exp(\hat{q})\,b_o \,=\, \left(\begin{array}{c} \cos(\arctan(\alpha))\\[0.1cm] \sin(\arctan(a))\end{array}\right) = 
   \left(\begin{array}{c} \mathrm{I}_p\\[0.1cm] a\end{array}\right)\left(\mathrm{I}_p + \alpha \right)^{-\frac{1}{2}}
\end{equation} 
where the second equality follows from the identities 
$$
\cos(\arctan(\alpha_{ii})) = (1+ \alpha^2_{ii})^{-\frac{1}{2}} \text{ and } \sin(\arctan(\alpha_{ii})) = \alpha_{ii}(1+ \alpha^2_{ii})^{-\frac{1}{2}}
$$ 
By (\ref{eq:proofgrass4}) and (\ref{eq:proofgrass5}), after ignoring the invertible matrix $\left(\mathrm{I}_p + \alpha \right)^{-\frac{1}{2}}$,
$$
\exp(\hat{q})\cdot o = \mathrm{Span}\left(\begin{array}{c} \mathrm{I}_p\\[0.1cm] a\end{array}\right) = \mathrm{Span}\left(b_o + b^{\scriptscriptstyle\perp}_oa\right)
$$
Replacing this into (\ref{eq:proofgrass33}), it follows that
\begin{equation} \label{eq:proofgrass6}
\mathrm{Exp}_x(\Phi_x(v)) = g\,k\cdot\mathrm{Span}\left(b_o + b^{\scriptscriptstyle\perp}_oa\right) = 
g\cdot \mathrm{Span}\left( k(b_o + b^{\scriptscriptstyle\perp}_oa)\right)
\end{equation}
and, by carrying out the matrix products, one may perform the simplification, 
$$
\mathrm{Span}\left( k(b_o + b^{\scriptscriptstyle\perp}_oa)\right) = 
\mathrm{Span}\left(b_o + b^{\scriptscriptstyle\perp}_ora\right) = 
\mathrm{Span}\left(b_o + b^{\scriptscriptstyle\perp}_o\omega\right)
$$ 
to obtain from (\ref{eq:proofgrass6}),
$$
\mathrm{Exp}_x(\Phi_x(v)) = g\cdot\mathrm{Span}\left(b_o + b^{\scriptscriptstyle\perp}_o\omega\right)
$$
which immediately yields (\ref{eq:phigrassproof1}), since $g\,b_o = b$ and  $g\,b^{\scriptscriptstyle \perp}_o = b^{\scriptscriptstyle \perp}$. \\[0.1cm]
\textbf{Proof of (\ref{eq:proofgrass5})\,:} write $q = (\kappa\,,0_{\scriptscriptstyle p\times q})^{\dagger\,}$, where $\kappa$ is $p \times p$ and diagonal. It is enough to show
\begin{equation} \label{eq:finalgrassproof1}
\exp(\hat{q}) = \left(\begin{array}{cc}\cos(\kappa) & -\sin(q)^\dagger \\[0.1cm] \sin(q) & \cos(\kappa)_{\scriptscriptstyle q\times q}\end{array}\right)
\end{equation}
wehre $\cos(\kappa)_{\scriptscriptstyle q\times q}$ is the $q \times q$ matrix,
$$
\cos(\kappa)_{\scriptscriptstyle q\times q} = \left(\begin{array}{cc} \cos(\kappa) & \\[0.1cm] & \mathrm{I}_{q-p}\!\end{array}\right)
$$
This follows by writing
$$
\hat{q} = \sum^p_{i=1}\,\kappa_{ii}\,\hat{f}_{i} \hspace{1cm} f_{i} = (\delta_{i}\,,0_{\scriptscriptstyle p\times q})^\dagger
$$
where $\delta_i$ is $p\times p$, diagonal, with its only non-zero element on the $i$-th line, and equal to $1$. Indeed, the matrices $\hat{f}_{i}$ commute with one another, so that 
\begin{equation} \label{eq:finalgrassproof3}
  \exp(\hat{q}) = \prod^p_{i=1}\,\exp(\kappa_{ii}\,\hat{f}_{i})
\end{equation}
and one readily checks $\hat{f}^{\scriptscriptstyle\hspace{0.03cm}2}_i = -e_{i\,}$, where $e_i$ is  $d\times d$, diagonal, with its only non-zero elements on the $i$-th and $(p+i)$-th lines, and equal to $1$. Therefore, 
\begin{equation} \label{eq:finalgrassproof2}
  \exp(t\,\hat{f}_i) = \mathrm{I}_d + (\cos(t)-1)\,e_i + \sin(t)\,\hat{f}_i
\end{equation}
Then, (\ref{eq:finalgrassproof1}) obtains after replacing (\ref{eq:finalgrassproof2}) into (\ref{eq:finalgrassproof3}), and using 
$$
\begin{array}{ccccc} 
e_i\,e_j = 0 && e_i\,\hat{f}_j = 0 & &\\[0.2cm]
\hat{f}_i\,e_j = 0 && \hat{f}_i\,\hat{f}_j = 0 & & \text{for $i \neq j$} 
\end{array}
$$
which may be shown by performing the matrix products. \\[0.1cm]
\textbf{Remark\,:} the above proof has a flavor of the structure theory of Riemannian symmetric spaces. In fact, $\mathrm{Gr}_{\scriptscriptstyle \mathbb{R}}(p\,,q) = \left. O(p+q) \middle/ O(p) \times O(q)\right.$ is a Riemannian symmetric space. The associated Cartan decomposition is
\begin{equation} \label{eq:grasscartan}
  \mathfrak{o}(p+q) = \mathfrak{k} \,+\,\mathfrak{p}
\end{equation} 
where $\mathfrak{k}$ is the Lie algebra of $K = O(p)\times O(q)$, and where $\mathfrak{p}$ was given in (\ref{eq:grassp}). Then,
\begin{equation} \label{eq:grassaa}
  \mathfrak{a} = \left\lbrace \hat{a}\,;\, a = (\alpha\,,0_{\scriptscriptstyle p\times q})^{\dagger\,}\,,\, \alpha \text{ is $p\times p$ diagonal}\right\rbrace
\end{equation}
is a maximal Abelian subspace of $\mathfrak{p}$. From~\cite{helgason} (Lemma 6.3, Chapter V), it follows that any $\hat{\omega} \in \mathfrak{p}$ is of the form $\hat{\omega} = \mathrm{Ad}(k)\,\hat{a}$ where $\mathrm{Ad}$ denotes the adjoint representation, $k \in K$ and $\hat{a} \in \mathfrak{a}$.\hfill\linebreak In the present context, this reads $\hat{\omega} = k\cdot \hat{a}$, which is indeed realised if $\omega$ has s.v.d. $\omega = ras^\dagger$, and $k$ is the same as in (\ref{eq:proofgrass3}). 

\section{The squared distance function} \label{sec:squareddistance}

\subsection{Cut locus} \label{subsec:cutlocus}
A Riemannian manifold $M$ becomes a metric space, when equipped with the distance function
\begin{equation} \label{eq:distance}
  d(x,y) = \inf \left\lbrace  L(c)\,;\, c \in C^1([0,1]\,,M):c(0) = x \,,\,c(1) = y\,\right\rbrace
\end{equation}
known as the \textit{Riemannian distance}.  Here, $L(c)$ is the length functional (\ref{eq:length}). When $M$ is geodesically complete, the infimum in (\ref{eq:distance}) is always achieved by some curve $c^*$, which is then said to be length-minimising. In addition, any length-minimising curve is a geodesic. 

This is not to say that all geodesics are length-minimising. A geodesic curve $c$, with $c(0) = x$, may reach a point $c(t) = y$, such that $L(\left.c\right|_{\scriptscriptstyle [0,t]}) \geq d(x,y)$. Roughly, this happens when $t$ is so large that $c$ becomes too long.  %\\[0.05cm]%

For $v \in T_xM$ with $\Vert v \Vert_x = 1$ ($\Vert\cdot \Vert_x$ is the norm given by the scalar product $\langle\cdot,\cdot\rangle_x$), define
\begin{equation} \label{eq:tv}
\mathrm{t}(v) = \sup\left\lbrace t\geq 0 : L(\left.c_{\scriptscriptstyle v}\right|_{\scriptscriptstyle [0,t]}) = d(x,c_{\scriptscriptstyle v}(t))\right\rbrace
\end{equation}
where $c_{\scriptscriptstyle v}$ denotes the geodesic curve with $\dot{c}_{\scriptscriptstyle v}(0) = v$. The following sets
\begin{equation} \label{eq:tangentcut}
   \mathrm{TC}(x) = \left\lbrace t\,v\,;\, t = \mathrm{t}(v)\,,\,\Vert v\Vert_x = 1 \right\rbrace \hspace{1cm}
  \mathrm{TD}(x) = \left\lbrace t\,v\,;\, t < \mathrm{t}(v)\,,\,\Vert v\Vert_x = 1 \right\rbrace
\end{equation}
are known as the tangent cut locus and tangent injectivity domain of $x$. The cut locus and injectivity domain of $x$ are the sets $\mathrm{Cut}(x) = \mathrm{Exp}\left( \mathrm{TC}(x)\right)$ and $\mathrm{D}(x) = \mathrm{Exp}\left( \mathrm{TD}(x)\right)$.

Since any two points $x$ and $y$ in $M$ are connected by a length-minimising geodesic $c^*$, 
\begin{equation} \label{eq:cutdecomp}
  M = \mathrm{D}(x) \,\cup\,\mathrm{Cut}(x)
\end{equation}
It is interesting to note that $\mathrm{Cut}(x)$ is a closed and negligible set. 

\subsection{Normal coordinates} \label{subsec:normalcoordinates}
The exponential map $\mathrm{Exp}_x$ is a diffeomorphism of $\mathrm{TD}(x)$ onto $\mathrm{D}(x)$ ($\mathrm{TD}(x)$ is the largest subset of $T_xM$ with this property).  Pick some orthonormal basis $(u_i)$ of $T_xM$, and define, for $y \in \mathrm{D}(x)$, 
\begin{equation} \label{eq:normalcoordinates}
  \tau^i(y) = \left\langle \mathrm{Exp}^{-1}_x(y)\hspace{0.02cm}, u_i\hspace{0.02cm}\right\rangle_x \hspace{1cm} i = 1\,,\ldots,\,n
\end{equation}
Then, $\tau^i:\mathrm{D}(x) \rightarrow \mathbb{R}$ are well-defined local coordinates, known as normal coordinates. These coordinates satisfy
\begin{equation} \label{eq:normal0}
  \tau^i(x) = 0 \hspace{0.4cm} g_{ij}(x) = \delta_{ij} \hspace{0.5cm} \Gamma^i_{jk}(x) = 0
\end{equation}
in the notation of (\ref{eq:lengthelement}) and (\ref{eq:christoffel}). Even more, (\ref{eq:normalcoordinates}) is equivalent to the property that geodesics through $x$ appear as straight lines through $0 \in \mathbb{R}^n$, in the normal coordinate map $\tau:\mathrm{D}(x) \rightarrow \mathbb{R}^n$. 

Now, the coordinate vector fields $\partial_i = \left.\partial\middle/\partial \tau^i\right.$ are given by 
\begin{equation} \label{eq:dexpnorm}
  \partial_i(y) \,=\,\mathrm{d\hspace{0.02cm}Exp}_x(v)(u_i) \hspace{1cm} \text{where } v = \tau^i(y)\hspace{0.02cm}u_i
\end{equation}
where $\mathrm{d\hspace{0.02cm}Exp}_x$ is the derivative of $\mathrm{Exp}_{x}:T_xM\rightarrow M$. This may be computed using Jacobi fields,
\begin{equation} \label{eq:dexpjacobi}
\mathrm{d\hspace{0.02cm}Exp}_x(tv)(tu) = J(t)  
\end{equation}
where $J$ is a vector field (Jacobi field) along the geodesic $c(t) = \mathrm{Exp}_x(tv)$, which solves the Jacobi equation
\begin{equation} \label{eq:jacobiequation}
  \nabla^2_{t\,}{J} - R(\dot{c},J)\dot{c} = 0
\end{equation}
where $J(0) = 0$, $\nabla_{t\,}J(0) = u$, and where $R$ denotes the Riemann curvature tensor. Of course, there do exist other means of computing the derivative $\mathrm{d\hspace{0.02cm}Exp}_x$ (\textit{e.g.} when $\mathrm{Exp}_x$ coincides with a matrix exponential). 
\vfill
\pagebreak
\subsection{Distance function} \label{subsec:rx}
For $x \in M$, consider the distance function $r_x(y)  = d(x,y)$. For $y \in \mathrm{D}(x)$, it is possible to show 
\begin{equation} \label{eq:localdistance}
  r_x(y) = \left( \sum^n_{i=1} \tau^i(y)^{\hspace{0.02cm}2}\right)^{\!\!\frac{1}{2}} \hspace{0.5cm}
%  f_x(y) = \frac{1}{2}\,\sum^n_{i=1} \tau^i(u)^{\hspace{0.02cm}2}
\end{equation}
in terms of the normal coordinates $\tau^i$. From (\ref{eq:localdistance}), the distance function $r_x$ is smooth on 
$\mathrm{U}_x = \mathrm{D}(x) - \lbrace x \rbrace$. %In the following, $r_x$ will simply be denoted by $r$. 

When $y \in \mathrm{U}_{x}$ is of the form $y = c_{\scriptscriptstyle v}(t)$, where $c_{\scriptscriptstyle v}$ is a geodesic with $\dot{c}_{\scriptscriptstyle v}(0) = v$ and $\Vert v \Vert_x = 1$, define $\partial_{\hspace{0.02cm}r}(y) = \dot{c}_{\scriptscriptstyle v}(t)$.  By the first variation of arc length formula (Theorem II.4.1 in~\cite{chavel}),
\begin{equation} \label{eq:gradr}
  \mathrm{grad}\,r_x(y) = \partial_{\hspace{0.02cm}r}(y) \hspace{1cm} \text{for } y \in \mathrm{U}_x 
\end{equation}
Introduce geodesic spherical coordinates $(r,\theta^{\scriptscriptstyle\,\alpha})$ on $\mathrm{U}_{x\,}$. If $y = c_{\scriptscriptstyle v}(t)$ these are given by $r = t$ and $(\theta^{\scriptscriptstyle\,\alpha}) = \theta(v)$, where $\theta$ identifies the unit sphere in $T_xM$ with the Euclidean unit sphere $S^{n-1}$.\hfill\linebreak In these coordinates, the metric is given by
\begin{equation} \label{eq:lengthspherical}
  g = dr \otimes dr \,+\, g^{\hspace{0.03cm} r}_{\scriptscriptstyle \alpha\beta}\;d\theta^{\scriptscriptstyle\,\alpha}\!\otimes\!d\theta^{\scriptscriptstyle\,\beta}
\end{equation}
reflecting the fact that $\partial_{\hspace{0.02cm}r}$ is orthogonal to constant $r_x$ surfaces, here parameterised by $(\theta^{\scriptscriptstyle\,\alpha})$.   
The coordinate vector fields $\partial_{\hspace{0.02cm}\scriptscriptstyle \alpha}$ are given by (\ref{eq:dexpjacobi})\,: $\partial_{\hspace{0.02cm}\alpha}(y) = J(r)$ for $y = c_{\scriptscriptstyle v}(r)$, where $J(0) = 0$ and $\nabla_{t\,}J(0) = u_{\hspace{0.02cm}\scriptscriptstyle \alpha}$ (where $u_{\hspace{0.02cm}\scriptscriptstyle \alpha} = \left.\partial\middle/\partial \theta^{\scriptscriptstyle\,\alpha}\right.$ are coordinate vector fields on the unit sphere in $T_xM$). In particular, if $A:T_xM \rightarrow T_yM$ solves the operator Jacobi equation (along the geodesic $c_{\scriptscriptstyle v}$)
\begin{equation} \label{eq:jacobioperator}
  \nabla^2_{t\,}{A} - R_{\dot{c}_{\scriptscriptstyle v}}\hspace{0.02cm}A = 0 \hspace{1cm} A(0) = 0 \,,\, \nabla_{t\,}A(0) = \mathrm{Id}_x
\end{equation}
where $R_{\hspace{0.02cm}\dot{c}_{\scriptscriptstyle v}}(\cdot) = R(\dot{c}_{\scriptscriptstyle v\hspace{0.02cm}},\cdot)\hspace{0.02cm}\dot{c}_{\scriptscriptstyle v\,}$, then $\partial_{\hspace{0.02cm}\alpha}(y) = A(r)\hspace{0.02cm} u_{\hspace{0.02cm}\scriptscriptstyle \alpha\,}$. Thus, if $\mathcal{A}(y):T_xM\rightarrow T_xM$ is given by $\mathcal{A}(y) = \Pi^{\scriptscriptstyle 0}_{r} \circ A(r)$, then $g^{\hspace{0.03cm} r}(y) = (\mathcal{A}(y))^*(h)$, the pullback under $\mathcal{A}(y)$ of the metric $h$ of the unit sphere in $T_xM$. It should be noted $\mathcal{A}(y)$ maps tangent spaces of this unit sphere to themselves. 

The Hessian of $r_x$ follows from (\ref{eq:hessian}) and (\ref{eq:gradr}), which yield (after using the fact that the vector fields $\partial_{\hspace{0.02cm}r}$ and $\partial_{\hspace{0.02cm}\scriptscriptstyle \alpha}$ commute)
$$
\mathrm{Hess}\,r_x\cdot \partial_{\hspace{0.02cm}r\hspace{0.02cm}} = 0 \hspace{0.3cm}\text{and}\hspace{0.3cm}\mathrm{Hess}\,r_x\cdot \partial_{\hspace{0.02cm}\alpha} =\nabla_{\partial_{\hspace{0.02cm}r}}\hspace{0.02cm}\partial_{\hspace{0.02cm}\alpha}
$$ 
Then, using the expression of the $\partial_{\hspace{0.02cm}\alpha}$ as Jacobi fields, 
\begin{equation} \label{eq:hessr}
  \mathrm{Hess}\,r_x(y) \,=\, \left.\nabla_{t\,}A(t)A^{-1}(t)\right|_{t=r}
\end{equation}
Taking the covariant derivative $\nabla_{t}$ of this formula yields the Ricatti equation
\begin{equation} \label{eq:ricatti}
  \nabla_{\partial_{\hspace{0.02cm}r}}  \mathrm{Hess}\,r_x \,=\, R_{\scriptscriptstyle \partial_{\hspace{0.02cm}r}} - \left(\mathrm{Hess}\,r_x\right)^{\hspace{0.02cm}2}
\end{equation}
The Jacobi equation (\ref{eq:jacobioperator}) and the Ricatti equation (\ref{eq:ricatti}) lead up to the comparison theorems\footnote{The inequalities (\ref{eq:metcomp}) and (\ref{eq:hescomp}) are in the sense of the usual Loewner order for self-adjoint operators.}.
\begin{theorem} \label{th:comp}
Assume the sectional curvatures of $M$ lie within the interval $[\kappa_{\min}\hspace{0.02cm},\kappa_{\max}\hspace{0.03cm}]\hspace{0.02cm}$. Then, 
\begin{equation} \label{eq:metcomp}
   \mathrm{sn}^2_{\kappa_{\max}}(r)\,h\,\leq\, g^{\hspace{0.03cm} r}(y) \,\leq\,    \mathrm{sn}^2_{\kappa_{\min}}(r)\,h   %g^{\hspace{0.03cm} r}(y)
\end{equation}  
\begin{equation} \label{eq:hescomp}
\mathrm{ct}_{\kappa_{\max}}(r)\,g^{\hspace{0.03cm} r}(y)\,\leq\, \mathrm{Hess}\,r_x(y)\,\leq\, 
\mathrm{ct}_{\kappa_{\min}}(r)\,g^{\hspace{0.03cm} r}(y)
\end{equation}
for $y \in \mathrm{U}_x\,$. Here, $\mathrm{sn}^{\prime\prime}_\kappa(r) + \kappa\,\mathrm{sn}_{\kappa}(r) = 0$ with $\mathrm{sn}_{\kappa}(0) = 0$ and $\mathrm{sn}^\prime_{\kappa}(0) = 1$, and $\mathrm{ct}_{\kappa} = \left. \mathrm{sn}^\prime_{\kappa}\middle/\mathrm{sn}_{\kappa}\right.$.
\end{theorem}
\noindent \textbf{Remark\,:} in addition to its singularity at $x$, the distance function $r_x$ is singular on $\mathrm{Cut}(x)$.\hfill\linebreak If $y\in\mathrm{Cut}(x)$, then either $y$ is a first conjugate point ($A(r)$ is singular, for the first time after $x$),\hfill\linebreak or there exist two distinct length-minimising geodesics connecting $x$ to $y$. In the first case, $\mathrm{Hess}\,r_x(y)$ has an eigenvalue equal to $-\infty$. In the second case, $\mathrm{grad}\,r_x$ is discontinuous at $y$.\hfill\linebreak The distributional Hessian of $r_x$ was studied in~\cite{distributional}. \\[0.1cm]
\noindent \textbf{Remark\,:} the reader may have noted, or recalled, that $y \in \mathrm{Cut}(x)$ if and only if $x \in \mathrm{Cut}(y)$.

\vfill
\pagebreak

\subsection{Squared distance} \label{ssec:sqd}
For $x \in M$, consider the squared distance function $f_x(y)  = d^{\hspace{0.03cm}2}(x,y)/2$. For $y \in \mathrm{D}(x)$,
\begin{equation} \label{eq:localfx}
f_x(y) = \frac{1}{2}\,\sum^n_{i=1} \tau^i(y)^{\hspace{0.02cm}2}
\end{equation}
in terms of the normal coordinates $\tau^i$. It follows that $f_x$ is smooth on $\mathrm{D}(x)$. Of course, $f_x = r^{\hspace{0.03cm}2}_x/2$. Therefore, applying the chain rule to (\ref{eq:gradr}), 
\begin{equation} \label{eq:gradfx}
  \mathrm{grad}\,f_x(y) \,=\, -\hspace{0.02cm}\mathrm{Exp}^{-1}_y(x) \hspace{0.5cm} \text{for } y \in \mathrm{D}(x)
\end{equation}
and, by another application of the chain rule,
\begin{equation} \label{eq:hessfx}
 \mathrm{Hess}\, f_x(y) = dr_x \otimes dr_x + r_x\,\mathrm{Hess}\,r_x 
\end{equation}
Just like $r_{x\,}$, $f_x$ is singular on $\mathrm{Cut}(x)$. If $y \in \mathrm{Cut}(x)$ is a first conjugate point, then $ \mathrm{Hess}\, f_x(y)$ has an eigenvalue equal to $-\infty$. 

The convexity of the function $f_x$ will play a significant rôle, in the following, especially when $M$ is a Hadamard manifold\,: a simply connected, geodesically complete Riemannian manifold of non-positive sectional curvature. When $M$ is a Hadamard manifold, the following properties hold\,:  any $x\hspace{0.02cm},y \in M$ are connected by a unique geodesic $c$\,;\,for all $x \in M$, $\mathrm{Cut}(x)$ is empty, and $f_x$ is smooth and $1/2$-strongly convex\,;\,all geodesic balls are convex (see the remarks below,\hfill\linebreak for the notions of convex set and function). 

Assume $M$ is a Hadamard manifold. In addition, assume that the sectional curvature of $M$ is bounded below by $\kappa_{\min} = -c^{\hspace{0.02cm}\scriptscriptstyle 2}$. Theorem \ref{th:comp} may be applied to (\ref{eq:hessfx}), after setting
$\kappa_{\max} = 0$. This yields
\begin{equation} \label{eq:hesfxcomp}
g(y) \,\leq\,\mathrm{Hess}\,f_x(y)\,\leq\, c\hspace{0.03cm} r_x(y)\coth(c\hspace{0.03cm} r_x(y))\,g(y)
\end{equation}
for $y \in M$. In addition to showing that $f_x$ is $1/2$-strongly convex, this shows that $\mathrm{Hess}\,f_x$ has, at most, linear growth
\begin{equation} \label{eq:hesfxcompbis}
\mathrm{Hess}\,f_x(y)\,\leq\, (1+c\hspace{0.03cm}r_x(y))\,g(y)
\end{equation}
since $x\coth(x) \leq 1 + x$ for $x \geq 0$. 
\\[0.1cm]
\noindent \textbf{Remark\,:} a subset $A \subset M$ is called convex (that is, strongly convex, in the terminology of~\cite{chavel}) if any $x\hspace{0.02cm},y \in A$ are connected by a unique length-minimising geodesic $c$,  and $c$ lies entirely in $A$.\hfill\linebreak A function $f:A \rightarrow \mathbb{R}$ is then called (strictly) convex if $f \circ c:\mathbb{R} \rightarrow \mathbb{R}$ is (strictly) convex, for any geodesic $c$ which lies in $A$. It is called $\alpha$-strongly convex (for some $\alpha > 0$) if $f \circ c:\mathbb{R} \rightarrow \mathbb{R}$ is $\alpha$-strongly convex, for any geodesic $c$ which lies in $A$, 
\begin{equation} \label{eq:strongconv}
(f\circ c)(p\hspace{0.02cm}s + q\hspace{0.02cm}t) \leq p\hspace{0.02cm}(f\circ c)(s) + q\hspace{0.02cm}(f\circ c)(t) - \alpha\hspace{0.02cm}p\hspace{0.02cm}q\,d^{\hspace{0.03cm}2}(c(s),c(t))
\end{equation}
whenever $p\hspace{0.02cm},q \geq 0$ and $p+q = 1$. For example, if $M$ is a sphere and $A$ is the open northern hemisphere, then $A$ is  convex. Then,  $f_x:A \rightarrow \mathbb{R}$, where $x$ denotes the north pole, is strictly convex, but not strongly convex. \\[0.1cm]
\textbf{Remark\,:} for $x \in M$, let $\mathrm{inj}(x) = d(x\hspace{0.02cm},\mathrm{Cut}(x))$ denote the injectivity radius at $x$. Then, let $\mathrm{inj}(M) = \inf_{x\in M}\mathrm{inj}(x)$, the injectivity radius of $M$. Assume all the sectional curvatures of $M$ are less than $\kappa_{\max} = c^{\hspace{0.02cm}\scriptscriptstyle 2}$. If $B(x,R)$ is a geodesic ball with radius $R \leq (1/2)\hspace{0.02cm}\min\lbrace \mathrm{inj}(M)\hspace{0.02cm},\pi\hspace{0.03cm}c^{\scriptscriptstyle -1}\rbrace$, then $B(x,R)$ is convex. Here, if $\kappa_{\max} = 0$, then $c^{\scriptscriptstyle -1}$ is understood to be $+\infty$. However, there do exist manifolds $M$ with negative sectional curvature, and with $\mathrm{inj}(M) = 0$ (\textit{e.g.} the quotient of the Poincar\'e upper half-plane, by a discrete group of translations). 

\section{Example\,: robust Riemannian barycentre} \label{sec:robust}
Let $M$ be a Hadamard manifold, with sectional curvatures bounded below by $\kappa_{\min} = -c^{\hspace{0.02cm}\scriptscriptstyle 2}$. Recall that $f_x$ is $1/2$-strongly convex, and $\mathrm{Hess}\,f_x$ has, at most, linear growth (as in (\ref{eq:hesfxcompbis})). On the other hand, consider the function
\begin{equation} \label{eq:hdist}
  V_x(y) =  \delta^{\hspace{0.02cm} \scriptscriptstyle 2}\left[\mathstrut 1+ \left(d(x,y)\middle/\delta\right)^{\scriptscriptstyle 2\,}\right]^{\scriptscriptstyle\frac{1}{\mathstrut 2}}\, -\,\delta^{\hspace{0.02cm} \scriptscriptstyle 2}
\end{equation}
where $\delta > 0$ is a cutoff parameter. Note that $V_x(y) \geq 0$, and $V_x(y) = 0$ if and only if $x = y$. Moreover, $V_x \sim f_{x\,}$, when $\left.d(x\hspace{0.02cm},y)\middle/\delta\right.$ is small, and $V_x \sim \delta\hspace{0.02cm}r_{x\,}$, when $\left.d(x\hspace{0.02cm},y)\middle/\delta\right.$ is large.
\begin{proposition} \label{prop:hdist}
Let $M$ be a Hadamard manifold, with sectional curvatures bounded below by $\kappa_{\min} = -c^{\hspace{0.02cm}\scriptscriptstyle 2}$. If $V_x:M\rightarrow \mathbb{R}$ is defined as in (\ref{eq:hdist}), then 
$V_x$ is smooth, strictly (but not strongly) convex, and $\mathrm{Hess}\,V_x$ is bounded by $1+\delta\hspace{0.02cm}c$.
\end{proposition}
Let $\pi$ be a probability measure on $M$, and consider the problem of minimising
\begin{equation} \label{eq:huberpi}
  V_{\pi}(y) = \int_M\,V_x(y)\,\pi(dx)
\end{equation}
A global minimiser of $V_\pi$ will be called a robust Riemannian barycentre of $\pi$. Here, the adjective ``robust"  comes from the field of robust statistics~\cite{huber}.
\begin{proposition} \label{prop:robustbcentre}
  Let $\pi$ be a probability distribution on a Hadamard manifold $M$. If $\pi$ has finite first-order moments, then the function $V_{\pi}$ is a proper, strictly convex function, with a unique global minimum $x^* \in \Theta$. Therefore, $\pi$ has a unique robust Riemannian barycentre $x^*$.
\end{proposition}
Recall that $\pi$ has finite first-order moments, if and only if there exists $y_o \in M$ with
\begin{equation} \label{eq:firstorder}
\int_M\,r_{x}(y_o)\,\pi(dx) \,<\,\infty
\end{equation}
and recall that $V_{\pi}$ is said to be proper if it takes on finite values.\\[0.1cm]
\noindent \textbf{Proof of Proposition \ref{prop:hdist}\,:} by applying the chain rule to (\ref{eq:hdist}), and using (\ref{eq:gradfx}),
\begin{equation} \label{eq:hgrad}
  \mathrm{grad}\,V_x(y) = -\,\frac{\mathrm{Exp}^{-1}_y(x)}{\mathstrut \left[\mathstrut 1+ \left(d(x,y)\middle/\delta\right)^{\scriptscriptstyle 2\,}\right]^{\scriptscriptstyle\frac{1}{\mathstrut 2}}}
\end{equation} 
Then, by applying (\ref{eq:hessian}),
\begin{equation} \label{eq:hhess}
\mathrm{Hess}\,V_x(y) = -\,
\frac{{\small \mathrm{Exp}^{-1}_y(x)\otimes \mathrm{Exp}^{-1}_y(x)}}{\mathstrut\delta^{\hspace{0.02cm} \scriptscriptstyle 2} \left[\mathstrut 1+ \left(d(x,y)\middle/\delta\right)^{\scriptscriptstyle 2\,}\right]^{\scriptscriptstyle \frac{3}{\mathstrut 2}}}\,-\, \frac{\nabla\,\mathrm{Exp}^{-1}_y(x)}{\mathstrut \left[\mathstrut 1+ \left(d(x,y)\middle/\delta\right)^{\scriptscriptstyle 2\,}\right]^{\scriptscriptstyle\frac{1}{\mathstrut 2}}}
\end{equation}
\vspace{0.2cm}
To conclude, it is enough to note the inequalities,
$$
0\,\leq\,\mathrm{Exp}^{-1}_y(x)\otimes \mathrm{Exp}^{-1}_y(x)\leq d^{\hspace{0.03cm}\scriptscriptstyle 2}(x,y)\hspace{0.03cm}g(y) %\hspace{0.4cm}\text{and}\hspace{0.4cm}
%1\,\leq\,-\,\nabla\,\mathrm{Exp}^{-1}_y(x)\,\leq\, (1+\kappa\hspace{0.02cm} r_x(y))
$$
which follows since $\mathrm{Exp}^{-1}_y(x)\otimes \mathrm{Exp}^{-1}_y(x)$ is a rank-one operator in $T_yM$, and
$$
g(y)\,\leq\,-\,\nabla\,\mathrm{Exp}^{-1}_y(x)\,\leq\, (1+c\hspace{0.03cm}r_x(y))\hspace{0.03cm}g(y)
$$
which is the same as (\ref{eq:hesfxcompbis}), and follows from (\ref{eq:hessian}) and (\ref{eq:gradfx}). Replacing these into (\ref{eq:hhess}), a direct calculation shows
\begin{equation} \label{eq:hdistproof}
0\,<\,\mathrm{Hess}\,V_x(y)\,\leq\, (1+\delta\hspace{0.02cm}c)\hspace{0.03cm}g(y)
\end{equation}
which completes the proof. \vfill\pagebreak
\noindent \textbf{Proof of Proposition \ref{prop:robustbcentre}\,:} using the sub-additivity of the square root, (\ref{eq:hdist}) and (\ref{eq:huberpi}) imply that for any $y \in M$,
$$
V_{\pi}(y) \,\leq\, \int_M\,r_{x}(y)\,\pi(dx) 
$$
But, by the triangle inequality, and (\ref{eq:firstorder}),
$$
\int_M\,r_{x}(y)\,\pi(dx) \leq d(y\hspace{0.02cm},y_o) + \int_M\,r_{x}(y_o)\,\pi(dx) \,<\infty
$$
Therefore, $V_\pi$ is proper. That $V_\pi$ is also strictly convex is an immediate result of Proposition \ref{prop:hdist}\,:\hfill\linebreak each function $V_x$ is strictly convex, and $V_\pi(y)$ is the expectation of $V_x(y)$ with respect to a random $x$ with distribution $\pi$. Now, to show that $V_\pi$ has a unique global minimum, it is enough to show that $V_\pi(y)$ goes to infinity as $y$ goes to infinity. Note that $\varphi(x) = (1+x^{\scriptscriptstyle 2})^{\scriptscriptstyle \frac{1}{2}}$ is convex. This implies (using the elementary fact that the graph of a convex function remains above any of its tangents),
$$
V_x(y) \geq (\sqrt{2} - 1)\hspace{0.04cm}\delta^{\hspace{0.02cm}\scriptscriptstyle 2}\,+\, \frac{\delta}{\sqrt{2}}\hspace{0.04cm}r_x(y)
$$
Taking the expectation with respect to $\pi$,
$$
V_\pi(y) \geq (\sqrt{2} - 1)\hspace{0.04cm}\delta^{\hspace{0.02cm}\scriptscriptstyle 2}\,+\, \frac{\delta}{\sqrt{2}}\hspace{0.04cm}\int_M\,r_x(y)\,\pi(dx)
$$
To see that $V_\pi(y)$ goes to infinity as $y$ goes to infinity, it is now enough to note, using the triangle inequality,
$$
\int_M\,r_x(y)\,\pi(dx)\,\geq\,d(y\hspace{0.02cm},y_o)\,-\int_M\,r_{x}(y_o)\,\pi(dx)
$$
where $d(y\hspace{0.02cm},y_o)$ goes to infinity as $y$ goes to infinity. \\[0.1cm]
\noindent \textbf{Remark\,:} the above Proposition \ref{prop:robustbcentre} only requires $M$ to be a Hadamard manifold, without the additional condition that it have sectional curvatures bounded below. Indeed, Proposition \ref{prop:robustbcentre} only relies on the fact that $V_x$ is strictly convex, and not on the fact that the Hessian of $V_x$ is bounded above by $1+\delta\hspace{0.02cm}c$.\\[0.1cm]
\noindent\textbf{Remark\,:} if a function $V:M\rightarrow \mathbb{R}$, on a Riemannian manifold $M$, has bounded Hessian, then it has Lipschitz-gradient. 
That is, if there exists $\ell \geq 0$ such that $\left|\mathrm{Hess}\,V(x)(u,u)\right|\leq \ell\hspace{0.02cm}g(u,u)$ for all $x \in M$ and $v \in T_xM$, then
\begin{equation} \label{eq:lipschitzgrad}
 \left\Vert \Pi^{\scriptscriptstyle 0}_{\scriptscriptstyle 1}\left(\mathrm{grad}\,V_{c(1)}\right) - \mathrm{grad}\,V_{c(0)}\hspace{0.03cm}\right\Vert_{c(0)}\,\leq\,\ell\hspace{0.02cm}L(c)
\end{equation}
for any smooth curve $c:[0,1]\rightarrow M$, where $L(c)$ is the length of $c$. This is due to the following. 
\begin{lemma} \label{lem:lipschitzlemma}
  Let $X$ be a vector field on a Riemannian manifold $M$. If the operator norm of the covariant derivative $\nabla X$ is bounded by $\ell \geq 0$, then
\begin{equation} \label{eq:lipschitzfield}
 \left\Vert \Pi^{\scriptscriptstyle 0}_{\scriptscriptstyle 1}\left(X_{c(1)}\right) - X_{c(0)}\hspace{0.03cm}\right\Vert_{c(0)}\,\leq\,\ell\hspace{0.02cm}L(c)
\end{equation}
for any smooth curve $c:[0,1]\rightarrow M$. 
\end{lemma}
\noindent \textbf{Sketch of proof\,:} let $u_i$ be a parallel orthonormal base along $c$ ($u_i$ are vector fields along $c$, with $u_i(t)$ an orthonormal basis of $T_{c(t)}M$, for each $t$). Let $X^i(t) = \langle X\hspace{0.02cm},u_i\rangle_{c(t)\,}$ and note 
$$
 \left\Vert \Pi^{\scriptscriptstyle 0}_{\scriptscriptstyle 1}\left(X_{c(1)}\right) - X_{c(0)}\hspace{0.03cm}\right\Vert^2_{c(0)} = 
 \sum^n_{i=1}\left(X^i(1) - X^i(0)\right)^{\!2} = \sum^n_{i=1}\left(\int^1_0\,\langle \nabla_{\dot{c}\,}X\hspace{0.02cm},u_{i\hspace{0.02cm}}\rangle_{c(t)\,}dt\right)^{\!\!2}
$$
the proof then follows by using Jensen's inequality, since $\Vert \nabla_{\dot{c}\,}X\Vert_{c(t)} \leq \ell\hspace{0.03cm}\Vert \dot{c}\Vert_{c(t)\,}$. 
%construction open problem
\section{Riemannian volume and integral formulae} \label{sec:volume}

\subsection{Elementary volume comparison} \label{subsec:volcomp}
If a Riemannian manifold $M$ is orientable, then $M$ admits a volume form, called the Riemannian volume form, to be denoted $\mathrm{vol}$, in the following. In terms of local coordinates $(x^i\,;i=1,\ldots,n)$
\begin{equation} \label{eq:volumeform}
  \mathrm{vol} = \det(g)^{\frac{1}{2}}\,dx^1 \wedge \ldots \wedge dx^n
\end{equation}
where $\det(g)$ is the determinant of the metric, which is equal the determinant of the matrix $(g_{ij})$, defined in  (\ref{eq:lengthelement}). Then, the integral of a continuous, compactly-supported function $f:M\rightarrow \mathbb{R}$, with respect to $\mathrm{vol}$, is the integral of the $n$-form $f\hspace{0.02cm}\mathrm{vol}$ over $M$. This is denoted $\int_M\, f(x)\,\mathrm{vol}(dx)$. There exists a unique measure $|\mathrm{vol}|$ on the Borel $\sigma$-algebra of $M$, such that~\cite{marsden} (Chapter 8), for continuous, compactly-supported $f$,
$$
\int_M\, f(x)\,\mathrm{vol}(dx) \,=\, \int_M\,f(x)\,|\mathrm{vol}|(dx)
$$
where the integral on the left is a Riemann integral, and the integral on the right is a Lebesgue integral. It is quite useful to study these integrals using geodesic spherical coordinates (which were introduced in \ref{subsec:rx}). Let $(r,\theta^{\scriptscriptstyle\,\alpha})$ be geodesic spherical coordinates, with origin at $x \in M$. Recall that these are defined on $\mathrm{U}_{x} = \mathrm{D}(x) - \lbrace x \rbrace$, where $\mathrm{D}(x)$ is the injectivity domain of $x$. Since $M$ can be decomposed as in (\ref{eq:cutdecomp}), and $\mathrm{Cut}(x)$ is negligible,
\begin{equation} \label{eq:integralux}
\int_M\,f(y)\,\mathrm{vol}(dy) \,=\,\int_{\,\mathrm{U}_{x}}f(y)\,\mathrm{vol}(dy)
\end{equation}
Using (\ref{eq:lengthspherical}) and (\ref{eq:volumeform}), $\mathrm{vol}(dy) = \det(\mathcal{A}(y))\,dr\wedge \omega_{n-1}(d\theta)$, where $\omega_{n-1}$ is the area measure on the unit sphere in $T_xM$ (as of now, this is identified with the Euclidean unit sphere $S^{n-1}$). Using (\ref{eq:tangentcut}) and $\mathrm{D}(x) = \mathrm{Exp}\left( \mathrm{TD}(x)\right)$, (\ref{eq:integralux}) yields
\begin{equation} \label{eq:integralspherical}
\int_M\,f(y)\,\mathrm{vol}(dy) \,=\, \int^{\mathrm{t}(\theta)}_0\!\!\!\int_{S^{n-1}}f(r,\theta)\,\det(\mathcal{A}(r,\theta))\,dr\hspace{0.03cm} \omega_{n-1}(d\theta)
\end{equation}
where $\mathrm{t}$ was defined in (\ref{eq:tv}). This formula expresses integrals, with respect to the Riemannian volume form, using geodesic spherical coordinates. 

Recall the Laplacian $\Delta\hspace{0.03cm} r_x = \mathrm{div}\hspace{0.03cm}\partial_{\hspace{0.02cm}r\,}$. By definition of the divergence, $\mathcal{L}_{\partial_{\hspace{0.02cm}r}}\hspace{0.02cm}\mathrm{vol} = (\mathrm{div}\hspace{0.02cm}\partial_{\hspace{0.02cm}r})\hspace{0.02cm}\mathrm{vol}$. Writing this in geodesic spherical coordinates,
\begin{equation} \label{eq:laplacelogdet}
 \Delta\hspace{0.03cm} r_x(r,\theta) \,=\, \partial_{\hspace{0.02cm}r}\hspace{0.02cm}\log\det(\mathcal{A}(r,\theta))
\end{equation}
Accordingly, the comparison theorems \ref{th:comp} can be used to obtain the volume comparison theorem.
\begin{theorem} \label{th:compvol}
Assume the sectional curvatures of $M$ lie within the interval $[\kappa_{\min}\hspace{0.02cm},\kappa_{\max}\hspace{0.03cm}]\hspace{0.02cm}$. Then, 
\begin{equation} \label{eq:volcomp}
   \mathrm{sn}^{n-1}_{\kappa_{\max}}(r)\,\leq\, \det(\mathcal{A}(r,\theta)) \,\leq\,    \mathrm{sn}^{n-1}_{\kappa_{\min}}(r)  %g^{\hspace{0.03cm} r}(y)
\end{equation}  
\begin{equation} \label{eq:laplacecomp}
(n-1)\hspace{0.02cm}\mathrm{ct}_{\kappa_{\max}}(r)\,\leq\, \partial_{\hspace{0.02cm}r}\hspace{0.02cm}\log\det(\mathcal{A}(r,\theta))\,\leq\, 
(n-1)\hspace{0.02cm}\mathrm{ct}_{\kappa_{\min}}(r)
\end{equation}
\end{theorem}
This volume comparison theorem is quite elementary, as stronger and deeper comparison results do exist\footnote{For example, Gromov's volume comparison theorem can be used to give a short proof of the famous ``sphere theorem", Theorem III.4.6 in~\cite{chavel}. }. Moreover, in this theorem, the lower bound on sectional curvature may be replaced by a lower bound on Ricci curvature, without any change to the conclusion.  

\noindent \textbf{Remark\,:} roughly, (\ref{eq:volcomp}) states that ``more curvature means less volume". If $f:M\rightarrow \mathbb{R}$ is a non-negative function of distance to $x$, so $f(y) = f(r)$ in terms of the coordinates $(r,\theta^{\scriptscriptstyle\,\alpha})$, then
\begin{equation} \label{eq:integralcomp}
\omega_{n-1}\,\int^{\scriptscriptstyle R}_{\scriptscriptstyle 0}\,f(r)\,\mathrm{sn}^{n-1}_{\kappa_{\max}}(r)\hspace{0.02cm}dr\leq\,\int_{\scriptscriptstyle B(x,R)}f(y)\,\mathrm{vol}(dy)\,\leq\,\omega_{n-1}\,
\int^{\scriptscriptstyle R}_{\scriptscriptstyle 0}\,f(r)\,\mathrm{sn}^{n-1}_{\kappa_{\min}}(r)\hspace{0.02cm}dr
\end{equation}
for any $R \leq \min\lbrace \mathrm{inj}(x)\hspace{0.03cm},\pi\hspace{0.03cm}c^{\scriptscriptstyle -1}\rbrace$. Here, $\mathrm{inj}(x)$ is the injectivity radius at $x$, $c = |\kappa_{\max}|^{\scriptscriptstyle 1/2\hspace{0.03cm}}$, and $\omega_{n-1}$ denotes the area of $S^{\hspace{0.02cm}n-1}$. In addition, if $\kappa_{\max} \leq 0$, then $c^{\scriptscriptstyle -1}$ is understood to be $+\infty$. \\[0.1cm]
\indent In general, it may be impossible to apply the integral formula (\ref{eq:integralspherical}), since $\mathrm{t}(\theta)$ may be unknown. Here are two examples where $\mathrm{t}(\theta)$ is known, and quite tractable (in fact, constant). \\[0.1cm]
\textbf{Example 1\,:} if $M$ is a Hadamard manifold, then for any choice of the origin $x$, and any $\theta \in S^{n-1}$, one has $\mathrm{t}(\theta) = \infty$, and (\ref{eq:integralspherical}) becomes
\begin{equation} \label{eq:integralsphericalhadamard}
\int_M\,f(y)\,\mathrm{vol}(dy) \,=\, \int^{\infty}_0\!\!\!\int_{S^{n-1}}f(r,\theta)\,\det(\mathcal{A}(r,\theta))\,dr\hspace{0.03cm} \omega_{n-1}(d\theta)
\end{equation}
\textbf{Example 2\,:} compact rank-one symmetric space are the following manifolds\,: spheres, real projective spaces, complex projective spaces,  quaternion projective spaces, or the Cayley plane. These are manifolds all of whose geodesics are closed (i.e. periodic) and isometric to one another
(see~\cite{besseclosed}, for a detailed account). Therefore, $\mathrm{t}(\theta)$ does not depend on $x$ nor on $\theta$, but is always equal to $l/2$, where $l$ is the length of a simple geodesic loop. Scaling the metric so the maximum sectional curvature is equal to $1$, it can be shown $l = \pi$ for real projective spaces, and $l = 2\pi$\hfill\linebreak in all other cases. Moreover (\ref{eq:integralspherical}) takes on the form (this may be found by looking up the solution of the Jacobi equation in~\cite{besseclosed}, Page 82),
\begin{equation} \label{eq:integralsphericalcross}
\int_M\,f(y)\,\mathrm{vol}(dy) \,=\, \int^{\frac{l}{\mathstrut 2}}_0\!\!\!\int_{S^{n-1}}f(r,\theta)\left(\sin(r)\right)^{k-1}\left(2\sin(r/2)\right)^{n-k}\,dr\hspace{0.03cm} \omega_{n-1}(d\theta)
\end{equation}
where $k= n$ for spheres and real projective spaces, and $k = 2$ or $4$ for complex or quaternion projective spaces, respectively. For the Cayley plane, $n = 16$ and $k = 8$. 

\subsection{Riemannian symmetric spaces} \label{ssec:sspace}
A Riemannian symmetric space is a Riemannian manifold $M$, such that, for each $x \in M$, there exists an isometry $s_x : M \rightarrow M$, with $s_x(x) = x$ and $d\hspace{0.02cm}s_x(x) = -\mathrm{Id}_x\,$. This isometry $s_x$ is called the geodesic symmetry at $x$. 

Let $G$ denote the identity component of the isometry goup of $M$, and $K = K_o$ be the stabiliser in $G$ of some point $o \in M$\footnote{According to the Myers-Steenrod theorem, $G$ is a connected Lie group, and $K$ a compact subgroup of $G$.}. Then, $M = G/K$ is a Riemannian homogeneous space. The mapping $\theta : G \rightarrow G$, where $\theta(g) = s_o\circ g \circ s_o$ is an involutive isomorphism of $G$.

Let $\mathfrak{g}$ denote the Lie algebra of $G$, and consider the Cartan decomposition, $\mathfrak{g} = \mathfrak{k} + \mathfrak{p}$, where $\mathfrak{k}$ is the $+1$ eigenspace of $d\hspace{0.02cm}\theta$ and $\mathfrak{p}$ is the $-1$ eigenspace of $d\hspace{0.02cm}\theta$. One clearly has the commutation relations,
\begin{equation} \label{eq:sscommute}
  [\mathfrak{k},\mathfrak{k}] \subset \mathfrak{k}\hspace{0.2cm};\hspace{0.2cm}
  [\mathfrak{k},\mathfrak{p}] \subset \mathfrak{p}\hspace{0.2cm};\hspace{0.2cm}
  [\mathfrak{p},\mathfrak{p}] \subset \mathfrak{k}
\end{equation}
In addition, it turns out that $\mathfrak{k}$ is the Lie algebra of $K$, and that $\mathfrak{p}$ may be identified with $T_oM$.

The Riemannian metric of $M$ may always be expressed in terms of an $\mathrm{Ad}(K)$-invariant scalar product $Q$ on $\mathfrak{g}$. If $x \in M$ is given by $x = g\cdot o$ for some $g \in G$ (where $g\cdot o = g(o)$), then
\begin{equation} \label{eq:ssmetric}
  \langle u,\!v\rangle_{\scriptscriptstyle x} = Q(g^{\scriptscriptstyle -1}\cdot u\hspace{0.02cm},g^{\scriptscriptstyle -1}\cdot v)
\end{equation} 
where the vectors $g^{\scriptscriptstyle -1}\cdot u$ and $g^{\scriptscriptstyle -1}\cdot v$, which belong to $T_oM$, are identified with elements of $\mathfrak{p}$. Here, by an abuse of notation, $d\hspace{0.02cm}g^{\scriptscriptstyle -1}\cdot u$ is denoted $g^{\scriptscriptstyle -1}\cdot u$. \vfill\pagebreak

Let $\exp:\mathfrak{g} \rightarrow G$ denote the Lie group exponential. If $v \in T_oM$, then the Riemannian exponential $\mathrm{Exp}_o(v)$ is given by
\begin{equation} \label{eq:ssexp1}
  \mathrm{Exp}_o(v) = \exp(v)\cdot o
\end{equation}
Moreover, if $\Pi^t_{\scriptscriptstyle 0}$ denotes parallel transport along the geodesic $c(t) = \mathrm{Exp}_o(tv)$, then 
\begin{equation} \label{eq:ssparallel}
  \Pi^t_{\scriptscriptstyle 0}(u) = \exp(tv)\cdot u
\end{equation}
for any $u \in T_o M$ (note that the identification $T_oM \simeq \mathfrak{p}$ is always made, implicitly). Using (\ref{eq:ssparallel}), one can derive the following expression for the Riemann curvature tensor at $o$,
\begin{equation} \label{eq:sscurvature}
 R_o(v,u)w = -[[v\hspace{0.03cm},u]\hspace{0.02cm},w] \hspace{1cm} v,u,w \in T_oM
\end{equation}
A fundamental property of symmetric spaces is that the curvature tensor is parallel\,:$\hspace{0.03cm}\nabla\,R = 0$. This is often used to solve the Jacobi equation (\ref{eq:jacobiequation}), and then express the derivative of the Riemannian exponential, using in (\ref{eq:dexpjacobi}),
\begin{equation} \label{eq:dexpss}
  \mathrm{d\hspace{0.02cm}Exp}_x(v)(u) \,=\, \exp(v)\cdot \mathrm{sh}(R_v)(u)
\end{equation}
where $\mathrm{sh}(R_v) = \sum^{\infty}_{n=0} (R_v)^n/(2n+1)!$ for the self-adjoint curvature operator $R_v(u) = [v\hspace{0.03cm},[v\hspace{0.02cm},u]]$. Since $\exp(v)$ is an isometry, the following expression of the Riemannian volume is immediate
%(of course, the star denotes pullback)
\begin{equation} \label{eq:ssvol}
  \mathrm{Exp}^*_o(\mathrm{vol}) = \left|\det(\mathrm{sh}(R_v))\right|\hspace{0.02cm}dv
\end{equation}
where $dv$ denotes the volume form on $T_oM$, associated with the restriction of $Q$ to $\mathfrak{p}$.

Expression (\ref{eq:ssvol}) yields applicable integral formulae, when $\mathfrak{g}$ is a reductive Lie algebra ($\mathfrak{g} = \mathfrak{z} + \mathfrak{g}_{ss}$\,: $\mathfrak{z}$ the centre of $\mathfrak{g}$ and $\mathfrak{g}_{ss}$ semisimple). If $\mathfrak{a}$ is a maximal Abelian subspace of $\mathfrak{p}$\footnote{Recall that the dimension of $\mathfrak{a}$ is known as the rank of $M$. In fact, $\mathrm{Exp}_o(\mathfrak{a})$ is a totally flat submanifold of $M$,\hfill\linebreak of maximal dimension, and the only such submanifold, up to isometry.}, any $v \in \mathfrak{p}$ is of the form $v = \mathrm{Ad}(k)\,a$ for some $k \in K$ and $a \in\mathfrak{a}$ (see~\cite{helgason}, Lemma 6.3, Chapter V).\hfill\linebreak Moreover, using the fact that $\mathrm{Ad}(k)$ is an isomorphism of $\mathfrak{g}$, 
\begin{equation} \label{eq:raeigen}
\mathrm{Ad}(k^{\scriptscriptstyle -1})\circ R_v \circ \mathrm{Ad}(k) = R_a = \sum_{\lambda \in \Delta_+} (\lambda(a))^2\;\Pi_{\lambda}
\end{equation}
where each $\lambda \in \Delta_+$ is a linear form $\lambda : \mathfrak{a} \rightarrow \mathbb{R}$, and $\Pi_{\lambda}$ is the orthogonal projectors onto the corresponding eigenspace of $R_{a\,}$. Here, $\Delta_+$ is the set of positive roots of $\mathfrak{g}$ with respect to $\mathfrak{a}$~\cite{helgason} (see Lemma 2.9, Chapter VII).

It is possible to use the diagonalisation (\ref{eq:raeigen}), in order to evaluate the determinant (\ref{eq:ssvol}). To obtain a regular parameterisation, let $S = K/K_{\mathfrak{a\,}}$, where $K_{\mathfrak{a}}$ is the centraliser of $\mathfrak{a}$ in $K$. Then, let $\varphi : S \times \mathfrak{a} \rightarrow M$ be given by $\varphi(s\hspace{0.02cm},a) = \mathrm{Exp}_o(\beta(s\hspace{0.02cm},a))$ where $\beta(s,a) = \mathrm{Ad}(s)\,a$. Now, by (\ref{eq:ssvol}) and (\ref{eq:raeigen}),
$$
\varphi^*(\mathrm{vol}) = \prod_{\lambda \in \Delta_+} \left| \frac{\sinh\hspace{0.02cm}\lambda(a)}{\lambda(a)}\right|^{m_\lambda}\hspace{0.03cm}\beta^*(dv)
$$
where $m_\lambda$ is the multiplicity of $\lambda$ (the rank of $\Pi_\lambda$). On the other hand, one may show that
\begin{equation} \label{eq:beta*}
\beta^*(dv) = \prod_{\lambda \in \Delta_+} |\lambda(a)|^{ m_\lambda}\hspace{0.05cm}da\,\omega(ds)
\end{equation}
where $da$ is the volume form on $\mathfrak{a}$, and $\omega$ is the invariant volume induced onto $S$ from $K$.

Finally, the Riemannian volume, in terms of the parameterisation $\varphi$, takes on the form
\begin{equation} \label{eq:ssvolka}
\varphi^*(\mathrm{vol}) = \prod_{\lambda \in \Delta_+} \left| \sinh\hspace{0.02cm}\lambda(a)\right|^{ m_\lambda}\hspace{0.03cm}da\,\omega(ds)
\end{equation}
 Using (\ref{eq:ssvolka}), it will be possible to write down integral formulae for Riemannian symmetric spaces, either non-compact or compact.

\subsubsection{The non-compact case}
This is the case were $\mathfrak{g}$ admits an $\mathrm{Ad}(G)$-invariant, non-degenerate, symmetric bilinear form $B$, such that $Q(u\hspace{0.02cm},z) = - B(u,d\hspace{0.02cm}\theta(z))$ is an $\mathrm{Ad}(K)$-invariant scalar product on $\mathfrak{g}$. In this case, $B$ is 
negative-definite on $\mathfrak{k}$ and positive-definite on $\mathfrak{p}$. Moreover, $\mathrm{ad}(z) = [z,\cdot]$ is skew-symmmetric or symmetric (with respect to $Q$), according to whether $z \in \mathfrak{k}$ or $z \in \mathfrak{p}$.

If $u_{\scriptscriptstyle 1\hspace{0.02cm}},u_{\scriptscriptstyle 2} \in \mathfrak{p}$ are orthonormal, the sectional curvature of $\mathrm{Span}(u_{\scriptscriptstyle 1\hspace{0.02cm}},u_{\scriptscriptstyle 2})$ is found from (\ref{eq:sscurvature}), $\kappa(u_{\scriptscriptstyle 1\hspace{0.02cm}},u_{\scriptscriptstyle 2}) = - \Vert [u_{\scriptscriptstyle 1\hspace{0.02cm}},u_{\scriptscriptstyle 2}] \Vert^2_o \leq 0$. Therefore, $M$ has  non-positive sectional curvature.

In fact, $M$ is a Hadamard manifold. It is geodesically complete by (\ref{eq:ssexp1}). It is moreover simply connected, because $\mathrm{Exp}_o:\mathfrak{p} \rightarrow M$ is a diffeomorphism~\cite{helgason} (Theorem 1.1, Chapter VI). Thus, (\ref{eq:ssvol}) yields a first integral formula,
\begin{equation} \label{eq:ssvolnc}
  \int_M\,f(x)\,\mathrm{vol}(dx) = \int_{\mathfrak{p}}\, f(\mathrm{Exp}_o(v))\hspace{0.02cm}\left|\det(\mathrm{sh}(R_v))\right|\hspace{0.02cm}dv
\end{equation}

To obtain an integral formula from (\ref{eq:ssvolka}), one should first note that $\beta : S \times \mathfrak{a} \rightarrow \mathfrak{p}$ is not regular, nor one-to-one. Recall the following\,:

\noindent $\bullet$ the hyperplanes $\lambda(a) = 0$, where $\lambda \in \Delta_+\,$, divide $\mathfrak{a}$ into finitely many connected components, which are open and convex sets, known as Weyl chambers. From (\ref{eq:beta*}), $\beta$ is regular on each Weyl chamber. 

\noindent $\bullet$ let $K^\prime_{\mathfrak{a}}$ denote the normaliser of $\mathfrak{a}$ in $K$. Then, $W = \left.K^\prime_{\mathfrak{a}}\middle/K_{\mathfrak{a}}\right.$ is a finite group of automorphisms of $\mathfrak{a}$, called the Weyl group, which acts freely transitively on the set of Weyl chambers~\cite{helgason} (Theorem 2.12, Chapter VII). 

Then, for each Weyl chamber $C$, $\beta$ is regular and one-to-one, from $S\times C$ onto its image in $\mathfrak{p}$. \hfill\linebreak Moreover, if $\mathfrak{a}_r$ is the union of the Weyl chambers ($a \in \mathfrak{a}_r$ if and only if $\lambda(a) \neq 0$ for any $\lambda \in \Delta_+$), then $\beta$ is regular and $|W|$-to-one from $S \times \mathfrak{a}_r$ onto its image in $\mathfrak{p}$. 

To obtain the desired integral formula, it only remains to note that $\varphi$ is a diffeomorphism from $S\times C$ onto its image in $M$. However, this image is the set $M_r$ of regular values of $\varphi$. By Sard's lemma, its complement is negligible~\cite{bogachev}.  
\begin{proposition} \label{prop:ssvolncka}
  Let $M = G/K$ be a Riemannian symmetric space, which belongs to the ``non-compact case", just described. Then, for any bounded continuous function $f:M\rightarrow \mathbb{R}$,
\begin{equation} \label{eq:ssvolncka1}
\int_M\,f(x)\,\mathrm{vol}(dx) \,=\,
\int_{C_+}\int_{S}f(\varphi(s,a))\hspace{0.03cm}\prod_{\lambda \in \Delta_+}\left( \sinh\hspace{0.02cm}\lambda(a)\right)^{ m_\lambda}\hspace{0.03cm}da\,\omega(ds)
\end{equation}
\begin{equation} \label{eq:ssvolncka2}
\phantom{\int_M\,f(x)\,\mathrm{vol}(dx)} \,=\,
\frac{1}{|W|}\,\int_{\mathfrak{a}}\int_{S}\,f(\varphi(s,a))\hspace{0.03cm}\prod_{\lambda \in \Delta_+}\left| \sinh\hspace{0.02cm}\lambda(a)\right|^{ m_\lambda}\hspace{0.03cm}da\,\omega(ds)
\end{equation}
Here, $C_+$ is the Weyl chamber $C_+ = \lbrace a \in \mathfrak{a}\,: \lambda \in \Delta_+ \Rightarrow \lambda(a) > 0\rbrace$.  
\end{proposition}
\noindent \textbf{Example 1\,:} consider $M = \mathrm{H}(N)$ the space of $N \times N$ Hermitian positive-definite matrices. Here, $G = \mathrm{GL}(N,\mathbb{C})$ and $K = U(N)$, the groups of $N \times N$, complex, invertible and unitary matrices. Moreover, $B(u,\!z) = \mathrm{Re}(\mathrm{tr}(uz))$ and $d\hspace{0.02cm}\theta(z) = -z^\dagger$. Thus, $\mathfrak{p}$ is the space of $N \times N$ Hermitian matrices, and one may choose $\mathfrak{a}$ the space of $N \times N$ real diagonal matrices. The positive roots are the linear maps $\lambda(a) = a_{ii} - a_{jj}$ where $i < j$, and each one has its multiplicity equal to $2$. Thus, $C_+$ is the cone of real diagonal matrices $a$ with $a_{\scriptscriptstyle 11} > \ldots > a_{\scriptscriptstyle NN} > 0$. The Weyl group $W$ is the group of permutation matrices in $U(N)$ (so $|W| = N!$). Finally, $S= U(N)/T_{\scriptscriptstyle N} \equiv S_{\scriptscriptstyle N\,}$, where $T_{\scriptscriptstyle N}$ is the torus of diagonal unitary matrices. By (\ref{eq:ssvolncka2}), 
\begin{equation} \label{eq:ssvolncka2hn}
\int_{\mathrm{H}(N)}\,f(x)\,\mathrm{vol}(dx) \,=\,
\frac{1}{N!}\,\int_{\mathfrak{a}}\int_{S_{\scriptscriptstyle N}}\,f\left(s\hspace{0.02cm}\exp(2a)\hspace{0.02cm}s^\dagger\right)\hspace{0.03cm}\prod_{i < j} \sinh^2(a_{ii} - a_{jj})\hspace{0.03cm}da\,\omega(ds)
\end{equation}
where $da = da_{\scriptscriptstyle 11}\ldots da_{\scriptscriptstyle NN\,}$. \vfill\pagebreak
\noindent \textbf{Example 2\,:} pursuing the previous example, assume $f$ is a class function\,: $f(k\cdot x) = f(x)$ for $k \in K$ and $x \in \mathrm{H}(N)$. That is, $f(x)$ depends only on the eigenvalues $x_i = e^{r_i}$ of $x$. By (\ref{eq:ssvolncka2hn}),
\begin{equation} \label{eq:integralhn}
\int_{\mathrm{H}(N)}\,f(x)\,\mathrm{vol}(dx) \,=\,
\frac{\omega(S_{\scriptscriptstyle N})}{2^{\scriptscriptstyle N}N!}\, \int_{\mathbb{R}^{\scriptscriptstyle N}}\,f\left(\exp(r)\right)\hspace{0.03cm}\prod_{i < j}\sinh^2((r_i - r_j)/{2})\hspace{0.03cm}dr
\end{equation}
or, by introducing the eigenvalues $x_i$ as integration variables,
\begin{equation} \label{eq:integralhnvmonde}
\int_{\mathrm{H}(N)}\,f(x)\,\mathrm{vol}(dx) \,=\,
\frac{\omega(S_{\scriptscriptstyle N})}{2^{\scriptscriptstyle N^2}N!} \,\int_{\mathbb{R}^{\scriptscriptstyle N}_+}\,f\left(x_{\scriptscriptstyle 1\,},\ldots,x_{\scriptscriptstyle N}\right)\hspace{0.03cm}|V(x)|^2\hspace{0.03cm}\prod^N_{i=1} x^{\scriptscriptstyle -N}_i\hspace{0.02cm}dx_i
\end{equation}
where $V(x) = \prod_{i<j} (x_j - x_i)$ is the Vandermonde determinant. Integrals of this form are well-known in random matrix theory~\cite{mehta}. %\\[0.1cm]
%\textbf{Remark\,:} here is a proof of the fact that $\mathrm{Exp}_o:\mathfrak{p} \rightarrow M$ is a diffeomorphism. It is enough to consider the case where $\mathfrak{g}$ is semisimple. %the cone !!

\subsubsection{The compact case} 
In this case, $\mathfrak{g}$ admits an $\mathrm{Ad}(G)$-invariant scalar product $Q$. Therefore, $\mathrm{ad}(z)$ is skew-symmmetric, with respect to $Q$, for each $z \in \mathfrak{g}$. Using (\ref{eq:sscurvature}), it follows that $M$ is compact, with non-negative sectional curvature.

In fact, the compact case may be obtained from the previous non-compact case by \textit{duality}. Denote $\mathfrak{g}_{\scriptscriptstyle\, \mathbb{C}}$ the complexification of $\mathfrak{g}$, and let $\mathfrak{g}^* = \mathfrak{k} + \mathfrak{p}_*$ where $\mathfrak{p}_* = i\hspace{0.02cm}\mathfrak{p}$. Then, $\mathfrak{g}^*$ is a compact real form of $\mathfrak{g}_{\scriptscriptstyle\, \mathbb{C}}$ (that is, $\mathfrak{g}^*$ is a compact Lie algebra, and its complexification is equal to $\mathfrak{g}_{\scriptscriptstyle\, \mathbb{C}}$). Denote $G^*$ the connected Lie group with Lie algebra $\mathfrak{g}^*$. %Of course, $G^*$ is compact, and $M^*$ is compact and has non-negative sectional curvature.  

If $M = G/K$ is a Riemannian symmetric space which belongs to the non-compact case, then $M^* = G^*\!/K$ is a Riemannian symmetric space which belongs to the compact case. Formally, to pass from the non-compact case to the compact case, all one has to do is replace $a$ by $i\hspace{0.02cm}a$. Applying this recipe to (\ref{eq:ssvolka}), one obtains
\begin{equation} \label{eq:ssvolkac}
\varphi^*(\mathrm{vol}) = \prod_{\lambda \in \Delta_+} \left| \sin\hspace{0.02cm}\lambda(a)\right|^{ m_\lambda}\hspace{0.03cm}da\,\omega(ds)
\end{equation}
where $da$ is the volume form on $\mathfrak{a}_* = i\hspace{0.02cm}\mathfrak{a}$, and $\omega$ is the invariant volume induced onto $S$ from $K$.
Note that the image under $\mathrm{Exp}_o$ of $\mathfrak{a}_*$ is the torus $T_* = \mathfrak{a}_*/\mathfrak{a}_{\scriptscriptstyle K\,}$,  where $\mathfrak{a}_{\scriptscriptstyle K}$ is the lattice given by $\mathfrak{a}_{\scriptscriptstyle K} = \lbrace a \in \mathfrak{a}_*: \mathrm{Exp}_o(a) = o \rbrace$. Recall the following\,:

\noindent $\bullet$ $\varphi(s,a)$ only depends on $t = \mathrm{Exp}_o(a)$. Thus, $\varphi$ may be considered as a map from $S \times T_*$ to $M$.

\noindent $\bullet$ if $a \in \mathfrak{a}_{\scriptscriptstyle K}$ then $\exp(2a) = e$ (the identity element in $G^*$). Thus, $\lambda(a) \in i\hspace{0.02cm}\pi\,\mathbb{Z}$ for all $\lambda \in \Delta_+$~\cite{helgason} (Page 383).  Therefore, there exists a function $D:T \rightarrow \mathbb{R}$, such that
$$
D(t) = \prod_{\lambda \in \Delta_+} \left| \sin\hspace{0.02cm}\lambda(a)\right|^{ m_\lambda} \hspace{0.5cm}
\text{whenever $t = \mathrm{Exp}_o(a)$}
$$
Now, $T_*$ is a totally flat submanifold of $M$. Therefore, $\mathrm{Exp}^*(dt) = da$, where $dt$ denotes the invariant volume induced onto $T_*$ from $M$. With a slight abuse of notation, (\ref{eq:ssvolkac}) now reads,
\begin{equation} \label{eq:ssvolkac1}
\varphi^*(\mathrm{vol}) = D(t)\hspace{0.03cm}dt\,\omega(ds)
\end{equation}
Denote $(T_*)_r$ the set of $t \in T_*$ such that $D(t) \neq 0$. By the same arguments as in the non-compact case, $\varphi$ is a regular $|W|$-to-one map from $S \times (T_*)_r$ onto $M_r\,$, the set of regular values of $\varphi$. 
\begin{proposition} \label{prop:ssvolkac}
  Let $M = G^*\!/K$ be a Riemannian symmetric space, which belongs to the ``compact case", just described. Then, for any bounded continuous function $f:M\rightarrow \mathbb{R}$,
\begin{equation} \label{eq:ssvolkac2}
\int_M\,f(x)\,\mathrm{vol}(dx) \,=\,
\frac{1}{|W|}\,\int_{T_*}\int_{S}\,f(\varphi(t,a))\hspace{0.03cm}D(t)\hspace{0.03cm}dt\,\omega(ds)
\end{equation}
\end{proposition}
\noindent \textbf{Example 1\,:} the dual of $\mathrm{H}(N)$ is the unitary group $U(N)$. Here, $G^* = U(N)\times U(N)$ and $K \simeq U(N)$, is the diagonal group $K = \lbrace (x\hspace{0.02cm},x)\,;x \in U(N)\rbrace$. The Riemannian metric is given by the trace scalar product $Q(u,\!z) = -\mathrm{tr}(uz)$. Moreover, $T_* = T_{\scriptscriptstyle N}$ and $S= S_{\scriptscriptstyle N}$ (this is $U(N)/T_{\scriptscriptstyle N}$). The positive roots are $\lambda(ia) = a_{ii} - a_{jj}$ where $i < j$ and where $a$ is $N \times N$, real and diagonal\footnote{Please do no confuse the imaginary number $i$ with the subscript $i$.}. By writing the integral over $T_{\scriptscriptstyle N}$ as a multiple integral, (\ref{eq:ssvolkac2}) reads,
\begin{equation} \label{eq:ssvolkacun}
\int_{U(N)}\,f(x)\,\mathrm{vol}(dx) \,=\,
\frac{1}{N!}\,\int_{[0\hspace{0.02cm},2\pi]^N}\int_{S_{\scriptscriptstyle N}}\,f\left(s\hspace{0.02cm}\exp(2ia)\hspace{0.02cm}s^\dagger\right)\hspace{0.03cm}\prod_{i < j} \sin^2(a_{ii} - a_{jj})\hspace{0.03cm}\omega(ds)\,da
\end{equation}
where $da = da_{\scriptscriptstyle 11}\ldots da_{\scriptscriptstyle NN\,}$. \\[0.1cm]
\textbf{Example 2\,:} assume $f$ is a class function.  That is, $f(x)$ depends only on eigenvalues $e^{i\theta_i}$ of $x$. Integrating out $s$, from (\ref{eq:ssvolkacun}), it follows,
\begin{equation} \label{eq:integralun}
\int_{U(N)}\,f(x)\,\mathrm{vol}(dx) \,=\,
\frac{\omega(S_{\scriptscriptstyle N})}{2^{\scriptscriptstyle N}N!}\, \int_{[0\hspace{0.02cm},2\pi]^N}\,f\left(\exp(i\theta)\right)\hspace{0.03cm}\prod_{i < j}\sin^2((\theta_i - \theta_j)/{2})\hspace{0.03cm}d\theta
\end{equation}
or, after an elementary manipulation, 
\begin{equation} \label{eq:integralunvmonde}
\int_{U(N)}\,f(x)\,\mathrm{vol}(dx) \,=\,
\frac{\omega(S_{\scriptscriptstyle N})}{2^{\scriptscriptstyle N^2}N!} \,\int_{[0\hspace{0.02cm},2\pi]^N}\,f\left(\theta_{\scriptscriptstyle 1\,},\ldots,\theta_{\scriptscriptstyle N}\right)\hspace{0.03cm}|V(e^{i\theta})|^2\hspace{0.03cm}d\theta_{\scriptscriptstyle 1}\ldots\theta_{\scriptscriptstyle N}
\end{equation}
where $V(e^{i\theta}) = \prod_{i<j} (e^{i\theta_j} - e^{i\theta_i})$ is the Vandermonde determinant. Integrals of this form are well-known in the random matrix theory of compact groups~\cite{meckes}. %\\[0.1cm]

%give the metrics
%rank

\section{Geodesics in symmetric spaces} \label{sec:geolemma}
Let $M = G/K$ be a Riemannian symmetric space, and assume $G$ has reductive Lie algebra $\mathfrak{g}$.\hfill\linebreak Let the metric of $M$ be given by an $\mathrm{Ad}(K)$-invariant scalar product $Q$ on $\mathfrak{g}$, according to (\ref{eq:ssmetric}). 
Recall the Cartan decomposition, $\mathfrak{g} = \mathfrak{k} + \mathfrak{p}$. Assume $\mathfrak{k}$ and $\mathfrak{p}$ are orthogonal, with respect to $Q$, and extend $Q$ to a left-invariant Riemannian metric $(\cdot,\cdot)$ on $G$.

Then, consider the natural projection $\pi:G \rightarrow M$,  which is given by $\pi(g) = g\cdot o$ for $g \in G$. Denote by $V_g$ the kernel of $d\pi(g)$, and by $H_g$ its orthogonal complement with respect to $(\cdot,\cdot)_{\scriptscriptstyle g\,}$.  Since $(\cdot,\cdot)$ is left-invariant, and $\mathfrak{k}$ and $\mathfrak{p}$ are orthogonal,
\begin{equation} \label{eq:kptovh}
  V_g = dL_g(\mathfrak{k}) \hspace{1cm}   H_g = dL_g(\mathfrak{p})
\end{equation}
where $L_g$ denotes left translation by $g$ (below, $R_g$ will denote right translation). 

For $v \in T_xM$, there exists a unique $v^{\scriptscriptstyle H}(g) \in H_g$ such that $d\pi\left( v^{\scriptscriptstyle H}(g)\right) = v$. By an abuse of notation, denote $v^{\scriptscriptstyle H}(e) = dL_{\scriptscriptstyle g^{-1}}\left( v^{\scriptscriptstyle H}(g)\right)$. Recall that for any $u\hspace{0.03cm},v \in T_xM$ (see~\cite{nomizu2}, Chapter X),
\begin{equation} \label{eq:vhtometric}
  \langle u\hspace{0.02cm},v\rangle_{\scriptscriptstyle x} \,=\, \left(u^{\scriptscriptstyle H}(g)\hspace{0.02cm},v^{\scriptscriptstyle H}(g)\right)_{\scriptscriptstyle g} \,=\, Q\left(u^{\scriptscriptstyle H}(e)\hspace{0.02cm},v^{\scriptscriptstyle H}(e)\right)
\end{equation}
For Propositions \ref{prop:geometriclemma} and \ref{prop:geodesiclemma}, consider the ``infinitesimal action",
\begin{equation} \label{eq:infaction}
\xi\cdot x = \left.\frac{d}{dt}\right|_{t=0} \exp(t\hspace{0.02cm}\xi)\cdot x \hspace{1cm} \xi \in \mathfrak{g}\text{ and } x \in M
\end{equation}
In addition, let $\mathrm{B}$ be the bilinear form on $\mathfrak{g}$, given by $\mathrm{B} = B$, if $M$ belongs to the non-compact case, and by $\mathrm{B} = Q$, if $M$ belongs to the compact case (these two cases were described in \ref{ssec:sspace}).
\begin{proposition} \label{prop:geometriclemma}
Let $M = G/K$ be a symmetric space, of the ``non-compact case", or of the ``compact case".  For $g \in G$, $x = \pi(g)$ and $v \in T_xM$, let $\omega_{\scriptscriptstyle v} = \mathrm{Ad}(g)\hspace{0.02cm}v^{\scriptscriptstyle H}(e)$. Then, $\omega_{\scriptscriptstyle v}\cdot x = v$ and $\langle u\hspace{0.02cm},v\rangle_{\scriptscriptstyle x} = \mathrm{B}(\xi\hspace{0.02cm},\omega_{\scriptscriptstyle v})$ whenever $u = \xi \cdot x$. 
\end{proposition}
\begin{proposition} \label{prop:geodesiclemma}
 In the notation of the previous proposition,
\begin{equation} \label{eq:geodesiclift}
\mathrm{Exp}_x(v) = \exp(\omega_{\scriptscriptstyle v})\cdot x
\end{equation}
for $x \in M$ and $v \in T_xM$.
\end{proposition}
Propositions \ref{prop:geometriclemma} and \ref{prop:geodesiclemma} offer a straightforward computational route to the Riemannian exponential map $\mathrm{Exp}$. To compute $\mathrm{Exp}_x(v)$, one begins by ``lifting" $v$ from $T_xM$ to $\mathfrak{g}$, under the form of $\omega_{\scriptscriptstyle v\,}$. Then, it is enough to compute the action of $\exp(\omega_{\scriptscriptstyle v})$, which is just a matrix exponential, in practice. \\[0.1cm]
\textbf{Example 1\,:} consider an example of the non-compact case, $M = \mathrm{H}(N)$, the space of $N \times N$ Hermitian positive-definite matrices. Here, $G = \mathrm{GL}(N,\mathbb{C})$ and $\pi(g) = gg^\dagger$ for $g \in G$. Then, 
$$
d\pi(g)\cdot h = h\hspace{0.02cm}g^\dagger + g\hspace{0.02cm}h^\dagger
$$
for $h \in T_gG$. For $x = \pi(g)$ and $v \in T_xM$, it follows that $v^{\scriptscriptstyle H}(g) = \left.v\hspace{0.03cm}\theta(g)\middle/2\right.$, where $\theta(g) = (g^\dagger)^{-1}$. By definition, $\omega_{\scriptscriptstyle v} = dR_{\scriptscriptstyle g^{-1}}(v^{\scriptscriptstyle H}(g))$. Since $gg^\dagger = x$, this gives $\omega_{\scriptscriptstyle v} = (v/2)\hspace{0.03cm}x^{-1}$. Therefore, using the fact that $g\cdot x = g\hspace{0.02cm}x\hspace{0.02cm}g^\dagger$, it follows
$$
\mathrm{Exp}_x(v) = \exp\left(v\hspace{0.03cm}x^{-1}\!\middle/2\right) \,x \,\exp\left(x^{-1}\hspace{0.03cm}v\middle/2\right) 
$$
Accordingly, by an elementary property of the matrix exponential\footnote{Matrix functions (powers, logarithms, \textit{etc.}) of Hermitian arguments should be understood as Hermitian matrix functions, obtained using the spectral decomposition --- see~\cite{higham}.}, 
\begin{equation} \label{eq:pennec}
\mathrm{Exp}_x(v) = x^{\frac{1}{2}}\exp\left(x^{-\frac{1}{2}}\hspace{0.04cm}v\hspace{0.04cm}x^{-\frac{1}{2}}\right)x^{\frac{1}{2}}
\end{equation}
which is the formula made popular by~\cite{pennec2006}. \\[0.2cm]
\textbf{Example 2\,:} let $M=G/K$ be a Riemannian symmetric space of the compact case. That is, the scalar product $Q$ on $\mathfrak{g}$ is $\mathrm{Ad}(G)$-invariant. Write $\mathfrak{g} = \mathfrak{k} + \mathfrak{p}$ the Cartan decomposition of $\mathfrak{g}$. For $x \in M$, denote $K_x$ the stabiliser of $x$ in $G$. If $x = \pi(g)$, this has Lie algebra $\mathfrak{k}_x =\mathrm{Ad}(g)(\mathfrak{k})$ (that is, the image under $\mathrm{Ad}(g)$ of $\mathfrak{k}$).  For $v \in T_xM$, by Proposition \ref{prop:geometriclemma}, its ``lift" $\omega_{\scriptscriptstyle v}$ should verify (note that, for the present example, $\mathrm{B} = Q$)
$$
\omega_{\scriptscriptstyle v}\cdot x = v \hspace{0.3cm}\text{and}\hspace{0.3cm} Q(\xi\hspace{0.02cm},\omega_{\scriptscriptstyle v}) = 0 \text{ for } \xi \in 
\mathfrak{k}_x
$$
where the second identity is because $\xi \cdot x = 0$ for $\xi \in \mathfrak{k}_x\,$. Because $Q$ is $\mathrm{Ad}(G)$-invariant, this second identity is equivalent to 
$Q(\kappa\hspace{0.02cm},\mathrm{Ad}(g^{\scriptscriptstyle -1})(\omega_{\scriptscriptstyle v})) = 0$ for $\kappa \in \mathfrak{k}$. That is, $\omega_{\scriptscriptstyle v} = \mathrm{Ad}(g)(\omega_{\scriptscriptstyle v}(o))$ for some $\omega_{\scriptscriptstyle v}(o) \in \mathfrak{p}$. This $\omega_{\scriptscriptstyle v}(o)$ is determined from $\omega_{\scriptscriptstyle v}\cdot x = v$, which yields $\omega_{\scriptscriptstyle v}(o) \cdot o = g^{\scriptscriptstyle -1}\cdot v$. However, the map $\omega \mapsto \omega \cdot o$ is an isomorphism from $\mathfrak{p}$ onto $T_oM$. Denoting its inverse by $\pi_o : T_oM \rightarrow \mathfrak{p}$, it follows that $\omega_{\scriptscriptstyle v}(o) = \pi_o(g^{\scriptscriptstyle -1}\cdot v)$. Finally,
\begin{equation} \label{eq:compactlift}
\omega_{\scriptscriptstyle v} = \mathrm{Ad}(g)\left( \pi_o(g^{\scriptscriptstyle -1}\cdot v)\right)
\end{equation}
A special case of this formula was used in (\ref{eq:grasslift}) of \ref{sec:pca}. \\[0.1cm]
\textbf{Proof of Proposition \ref{prop:geometriclemma}\,:} to begin, one must prove
\begin{equation} \label{eq:proofgeolemma1}
  \omega_{\scriptscriptstyle v}\cdot x = v
\end{equation}
From the definition of $\omega_{\scriptscriptstyle v}$ and $v^{\scriptscriptstyle H}(e)$, it is clear $\omega_{\scriptscriptstyle v} = dR_{g^{\scriptscriptstyle -1}}\hspace{0.02cm}(v^{\scriptscriptstyle H}(g))$. Replacing this into (\ref{eq:infaction}), the left-hand side of (\ref{eq:proofgeolemma1}) becomes, 
$$
  \omega_{\scriptscriptstyle v}\cdot x = \left.\frac{d}{dt}\right|_{t=0}\exp(t\,dR_{g^{\scriptscriptstyle -1}}\hspace{0.03cm}v^{\scriptscriptstyle H}(g))\cdot x = \left.\frac{d}{dt}\right|_{t=0} \left(\gamma(t)\hspace{0.04cm}g^{\scriptscriptstyle-1}\right)\cdot x
$$
where $\gamma$ is any curve in $G$, through $g$ with $\dot{\gamma}(0) = v^{\scriptscriptstyle H}(g)$. Therefore,
$$
  \omega_{\scriptscriptstyle v}\cdot x = \left.\frac{d}{dt}\right|_{t=0} \gamma(t)\cdot o = d\pi\left( v^{\scriptscriptstyle H}(g)\right) = v
$$
from the definition of $v^{\scriptscriptstyle H}(g)$. This proves (\ref{eq:proofgeolemma1}). It remains to show,
\begin{equation} \label{eq:proofgeolemma2}
 \langle u\hspace{0.02cm},v\rangle_{\scriptscriptstyle x} = Q(\xi\hspace{0.02cm},\omega_{\scriptscriptstyle v}) \hspace{1cm}
\text{for } u = \xi \cdot x
\end{equation}
The proof is separated into two cases. \\[0.1cm]
\textbf{non-compact case\,:} in this case, $Q(\xi\hspace{0.02cm},\omega) = - B(\xi,d\hspace{0.02cm}\theta(\omega))$, where $B$ is an $\mathrm{Ad}(G)$-invariant, non-degenerate, symmetric bilinear form. To prove (\ref{eq:proofgeolemma2}), note that
$$
d\pi(g)\left( dR_{\scriptscriptstyle g}(\xi)\right) = \left.\frac{d}{dt}\right|_{t=0} \left(\exp(t\xi)\,g\right)\cdot o = \left.\frac{d}{dt}\right|_{t=0} \exp(t\xi) \cdot x = u
$$
Therefore, $dR_{\scriptscriptstyle g}(\xi) = u^{\scriptscriptstyle H}(g) + w$ where $w \in V_g$. From (\ref{eq:vhtometric}), using left-invariance of $(\cdot,\cdot)$,
\begin{equation} \label{eq:proofgeolemma21}
  \langle u\hspace{0.02cm},v\rangle_{\scriptscriptstyle x} \,=\, \left(dR_{\scriptscriptstyle g}(\xi)\hspace{0.02cm},v^{\scriptscriptstyle H}(g)\right)_{\scriptscriptstyle g} \,=\, Q\left(\mathrm{Ad}(g^{\scriptscriptstyle -1})(\xi)\hspace{0.02cm},v^{\scriptscriptstyle H}(e)\right)
\end{equation}
Thus, using the definition of $Q$, and the fact that $v^{\scriptscriptstyle H}(e) \in \mathfrak{p}$,
$$
  \langle u\hspace{0.02cm},v\rangle_{\scriptscriptstyle x} = - B(\mathrm{Ad}(g^{\scriptscriptstyle -1})(\xi),d\hspace{0.02cm}\theta(v^{\scriptscriptstyle H}(e))) =
B(\mathrm{Ad}(g^{\scriptscriptstyle -1})(\xi),v^{\scriptscriptstyle H}(e))
$$
Finally, since $B$ is $\mathrm{Ad}(G)$-invariant, 
$$
  \langle u\hspace{0.02cm},v\rangle_{\scriptscriptstyle x} = B(\mathrm{Ad}(g^{\scriptscriptstyle -1})(\xi),v^{\scriptscriptstyle H}(e)) = 
B(\xi,\mathrm{Ad}(g)(v^{\scriptscriptstyle H}(e)))
$$
which is the same as (\ref{eq:proofgeolemma2}), by the definition of  $\omega_{\scriptscriptstyle v\,}$. Indeed, in the present case, $\mathrm{B} = B$.\\[0.1cm]
\textbf{compact case\,:} this follows from (\ref{eq:proofgeolemma21}), using the fact that $Q$ is $\mathrm{Ad}(G)$-invariant. Indeed, in the present case, $\mathrm{B} = Q$. \\[0.1cm]
\textbf{Proof of Proposition \ref{prop:geodesiclemma}\,:} for $\xi \in \mathfrak{g}$, introduce the corresponding vector fields $X_\xi$ on $M$, given by $
X_\xi(x) = \xi\cdot x$. Since this is a Killing vector field~\cite{nomizu2}, if $c:\mathbb{R} \rightarrow M$ is a geodesic curve in $M$, then $\ell(\xi) = \langle X_\xi\hspace{0.03cm},\dot{c}\rangle_{\scriptscriptstyle c(t)}$ is a constant, (a law of conservation, really due to Noether's theorem!). Now, in the notation of Proposition \ref{prop:geometriclemma}, let $\omega(t) = \omega_{\scriptscriptstyle \dot{c}(t)\,}$. By Proposition \ref{prop:geometriclemma}, 
$$
\mathrm{B}(\omega(t),\xi) = \ell(\xi)
$$
Since this is a constant, and since $\mathrm{B}$ is non-degenerate, it follows that $\omega(t) = \omega$ is a constant. Proposition \ref{prop:geometriclemma} also implies that $c$ satisfies the ordinary differential equation
$$
\dot{c} = \omega \cdot c
$$
But this differential equation is also satisfied  by $c(t) = \exp(t\hspace{0.03cm}\omega)\cdot\hspace{0.02cm} c(0)$, as one may see from (\ref{eq:infaction}).\hfill\linebreak By uniqueness of the solution, for given initial conditions, 
$$
c(t) = \exp(t\hspace{0.03cm}\omega_{\scriptscriptstyle \dot{c}(0)})\cdot c(0)
$$
This immediately implies (\ref{eq:geodesiclift}), by setting $t = 1$, $c(0) = x$ and $\dot{c}(0) = v$.

\chapter{The barycentre problem} \label{barycentre}

\minitoc
\vspace{0.1cm}

{\small State-of-the art results establish the existence and uniqueness of the Riemannian barycentre of a probability distribution which is supported inside a compact convex geodesic ball. What happens for a probability distribution which is not supported, but concentrated, inside a convex geodesic ball\,?\hfill\linebreak This question raises new difficulties that cannot be resolved by using the tools applicable to distributions which have compact convex support. The present chapter develops new tools, able to deal with these difficulties (at least in part), following the approach in~\cite{salemgibbs}.
\begin{itemize}
  \item \ref{sec:frechet} and \ref{sec:barycentrebackground} review some of the major contributions to the study of Riemannian barycentres, due to Fréchet, Emery, Kendall, Afsari, and Arnaudon. 
 \item \ref{sec:gibbs} introduces the main problem treated in the following, which is to study the existence and uniqueness of Riemannian barycentres of Gibbs distributions on compact Riemannian symmetric spaces.
 \item \ref{sec:concentration} -- \ref{sec:unique} lead up to the following conclusion\,: let $\pi_{\scriptscriptstyle T} \propto \exp(-U/T)$ be a Gibbs distribution on a simply connected compact Riemannian symmetric space $M$, such that the potential function $U$ has a unique global minimum at $x^* \in M$. If $M$ is simply connected, then for each $\delta < r_{\scriptscriptstyle cx}/2$ (where $r_{\scriptscriptstyle cx}$ is the convexity radius of $M$), there exists a critical temperature $T_{\scriptscriptstyle \delta}$ such that $T < T_{\scriptscriptstyle \delta}$ implies that $\pi_{\scriptscriptstyle T}$ has a unique Riemannian barycentre $\hat{x}_{\scriptscriptstyle T}$ and this $\hat{x}_{\scriptscriptstyle T}$ belongs to the geodesic ball $B(x^*\!,\delta)$. The assumption that $M$ is simply connected cannot be removed (see Lemma \ref{lem:hddim} and the following remark). 
\item \ref{sec:twtd} provides expressions which can be used to analytically compute the critical temperature $T_{\scriptscriptstyle \delta\hspace{0.03cm}}$.
\item \ref{sec:ssbis} introduces additional background material, on the geometry of compact symmetric spaces, which is required for the proofs of the results in \ref{sec:concentration} -- \ref{sec:unique}.
\item \ref{sec:bcentreproofs} details the proofs of these results, concerning the concentration, differentiability, convexity,\hfill\linebreak existence and uniqueness of Riemannian barycentres of Gibbs distributions on compact symmetric spaces.
\end{itemize}
}
\vfill
\pagebreak

%repetition of paper, with many additional typing mistakes

\section{Fréchet's fruitful idea} \label{sec:frechet}
In 1948, Maurice Fréchet proposed a generalisation of the concept of mean value, from Euclidean spaces to general metric spaces~\cite{frechet}. 
Today, this generalisation is known as the Fréchet mean. Precisely, a Fréchet mean, of a probability distribution $\pi$ on a metric space $M$, is any global minimum of the so-called variance function
\begin{equation} \label{eq:frechet}
\mathcal{E}_{\pi}(y) \,=\,
\frac{1}{2}\hspace{0.03cm}
\int_M\,d^{\hspace{0.03cm}\scriptscriptstyle 2}(y\hspace{0.02cm},x)\hspace{0.03cm}\pi(dx)
%ONE HALF
%mathcal{E} separately
\end{equation}
where $d(x\hspace{0.02cm},y)$ denotes the distance between $x$ and $y$ in $M$. In the following, the focus will be on the case where $M$ is a Riemannian manifold. Then, a Fréchet mean of $\pi$ will be called a Riemannian barycentre, or just a barycentre, of $\pi$.

If $\mathcal{E}_{\pi}(y)$ takes on finite values (in fact, if it is finite for just one $y = y_o$), then $\pi$ has at least one Fréchet mean. In particular, if $M$ is a Euclidean space, then this Fréchet mean is always unique, and equal to the mean value (expectation) of $\pi$. In general, the Fréchet mean of a probability distribution $\pi$ is not unique, and one may think of the Fréchet mean of $\pi$ as the set $F(\pi)$, of all global minima of its variance function $\mathcal{E}_{\pi\hspace{0.03cm}}$. \\[0.1cm]%In contrast to the case where $M$ is a Euclidean space, if $M$ is a circle, with its canonical metric, and $\pi$ the uniform distribution on this circle, then $F(\pi) = M$. 
\textbf{Example 1\,:} if $M = S^1$, the unit circle, and $\pi$ is the uniform distribution (\textit{i.e.} Haar measure), on $S^1$, then $F(\pi) = S^1$. Any point on the circle is a barycentre of the uniform distribution.\\[0.1cm]
\indent If $x_{\scriptscriptstyle 1},\ldots,x_{\scriptscriptstyle N} \in M$, then an empirical Fréchet mean of $(x_{\scriptscriptstyle 1},\ldots,x_{\scriptscriptstyle N})$ is any Fréchet mean of the empirical distribution $(\delta_{x_{\scriptscriptstyle 1}}+\ldots+\delta_{x_{\scriptscriptstyle N}})/N$ ($\delta_x$ denotes the Dirac distribution concentrated at $x$). In other words, an empirical Fréchet mean of $(x_{\scriptscriptstyle 1},\ldots,x_{\scriptscriptstyle N})$ is any global minimum of the empirical variance function
\begin{equation} \label{eq:empiricalfrechet}
\mathcal{E}_{\scriptscriptstyle N}(y) \,=\,\frac{1}{2N}\sum^N_{n=1}d^{\hspace{0.03cm}2}(y\hspace{0.02cm},x_n)
%ONE HALF
%mathcal{E}_N separately
\end{equation}
When $M$ is a Riemannian manifold, the term ``empirical Fréchet mean" will be replaced by the term ``empirical barycentre". \\[0.1cm]
\textbf{Example 2\,:} if $M = S^1$, and $x_{\scriptscriptstyle 1\hspace{0.02cm}},x_{\scriptscriptstyle 2}$ are two opposite points on $S^1$, then the empirical barycentre of $(x_{\scriptscriptstyle 1\hspace{0.02cm}},x_{\scriptscriptstyle 2})$ is a two-point set. For example, if $x_{\scriptscriptstyle 1} = 1$ and $x_{\scriptscriptstyle 2} = -1$, then the empirical barycentre is the set $\lbrace i,-i\rbrace$ ($i$ being the square root of $-1$). \\[0.1cm]
Two important problems, in relation to the concept of Fréchet mean, are establishing the uniqueness of the Fréchet mean of some probability distribution, and effectively computing this Fréchet mean, (or the set of Fréchet means, in case uniqueness does not hold). 

Another type of problem is related to the large-sample theory of the Fréchet mean, and was treated in~\cite{bhatta1}\cite{bhatta2}.  Assume $(x_n\,;n\geq 1)$ are independent samples from the distribution $\pi$. If $F_{\scriptscriptstyle N}$ is the set of empirical Fréchet means of $(x_{\scriptscriptstyle 1},\ldots,x_{\scriptscriptstyle N})$, then one is interested in using $F_{\scriptscriptstyle N}$ to somehow approximate $F(\pi)$.  

In~\cite{bhatta1}, it was shown that, if the metric space $M$ is such that any closed and bounded subset of $M$ is compact\footnote{That is, for any $x \in M$, the function $y \mapsto d(x\hspace{0.02cm},y)$ is a proper function, meaning it has compact sublevel sets.}, then for any $\epsilon > 0$, the set $F_{\scriptscriptstyle N}$ almost-surely belongs to the $\epsilon$-neighborhood of the set $F(\pi)$, when $N$ is sufficiently large. 

Moreover, if $\pi$ has a unique Fréchet mean, say $F(\pi) = \lbrace\hat{x}_{\pi}\rbrace$, then any sequence of empirical Fréchet means, $\bar{x}_{\scriptscriptstyle N} \in F_{\scriptscriptstyle N}$ converges almost-surely to $\hat{x}_{\pi}$ (an extension of this last result, from independent to Markovian samples, is obtained in \ref{ssec:recursiveb}). 

In~\cite{bhatta2}, a central limit theorem was added to this last convergence result. Specifically, if $M$ is a Riemannian manifold, the distribution of $N^{\scriptscriptstyle \frac{1}{2}}\hspace{0.03cm}\mathrm{Exp}_{\hat{x}_\pi}(\bar{x}_{\scriptscriptstyle N})$ converges to a multivariate normal distribution (in the tangent space at $\hat{x}_\pi$). This ``central limit theorem" requires several technical conditions, in order to hold true, and should therefore only be applied after due verification.
%Establishing the existence and uniqueness, of the Fréchet mean of a probability distribution, an
%Here, $\mathcal{E}_{\scriptscriptstyle N}$ is the variance function of the empirical distribution $\hat{\pi}_{}$

%distance is proper
%bhattacharya
%hausdorff distance

\vfill
\pagebreak

\section{Existence and uniqueness} \label{sec:barycentrebackground}
The problem of the existence and uniqueness of Riemannian barycentres has generated a rich literature, with ramifications in stochastic analysis on manifolds, Riemannian geometry, and probability theory. The present section attempts a quick, non-exhaustive summary of some famous results from this literature.

\subsection{Emery and Kendall}  \label{subsec:emery}
The works of Emery and Kendall~\cite{kendall}, later expanded upon by Afsari~\cite{afsari}, are related to the existence and uniqueness of the Riemannian barycentre of a probability distribution $\pi$, supported inside some geodesic ball $B(x^*\!,\delta)$, in a Riemannian manifold $M$.

Emery and Kendall, among others, considered the so-called Karcher mean of $\pi$. This is a local minimum of the variance function $\mathcal{E}_\pi$ in (\ref{eq:frechet}). In~\cite{kendall}, $\pi$ is assumed to have compact support, inside a so-called regular geodesic ball $B(x^*\!,\delta)$. Here, ``regular geodesic ball" means  \\[0.1cm]
$\bullet$ $\delta < \frac{\pi}{2}c^{\scriptscriptstyle -1}$, where all sectional curvatures of $M$ are less than $\kappa_{\max} = c^{\hspace{0.02cm}\scriptscriptstyle 2}$.\\
$\bullet$ the cut locus of $x^*$ does not intersect $B(x^*\!,\delta)$ (that is $\delta < \mathrm{inj}(x^*)$).\\[0.1cm]
These two conditions guarantee that the closed ball $\bar{B}(x^*\!,\delta)$ is weakly convex, and that it has convex geometry.

Weakly convex means for any $x\hspace{0.02cm},y \in  \bar{B}(x^*\!,\delta)$ there exists a unique geodesic $\gamma:[0,1]\rightarrow M$, such that $\gamma(0) = x$, $\gamma(1) = y$ and $\gamma(t) \in \bar{B}(x^*\!,\delta)$ for all $t \in [0,1]$ (this is equivalent to the terminology of~\cite{chavel}\footnote{This geodesic $\gamma$ is the unique length-minimising curve, among all curves which connect $x$ to $y$ and lie in $\bar{B}(x^*\!,\delta)$. See the proof of Theorem IX.6.2, Page 405 in~\cite{chavel}.}). 

Convex geometry means there exists a positive, bounded, continuous, and convex function $\Psi$, defined on $\bar{B}(x^*\!,\delta) \times \bar{B}(x^*\!,\delta)$, such that $\Psi(x\hspace{0.02cm},y) = 0$ if and only if $x = y$.  

When $\pi$ is supported inside $B(x^*\!,\delta)$, the function $\mathcal{E}_\pi$ takes on finite values, and therefore has a global minimum $\hat{x}_\pi\hspace{0.03cm}$. However, it is not immediately clear this $\hat{x}_\pi$ should lie within $B(x^*\!,\delta)$. In~\cite{kendall}, the existence of a local minimum, \textit{i.e.} Karcher mean, within $B(x^*\!,\delta)$ is guaranteed, subject to interpreting the distance in (\ref{eq:frechet}) as geodesic distance within $B(x^*\!,\delta)$.

If $\hat{x}_\pi$ is a local minimum of $\mathcal{E}_\pi$ in $B(x^*\!,\delta)$, then the convex geometry property of the closed ball $\bar{B}(x^*\!,\delta)$ guarantees this local minimum is unique. This follows by using a general form of Jensen's inequality, due to Emery. Specifically, if $\hat{x}_{\scriptscriptstyle 1}$ and $\hat{x}_{\scriptscriptstyle 2}$ are Karcher means in $B(x^*\!,\delta)$, then $(\hat{x}_{\scriptscriptstyle 1\hspace{0.02cm}},\hat{x}_{\scriptscriptstyle 2})$ is a Karcher mean of the image distribution $\delta^*\pi$ of $\pi$, under the map $\delta(x) = (x,x)$. Then, applying Jensen's inequality to the convex function $\Psi$, it follows
$$
\Psi(\hat{x}_{\scriptscriptstyle 1\hspace{0.02cm}},\hat{x}_{\scriptscriptstyle 2}) \leq \int_{\scriptscriptstyle \bar{B}(x^*\!,\delta)}\Psi(x,x)\hspace{0.03cm}\pi(dx) 
$$
so $\Psi(\hat{x}_{\scriptscriptstyle 1\hspace{0.02cm}},\hat{x}_{\scriptscriptstyle 2}) = 0$, and therefore $\hat{x}_{\scriptscriptstyle 1} = \hat{x}_{\scriptscriptstyle 2\hspace{0.03cm}}$. \\[0.1cm]
\textbf{Remark\,:} it was conjectured by Emery that any weakly convex geodesic ball should also have convex geometry. A counterexample to this conjecture was provided by Kendall, in the form of his ``propeller"~\cite{propeller}.

\subsection{Afsari's contribution} \label{subsec:afsari}
Afsari's seminal work on Riemannian barycentres was published ten years ago~\cite{afsari}. It provided the following statement\,: if $\pi$ is supported inside a geodesic ball $B(x^*\!,\delta)$, then $\pi$ has a unique Riemannian barycentre $\hat{x}_{\pi}$ and $\hat{x}_{\pi} \in B(x^*\!,\delta)$, as soon as
\begin{equation} \label{eq:afsari1}
  \delta < \frac{1}{2}\hspace{0.02cm}\min\left\lbrace \pi c^{\scriptscriptstyle -1},\mathrm{inj}(M)\right\rbrace
\end{equation}
Here, $c$ is such that all sectional curvatures of $M$ are less than $\kappa_{\max} = c^{\hspace{0.02cm}\scriptscriptstyle 2}$ (if $M$ has negative sectional curvatures, $c^{\scriptscriptstyle -1}$ is understood to be $+\infty$), and $\mathrm{inj}(M)$ is the injectivity radius of $M$. 

Condition (\ref{eq:afsari1}) ensures the geodesic ball $B(x^*\!,\delta)$ is convex, in the sense of \ref{ssec:sqd} (strongly convex, in the terminology of~\cite{chavel}), rather than just weakly convex as in \ref{subsec:emery}. This stronger condition is required, because the Riemannian barycentre (Fréchet mean) is considered, rather than just the Karcher mean. In fact, Afsari extended his results beyond Riemannian barycentres to $L^{\scriptscriptstyle p}$ Riemannian barycentres, which are obtained by replacing the squared distance in (\ref{eq:frechet}) with a distance elevated to the power $p$, where $p\geq 1$.

In~\cite{afsari}, the following approach is used, for the proof of existence and uniqueness. First, it is shown that any global minimum of $\mathcal{E}_\pi$ must lie inside $B(x^*\!,\delta)$. This is done using the Alexandrov-Toponogov comparison theorem, under its form stated in~\cite{chavel} (Page 420). Then, the Poincaré-Hopf theorem is employed, in order to prove uniqueness of local minima, inside the geodesic ball $B(x^*\!,\delta)$. 

Specifically, $\mathcal{E}_\pi$ is differentiable at any point $y$ which belongs to the closed ball $\bar{B}(x^*\!,\delta)$, and 
$$
\mathrm{grad}\,\mathcal{E}_\pi(y) \,=\, -\int_{M}\,\mathrm{Exp}^{-1}_y(x)\hspace{0.03cm}\pi(dx) \hspace{0.5cm} \text{for } y \in \bar{B}(x^*\!,\delta)
$$
Then, it is shown that, if $y \in B(x^*\!,\delta)$ and $\mathrm{grad}\,\mathcal{E}_\pi(y) = 0$, then $\mathrm{Hess}\,\mathcal{E}_\pi(y)$ is positive-definite. In other words, the singular point $y$ of the gradient vector field $\mathrm{grad}\,\mathcal{E}_\pi$ has its index equal to $1$.\hfill\linebreak Since this vector field is outward pointing on the boundary of $\bar{B}(x^*\!,\delta)$, the Poincaré-Hopf theorem implies the sum of the indices of all its singular points in $B(x^*\!,\delta)$ is equal to the Euler-Poincaré characteristic of $\bar{B}(x^*\!,\delta)$, which is equal to $1$ (since $\bar{B}(x^*\!,\delta)$ is homeomorphic to a closed ball in $\mathbb{R}^n$). \\[0.1cm]
\textbf{Remark\,:} the argument just summarised not only shows that $\mathcal{E}_\pi$ has a unique local minimum in $B(x^*\!,\delta)$, but that it has a unique stationary point in $B(x^*\!,\delta)$. Moreover, the advantage of this argument, over the ``convex geometry" uniqueness argument, (summarised in \ref{subsec:emery}),\hfill\linebreak is that it can be used to show the uniqueness of $L^{\scriptscriptstyle p}$ Riemannian barycentres, for general $p > 1$.

\subsection{Hadamard manifolds} \label{subsec:frechhad}
Existence and uniqueness of Riemannian barycentres hold under quite general conditions, when the underlying Riemannian manifold $M$ is a Hadamard manifold (recall definition from \ref{ssec:sqd}). Mostly, these existence and uniqueness properties are just special cases of the properties of Fréchet means in metric spaces of non-positive curvature, which were developed by Sturm~\cite{sturm}.

Let $\pi$ be a probability distribution on a Hadamard manifold $M$. As already mentioned in \ref{sec:frechet}, if the variance function $\mathcal{E}_\pi$ in (\ref{eq:frechet}) takes on finite values, then $\pi$ has at least one Riemannian barycentre, say $\hat{x}_{\pi\hspace{0.03cm}}$. For this, it is enough that
$\mathcal{E}_\pi(y_o) < \infty$, for just one $y_o \in M$. In other words, it is enough that $\pi$ should have a finite second-order moment
\begin{equation} \label{eq:secondordermoment}
\int_M\,d^{\hspace{0.03cm} 2}(y_o\hspace{0.03cm},x)\,\pi(dx) \,<\,\infty
\end{equation}
Indeed, if (\ref{eq:secondordermoment}) is verified, then a straightforward application of the triangle inequality implies that $\mathcal{E}_\pi(y) < \infty$ for all $y \in M$. 

When $M$ is a Hadamard manifold, existence of a Riemannian barycentre automatically implies its uniqueness. This can be shown using the ``convex geometry" uniqueness argument, discussed in \ref{subsec:emery}. Indeed, if $M$ is a Hadamard manifold, then $\Psi:M\times M \rightarrow \mathbb{R}$, where $\Psi(x\hspace{0.02cm},y) = d(x\hspace{0.02cm},y)$ is convex, and $\Psi(x\hspace{0.02cm},y) = 0$ if and only if $x = y$. Alternatively, uniqueness of the Riemannian barycentre follows from the strong convexity of the variance function $\mathcal{E}_{\pi\hspace{0.03cm}}$. Recall from \ref{ssec:sqd} that $f_x(y) = d^{\hspace{0.03cm} 2}(x\hspace{0.02cm},y)/2$ is a $1/2$-strongly convex function, for each $x \in M$. Then, (\ref{eq:frechet}) says that $\mathcal{E}_\pi$ is an expectation of $1/2$-strongly convex functions, and is therefore $1/2$-strongly convex. In turn, this implies that $\mathcal{E}_\pi$ has a unique global minimum, $\hat{x}_\pi \in M$. 

When $M$ is a Hadamard manifold, it should also be noted that $\mathcal{E}_\pi$ is smooth throughout $M$, and that its gradient is given by
\begin{equation} \label{eq:gradepsilonhadamard}
\mathrm{grad}\,\mathcal{E}_\pi(y) = -\int_M\,\mathrm{Exp}^{-1}_y(x)\hspace{0.03cm}\pi(dx)
\end{equation}
as can be found by applying (\ref{eq:gradfx}) under the integral in (\ref{eq:frechet}). Strong convexity of $\mathcal{E}_\pi$ implies its global minimum $\hat{x}_\pi$ is also its unique stationary point in $M$ (\textit{i.e.} the unique point where $\mathrm{grad}\,\mathcal{E}_\pi$ is equal to zero). 

\subsection{Generic uniqueness} \label{subsec:generic}
The empirical barycentre of the points $(x_{\scriptscriptstyle 1},\ldots,x_{\scriptscriptstyle N})$, in any complete Riemannian manifold  $M$,\hfill\linebreak is generically unique. This means that this empirical barycentre is unique, for almost all $(x_{\scriptscriptstyle 1},\ldots,x_{\scriptscriptstyle N})$ in the product Riemannian manifold $M^{\scriptscriptstyle N} = M \times \ldots \times M$, equipped with its Riemannian volume measure. This interesting result was obtained by Arnaudon and Miclo~\cite{arnaumic}. In particular, it implies that when $(x_{\scriptscriptstyle 1},\ldots,x_{\scriptscriptstyle N})$ are independent samples, from a distribution $\pi$, which has a probability density with respect to the Riemannian volume of $M$, then their empirical barycentre $\bar{x}_{\scriptscriptstyle N}$ is almost-surely unique. 

\section{Gibbs distributions\,: an open problem} \label{sec:gibbs}
Throughout the following, $M$ will be a compact, orientable Riemannian manifold, with positive sectional curvatures, all less than $\kappa_{\max} = c^{\hspace{0.02cm}\scriptscriptstyle 2}$. Afsari's statement, recalled in \ref{subsec:afsari}, says that if  $\pi$ is a probability distribution on $M$, supported inside a convex geodesic ball $B(x^*\!,\delta)$, then $\pi$ has a unique Riemannian barycentre $\hat{x}_\pi\hspace{0.03cm}$, as soon as  
\begin{equation} \label{eq:afsaribis}
  \delta < \frac{1}{2}\hspace{0.02cm}\min\left\lbrace \pi c^{\scriptscriptstyle -1},\mathrm{inj}(M)\right\rbrace
\end{equation}
where $\mathrm{inj}(M)$ denotes the injectivty radius of $M$. 

Inequality (\ref{eq:afsaribis}) is optimal. Indeed, it is easy to think of examples which show that, if it is replaced by an equality, then $\hat{x}_\pi$ will immediately fail to be unique. On the other hand, this inequality does not tell us what happens in the important case where $\pi = \pi_{\scriptscriptstyle T}$ is a Gibbs distribution,
\begin{equation} \label{eq:gibbs}
  \pi_{\scriptscriptstyle T}(dx) = \left(Z(T)\right)^{\scriptscriptstyle -1}\exp\left[-\frac{U(x)}{T}\right]\hspace{0.03cm}\mathrm{vol}(dx) \\[0.12cm]
\end{equation}
for some temperature $T$, and potential function $U:M\rightarrow \mathbb{R}$, where $Z(T)$ is a normalising constant ($\mathrm{vol}$ denotes the Riemannian volume form). 

The present chapter will introduce several results, which deal with this case. These are concerned with the concentration, differentiability, convexity, and uniqueness properties, of the Riemannian barycentre $\hat{x}_{\scriptscriptstyle T}$ of the Gibbs distribution $\pi_{\scriptscriptstyle T}$.

The starting assumption for these results is that the potential function $U$ has a unique global minimum at $x^* \in M$. Under this assumption, while $\pi_{\scriptscriptstyle T}$ is not supported inside any convex geodesic ball $B(x^*\!,\delta)$, it is still \textit{concentrated} on any such ball, provided the temperature $T$ is sufficiently small. Then, the aim is to know exactly how small $T$ should be made, in order to ensure the required properties of $\hat{x}_{\scriptscriptstyle T\hspace{0.02cm}}$. This aim can be fully achieved, under the further assumption that $M$ is a simply connected compact Riemannian symmetric space. 

Given these two assumptions, the following conclusion will be obtained\,: for each $\delta < \frac{1}{2}r_{\scriptscriptstyle cx}$ ($r_{\scriptscriptstyle cx}$ denotes the convexity radius of $M$), there exists a critical temperature $T_{\scriptscriptstyle \delta}$ such that $T < T_{\scriptscriptstyle \delta}$ implies $\pi_{\scriptscriptstyle T}$ has a unique Riemannian barycentre $\hat{x}_{\scriptscriptstyle T}$ and this $\hat{x}_{\scriptscriptstyle T}$ belongs to the geodesic ball $B(x^*\!,\delta)$. Moreover, if $U$ is invariant by geodesic symmetry about $x^*$, then $\hat{x}_{\scriptscriptstyle T} = x^*$. \\[0.1cm]
\textbf{Remark\,:} if $M$ is a Riemannian manifold, the convexity radius $r_{\scriptscriptstyle cx}(x)$ of $x \in M$ is the supremum of $R > 0$ such that the geodesic ball $B(x\hspace{0.02cm},R)$ is convex (this is strictly positive, for any $x \in M$). The convexity radius $r_{\scriptscriptstyle cx}(M)$ of $M$ is the infimum of $r_{\scriptscriptstyle cx}(x)$, over all $x \in M$ (if $M$ is compact, this is strictly positive). Here, $r_{\scriptscriptstyle cx}(M)$ is just denoted $r_{\scriptscriptstyle cx\hspace{0.03cm}}$.

%%state assumptions
%%add arnaudon and miclo to bibliography

%%optimality of Afsari

\section{Concentration of barycentres} \label{sec:concentration}
Denote the variance function of the Gibbs distribution $\pi_{\scriptscriptstyle T}$ in (\ref{eq:gibbs}) by $\mathcal{E}_{\scriptscriptstyle T\hspace{0.03cm}}$. According to (\ref{eq:frechet}),
\begin{equation} \label{eq:ETT}
\mathcal{E}_{\scriptscriptstyle T}(y) = 
\frac{1}{2}\hspace{0.03cm}
\int_M\,d^{\hspace{0.03cm}\scriptscriptstyle 2}(y\hspace{0.02cm},x)\hspace{0.03cm}\pi_{\scriptscriptstyle T}(dx)
\end{equation}
Throughout the following, it will be assumed that the potential function $U$, which appears in (\ref{eq:gibbs}), has a unique global minimum at $x^* \in M$. While $U$ is not required to be smooth, it is required to be well-behaved near $x^*$, in the sense that there exist $\mu_{\min}\hspace{0.02cm},\mu_{\max} > 0$ and $\rho > 0$ such that
\begin{equation} \label{eq:wellbehaved}
 \mu_{\min}\hspace{0.03cm}d^{\hspace{0.03cm}\scriptscriptstyle 2}(x\hspace{0.02cm},x^*) \,\leq\, 2(U(x) - U(x^*))\,\leq
 \mu_{\max}\hspace{0.03cm}d^{\hspace{0.03cm}\scriptscriptstyle 2}(x\hspace{0.02cm},x^*)
\end{equation}
whenever $d(x\hspace{0.02cm},x^*) \leq \rho$. This is always verified if $U$ is twice differentiable at $x^*$, and the spectrum of $\mathrm{Hess}\,U(x^*)$ is contained in the open interval $(\mu_{\min}\hspace{0.02cm},\mu_{\max})$.

The following Proposition \ref{prop:concentration} establishes the concentration property of the Riemannian barycentres of $\pi_{\scriptscriptstyle T}$ as the temperature $T$ is made small. In this proposition, $W$ denotes the Kantorovich ($L^{\scriptscriptstyle 1}$-Wasserstein) distance, and $\delta_{x^*}$ the Dirac distribution concentrated at $x^*$.
\begin{proposition} \label{prop:concentration}
 Let $M$ be a compact, orientable Riemannian manifold, with positive sectional curvatures, and dimension equal to $n$. \\[0.1cm]
(i) Let $\eta > 0$. For any Riemannian barycentre $\hat{x}_{\scriptscriptstyle T}$ of $\pi_{\scriptscriptstyle T}$
\begin{equation} \label{eq:concentration1}
  W(\pi_{\scriptscriptstyle T}\hspace{0.02cm},\delta_{x^*}) < \frac{\eta^2}{4\hspace{0.02cm}\mathrm{diam}\,M} \;\Longrightarrow\; d(\hat{x}_{\scriptscriptstyle T}\hspace{0.02cm},x^*) < \eta
\end{equation}
where $\mathrm{diam}\,M$ is the diameter of $M$. \\
(ii) There exists a temperature $T_{\scriptscriptstyle W}$ such that $T \leq T_{\scriptscriptstyle W}$ implies
\begin{equation} \label{eq:concentration2}
W(\pi_{\scriptscriptstyle T}\hspace{0.02cm},\delta_{x^*}) \leq (8\pi)^{\!\frac{1}{2}}\hspace{0.02cm}B^{-1}_n\left(\frac{\pi}{2}\right)^{\!n-1}
\left(\frac{\mu_{\max}}{\mu_{\min}}\right)^{\!\!\frac{n}{2}}\left(\frac{T}{\mu_{\min}}\right)^{\!\!\frac{1}{2}}
\end{equation}
where $B_n = B(1/2\hspace{0.02cm},n/2)$ in terms of the Euler Beta function.
\end{proposition}
Proposition \ref{prop:concentration} shows exactly how small $T$ should be made, in order to ensure that all the Riemannian barycentres $\hat{x}_{\scriptscriptstyle T}$ concentrate within an open ball $B(x^*\!,\eta)$. Roughly, (i) states that, if $\pi_{\scriptscriptstyle T}$ is close to $\delta_{x^*}\hspace{0.03cm}$, then all $\hat{x}_{\scriptscriptstyle T}$ will be close to $x^*$. On the other hand, (ii) bounds the distance between $\pi_{\scriptscriptstyle T}$ and $\delta_{x^*}\hspace{0.03cm}$, as a function of $T$. The temperature $T_{\scriptscriptstyle W}$ mentioned in (ii) will be expressed explicitly  in \ref{sec:twtd}, below.

Here, two things should be noted, concerning (\ref{eq:concentration2}). First, this inequality is both optimal and explicit. It is optimal because the dependence on $T^{\frac{1}{2}}$ in its right-hand side cannot be improved. Indeed, the multi-dimensional Laplace approximation (for example, see~\cite{wong}), shows the left-hand side is equivalent to $\mathrm{L}\cdot T^{\frac{1}{2}}$ when $T \rightarrow 0$. While this constant $\mathrm{L}$ is not tractable, the constants appearing in (\ref{eq:concentration2}) depend explicitly on the manifold $M$ and the function $U$. In fact, (\ref{eq:concentration2}) does not follow from the multi-dimensional Laplace approximation, but rather from the volume comparison theorems, in \ref{subsec:volcomp}. 

Second, in spite of these nice properties, (\ref{eq:concentration2}) does not escape the curse of dimensionality. Indeed, for fixed $T$, its right-hand side increases exponentially with the dimension $n$ of $M$ (note that $B_n$ decreases like $n^{\scriptscriptstyle -\frac{1}{2}}$). In fact, the temperature $T_{\scriptscriptstyle W}$ also depends on $n$, but it is typically much less affected by it, and decreases slower than $n^{\scriptscriptstyle -1}$ as $n$ increases.

%highlight conclusion

\section{Differentiability of the variance function} \label{sec:diff}
Assume that $M$ is a simply connected compact Riemannian symmetric space. Under this assumption, it turns out that the variance function $\mathcal{E}_{\scriptscriptstyle T}(y)$ is $C^2$ throughout $M$, for any value $T > 0$ of the temperature $T$. This surprising result is contained in Proposition \ref{prop:differentiability}. 

To state  Proposition \ref{prop:differentiability}, consider for $x \in M$ the function $f_x(y) = d^{\hspace{0.03cm}2}(x,y)/2$. Recall from \ref{ssec:sqd} that this function is $C^2$ on the open set $\mathrm{D}(x) = M - \mathrm{Cut}(x)$. When $y \in \mathrm{D}(x)$, denote $G_y(x)$ and $H_y(x)$ the gradient and Hessian of $f_x(y)$.

With this notation, for any $x \in M$, the gradient $G_y(x)$ belongs to $T_yM$, and the Hessian $H_y(x)$ defines a symmetric bilinear form on $T_yM$. However (recall the remarks in \ref{subsec:rx}), both $G_y(x)$ and $H_y(x)$ are singular on $\mathrm{Cut}(x)$, where $H_y(x)$ will even blow up, as it has an eigenvalue equal to $-\infty$.
\begin{proposition} \label{prop:differentiability}
Let $M$ be a simply connected compact Riemannian symmetric space. \\[0.1cm]
(i) The following integrals converge for any temperature $T > 0$
\begin{equation} \label{eq:GH}
  G_y = \int_{\mathrm{D}(y)}G_y(x)\hspace{0.03cm}\pi_{\scriptscriptstyle T}(dx) \hspace{0.2cm};\hspace{0.2cm}
  H_y = \int_{\mathrm{D}(y)}H_y(x)\hspace{0.03cm}\pi_{\scriptscriptstyle T}(dx)
\end{equation}
and both depend continuously on $y$. \\
\noindent (ii)  The gradient and Hessian of the variance function $\mathcal{E}_{\scriptscriptstyle T}(y)$ are given by
\begin{equation} \label{eq:derivatives}
 \mathrm{grad}\,\mathcal{E}_{\scriptscriptstyle T}(y) = G_y \hspace{0.2cm};\hspace{0.2cm}
 \mathrm{Hess}\,\mathcal{E}_{\scriptscriptstyle T}(y) = H_y
\end{equation}
so that $\mathcal{E}_{\scriptscriptstyle T}(y)$ is $C^2$ throughout $M$.
\end{proposition}
The proof of Proposition \ref{prop:differentiability} relies on the following lemma. 
\begin{lemma} \label{lem:hddim}
  Assume $M$ is a simply connected compact Riemannian symmetric space. Let $\gamma : I \rightarrow M$ be a geodesic defined on a compact interval $I$. Denote by $\mathrm{Cut}(\gamma)$ the union of all cut loci $\mathrm{Cut}(\gamma(t))$ for $t \in I$. Then, the Hausdorff dimension of $\mathrm{Cut}(\gamma)$ is strictly less than the dimension of $M$. In particular, $\mathrm{Cut}(\gamma)$ is a set with Riemannian volume equal to zero.  
\end{lemma}
\noindent \textbf{Remark\,:} the assumption that $M$ is simply connected cannot be removed. For example, the conclusion of Lemma \ref{lem:hddim} does not hold if $M$ is a real projective space. \\[0.1cm]
\noindent The proof of Lemma \ref{lem:hddim} uses the structure of Riemannian symmetric spaces, as well as some results from dimension theory, found in~\cite{hausdorff}. The notion of Hausdorff dimension is needed, because $\mathrm{Cut}(\gamma)$ may fail to be a manifold.

Lemma \ref{lem:hddim} is crucial to Proposition \ref{prop:differentiability}, because it leads to the following expression,
$$
\mathcal{E}_{\scriptscriptstyle T}(\gamma(t)) = \int_M\,f_x(\gamma(t))\hspace{0.03cm}\pi_{\scriptscriptstyle T}(dx) =
\int_{\mathrm{D}(\gamma)}f_x(\gamma(t))\hspace{0.03cm}\pi_{\scriptscriptstyle T}(dx) \hspace{0.5cm} \text{for all $t \in I$}
$$
where $\mathrm{D}(\gamma) = M - \mathrm{Cut}(\gamma)$, and the second inequality follows since $\mathrm{Cut}(\gamma)$ has Riemannian volume equal to zero. Then, recalling that $x \in \mathrm{Cut}(\gamma(t))$ if and only if $\gamma(t) \in \mathrm{Cut}(x)$, it becomes possible to differentiate $f_x(\gamma(t))$ under the integral. This leads to the proof of (ii).
\section{Uniqueness of the barycentre} \label{sec:unique}
The following Proposition \ref{prop:unique} establishes the uniqueness of $\hat{x}_{\scriptscriptstyle T}$ as the temperature $T$ is made small. As in the previous Proposition \ref{prop:differentiability}, $M$ is a simply connected compact Riemannian symmetric space. The convexity radius of $M$ is denoted $r_{\scriptscriptstyle cx\hspace{0.03cm}}$. This is given by $r_{\scriptscriptstyle cx} = \frac{\pi}{2}\hspace{0.03cm}c^{\scriptscriptstyle -1}$ (see \ref{sec:ssbis}, below).  

Recall the definition (\ref{eq:gibbs}) of the Gibbs distribution $\pi_{\scriptscriptstyle T\hspace{0.03cm}}$, where the potential function $U$ has a unique global minimum at $x^* \in M$. Let $s_{x^*}$ denote the geodesic symmetry at $x^*$ (recall definition from \ref{ssec:sspace}). The potential function $U$ is said to be invariant by geodesic symmetry about $x^*$, if $U \circ s_{x^*} = U$.
\begin{proposition} \label{prop:unique}
 Let $M$ be a simply connected compact Riemannian symmetric space, with convexity radius $r_{\scriptscriptstyle cx\hspace{0.03cm}}$. For $\delta < \frac{1}{2}r_{\scriptscriptstyle cx\hspace{0.03cm}}$, there exists a critical temperature $T_{\scriptscriptstyle \delta}$ such that  \\[0.1cm]
(i) When $T < T_{\scriptscriptstyle \delta\hspace{0.03cm}}$, the Riemannian barycentre $\hat{x}_{\scriptscriptstyle T}$ of $\pi_{\scriptscriptstyle T}$ is unique and $\hat{x}_{\scriptscriptstyle T} \in B(x^*\!,\delta)$. \\
\noindent
(ii) If, in addition, $U$ is invariant by geodesic symmetry about $x^*$, then $\hat{x}_{\scriptscriptstyle T} = x^*$.
\end{proposition}
Proposition \ref{prop:unique} shows exactly how small $T$ should be made, in order to ensure that the Riemannian barycentre $\hat{x}_{\scriptscriptstyle T}$ is unique. In turn, this uniqueness of $\hat{x}_{\scriptscriptstyle T}$ follows from the convexity of the variance function $\mathcal{E}_{\scriptscriptstyle T}(y)$, obtained in the following Lemma \ref{lem:etconv}.

To state this lemma, consider the function $f(T)$ of the temperature $T$
\begin{equation} \label{eq:fT}
  f(T) = \left(\frac{2}{\pi}\right)\left(\frac{\pi}{8}\right)^{\!n-1}\left(\frac{\mu_{\max}}{T}\right)^{\!\!\frac{n}{2}}\exp\left(-\frac{U_\delta}{T}\right)
\end{equation}
for any given $\delta$, where $U_\delta = \inf\lbrace U(x) - U(x^*)\,; x \notin B(x^*\!,\delta)\rbrace$. Note that $f(T)$ decreases to zero as $T$ is made arbitrarily small.
\begin{lemma} \label{lem:etconv}
Under the same assumptions as Proposition \ref{prop:unique}, let $\delta < \frac{1}{2}r_{\scriptscriptstyle cx\hspace{0.03cm}}$. \\[0.1cm]
(i) For all $y \in B(x^*\!,\delta)$,
\begin{equation} \label{eq:ethlower}
  \mathrm{Hess}\,\mathcal{E}_{\scriptscriptstyle T}(y) \,\geq\, \mathrm{Ct}(2\delta)\hspace{0.03cm}[1 - \mathrm{vol}(M)f(T)] - \pi A_{\scriptscriptstyle M}\hspace{0.03cm}f(T)
\end{equation}
where $\mathrm{Ct}(2\delta) = 2c\delta\hspace{0.02cm}\cot(2c\delta)$ and $A_{\scriptscriptstyle M} > 0$ is a constant which depends only on the symmetric space $M$. \\
\noindent (ii) There exists a critical temperature $T_{\scriptscriptstyle \delta}$ such that $T < T_{\scriptscriptstyle \delta}$ implies the variance function $\mathcal{E}_{\scriptscriptstyle T}(y)$ is strongly convex on $B(x^*\!,\delta)$.
\end{lemma}
The inequality in (\ref{eq:ethlower}) should be understood as saying all the eigenvalues of $\mathrm{Hess}\,\mathcal{E}_{\scriptscriptstyle T}(y)$ are greater than the right-hand side (of course, this is an abuse of notation). The critical temperature $T_{\scriptscriptstyle \delta}$ will be expressed in the following section.

\section{Finding $T_{\scriptscriptstyle W}$ and $T_{\scriptscriptstyle \delta}$} \label{sec:twtd}
The present paragraph provides expressions of the temperatures $T_{\scriptscriptstyle W}$ and $T_{\scriptscriptstyle \delta\hspace{0.03cm}}$, which appear in Propositions \ref{prop:concentration} and \ref{prop:unique}. These are expressions (\ref{eq:tw}) and (\ref{eq:td}) below, which should be considered as part of Propositions \ref{prop:concentration} and \ref{prop:unique}, and will accordingly be proved in \ref{sec:bcentreproofs}.

Expressions (\ref{eq:tw}) and (\ref{eq:td}) allow $T_{\scriptscriptstyle W}$ and $T_{\scriptscriptstyle \delta}$ to be computed as solutions of scalar non-linear equations, which depend on Condition (\ref{eq:wellbehaved}) and on the Riemannian symmetric space $M$.\hfill\linebreak In order to state them, write
\begin{equation} \label{eq:fTm}
  f(T,m,\rho) = \left(\frac{2}{\pi}\right)^{\!\!\frac{1}{2}}\left(\frac{\mu_{\max}}{T}\right)^{\!\!\frac{m}{2}}
\exp\left(-\frac{U_\rho}{T}\right)
\end{equation}
in terms of the temperature $T$ and positive $m$ and $\rho$, where $U_\rho$ is defined as in (\ref{eq:fT}). It should be noted that $f(T,m,\rho)$ decreases to $0$ as $T$ is made arbitrarily small, for fixed $m$ and $\rho$. The following expression holds for $T_{\scriptscriptstyle W}\hspace{0.03cm}$,
\begin{equation} \label{eq:tw}
 T_{\scriptscriptstyle W} \,=\,\min\hspace{0.02cm}\lbrace T^1_{\scriptscriptstyle W}\hspace{0.02cm},T^2_{\scriptscriptstyle W}\rbrace
\end{equation}
with $T^1_{\scriptscriptstyle W}$ and $T^2_{\scriptscriptstyle W}$ given by
$$
\begin{array}{l}
T^1_{\scriptscriptstyle W} = \inf\hspace{0.02cm}\lbrace T>0: f(T,n-2,\rho) > \rho^{2-n}\hspace{0.02cm} A_{n-1}\rbrace \\[0.4cm]
T^2_{\scriptscriptstyle W} = \inf\hspace{0.02cm}\left\lbrace T>0: f(T,n+1,\rho) > \left(\mu_{\max}\middle/\mu_{\min}\right)^{\!\scriptscriptstyle \frac{n}{2}}\hspace{0.01cm}C_n\right\rbrace
\end{array}
$$
where $A_n$ is the $n$-th absolute moment of a standard normal random variable ($A_n = \mathbb{E}|X|^n$ where $X \sim N(0,1)$), and $C_n = \left.(\omega_{n-1}\hspace{0.02cm}A_n)\middle/(\mathrm{diam}\,M\times\mathrm{vol}\,M)\right.$, where $\omega_{n-1}$ is the area of the unit sphere $S^{n-1} \subset \mathbb{R}^n$. Moreover, for $T_{\scriptscriptstyle \delta}\hspace{0.03cm}$,
\begin{equation} \label{eq:td}
 T_{\scriptscriptstyle \delta} \,=\,\min\hspace{0.02cm}\lbrace T^1_{\scriptscriptstyle \delta}\hspace{0.02cm},T^2_{\scriptscriptstyle \delta}\rbrace
\end{equation}
where, in the notation of (\ref{eq:concentration2}) and (\ref{eq:ethlower}),
$$
\begin{array}{l}
T^1_{\scriptscriptstyle \delta} = \inf\hspace{0.02cm}\left\lbrace T\leq T_{\scriptscriptstyle W}: \left(2\pi T\middle/\mu_{\min}\right)^{\!\frac{1}{2}} > \delta^{\scriptscriptstyle 2}\hspace{0.03cm}\left(\mu_{\min}\middle/\mu_{\max}\right)^{\!\frac{n}{2}}\hspace{0.03cm}D_n\right\rbrace \\[0.4cm]
T^2_{\scriptscriptstyle \delta} = \inf\hspace{0.02cm}\left\lbrace T\leq T_{\scriptscriptstyle W}: f(T) > \mathrm{Ct}(2\delta)[\mathrm{Ct}(2\delta)\mathrm{vol}(M) + \pi A_{\scriptscriptstyle M}]^{-1}\right\rbrace
\end{array}
$$
with $D_n = (2/\pi)^{n-1}\!\left.B_n\middle/(4\hspace{0.02cm}\mathrm{diam}\,M)\right.$. \\[0.1cm]
\noindent \textbf{Remark\,:} the following formulae for $A_n$ and $\omega_{n-1}$ will be useful in \ref{sec:bcentreproofs},
\begin{equation} \label{eq:gammastuff}
  A_n = \pi^{\scriptscriptstyle -\frac{1}{2}}2^{\scriptscriptstyle \frac{n}{2}}\hspace{0.03cm}\Gamma((n+1)/2) \hspace{0.2cm};\hspace{0.2cm}
  \omega_{n-1} = \frac{2\hspace{0.02cm}\pi^{\scriptscriptstyle \frac{n}{2}}}{\Gamma(n/2)}
\end{equation}
These are well-known, and follow easily from the definition of the Euler Gamma function~\cite{watson}. 

\section{Compact symmetric spaces} \label{sec:ssbis}
Compact Riemannian symmetric spaces belong to the ``compact case", already treated in \ref{ssec:sspace}. Some additional material, on these spaces, is needed for the proofs of Propositions \ref{prop:differentiability} and \ref{prop:unique}. 

\subsection{Roots and the Jacobi equation} \label{subsec:rootsjacobis}
As of now, let $M = G/K$ be a symmetric space, where $G$ is semisimple and compact, and $K = K_y$ the stabiliser in $G$ of some point $y \in M$. Recall the Cartan decomposition $\mathfrak{g} = \mathfrak{k} + \mathfrak{p}$, where $\mathfrak{g}$ and $\mathfrak{k}$ are the Lie algebras of $G$ and $K$, respectively. Moreover, let $\mathfrak{a}$ be a maximal Abelian subspace of $\mathfrak{p}$, and denote $\Delta_+$ the corresponding set of positive roots $\lambda : \mathfrak{a} \rightarrow \mathbb{R}$.

Then, $\mathfrak{p}$ may be identified with $T_yM$, and any $v \in \mathfrak{p}$ can be written $v = \mathrm{Ad}(k)\,a$ for some $k \in K$ and $a \in\mathfrak{a}$. Accordinly, the self-adjoint curvature operator, $R_v$ (given by $R_v(u) = [v\hspace{0.03cm},[v\hspace{0.02cm},u]]$ for $u \in T_yM$), can be diagonalised (the reader may wish to note (\ref{eq:raeigenbis}) differs from (\ref{eq:raeigen}) by a minus sign, since the space here denoted $\mathfrak{p}$ would have been $\mathfrak{p}_* = i\hspace{0.02cm}\mathfrak{p}$, in Chapter \ref{introduction})
\begin{equation} \label{eq:raeigenbis}
\mathrm{Ad}(k^{\scriptscriptstyle -1})\circ R_v \circ \mathrm{Ad}(k) = R_a \hspace{0.4cm} \text{where } R_a= - \sum_{\lambda \in \Delta_+} (\lambda(a))^2\;\Pi_{\lambda}
\end{equation}
and where $\Pi_{\lambda}$ is the orthogonal projector onto the eigenspace of $R_a$ which corresponds to the eigenvalue $-(\lambda(a))^2$. The rank of $\Pi_\lambda$ is denoted $m_\lambda$ and called the multiplicity of $\lambda$.

Recall that the curvature tensor of a symmetric space is parallel\,:$\hspace{0.03cm}\nabla\,R = 0$. This property, when combined with the diagonalisation (\ref{eq:raeigenbis}), yields the solutions of the operator Jacobi equation (\ref{eq:jacobioperator}), and of the Ricatti equation (\ref{eq:ricatti}). 

Alternatively, if $A(t)$ solves (\ref{eq:jacobioperator}) and $\mathcal{A}(t) = \Pi^{\scriptscriptstyle 0}_{t} \circ A(t)$, where $\Pi^{\scriptscriptstyle 0}_{t}$ denotes parallel transport, along the geodesic $c_{\scriptscriptstyle v}$ with $c_{\scriptscriptstyle v}(0) = y$ and $\dot{c}_{\scriptscriptstyle v}(0) = v$, then $\mathcal{A}(t)$ solves the differential equation
\begin{equation} \label{eq:jacobisss}
  \mathcal{A}^{\prime\prime} - R_v\hspace{0.02cm}\mathcal{A} = 0 \hspace{1cm} \mathcal{A}(0) = 0 \,,\, \mathcal{A}^\prime(0) = \mathrm{Id}_{y}
\end{equation}
where the prime denotes differentiation with respect to $t$. Using (\ref{eq:raeigenbis}), it follows that
\begin{equation} \label{eq:compactAA}
  \mathcal{A}(t) \,=\, \Pi^{k}_{\mathfrak{a}} \,+\, \sum_{\lambda \in \Delta_+}\left(\sin(\lambda(a)t)\middle/\lambda(a)\right)\hspace{0.02cm}
\Pi^{k}_{\mathfrak{\lambda}}
\end{equation} 
where $\Pi^{k}_{\mathfrak{a}} = \mathrm{Ad}(k)\circ \Pi_{\mathfrak{a}} \circ \mathrm{Ad}(k^{\scriptscriptstyle -1})$ and 
$\Pi^{k}_{\lambda} = \mathrm{Ad}(k)\circ \Pi_{\lambda} \circ \mathrm{Ad}(k^{\scriptscriptstyle -1})$, with $\Pi_{\mathfrak{a}}$ the orthogonal projector onto $\mathfrak{a}$.

\subsection{The cut locus} \label{subsec:ssbiscut}
Let $M$ be a compact Riemannian symmetric space, as above. Assume, as in Propositions \ref{prop:differentiability} and \ref{prop:unique}, that $M$ is simply connected. In ths case, the following important property holds~\cite{sakai}\,:\hfill\linebreak the cut locus of any point $y \in M$ is identical to the first conjugate locus of this point. 

Accordingly, if $v$ is a unit vector in $\mathfrak{p} \simeq T_yM$, the geodesic $c_{\scriptscriptstyle v}$ will meet the cut locus of $y$\hfill\linebreak for the first time, when $\det(\mathcal{A}(t)) = 0$ for the first time after $t = 0$. But, as seen from (\ref{eq:compactAA}),\hfill\linebreak if $v = \mathrm{Ad}(k)\,a$, then this happens when $t = \mathrm{t}(v)$ given by
\begin{equation} \label{eq:tc1}
  \mathrm{t}(v) = \min_{\lambda \in \Delta_+}\,\frac{\pi}{|\lambda(a)|} = 
                                      \min_{\lambda \in \Delta_+}\,\frac{\pi}{\lambda(a)}
\end{equation}
where the absolute value can be dropped because it is always possible to assume $a$ belongs to $\bar{C}_+\hspace{0.02cm}$, the closure of the Weyl chamber $C_+$ (the set of $a \in \mathfrak{a}$ such that $\lambda(a) > 0$ for each $\lambda \in \Delta_+$). If $M$ is an irreducible symmetric space, then there exists a maximal root $c \in \Delta_+\hspace{0.02cm}$, so that $c(a) \geq \lambda(a)$ for all $\lambda \in \Delta_+$ and $a \in \bar{C}_+$~\cite{helgason}. In this case, $\mathrm{t}(v) = \pi/c(a)$. On the other hand, if $M$ is not irreducible, it is a product of irreducible compact Riemannian symmetric spaces, say $M = M_{\scriptscriptstyle 1}\times\ldots\times M_{ s\hspace{0.03cm}}$. If $c_{\scriptscriptstyle 1},\ldots,c_s$ are the corresponding maximal roots,
\begin{equation}\label{eq:tc2}
  \mathrm{t}(v) = \min_{\ell= 1,\ldots,\hspace{0.02cm}s}\,\frac{\pi}{c_{\ell}(a)}
\end{equation}
The cut locus of $y$ is the set of all points $c_{\scriptscriptstyle v}(\mathrm{t}(v))$ where $v$ is a unit vector in $T_yM$. Then, the injectivity radius $\mathrm{inj}(y)$ of $y$ is equal to the minimum of $\mathrm{t}(v)$, taken over all unit vectors $v$.\hfill\linebreak From (\ref{eq:tc2}), this is equal to $\pi\hspace{0.03cm}c^{\scriptscriptstyle -1}$ where $c = \max_{\ell= 1,\ldots,\hspace{0.02cm}s} \Vert c_{\ell}\Vert$ and $\Vert c_{\ell}\Vert$ denotes the norm of $c_\ell \in\mathfrak{a}^*$ (the dual space of $\mathfrak{a}$). Since $M$ is a homogeneous space, the injectivity radius of $M$ is also equal to $\pi\hspace{0.03cm}c^{\scriptscriptstyle -1}$, since it is equal to the injectivity radius of any point $y$ in $M$. Incidentally, $c^{\hspace{0.02cm}\scriptscriptstyle 2}$ is the maximum sectional curvature of $M$.

With a bit of additional work, the above description of the cut locus of $y$ can be strengthened, to yield the following statements. Let $S = K/K_{\mathfrak{a}}$ where $K_{\mathfrak{a}}$ is the centraliser of $\mathfrak{a}$ in $K$. Moreover, denote $Q_+$ the set of $a \in \mathfrak{a}$ such that $\lambda(a) \in (0,\pi)$ for each $\lambda \in \Delta_+\hspace{0.02cm}$. Then, consider the mapping
\begin{equation} \label{eq:varphibis}
  \varphi(s\hspace{0.02cm},a) = \mathrm{Exp}_y(\beta(s\hspace{0.02cm},a)) \hspace{1cm} (s\hspace{0.02cm},a) \in S \times \bar{Q}_+
\end{equation}
where $\beta(s\hspace{0.02cm},a) = \mathrm{Ad}(s)\,a$ and $\bar{Q}_+$ is the closure of $Q_+\hspace{0.02cm}$. This mapping $\varphi$ is onto $M$, and is a diffeomorphism of $S \times Q_+$ onto its image $M_r\hspace{0.03cm}$, which is also the set of regular values of $\varphi$. Finally,
\begin{equation} \label{eq:cutysccss}
\mathrm{Cut}(y) = \varphi(S\times \bar{Q}_{\pi})  \hspace{0.5cm} \text{where } \bar{Q}_{\pi} = \bar{Q}_+ \,\cap\,(\cup_{\ell}\,\lbrace a: c_{\ell}(a) = \pi \rbrace)
\end{equation} 

\subsection{The squared distance function} \label{subsec:ssbishess}
For $x \in M$, consider the squared distance function $f_x(y) =  d^{\hspace{0.03cm}2}(x,y)/2$. If $x \notin \mathrm{Cut}(y)$, then $f_x$ is $C^2$ near $y$ (this is because $y \in \mathrm{Cut}(x)$ if and only if $x \in \mathrm{Cut}(y)$). %Since $M$ is simply connected, $\mathrm{Cut}(y)$ is given by (\ref{eq:cutysccss}). 

In this case, write $x = \varphi(s\hspace{0.02cm},a)$, where the map $\varphi$ was defined in (\ref{eq:varphibis}). Let $G_y(x)$ and $H_y(x)$ denote the gradient and Hessian of $f_x$ at $y$. These are given by
\begin{equation} \label{eq:ssgradfy}
  G_y(x) = - \beta(s,a) \hspace{3.3cm}
\end{equation}
\begin{equation} \label{eq:sshessfy}
  H_y(x) = \Pi^s_\mathfrak{a} \,+\, \sum_{\lambda \in \Delta_+} \lambda(a)\cot\lambda(a)\;\Pi^s_{\lambda}
\end{equation}
in the notation of (\ref{eq:compactAA}). Here, (\ref{eq:ssgradfy}) follows from (\ref{eq:gradfx}), since $x = \mathrm{Exp}_y\hspace{0.03cm}(\beta(s\hspace{0.02cm},a))$,  and (\ref{eq:sshessfy}) follows from the solution of the Ricatti equation (\ref{eq:ricatti}), discussed in \ref{subsec:rootsjacobis}.

If $M$ is simply connected, then $\mathrm{Cut}(y)$ is given by (\ref{eq:cutysccss}). Now, if $x \in \mathrm{Cut}(y)$ is written $x = \varphi(s\hspace{0.02cm},a)$, then $\lambda(a) = \pi$ for some $\lambda \in \Delta_+\hspace{0.02cm}$ ($\lambda = c_{\ell}$ which achieves the minimum in (\ref{eq:tc2})). By (\ref{eq:sshessfy}), $H_y(x)$ then has an eigenvalue equal to $-\infty$. In other words, $H_y(x)$ blows up when $x$ approaches $\mathrm{Cut}(y)$.

The convexity radius of a simply connected compact Riemannian symmetric space $M$ \hfill\linebreak is equal to half its injectivity radius. Accordingly, the convexity radius of $M$ is $r_{\scriptscriptstyle cx} = (\pi/2)\hspace{0.03cm}c^{\scriptscriptstyle -1}$. The proof of this statement may be summarised in the following way\,:

If $\delta < r_{\scriptscriptstyle cx\hspace{0.02cm}}$, then any $y_{\scriptscriptstyle 1\hspace{0.02cm}},y_{\scriptscriptstyle 2}$ in $B(x\hspace{0.02cm},\delta)$ must have $d(y_{\scriptscriptstyle 1\hspace{0.02cm}},y_{\scriptscriptstyle 2}) <  \pi\hspace{0.03cm}c^{\scriptscriptstyle -1}$, the injectivity radius of $M$,\hfill\linebreak and are therefore connected by a unique length-minimising geodesic curve $\gamma$. But, by (\ref{eq:sshessfy}), the squared distance function $f_x$ is convex on $B(x\hspace{0.02cm},\delta)$, where all eigenvalues of its Hessian are greater than $c\delta\hspace{0.02cm}\cot(c\delta) > 0$. This can be used to show that the geodesic $\gamma$ lies entirely in $B(x\hspace{0.02cm},\delta)$~\cite{petersen} (Page 177). In other words, the geodesic ball $B(x\hspace{0.02cm},\delta)$ is convex. On the other hand~\cite{sakai}, if $\delta = r_{\scriptscriptstyle cx}$ then there exists a closed (i.e. periodic) geodesic, of length $2\pi\hspace{0.03cm}c^{\scriptscriptstyle -1}$, contained in $B(x\hspace{0.02cm},\delta)$, so that this geodesic ball cannot be convex.

\subsection{A different integral formula} \label{subsec:ssbisintegral}
Consider again the map $\varphi$, defined in (\ref{eq:varphibis}). Let $M_r$ denote the set of regular values of $\varphi$. By Sard's lemma~\cite{bogachev}, the complement of $M_r$ in $M$ has zero Riemannian volume. Therefore,\hfill\linebreak if $f : M\rightarrow \mathbb{R}$ is a measurable function,
\begin{equation}\label{eq:ssbisintegral1}
\int_M\,f(x)\hspace{0.03cm}\mathrm{vol}(dx) = \int_{M_r}\,f(x)\hspace{0.03cm}\mathrm{vol}(dx)
\end{equation}
However, it was seen in \ref{subsec:ssbiscut} that $\varphi$ is a diffeomorphism of $S \times Q_+$ onto $M_r\hspace{0.02cm}$. Then, performing a ``change of variables", it follows that
\begin{equation} \label{eq:ssbisintegral2}
\int_M\,f(x)\hspace{0.03cm}\mathrm{vol}(dx) =
\int_{Q_+}\!\int_S\,f(s\hspace{0.02cm},a)\hspace{0.02cm}D(a)\hspace{0.03cm}da\hspace{0.02cm}\omega(ds)
\end{equation}
where $f(s\hspace{0.02cm},a) = f(\varphi(s\hspace{0.02cm},a))$ and $\varphi^*(\mathrm{vol}) = D(a)\hspace{0.03cm}da\hspace{0.02cm}\omega(ds)$. In particular, the ``volume density" $D(a)$ can be read from (\ref{eq:ssvolkac}),
\begin{equation} \label{eq:Dbis}
  D(a) = \prod_{\lambda \in \Delta_+} \left| \sin\hspace{0.02cm}\lambda(a)\right|^{ m_\lambda} \hspace{2.4cm}
\end{equation}
where the absolute value may be dropped, whenever $a \in Q_+$ is understood from the context. \\[0.1cm]
\textbf{Remark\,:} the integral formula (\ref{eq:ssbisintegral2}) is somewhat similar to (\ref{eq:ssvolkac2}). Roughly, both formulae involve the same change of variables, but (\ref{eq:ssbisintegral2}) takes advantage of the the description of the cut locus of $y$ in (\ref{eq:cutysccss}). Of course, (\ref{eq:ssbisintegral2}) only works when the compact symmetric space $M$ is simply connected.

%%%what is \Vert

\vfill
\pagebreak

\section{All the proofs} \label{sec:bcentreproofs}
Throughout the following proofs, it will be assumed that $U(x^*) = 0$. There is no loss of generality in making this assumption. Indeed, looking back at the definition (\ref{eq:gibbs}) of the Gibbs distribution $\pi_{\scriptscriptstyle T}\hspace{0.03cm}$, it is clear that a factor $\exp(-U(x^*)/T)$ may always be absorbed into $Z(T)$.
%\pi_{\scriptscriptstyle T}
%d^{\hspace{0.03cm}2}(y\hspace{0.02cm},x)
%  \mathcal{E}_{\scriptscriptstyle T}(y)
% \hat{x}_{\scriptscriptstyle T}
\subsection{Proof of Proposition \ref{prop:concentration}}

\subsubsection{Proof of (i)} 
For each $y \in M$, let $f_y(x) = d^{\hspace{0.03cm}2}(y\hspace{0.02cm},x)/2$. It follows from (\ref{eq:ETT}) that
\begin{equation} \label{proofconcentration11}
  \mathcal{E}_{\scriptscriptstyle T}(y) = \int_M\, f_y(x)\hspace{0.03cm}\pi_{\scriptscriptstyle T}(dx)
\end{equation}
On the other hand, consider the function $\mathcal{E}_{\scriptscriptstyle 0}(y)$,
\begin{equation} \label{eq:E0}
\mathcal{E}_{\scriptscriptstyle 0}(y) = \int_M\, f_y(x)\hspace{0.03cm}\delta_{x^*}(dx) = d^{\hspace{0.03cm}2}(y\hspace{0.02cm},x^*)/2
\end{equation}
For any $y \in M$, it is elementary that $f_y(x)$ is a Lipschitz function of $x$, with Lipschitz constant $\mathrm{diam}\,M$. Then, from the Kantorovich-Rubinshtein formula~\cite{kantorovich} (see VIII.4)
\begin{equation} \label{eq:proofconcentration12}
  \left| \mathcal{E}_{\scriptscriptstyle T}(y) - \mathcal{E}_{\scriptscriptstyle 0}(y)\right| \leq (\mathrm{diam}\,M)\hspace{0.02cm}W(\pi_{\scriptscriptstyle T}\hspace{0.02cm},\delta_{x^*})
\end{equation}
a uniform bound in $y \in M$. It now follows that, for any $\eta > 0$, 
\begin{equation} \label{eq:proofconcentration13}
  \inf_{y \in B(x^*\!,\eta)}   \mathcal{E}_{\scriptscriptstyle T}(y) -   \inf_{y \in B(x^*\!,\eta)}  \mathcal{E}_{\scriptscriptstyle 0}(y) \leq (\mathrm{diam}\,M)\hspace{0.02cm}W(\pi_{\scriptscriptstyle T}\hspace{0.02cm},\delta_{x^*}) 
\end{equation}
\begin{equation} \label{eq:proofconcentration14}
  \inf_{y \notin B(x^*\!,\eta)}   \mathcal{E}_{\scriptscriptstyle 0}(y) -   \inf_{y \notin B(x^*\!,\eta)}  \mathcal{E}_{\scriptscriptstyle T}(y) \leq (\mathrm{diam}\,M)\hspace{0.02cm}W(\pi_{\scriptscriptstyle T}\hspace{0.02cm},\delta_{x^*})
\end{equation}
However, from (\ref{eq:E0}), it is clear that
$$
\inf_{y \in B(x^*\!,\eta)}  \mathcal{E}_{\scriptscriptstyle 0}(y)  = 0 \hspace{0.25cm}\text{and}\hspace{0.2cm}
\inf_{y \notin B(x^*\!,\eta)}   \mathcal{E}_{\scriptscriptstyle 0}(y) = \frac{\eta^2}{2}
$$
To complete the proof, replace these into (\ref{eq:proofconcentration13}) and (\ref{eq:proofconcentration14}), and assume the condition in (\ref{eq:concentration1}) is verified. It then follows that
\begin{equation} \label{eq:proofconcentration15}
  \inf_{y \in B(x^*\!,\eta)}   \mathcal{E}_{\scriptscriptstyle T}(y) < \frac{\eta^2}{4} <   \inf_{y \notin B(x^*\!,\eta)}   \mathcal{E}_{\scriptscriptstyle T}(y)
\end{equation}
However, this means any global minimum of $\mathcal{E}_{\scriptscriptstyle T}(y)$ must belong to $B(x^*\!,\eta)$. Equivalently, any Riemannian barycentre $\hat{x}_{\scriptscriptstyle T}$ of $\pi_{\scriptscriptstyle T}$ must verify $d(\hat{x}_{\scriptscriptstyle T}\hspace{0.02cm},x^*) < \eta$. Thus, the conclusion in (\ref{eq:concentration1}) holds.

\subsubsection{Proof of (ii)}
Recall the condition in (\ref{eq:wellbehaved}), which holds for $d(x\hspace{0.02cm},x^*) \leq \rho$. By choosing $\rho < \min\lbrace \mathrm{inj}(x^*),\frac{\pi}{2}c^{\scriptscriptstyle -1}\rbrace$, it will be possible to apply (\ref{eq:integralcomp}) from \ref{subsec:volcomp}, in the remainder of the proof. Consider the truncated distribution
\begin{equation} \label{eq:pitruncate}
 \pi^{\scriptscriptstyle \rho}_{\scriptscriptstyle T}(dx) \,= \frac{\mathbf{1}_{\scriptscriptstyle B_\rho}(x)}{\pi_{\scriptscriptstyle T}(B_{\scriptscriptstyle \rho})}\hspace{0.03cm}   \pi^{\phantom{\scriptscriptstyle \rho}}_{\scriptscriptstyle T}(dx)
\end{equation}
where $\mathbf{1}$ denotes the indicator function, and $B_{\scriptscriptstyle \rho}$ denotes the open ball $B(x^*\!,\rho)$.
Of course, by the triangle inequality
\begin{equation} \label{eq:kantortriangle}
  W(\pi_{\scriptscriptstyle T}\hspace{0.02cm},\delta_{x^*}) \leq 
W(\pi^{\phantom{\scriptscriptstyle \rho}}_{\scriptscriptstyle T}\hspace{0.02cm},\pi^{\scriptscriptstyle \rho}_{\scriptscriptstyle T}) + 
W(\pi^{\scriptscriptstyle \rho}_{\scriptscriptstyle T}\hspace{0.02cm},\delta_{x^*})
\end{equation}
Now, the proof relies on the following estimates, which use the notation of \ref{sec:twtd}. \\[0.1cm]
\textbf{-- first estimate\,:} if $T\leq T^1_{\scriptscriptstyle W\hspace{0.03cm}}$, then
\begin{equation} \label{eq:estimate1}
 W(\pi^{\phantom{\scriptscriptstyle \rho}}_{\scriptscriptstyle T}\hspace{0.02cm},\pi^{\scriptscriptstyle \rho}_{\scriptscriptstyle T}) \leq (\mathrm{diam}\,M\times\mathrm{vol}\,M)\left(\frac{2}{\pi}\right)\left(\frac{\pi}{8}\right)^{\!\frac{n}{2}}\left(\frac{\mu_{\max}}{T}\right)^{\!\frac{n}{2}}\exp\left(-\frac{U_\rho}{T}\right)
\end{equation} 
\textbf{-- second estimate\,:}  if $T\leq T^1_{\scriptscriptstyle W\hspace{0.03cm}}$, then
\begin{equation} \label{eq:estimate2}
W(\pi^{\scriptscriptstyle \rho}_{\scriptscriptstyle T}\hspace{0.02cm},\delta_{x^*}) \leq 
(2\pi)^{\!\frac{1}{2}}
\hspace{0.02cm}B^{-1}_n
\left(\frac{\pi}{2}\right)^{\!n-1}
\left(\frac{\mu_{\max}}{\mu_{\min}}\right)^{\!\!\frac{n}{2}}
\left(\frac{T}{\mu_{\min}}\right)^{\!\!\frac{1}{2}}
\end{equation}
These two estimates will be proved below. To obtain (\ref{eq:concentration2}), assume that they hold, and that $T\leq T_{\scriptscriptstyle W\hspace{0.03cm}}$. Then, $T \leq T^2_{\scriptscriptstyle W}$ and the definition of $T^2_{\scriptscriptstyle W}$ implies
$$
f(T,n+1,\rho) \leq \left(\mu_{\max}\middle/\mu_{\min}\right)^{\!\scriptscriptstyle \frac{n}{2}}\hspace{0.01cm}C_n
$$
Using the definition of $C_n$ and formulae (\ref{eq:gammastuff}), this inequality reads
$$
(\mathrm{diam}\,M\times\mathrm{vol}\,M)\hspace{0.02cm}f(T,n+1,\rho) 
\leq 2\hspace{0.03cm}
(2\pi)^{\!\frac{n}{2}}\hspace{0.02cm}B^{-1}_n
\left(\mu_{\max}\middle/\mu_{\min}\right)^{\!\scriptscriptstyle \frac{n}{2}}
$$
This is the same as
$$
(\mathrm{diam}\,M\times\mathrm{vol}\,M)\hspace{0.02cm}\pi^{\scriptscriptstyle -1}
\left(\frac{\pi}{8}\right)^{\!\frac{n}{2}}f(T,n+1,\rho) \leq 
\left(\frac{\pi}{2}\right)^{\!n-1}
B^{-1}_n
\left(\mu_{\max}\middle/\mu_{\min}\right)^{\!\scriptscriptstyle \frac{n}{2}}
$$
From the definition of $f(T,n+1,\rho)$, it then follows that the right-hand side of (\ref{eq:estimate1}) is less than half the right-hand side of (\ref{eq:estimate2}). Since this is the case, (\ref{eq:concentration2}) follows from the triangle inequality (\ref{eq:kantortriangle}). \\[0.1cm]
\textbf{-- proof of first estimate\,:} consider the probability distribution $K$ on $M \times M$,
\begin{equation} \label{eq:coupling}
 K(dx_{\scriptscriptstyle 1}\times
   dx_{\scriptscriptstyle 2}) =\pi^{\scriptscriptstyle \rho}_{\scriptscriptstyle T}(dx_{\scriptscriptstyle 1})\left[ \pi_{\scriptscriptstyle T}(B_{\scriptscriptstyle \rho})\hspace{0.02cm}\delta_{x_{\scriptscriptstyle 1}}(dx_{\scriptscriptstyle 2}) + \mathbf{1}_{\scriptscriptstyle B^{\scriptscriptstyle c}_{\scriptscriptstyle \rho}}(x_{\scriptscriptstyle 2})\pi_{\scriptscriptstyle T}(dx_{\scriptscriptstyle 2})\right]
\end{equation}
where $B^{\scriptscriptstyle c}_{\scriptscriptstyle \rho}$ denotes the complement of $B_{\scriptscriptstyle \rho}$ in $M$. This distribution $K$ provides a coupling between $\pi^{\phantom{\scriptscriptstyle \rho}}_{\scriptscriptstyle T}$ and $\pi^{\scriptscriptstyle \rho}_{\scriptscriptstyle T\hspace{0.03cm}}$. Therefore, replacing (\ref{eq:coupling}) into the definition of the Kantorovich distance, it follows 
\begin{equation} \label{eq:proofestimate11}
W(\pi_{\scriptscriptstyle T}\hspace{0.02cm},\pi^{\scriptscriptstyle \rho}_{\scriptscriptstyle T}) \leq (\mathrm{diam}\,M)\hspace{0.02cm} 
\pi_{\scriptscriptstyle T}(B^{\scriptscriptstyle c}_{\scriptscriptstyle \rho})
\end{equation}
However, the definition (\ref{eq:gibbs}) of $\pi_{\scriptscriptstyle T}$ implies
\begin{equation} \label{eq:proofestimate12}
 \pi_{\scriptscriptstyle T}(B^{\scriptscriptstyle c}_{\scriptscriptstyle \rho}) \leq \left(Z(T)\right)^{-1}(\mathrm{vol}\,M)\exp\left(-\frac{U_\rho}{T}\right)
\end{equation}
Now, (\ref{eq:estimate1}) follows directly from (\ref{eq:proofestimate11}) and (\ref{eq:proofestimate12}), if the following lower bound on $Z(T)$ can be proved
\begin{equation} \label{eq:lowerz}
  Z(T) \geq \left(\frac{\pi}{2}\right)\left(\frac{8}{\pi}\right)^{\!\frac{n}{2}}\left(\frac{T}{\mu_{\max}}\right)^{\!\frac{n}{2}} \hspace{1cm} \text{for } T \leq T^1_{\scriptscriptstyle W}
\end{equation}
To prove this lower bound, note that
$$
Z(T) = \int_M\,\exp\left(-\frac{U(x)}{T}\right)\mathrm{vol}(dx) \,\geq \int_{B_{\scriptscriptstyle \rho}}\,\exp\left(-\frac{U(x)}{T}\right)\mathrm{vol}(dx)
$$
Replacing (\ref{eq:wellbehaved}) into this last inequality, it is possible to write
\begin{equation} \label{eq:prooflowerz1}
Z(T) \geq \int_{B_{\scriptscriptstyle \rho}}\,\exp\left(-\frac{U(x)}{T}\right)\mathrm{vol}(dx) \geq
\int_{B_{\scriptscriptstyle \rho}}\,\exp\left(-\frac{\mu_{\max}}{2T}d^{\hspace{0.03cm}2}(x,x^*)\right)\mathrm{vol}(dx)
\end{equation}
Since $\rho < \min\lbrace \mathrm{inj}(x^*),\frac{\pi}{2}c^{\scriptscriptstyle -1}\rbrace$, it is possible to apply (\ref{eq:integralcomp}) from \ref{subsec:volcomp}, to (\ref{eq:prooflowerz1}). Specifically, the lower bound in (\ref{eq:integralcomp}) yields,
\begin{equation} \label{eq:prooflowerz2}
  Z(T) \geq 
\omega_{n-1}\,\int^{\rho}_{\scriptscriptstyle 0}\,e^{-\frac{\mu_{\max}}{2T}r^2}\left(c^{\scriptscriptstyle -1}\!\sin(c\hspace{0.02cm}r)\right)^{n-1}\hspace{0.02cm}dr \geq
\omega_{n-1}(2/\pi)^{n-1}\,\int^{\rho}_{\scriptscriptstyle 0}\,e^{-\frac{\mu_{\max}}{2T}r^2}\hspace{0.02cm}r^{n-1}\hspace{0.02cm}dr
\end{equation}
where the second inequality follows since $\sin(t)$ is a concave function of $t \in [0,\pi/2]$, so that $
\sin(c\hspace{0.02cm}r) \geq (2/\pi)c\hspace{0.02cm}r$ for $r \in [0,\rho]$. Now, the required bound (\ref{eq:lowerz}) follows from (\ref{eq:prooflowerz2}) by noting 
$$
\int^{\rho}_{\scriptscriptstyle 0}\,e^{-\frac{\mu_{\max}}{2T}r^2}r^{n-1}\hspace{0.02cm}dr = 
(2\pi)^{\frac{1}{2}}\left(\frac{T}{\mu_{\max}}\right)^{\!\frac{n}{2}}A_{n-1} - 
\int^{\scriptscriptstyle \infty}_{\rho}\,e^{-\frac{\mu_{\max}}{2T}r^2}r^{n-1}\hspace{0.02cm}dr
$$
where $A_n = \mathbb{E}|X|^n$ for $X \sim N(0,1)$, and that
$$
\int^{\scriptscriptstyle \infty}_{\rho}\,e^{-\frac{\mu_{\max}}{2T}r^2}r^{n-1}\hspace{0.02cm}dr \leq \frac{\rho^{n-2\hspace{0.02cm}}T}{\mu_{\max}}
\,e^{-\frac{\mu_{\max}}{2T}\rho^2}\leq 
\frac{\rho^{n-2\hspace{0.02cm}}T}{\mu_{\max}}
\,e^{-\frac{U_\rho}{T}}
$$
Indeed, taken together, these give
\begin{equation} \label{eq:prooflowerz3}
  Z(T) \geq \omega_{n-1}(2/\pi)^{n-1}\,\left[ 
(2\pi)^{\frac{1}{2}}\left(\frac{T}{\mu_{\max}}\right)^{\!\frac{n}{2}}A_{n-1} -
\frac{\rho^{n-2\hspace{0.02cm}}T}{\mu_{\max}}
\,e^{-\frac{U_\rho}{T}}
\right]
\end{equation}
Then, (\ref{eq:lowerz}) can be obtained by noting that the second term in square brackets is negligible in comparison to the first as  $T\rightarrow 0$, and using formulae (\ref{eq:gammastuff}) for $A_{n-1}$ and $\omega_{n-1\hspace{0.03cm}}$. \\[0.1cm]
\textbf{-- proof of second estimate\,:} the Kantorovich distance $W(\pi^{\scriptscriptstyle \rho}_{\scriptscriptstyle T}\hspace{0.02cm},\delta_{x^*})$ between $\pi^{\scriptscriptstyle \rho}_{\scriptscriptstyle T}$ and $\delta_{x^*}$ is equal to the first-order moment
$$
\int_M\,d(x\hspace{0.02cm},x^*)\hspace{0.02cm}
\pi^{\scriptscriptstyle \rho}_{\scriptscriptstyle T}(dx)
$$
According to (\ref{eq:gibbs}) and (\ref{eq:pitruncate}), this means
$$
W(\pi^{\scriptscriptstyle \rho}_{\scriptscriptstyle T}\hspace{0.02cm},\delta_{x^*}) =
\left(\pi_{\scriptscriptstyle T}(B_{\scriptscriptstyle \rho})Z(T)\right)^{-1}\,\int_{B_{\scriptscriptstyle \rho}}
d(x\hspace{0.02cm},x^*)\hspace{0.02cm}\exp\left(-\frac{U(x)}{T}\right)\mathrm{vol}(dx)
$$
Using (\ref{eq:gibbs}) to express the probability in parentheses, this becomes
\begin{equation} \label{eq:proofestimate21}
W(\pi^{\scriptscriptstyle \rho}_{\scriptscriptstyle T}\hspace{0.02cm},\delta_{x^*}) = \frac{\int_{B_{\scriptscriptstyle \rho}}
d(x\hspace{0.02cm},x^*)\hspace{0.02cm}\exp\left(-\frac{U(x)}{T}\right)\mathrm{vol}(dx)}{
\int_{B_{\scriptscriptstyle \rho}}\exp\left(-\frac{U(x)}{T}\right)\mathrm{vol}(dx)}
\end{equation}
A lower bound on the denominator can be found from (\ref{eq:prooflowerz1}) and the subsequent inequalities, which were used to prove (\ref{eq:lowerz}). These inequalities provide
\begin{equation} \label{eq:proofestimate22}
\int_{B_{\scriptscriptstyle \rho}}\exp\left(-\frac{U(x)}{T}\right)\mathrm{vol}(dx) \geq
 \omega_{n-1}\hspace{0.2cm}
(2/\pi)^{n-1}
(2\pi)^{\!\frac{1}{2}}
\hspace{0.02cm}A_{n-1}
\left(T\middle/\mu_{\max}\right)^{\!\frac{n}{2}}
\end{equation}
whenever $T \leq T^1_{\scriptscriptstyle W}$. For the numerator in (\ref{eq:proofestimate21}), it will be shown that
\begin{equation} \label{eq:proofestimate23}
\int_{B_{\scriptscriptstyle \rho}}
d(x\hspace{0.02cm},x^*)\hspace{0.02cm}\exp\left(-\frac{U(x)}{T}\right)\mathrm{vol}(dx)
 \leq
 \omega_{n-1}
(2\pi)^{\!\frac{1}{2}}
\hspace{0.02cm}A_{n}
\left(T\middle/\mu_{\min}\right)^{\!\frac{n+1}{2}}
\end{equation} 
Then, (\ref{eq:estimate2}) follows by dividing (\ref{eq:proofestimate23}) by (\ref{eq:proofestimate22}), and replacing in (\ref{eq:proofestimate21}), using the fact that $A_n/A_{n-1} = (2\pi)^{\!\frac{1}{2}}B^{\scriptscriptstyle -1}_n\hspace{0.03cm}$, which can be found from formulae (\ref{eq:gammastuff}). Now, it only remails to prove (\ref{eq:proofestimate23}). This is done by noting, from (\ref{eq:wellbehaved}),
$$
\int_{B_{\scriptscriptstyle \rho}}
d(x\hspace{0.02cm},x^*)\hspace{0.02cm}\exp\left(-\frac{U(x)}{T}\right)\mathrm{vol}(dx)
 \leq
\int_{B_{\scriptscriptstyle \rho}}
d(x\hspace{0.02cm},x^*)\hspace{0.02cm}
\exp\left(-\frac{\mu_{\min}}{2T}d^{\hspace{0.03cm}2}(x,x^*)\right)\mathrm{vol}(dx)
$$
Applying the upper bound in (\ref{eq:integralcomp}) (with $\kappa_{\min} = 0$), to the last integral, it follows that 
$$
\int_{B_{\scriptscriptstyle \rho}}
d(x\hspace{0.02cm},x^*)\hspace{0.02cm}
\exp\left(-\frac{\mu_{\min}}{2T}d^{\hspace{0.03cm}2}(x,x^*)\right)\mathrm{vol}(dx) \leq
\omega_{n-1}\int^{\rho}_{\scriptscriptstyle 0}
e^{-\frac{\mu_{\min}}{2T}r^2}\hspace{0.02cm}r^n\hspace{0.02cm}dr \leq 
\omega_{n-1}\int^{\scriptscriptstyle \infty}_{\scriptscriptstyle 0}
e^{-\frac{\mu_{\min}}{2T}r^2}\hspace{0.02cm}r^n\hspace{0.02cm}dr
$$ 
The integral on the right-hand side is half the $n$-th absolute moment of a normal distribution. By expressing it in terms of $A_n\hspace{0.03cm}$, it is possible to directly recover (\ref{eq:proofestimate23}).

\subsection{Proof of Proposition \ref{prop:differentiability}}

\subsubsection{Proof of  (i)}
Under the integrals in (\ref{eq:GH}), $G_y(x)$ and $H_y(x)$ are given by (\ref{eq:ssgradfy}) and (\ref{eq:sshessfy}), for any $x \in \mathrm{D}(y)$. Furthermore, by (\ref{eq:ssgradfy}) and (\ref{eq:sshessfy}), both integrands $G_y(x)$ and $H_y(x)$ are continuous on the domaine of integration  $\mathrm{D}(y)$.

The integral $G_y$ converges, because $G_y(x)$ is uniformly bounded on $\mathrm{D}(y)$. Indeed, from (\ref{eq:ssgradfy}),
$$
\Vert G_y(x)\Vert_y = \Vert \beta(s,a)\Vert_y = d(y\hspace{0.02cm},x)
$$
where the second equality follows from the fact that $x = \mathrm{Exp}_y(\beta(s,a))$. Of course, $d(y\hspace{0.02cm},x)$ is always less than $\mathrm{diam}\,M$. 

The integral $H_y$ is an improper integral, since $H_y(x)$ blows up when $x$ approaches $\mathrm{Cut}(y)$, as explained in \ref{subsec:ssbishess}. Nonetheless, this integral converges absolutely, as shall be seen from the material in \ref{subsec:ssbisintegral}.

Precisely, recall the mapping $\varphi$ defined in (\ref{eq:varphibis}). Because $M$ is simply connected, $\mathrm{Cut}(y)$ is identical to the first conjugate locus of $y$. This means that $\mathrm{Cut}(y)$ is contained in the set of critical values of $\mathrm{Exp}_y\hspace{0.03cm}$, and therefore also in the set of critical values of $\varphi$. Equivalently, $\mathrm{D}(y)$ contains the set of regular values of $\varphi$, denoted $M_r$ in (\ref{eq:ssbisintegral1}). It then follows, as in (\ref{eq:ssbisintegral2}),
\begin{equation} \label{eq:Hypolar}
  H_y = \int_{Q_+}\!\int_S\,H_y(s\hspace{0.02cm},a)\hspace{0.02cm}p_{\scriptscriptstyle T}(s\hspace{0.02cm},a)\hspace{0.02cm}D(a)\hspace{0.03cm}da\hspace{0.02cm}\omega(ds)
\end{equation}
where $p_{\scriptscriptstyle T}$ denotes the density of $\pi_{\scriptscriptstyle T}$ with respect to the Riemannian volume. 

To prove that $H_y$ converges absolutely, it is enough to prove the integrand in (\ref{eq:Hypolar}) is uniformly bounded. However, the density $p_{\scriptscriptstyle T}$ is bounded, since it is continuous and $M$ is compact. Moreover, it is clear from (\ref{eq:sshessfy}) and (\ref{eq:Dbis}) that
\begin{equation} \label{eq:HypolarH}
  H_y(s\hspace{0.02cm},a) = \Pi^s_\mathfrak{a} \,+\, \sum_{\lambda \in \Delta_+} \lambda(a)\cot\lambda(a)\;\Pi^s_{\lambda}
\end{equation}
\begin{equation} \label{eq:HypolarD}
  D(a) = \prod_{\lambda \in \Delta_+} \left| \sin\hspace{0.02cm}\lambda(a)\right|^{ m_\lambda} \hspace{2.4cm}
\end{equation}
The product of these two expressions is uniformly bounded, because $\lambda(a) \in (0,\pi)$ on $Q_+$.

Thus, the integrals $G_y$ and $H_y$ converge, and it is clear from the above that this is true for any temperature $T > 0$. The fact that both $G_y$ and $H_y$ depend continuously on $y$ will be clear from the arguments in the proof of (ii).

\subsubsection{Proof of (ii)}
The proof relies in a crucial way on Lemma \ref{lem:hddim}, which is proved in \ref{subsec:proofhddim}, below. To compute the gradient and Hessian of $\mathcal{E}_{\scriptscriptstyle T}$ at $y \in M$, consider any geodesic $\gamma:I\rightarrow M$, defined on a compact interval $I = [-\tau,\tau]$, with $\gamma(0) = y$. For each $t \in I$, it is immediate from (\ref{eq:ETT}) that
\begin{equation} \label{eq:proofderivatives1}
\mathcal{E}_{\scriptscriptstyle T}(\gamma(t)) = \int_M\,f_x(\gamma(t))\hspace{0.03cm}\pi_{\scriptscriptstyle T}(dx) 
\end{equation}
However, Lemma \ref{lem:hddim} states that the set 
$$
\mathrm{Cut}(\gamma) = \bigcup_{t \in I}\,\mathrm{Cut}(\gamma(t))
$$
has Riemannian volume equal to zero. Thus, since $\pi_{\scriptscriptstyle T}$ is uniformly continuous with respect to Riemannian volume, $\mathrm{Cut}(\gamma)$ can be removed from the domain of integration in (\ref{eq:proofderivatives1}), to obtain
\begin{equation} \label{eq:proofderivatives2}
\mathcal{E}_{\scriptscriptstyle T}(\gamma(t)) =
\int_{\mathrm{D}(\gamma)}f_x(\gamma(t))\hspace{0.03cm}\pi_{\scriptscriptstyle T}(dx) \hspace{0.5cm} \text{for all $t \in I$}
\end{equation}
where $\mathrm{D}(\gamma) = M - \mathrm{Cut}(\gamma)$. Now, if $x \in \mathrm{D}(\gamma)$, then $x \notin \mathrm{Cut}(\gamma(t))$ for any $t \in I$. According to \ref{subsec:ssbishess}, this implies that $f_x(\gamma(t))$ is $C^2$ near each $t \in I$. In other words, $f_x(t) = f_x(\gamma(t))$ is $C^2$ inside the interval $I$. Then, the first and second derivatives of this function are given by
\begin{equation} \label{eq:fprimeprimef}
  f^\prime_x(t) =  \left\langle G_{\gamma(t)}(x)\hspace{0.02cm},\dot{\gamma}\right\rangle_{\gamma(t)} \hspace{0.3cm};\hspace{0.3cm} f^{\prime\prime}_x(t) = \left\langle H_{\gamma(t)}(x)\cdot \dot{\gamma}\hspace{0.02cm},\dot{\gamma}\right\rangle_{\gamma(t)}
\end{equation}
as in Proposition \ref{prop:geodesic}. Formally, (\ref{eq:derivatives}) follows by differentiating under the integral sign in (\ref{eq:proofderivatives2}), replacing from (\ref{eq:fprimeprimef}), and then putting $t = 0$. This differentiation under the integral sign is justified, as soon as it is shown that the families of functions,
$$
\left\lbrace x \mapsto G_{\gamma(t)}(x)\,;t\in I\right\rbrace \hspace{0.3cm};\hspace{0.3cm}
\left\lbrace x \mapsto H_{\gamma(t)}(x)\,;t\in I\right\rbrace
$$
which all have common domain of definition $\mathrm{D}(\gamma)$, are uniformly integrable with respect to the probability distribution $\pi_{\scriptscriptstyle T}$ (precisely, with respect to the restriction of $\pi_{\scriptscriptstyle T}$ to  $\mathrm{D}(\gamma)$). 

Roughly (for the exact definition, see~\cite{bogachev}), uniform integrability means that the rate of absolute convergence of the following integrals 
\begin{equation} \label{eq:GHt}
G_{\gamma(t)} = \int_{\mathrm{D}(\gamma)}G_{\gamma(t)}(x)\hspace{0.03cm}\pi_{\scriptscriptstyle T}(dx) \hspace{0.2cm};\hspace{0.2cm}
  H_{\gamma(t)} = \int_{\mathrm{D}(\gamma)}H_{\gamma(t)}(x)\hspace{0.03cm}\pi_{\scriptscriptstyle T}(dx)
\end{equation} 
does not depend on $t \in I$. This is clear for the integrals $G_{\gamma(t)}$ because, as in the proof of (i),
$$
\Vert G_{\gamma(t)}(x)\Vert_{\gamma(t)} = d(\gamma(t)\hspace{0.02cm},x)
$$
and this is bounded by $\mathrm{diam}\,M$, independently of $x$ and $t$.

Then, consider the integral $H_{\gamma(0)} = H_y\hspace{0.03cm}$, and recall Formulae (\ref{eq:Hypolar})--(\ref{eq:HypolarH}). For simplicity, assume that $M$ is an irreducible symmetric space (see Chapter VIII of~\cite{helgason}, Page 307). In this case, according to \ref{subsec:ssbiscut}, there exists a maximal root $c \in \Delta_+\hspace{0.03cm}$, so that $c(a) \geq \lambda(a)$ for all $\lambda \in \Delta_+$ and $a \in Q_+\hspace{0.02cm}$. Therefore, it follows from (\ref{eq:HypolarH}) that
\begin{equation} \label{eq:proofui1}
 \Vert H_y(x)\Vert_{\scriptscriptstyle F} \leq (\dim\,M)^{\!\frac{1}{2}}\hspace{0.03cm}\max\lbrace 1,|c(a)\cot c(a)|\rbrace
\end{equation}
where $\Vert_\cdot\Vert_{\scriptscriptstyle F}$ denotes the Frobenius norm given by the Riemannian metric of $M$. Now, the required uniform integrability is equivalent to the statement that
\begin{equation} \label{eq:uicondition}
  \lim_{K\rightarrow \infty}\,\int_{\mathrm{D}(\gamma)}\, \Vert H_y(x)\Vert_{\scriptscriptstyle F}\hspace{0.03cm} \mathbf{1}\lbrace  \Vert H_y(x)\Vert_{\scriptscriptstyle F}  > K\rbrace\hspace{0.03cm}\pi_{\scriptscriptstyle T}(dx) = 0 \;\;\text{uniformly in $y$}
\end{equation}
But from the inequality in (\ref{eq:proofui1}), if $K > 1$ there exists $\epsilon > 0$ such that
$$
\lbrace  \Vert H_y(x)\Vert_{\scriptscriptstyle F}  > K\rbrace \subset \lbrace c(a) > \pi - \epsilon\rbrace
$$
and $\epsilon \rightarrow 0$ as $K \rightarrow \infty$. In this case, the integral in (\ref{eq:uicondition}) is less than
\begin{equation} \label{eq:proofui2}
 (\dim\,M)^{\!\frac{1}{2}}\left( \sup\hspace{0.02cm} p_{\scriptscriptstyle T}(x)\right)\int_{\mathrm{D}(\gamma)}\,
|c(a)\cot c(a)|\hspace{0.03cm}\mathbf{1}\lbrace c(a) > \pi - \epsilon\rbrace\hspace{0.03cm}\mathrm{vol}(dx)
\end{equation}
By expressing this integral as in (\ref{eq:Hypolar}), it is seen to be equal to 
$$
\begin{array}{l}
\int_{Q_+}\!\int_S\,|c(a)\cot c(a)|\hspace{0.03cm}\mathbf{1}\lbrace c(a) > \pi - \epsilon\rbrace\hspace{0.02cm}D(a)\hspace{0.03cm}da\hspace{0.02cm}\omega(ds) \,= \\[0.2cm]
\omega(S)\,\int_{Q_+}\left[ |c(a)\cot c(a)|\hspace{0.02cm}D(a)\right]\mathbf{1}\lbrace c(a) > \pi - \epsilon\rbrace\hspace{0.02cm}da
\end{array}
$$
In view, of (\ref{eq:HypolarD}), since $c \in \Delta_+\hspace{0.02cm}$, the function in square brackets is bounded on the closure of $Q_+$ by $c^{\hspace{0.02cm}\scriptscriptstyle 2} = \Vert c \Vert^2$ (incidentally, this is the maximum sectional curvature of $M$, as explained in \ref{subsec:ssbiscut}). Finally, by (\ref{eq:proofui2}), the integral in (\ref{eq:uicondition}) is less than
$$
(\dim\,M)^{\!\frac{1}{2}}
\left( \sup\hspace{0.02cm} p_{\scriptscriptstyle T}(x)\right)
\omega(S)\hspace{0.02cm}c^{\hspace{0.02cm} \scriptscriptstyle 2}\,
\int_{Q_+}\mathbf{1}\lbrace c(a) > \pi - \epsilon\rbrace\hspace{0.02cm}da
$$
Recall that $c(a) \in (0,\pi)$ for $a \in Q_+\hspace{0.03cm}$. It is then clear this last integral converges to $0$ as $\epsilon \rightarrow 0$, at a rate which does not depend on $y$. This proves the required uniform integrability, so the proof is now complete, at least in the case where $M$ is an irreducible symmetric space.

In the general case, where $M$ is not irreducible, it is enough to note that, according to \ref{subsec:ssbiscut}, $M$ is a product of irreducible Riemannian symmetric spaces, $M = M_{\scriptscriptstyle 1}\times\ldots\times M_{ s\hspace{0.03cm}}$. Then, the proof boils down to the special case where $M$ is irreducible, as treated above.

\subsection{Proof of Lemma \ref{lem:hddim}} \label{subsec:proofhddim}
The proof uses the following general remark.\\[0.1cm]
\textbf{Remark\,:} let $M$ be a Riemannian manifold, and $g:M\rightarrow M$ be an isometry. Recall that $g\cdot y$ is used to denote $g(y)$, for $y \in M$. Similarly, if $A \subset M$, let $g\cdot A$ denote the image of $A$ under $g$.\hfill\linebreak Then, for any $y \in M$, $\mathrm{Cut}(g\cdot y) = g\cdot \mathrm{Cut}(y)$.
This is because a point $x \in M$ belongs to $\mathrm{Cut}(y)$, if and only if $x$ is a first conjugate point to $y$ along some geodesic, or there exist two different length-minimising geodesics connecting $y$ to $x$, and because both of these properties are preserved by any isometry $g$.\\[0.1cm]
\indent Assume $M$ is a simply connected compact Riemannian symmetric space. In the notation of \ref{sec:ssbis}, $M \simeq G/K$. Recall (by Proposition \ref{prop:geodesiclemma} of \ref{sec:geolemma}) any geodesic $\gamma:I\rightarrow M$ is given by
\begin{equation} \label{eq:prooflemgamma}
  \gamma(t) = \exp(t\hspace{0.02cm}\omega)\cdot y
\end{equation}
for some $y \in M$ and $\omega \in \mathfrak{p}$, where $\exp$ denotes the Lie group exponential. From the above remark, for each $t \in I$, the cut locus $\mathrm{Cut}(\gamma(t))$ of $\gamma(t)$ is given by
\begin{equation} \label{eq:cutgammat}
  \mathrm{Cut}(\gamma(t)) = \exp(t\hspace{0.02cm}\omega)\cdot\mathrm{Cut}(y)
\end{equation}
However, $\mathrm{Cut}(y)$ is described by (\ref{eq:cutysccss}) in \ref{subsec:ssbiscut}, which reads
\begin{equation} \label{eq:cutysccss1}
  \mathrm{Cut}(y) = \varphi(S\times \bar{Q}_{\pi})  
\end{equation}
in terms of the mapping $\varphi$ defined in (\ref{eq:varphibis}). It follows from (\ref{eq:cutgammat}) and (\ref{eq:cutysccss1}) that
\begin{equation} \label{eq:CUTgamma}
  \mathrm{Cut}(\gamma) = \Phi(I\times S\times \bar{Q}_{\pi}) \hspace{0.5cm} \Phi(t,s,a) = \exp(t\hspace{0.02cm}\omega)\cdot\varphi(s\hspace{0.02cm},a)
\end{equation} 
The aim is to show that this set has Hausdorff dimension strictly less than $\dim\,M$. This is done using results from dimension theory~\cite{hausdorff}. Precisely, note from (\ref{eq:cutysccss}) that
$$
\bar{Q}_{\pi} = \cup_{\ell}\, \bar{Q}_{\ell} \hspace{0.4cm} \text{where } \bar{Q}_{\ell} = \bar{Q}_+\,\cap\,\lbrace a:c_{\scriptscriptstyle \ell}(a) = \pi\rbrace
$$
Therefore, it is clear that
\begin{equation} \label{eq:CUTgammaunion}
 \mathrm{Cut}(\gamma) = \bigcup_{\ell} \Phi(I\times S\times \bar{Q}_{\ell})
\end{equation}
Then, it follows from~\cite{hausdorff} (Item (2) of Theorem 2) that
\begin{equation} \label{eq:CUTgammauniondim}
  \mathrm{dim}_{\scriptscriptstyle H}\,\mathrm{Cut}(\gamma) \leq \max_{\ell}\,\mathrm{dim}_{\scriptscriptstyle H}\,\Phi(I\times S\times \bar{Q}_{\ell})
\end{equation}
where $\mathrm{dim}_{\scriptscriptstyle H}$ is used to denote the Hausdorff dimension. Now, for each $\ell$,
$$
\Phi(I\times S\times \bar{Q}_{\ell}) = \Phi(I\times S_{\ell}\times \bar{Q}_{\ell}) \subset
\Phi(\mathbb{R} \times S_{\ell} \times \lbrace c_{\scriptscriptstyle \ell}(a) = \pi\rbrace)
$$
where $S_{\ell} = K/K_{\ell}$ with $K_{\ell}$ the centraliser of $\lbrace c_{\scriptscriptstyle \ell}(a) = \pi\rbrace$ in $K$. This last inclusion implies~\cite{hausdorff} (Item (1) of Theorem 2)
\begin{equation} \label{eq:diminclusion}
  \mathrm{dim}_{\scriptscriptstyle H}\,\Phi(I\times S\times \bar{Q}_{\ell}) \leq  \mathrm{dim}_{\scriptscriptstyle H}\,
\Phi(\mathbb{R} \times S_{\ell} \times \lbrace c_{\scriptscriptstyle \ell}(a) = \pi\rbrace)
\end{equation}
To conclude, note that the set $\mathbb{R} \times S_{\ell} \times \lbrace c_{\scriptscriptstyle \ell}(a) = \pi\rbrace$ is a differentiable manifold. Let this differentiable manifold be equipped with a product Riemannian metric (arising from flat metrics on $\mathbb{R}$ and $\lbrace c_{\scriptscriptstyle \ell}(a) = \pi\rbrace$, and from the invariant metric induced onto $S_{\ell}$ from $K$). It is clear from (\ref{eq:CUTgamma}) that $\Phi$ is smooth, and therefore locally Lipschitz. Then~\cite{hausdorff} (Item (5) of Theorem 2), 
\begin{equation} \label{eq:dimlipsch}
 \mathrm{dim}_{\scriptscriptstyle H}\,
\Phi(\mathbb{R} \times S_{\ell} \times \lbrace c_{\scriptscriptstyle \ell}(a) = \pi\rbrace) \leq
 \mathrm{dim}_{\scriptscriptstyle H}\,
(\mathbb{R} \times S_{\ell} \times \lbrace c_{\scriptscriptstyle \ell}(a) = \pi\rbrace)
\end{equation}
But the Hausdorff dimension of a Riemannian manifold is the same as its (usual) dimension. Accordingly, 
$$
 \mathrm{dim}_{\scriptscriptstyle H}\,
(\mathbb{R} \times S_{\ell} \times \lbrace c_{\scriptscriptstyle \ell}(a) = \pi\rbrace)
= 1 + \dim\,S_{\ell} + (\dim\,\mathfrak{a} - 1)
$$
since the dimension of a hyperplane in $\mathfrak{a}$ is $\dim\,\mathfrak{a} - 1$. In addition, from~\cite{helgason} (Page 253), $\dim\,S_{\ell} < \dim\,S$. Therefore,
$$
 \mathrm{dim}_{\scriptscriptstyle H}\,
(\mathbb{R} \times S_{\ell} \times \lbrace c_{\scriptscriptstyle \ell}(a) = \pi\rbrace) = \dim\,S_{\ell} + \dim\,\mathfrak{a} < \dim\,M
$$
since $\dim\,M = \dim\,S + \dim\,\mathfrak{a}$, as can be seen from (\ref{eq:varphibis}). Replacing this into (\ref{eq:dimlipsch}), it follows from (\ref{eq:CUTgammauniondim}) and (\ref{eq:diminclusion}) that $\dim\,\mathrm{Cut}(\gamma) < \dim\,M$. The lemma has therefore been proved.

\subsection{Proof of Proposition \ref{prop:unique}}
Assume that Lemma \ref{lem:etconv} is true. This lemma is proved in \ref{subsec:proofetconv}, below.

\subsubsection{Proof of (i)}
For $\delta < \frac{1}{2}r_{\scriptscriptstyle cx\hspace{0.03cm}}$, let $T_{\scriptscriptstyle \delta}$ be given by (\ref{eq:td}). By (ii) of Lemma \ref{lem:etconv}, $T < T_{\scriptscriptstyle \delta}$ implies the variance function $\mathcal{E}_{\scriptscriptstyle T}(y)$ is strongly convex on $B(x^*\!,\delta)$. It will be proved that any Riemannian barycentre $\hat{x}_{\scriptscriptstyle T}$ of $\pi_{\scriptscriptstyle T}$ belongs to $B(x^*\!,\delta)$. Then, since $\hat{x}_{\scriptscriptstyle T}$ is a minimum of $\mathcal{E}_{\scriptscriptstyle T}(y)$ in $B(x^*\!,\delta)$, it follows that $\hat{x}_{\scriptscriptstyle T}$ is unique (thanks to the strong convexity of $\mathcal{E}_{\scriptscriptstyle T}(y)$).

By (i) of Proposition \ref{prop:concentration}, to prove that any $\hat{x}_{\scriptscriptstyle T}$ belongs to $B(x^*\!,\delta)$, it is enough to prove
\begin{equation} \label{eq:bcproofunique1}
 W(\pi_{\scriptscriptstyle T\hspace{0.02cm}},\delta_{x^*}) < \frac{\delta^2}{4\hspace{0.02cm}\mathrm{diam}\,M}
\end{equation}
However, if $T < T_{\scriptscriptstyle \delta}$ then $T < T_{\scriptscriptstyle W}$ and, by (ii) of Proposition \ref{prop:concentration}, 
$W(\pi_{\scriptscriptstyle T\hspace{0.02cm}},\delta_{x^*})$ satisfies inequality (\ref{eq:concentration2}). In addition (from the definition of $T^1_{\scriptscriptstyle \delta}$ and $T^2_{\scriptscriptstyle \delta}$) one has $T < T^1_{\scriptscriptstyle \delta}$ and
$$
(2\pi)^{\!\frac{1}{2}}\hspace{0.02cm}(T/\mu_{\min})^{\!\frac{1}{2}} < \delta^2\hspace{0.03cm}(\mu_{\min}/\mu_{\max})^{\!\frac{n}{2}}\hspace{0.03cm}D_n
$$
By replacing the expression of $D_n$ and simplifying, this is the same as 
\begin{equation} \label{eq:bcproofunique2}
(2\pi)^{\!\frac{1}{2}}\hspace{0.02cm}
B^{-1}_n\hspace{0.02cm}(\pi/2)^{n-1}  
\hspace{0.03cm}(\mu_{\max}/\mu_{\min})^{\!\frac{n}{2}}
(T/\mu_{\min})^{\!\frac{1}{2}}
< \frac{\delta^2}{4\hspace{0.02cm}\mathrm{diam}\,M}
\end{equation}
Now, (\ref{eq:bcproofunique1}) follows from (\ref{eq:concentration2}) and (\ref{eq:bcproofunique2}).

\subsubsection{Proof of (ii)}
From the proof of (i), $\mathcal{E}_{\scriptscriptstyle T}(y)$ is strongly convex on $B(x^*\!,\delta)$, and $\hat{x}_{\scriptscriptstyle T}$ is the minimum of 
$\mathcal{E}_{\scriptscriptstyle T}(y)$ in $B(x^*\!,\delta)$. To prove that $\hat{x}_{\scriptscriptstyle T} = x^*$, it is then enough to prove that $x^*$ is a stationary point of $\mathcal{E}_{\scriptscriptstyle T}(y)$. However, the fact that $U$ is invariant by the geodesic symmetry $s_{x^*}$  will be seen to imply
\begin{equation} \label{eq:bcproofsx}
  d\hspace{0.02cm}s_{x^*}\cdot G_{x^*} = G_{x^*}
\end{equation} 
which is equivalent to $G_{x^*} = 0$, since $d\hspace{0.02cm}s_{x^*}$ is equal to minus the identity on $T_{x^*}M$ (see \ref{ssec:sspace}). Then, by (\ref{eq:derivatives}) in Proposition \ref{prop:differentiability}, $\mathrm{grad}\,\mathcal{E}_{\scriptscriptstyle T}(x^*) = 0$, so $x^*$ is indeed a stationary point of $\mathcal{E}_{\scriptscriptstyle T}(y)$.

To obtain (\ref{eq:bcproofsx}), it is possible to write from (\ref{eq:GH}),
\begin{equation} \label{eq:bcproofsx1}
  d\hspace{0.02cm}s_{x^*}\cdot G_{x^*} = d\hspace{0.02cm}s_{x^*}\cdot \int_{\mathrm{D}(x^*)}\,G_{x^*}(x)
\hspace{0.03cm}\pi_{\scriptscriptstyle T}(dx)
\end{equation}
But, it follows from (\ref{eq:varphibis}) and (\ref{eq:ssgradfy}) that $x = \mathrm{Exp}_{x^*}(-G_{x^*}(x))$. Then, since $s_{x^*}$ reverses geodesics through $x^*$,
$$
d\hspace{0.02cm}s_{x^*}\cdot G_{x^*}(x) = G_{x^*}(s_{x^*}\cdot x)
$$
Replacing this into (\ref{eq:bcproofsx2}), and introducing the new variable of integration $z = s_{x^*}\cdot x$\footnote{$\pi_{\scriptscriptstyle T}\circ s_{x^*}$ is the image of the distribution $\pi_{\scriptscriptstyle T}$ under the map $s_{x^*}:M\rightarrow M$. In other places of this thesis, this would be noted $(s_{x^*})^*\pi_{\scriptscriptstyle T\hspace{0.03cm}}$, but this notation seems kind of clumsy, in the present case.}, 
\begin{equation} \label{eq:bcproofsx2}
  d\hspace{0.02cm}s_{x^*}\cdot G_{x^*} = \int_{\mathrm{D}(x^*)}\,G_{x^*}(z)\hspace{0.02cm}(\pi_{\scriptscriptstyle T}\circ s_{x^*})(dz)  
\end{equation}
since $s^{\scriptscriptstyle -1}_{x^*} = s_{x^*}$ and $s_{x^*}$ maps $\mathrm{D}(x^*)$ onto itself. Now, note that $\pi_{\scriptscriptstyle T}\circ s_{x^*} = \pi_{\scriptscriptstyle T\hspace{0.03cm}}$. This is clear, since from (\ref{eq:gibbs}),
$$
(\pi_{\scriptscriptstyle T}\circ s_{x^*})(dz) = \left(Z(T)\right)^{-1}\exp\left[-\frac{(U\circ s_{x^*})(z)}{T}\right](\mathrm{vol}\circ s_{x^*})(dz)
$$
However, by assumption, $(U\circ s_{x^*})(z) = U(z)$,  since $U$ is invariant by geodesic symmetry about $x^*$. Moreover, because $s_{x^*}$ is an isometry, it must preserve Riemannian volume, so 
$(\mathrm{vol}\circ s_{x^*})(dz) = \mathrm{vol}(dz)$. Thus, (\ref{eq:bcproofsx2}) reads
$$
d\hspace{0.02cm}s_{x^*}\cdot G_{x^*} = 
\int_{\mathrm{D}(x^*)}\,G_{x^*}(z)\hspace{0.03cm}\pi_{\scriptscriptstyle T}(dz)
$$
From (\ref{eq:GH}), the right-hand side is $G_{x^*\hspace{0.03cm}}$, so (\ref{eq:bcproofsx}) is obtained. From (\ref{eq:bcproofsx}), since $d\hspace{0.02cm}s_{x^*} = -\mathrm{Id}_{x^*}\hspace{0.03cm}$, and $G_{x^*}$ belongs to $T_{x^*}M$,
$$
G_{x^*} = -\,G_{x^*}
$$
Of course, this means $G_{x^*} = 0$, as required.

\subsection{Proof of Lemma \ref{lem:etconv}} \label{subsec:proofetconv}

\subsubsection{Proof of (i)}
Let $y \in B(x^*\!,\delta)$ where $\delta < \frac{1}{2}r_{\scriptscriptstyle cx\hspace{0.03cm}}$. Now, recall that $\mathrm{Hess}\,\mathcal{E}_{\scriptscriptstyle T}(y) = H_y$ for all $y \in B(x^*\!,\delta)$. Then, from (\ref{eq:GH}), it is possible to write
\begin{equation} \label{eq:proofetconv1}
  H_y \,= \int_{B(y\hspace{0.02cm},r_{\scriptscriptstyle cx})}\,H_y(x)\hspace{0.03cm}\pi_{\scriptscriptstyle T}(dx) +
             \int_{\mathrm{D}(y)-B(y\hspace{0.02cm},r_{\scriptscriptstyle cx})}\,H_y(x)\hspace{0.03cm}\pi_{\scriptscriptstyle T}(dx)
\end{equation} 
Indeed, $B(y\hspace{0.02cm},r_{\scriptscriptstyle cx}) \subset \mathrm{D}(y)$, since the injectivity radius of $M$ is $2r_{\scriptscriptstyle cx}$ as given in \ref{subsec:ssbiscut}. The first integral in (\ref{eq:proofetconv1}) will be denoted $I_{\scriptscriptstyle 1}$ and the second integral $I_{\scriptscriptstyle 2\hspace{0.03cm}}$.

With regard to $I_{\scriptscriptstyle 1\hspace{0.03cm}}$, note the inclusions $B(x^*\!,\delta) \subset B(y\hspace{0.02cm},2\delta) \subset B(y\hspace{0.02cm},r_{\scriptscriptstyle cx})$, which follow from the triangle inequality. In addition, note that $H_y(x) \geq 0$ for $x \in B(y\hspace{0.02cm},r_{\scriptscriptstyle cx})$. Therefore,
\begin{equation} \label{eq:proofetconv2}
 I_{\scriptscriptstyle 1} \geq \int_{B(x^*\!,\delta)}\,H_y(x)\hspace{0.03cm}\pi_{\scriptscriptstyle T}(dx)
\end{equation}
where $H_y(x)$ is given by (\ref{eq:sshessfy}). But, from (\ref{eq:sshessfy}); the eigenvalues of $H_y(x)$ are
\begin{equation} \label{eq:etconveigen}
   \lambda(a)\cot\lambda(a) \geq \min_{\ell}\,c_{\scriptscriptstyle \ell}(a)\cot c_{\scriptscriptstyle \ell}(a)
\end{equation}
where the maximal roots $c_{\ell}$ were introduced before (\ref{eq:tc2}). By the Cauchy-Schwarz inequality, $c_{\scriptscriptstyle \ell}(a) \leq \Vert c_{\scriptscriptstyle \ell}\Vert\hspace{0.02cm}\Vert a \Vert \leq c \Vert a\Vert$, where $c^{\hspace{0.02cm}\scriptscriptstyle 2}$ denotes the maximum sectional curvature of $M$, whose expression was recalled in \ref{subsec:ssbiscut}. Now, if $x \in B(y\hspace{0.02cm},2\delta)$, then $\Vert a \Vert = d(y\hspace{0.02cm},x) < 2\delta$, and it follows from (\ref{eq:etconveigen}) that
\begin{equation} \label{eq:proofetconv3}
  H_y(x) \geq \min_{\ell}\,c_{\scriptscriptstyle \ell}(a)\cot c_{\scriptscriptstyle \ell}(a) \geq
2c\delta\hspace{0.02cm}\cot(2c\delta) = \mathrm{Ct}(2\delta) > 0
\end{equation}
where the last inequality is because $2c\delta < \frac{\pi}{2}$. Replacing in (\ref{eq:proofetconv2}) gives
$$
I_{\scriptscriptstyle 1} \geq \mathrm{Ct}(2\delta)\hspace{0.03cm}\pi_{\scriptscriptstyle T}(B(x^*\!,\delta)) = 
\mathrm{Ct}(2\delta)[1 - \pi_{\scriptscriptstyle T}(B^{\scriptscriptstyle c}(x^*\!,\delta))]
$$
Finally, (\ref{eq:proofestimate12}) and (\ref{eq:lowerz}) imply that $\pi_{\scriptscriptstyle T}(B^{\scriptscriptstyle c}(x^*\!,\delta)) \leq \mathrm{vol}(M)\hspace{0.02cm} f(T)$, where $f(T)$ was defined in (\ref{eq:fT}) (precisely, this follows after replacing $\rho$ by $\delta$ in (\ref{eq:proofestimate12})).
Thus,
\begin{equation} \label{eq:proofetconv4}
I_{\scriptscriptstyle 1} \geq   \mathrm{Ct}(2\delta)[1 -\mathrm{vol}(M)\hspace{0.02cm} f(T)]
\end{equation}
The proof of (\ref{eq:ethlower}) will be completed by showing
\begin{equation} \label{eq:AMM}
I_{\scriptscriptstyle 2} \geq - \pi A_{\scriptscriptstyle M}\hspace{0.03cm}f(T)
\end{equation}
To do so, introduce the function
\begin{equation} \label{eq:funka}
k(a) = \min_{\ell}\,c_{\scriptscriptstyle \ell}(a)\cot c_{\scriptscriptstyle \ell}(a)  \hspace{0.5cm} \text{for $a \in Q_+$}
\end{equation}
and note using (\ref{eq:etconveigen}) that
\begin{equation} \label{eq:AMM1}
I_{\scriptscriptstyle 2} \geq \int_{\mathrm{D}(y) - B(y\hspace{0.02cm},r_{\scriptscriptstyle cx})}\,k(a)\hspace{0.03cm}\pi_{\scriptscriptstyle T}(dx) \geq
\int_{\mathrm{D}(y)}\,\mathbf{1}\lbrace k(a) \leq 0\rbrace k(a)\hspace{0.03cm}\pi_{\scriptscriptstyle T}(dx)
\end{equation}
Indeed, the set of $a$ such that $k(a) \leq 0$ is contained in $\mathrm{D}(y) - B(y\hspace{0.02cm},r_{\scriptscriptstyle cx})$, because
\begin{equation} \label{eq:AMM2}
\lbrace k(a) \leq 0\rbrace = \cup_{\ell}\hspace{0.03cm} \lbrace c_{\scriptscriptstyle \ell}(a)\cot c_{\scriptscriptstyle \ell}(a) \leq 0 \rbrace =
\cup_{\ell}\hspace{0.03cm} \lbrace c_{\scriptscriptstyle \ell}(a) \geq \pi/2 \rbrace
\end{equation}
and $c_{\scriptscriptstyle \ell}(a) \geq \pi/2$ implies $d(y\hspace{0.02cm},x) = \Vert a \Vert \geq \frac{\pi}{2}\hspace{0.03cm}c^{\scriptscriptstyle -1} = r_{\scriptscriptstyle cx}$ (by Cauchy-Schwarz). By expressing the last integral in (\ref{eq:AMM1}) as in (\ref{eq:Hypolar}), it is seen to be equal to 
$$
\begin{array}{l}
  \int_{Q_+}\!\int_S\,
\mathbf{1}\lbrace k(a) \leq 0\rbrace k(a)\hspace{0.02cm}
p_{\scriptscriptstyle T}(s\hspace{0.02cm},a)\hspace{0.02cm}D(a)\hspace{0.03cm}da\hspace{0.02cm}\omega(ds) \geq \\[0.3cm]
-\pi\,
  \int_{Q_+}\!\int_S\,
\mathbf{1}\lbrace k(a) \leq 0\rbrace\hspace{0.02cm}
p_{\scriptscriptstyle T}(s\hspace{0.02cm},a)\hspace{0.03cm}da\hspace{0.02cm}\omega(ds)
\end{array}
$$
Indeed, it follows from (\ref{eq:HypolarD}) and (\ref{eq:funka}) that $k(a)\hspace{0.02cm}D(a) \geq \min_{\ell} -c_{\scriptscriptstyle \ell}(a)$, and this is greater than $-\pi$ because $c_{\scriptscriptstyle \ell}(a) \in (0,\pi)$ for $a \in Q_{+}\hspace{0.02cm}$. Now, (\ref{eq:AMM1}) implies
\begin{equation} \label{eq:AMM3}
I_{\scriptscriptstyle 2} \geq -\pi\,
  \int_{Q_+}\!\int_S\,
\mathbf{1}\lbrace k(a) \leq 0\rbrace\hspace{0.02cm}
p_{\scriptscriptstyle T}(s\hspace{0.02cm},a)\hspace{0.03cm}da\hspace{0.02cm}\omega(ds)
\end{equation}
As seen from (\ref{eq:AMM2}), the set $\lbrace k(a) \leq 0\rbrace \subset B^{\scriptscriptstyle c}(y\hspace{0.02cm},r_{\scriptscriptstyle cx}) \subset B^{\scriptscriptstyle c}(x^*\!,\delta)$. On the other hand,  $p_{\scriptscriptstyle T}(x) \leq f(T)$ for $x \in B^{\scriptscriptstyle c}(x^*\!,\delta)$. Replacing in (\ref{eq:AMM3}),
\begin{equation} \label{eq:AMM4}
I_{\scriptscriptstyle 2} \geq -\pi f(T)\,\int_{Q_+}\!\int_S\,da\hspace{0.02cm}\omega(ds)
\end{equation}
The double integral on the right-hand side is a positive constant which depends only on the symmetric space $M$. Denoting this by $A_{\scriptscriptstyle M}$ yields (\ref{eq:AMM}).

\subsubsection{Proof of (ii)}
Let $\delta < \frac{1}{2}r_{\scriptscriptstyle cx\hspace{0.03cm}}$. According to (\ref{eq:ethlower}), which has just been proved
\begin{equation} \label{eq:chap2fin1}
   \mathrm{Hess}\,\mathcal{E}_{\scriptscriptstyle T}(y) \,\geq\, \mathrm{Ct}(2\delta)\hspace{0.03cm}[1 - \mathrm{vol}(M)f(T)] - \pi A_{\scriptscriptstyle M}\hspace{0.03cm}f(T)
\end{equation}
for all $y \in B(x^*\!,\delta)$. Now, let $T_{\scriptscriptstyle \delta}$ be given by (\ref{eq:td}). It follows from the definition of $T^2_{\scriptscriptstyle \delta}$ that $T < T_{\scriptscriptstyle \delta}$ implies 
\begin{equation} \label{eq:chap2fin2}
 f(T) < \frac{\mathrm{Ct}(2\delta)}{\mathrm{Ct}(2\delta)\hspace{0.02cm}\mathrm{vol}(M) + \pi\hspace{0.02cm}A_{\scriptscriptstyle M}}
\end{equation}
This amounts to saying the right-hand side of (\ref{eq:chap2fin1}) is strictly positive. Since this is independent of $y$, it is clear that the variance function $\mathcal{E}_{\scriptscriptstyle T}(y)$ is indeed strongly convex on $B(x^*\!,\delta)$.

%%REVIEW REFS TO SEC:SSBIS

\chapter{Gaussian distributions and RMT} \label{gaussian}

\minitoc
\vspace{0.1cm}
{\small
Gaussian distributions on Riemannian symmetric spaces were introduced in~\cite{said2}. The present chapter expands on this work in several ways. In particular, it uncovers and exploits the connection between Gaussian distributions on Riemannian symmetric spaces and random matrix theory (RMT). 
\begin{itemize}
 \item \ref{sec:gausshistory} attempts to answer the seemingly easy question \textit{what is a Gaussian distribution\,?}, by adopting a historical perspective (the source material used is from~\cite{historynormal}\cite{borel}\cite{perrin}). 
 \item \ref{sec:rgd} defines Gaussian distributions as a family of distributions on a Riemannian manifold, for which \textit{maximum-likelihood estimation is equivalent to the Riemannian barycentre problem}. 
\item \ref{sec:z} expresses the normalising factor of a Gaussian distribution on a Riemannian symmetric space, which belongs to the non-compact case, under the form of a multiple integral. It also discusses analytic and numerical evaluation of this multiple integral. 
\item \ref{sec:gaussmle} proves existence and uniqueness of maximum-likelihood estimates for Gaussian distributions, defined on a Hadamard manifold which is also a homogeneous space. It also shows that these distributions maximise Shannon entropy, for given barycentre and dispersion. 
\item \ref{sec:gaussschur} describes the Riemannian barycentre and the covariance tensor of a Gaussian distribution.
\item \ref{sec:orthpo} and draws on random matrix theory to obtain an analytic formula for the normalising factor of a Gaussian distribution, on the space $\mathrm{H}(N)$ of $N \times N$ Hermitian positive-definite matrices.
\item \ref{sec:rmt} and \ref{sec:rmtbis} study large $N$ asymptotics of Gaussian distributions on $\mathrm{H}(N)$. In particular, \ref{sec:rmt} provides an asymptotic expression of the normalising factor. 
\item \ref{sec:THETA} introduces $\Theta$ distributions, a family of distributions, on the unitary group $U(N)$ (which is the dual symmetric space of $\mathrm{H}(N))$, with a remarkable connection to Gaussian distributions on $\mathrm{H}(N)$. 
%\item \ref{sec:warped} shows that the Fisher information metric, of the family of Gaussian distributions, on an irreducible Riemannian symmetric space, is a so-called warped Riemannian metric, and describes its geodesic curves. 
\end{itemize}
}

\vfill
\pagebreak

\section{From Gauss to Shannon} \label{sec:gausshistory}
The story of Gaussian distributions is a story of discovery and re-discovery. Different scientists, at different times, were repeatedly lead to these distributions, through different routes. 

In 1801, on New Year's day, Giuseppe Piazzi sighted a heavenly body (in fact, the asteroid Ceres), which he thought to be a new planet. Less than six weeks later, this ``new planet" disappeared behind the sun. Using a method of least squares, Gauss predicted the area in the sky, where it re-appeared one year later. His justification of this method of least squares (cast in modern language) is that measurement errors follow a family of distributions, which satisfies \\
\textbf{property 1\,:} maximum-likelihood estimation is equivalent to the least-squares problem. \\
\indent In an 1809 paper, he used this property to show that the distribution of measurement errors is (again, in modern language) a Gaussian distribution.

In 1810, Laplace studied the distribution of a quantity, which is the aggregate of a great number of elementary observations. He was lead in this (completely different) way, to the same distribution discovered by Gauss. Laplace was among the first scientists to show\\
\textbf{property 2\,:} the distribution of the sum of a large number of elementary observations is (asymptotically) a Gaussian distribution. \\
\indent Around 1860, Maxwell rediscovered Gaussian distributions, through his investigation of the velocity distribution of particles in an ideal gas (which he viewed as freely colliding perfect elastic spheres). Essentially, he showed that \\
\textbf{property 3\,:} the distribution of a rotationally-invariant random vector, which has independent components, is a Gaussian distribution.\footnote{A deeper version of Maxwell's idea was obtained by Poincaré and Borel, around 1912, who showed that\,: if $v = (v_n\,;n=1,\ldots,N)$ is uniformly distributed, on an $(N-1)$-dimensional sphere of radius $N^{\frac{1}{2}}$, then the distribution of $v_{\scriptscriptstyle 1}$ is (asymptotically) a Gaussian distribution. This is Poincaré's model of the one-dimensional ideal gas, with $N$ particles.}\\
\indent Kinetic theory lead to another fascinating development, related to Gaussian distributions. Around 1905, Einstein (and, independently, Smoluchowsky) showed that \\
\textbf{property 4\,:} the distribution of the position of a particle, which is undergoing a Brownian motion, is a Gaussian distribution. \\
\indent In addition to kinetic theory, alternative routes to Gaussian distributions have been found in quantum mechanics, information theory, and other fields.
In quantum mechanics, a Gaussian distribution is a position distribution with minimum uncertainty. That is, it achieves equality in Heisenberg's inequality (this is because a Gaussian function is proportional to its own Fourier transform). In information theory, one may attribute to Shannon the following maximum\hfill\linebreak entropy characterisation \\
\textbf{property 5\,:} a probability distribution with maximum entropy, among all distributions with a given mean and variance, is a Gaussian distribution. \\
\indent The above list of re-discoveries of Gaussian distributions, by means of different definitions, may be extended much longer. However, the main point is the following. In Euclidean space, any one of the above five properties leads to the same famous expression of a Gaussian distribution,
$$
P(dx|\bar{x},\sigma) = \left(2\pi\sigma^2\right)^{\!-\frac{n}{2}}\exp\left[ -\frac{(x-\bar{x})^2}{2\sigma^2}\right]dx
$$
In non-Euclidean space, each one of these properties may lead to a different kind of distribution, which may then be called a Gaussian distribution, but only from a more restricted point of view. People interested in Brownian motion may call the heat kernel of a Riemannian manifold a Gaussian distribution on that manifold. However, statisticians will not like this definition, since it will (in general) fail to have a straightforward connection to maximum-likelihood estimation.

Going from Euclidean to non-Euclidean spaces, the concept of Gaussian distribution breaks down into several different concepts. At best, one may define Gaussian distributions based on a practical motivation, which makes one (or more) of their classical properties seem more advantageous than the others. 
\section{The ``right" Gaussian} \label{sec:rgd}
As of now, the following definition of Gaussian distributions is chosen. Gaussian distributions, on a Riemannian manifold $M$, are a family of distributions $P(\bar{x}\hspace{0.03cm},\sigma)$, parameterised by $\bar{x} \in M$ and $\sigma > 0$, such that\,: a maximum-likelihood estimate $\hat{x}_{\scriptscriptstyle N}$ of $\bar{x}$, based on independent samples $(x_n\,;n=1,\ldots,N)$ from $P(\bar{x}\hspace{0.02cm},\sigma)$, is a solution of the least-squares problem 
$$
\text{minimise over $x \in M$} \hspace{0.3cm} \mathcal{E}_{\scriptscriptstyle N}(x) = \sum^N_{n=1}d^{\hspace{0.03cm}2}(x_n\hspace{0.02cm},x)
$$
Of course, this is the same least-squares problem as (\ref{eq:empiricalfrechet}), so $\hat{x}_{\scriptscriptstyle N}$ is an empirical barycentre of the samples $(x_n)$. Therefore (as discussed in \ref{subsec:generic}), $\hat{x}_{\scriptscriptstyle N}$ is almost-surely unique, if $P(\bar{x}\hspace{0.02cm},\sigma)$ has a probability density with respect to the Riemannian volume of $M$ (this will indeed be the case). 
Now, consider the density profile 
\begin{equation} \label{eq:gaussprofile}
  f(x|\bar{x}\hspace{0.02cm},\sigma) \,=\, \exp\left[ -\frac{d^{\hspace{0.03cm}2}(x,\bar{x})}{2\sigma^2}\right]
\end{equation}
and the normalising factor,
\begin{equation} \label{eq:gaussnorm}
  Z(\bar{x}\hspace{0.02cm},\sigma) = \int_M\,   f(x|\bar{x}\hspace{0.02cm},\sigma)\,\mathrm{vol}(dx)
\end{equation}
If this is finite, then 
\begin{equation} \label{eq:pregaussdensity}
 P(dx|\bar{x},\sigma) \,=\, \left(Z(\bar{x}\hspace{0.02cm},\sigma)\right)^{-1}\hspace{0.03cm}f(x|\bar{x}\hspace{0.02cm},\sigma)\hspace{0.04cm}\mathrm{vol}(dx)
\end{equation}
is a well-defined probability distribution on $M$. In \ref{sec:gaussmle}, below, it will be shown that $P(\bar{x}\hspace{0.02cm},\sigma)$, as defined by (\ref{eq:pregaussdensity}), is indeed a Gaussian distribution, if $M$ is a Hadamard manifold, and also a homogeneous space. The following propositions will then be helpful.
\begin{proposition} \label{prop:zhadamardk}
 Let $M$ be a Hadamard manifold, whose sectional curvatures lie in $[\kappa\hspace{0.03cm},0]$, where $\kappa = -c^{\hspace{0.02cm}\scriptscriptstyle 2}$. Then, for any $\bar{x} \in M$ and $\sigma > 0$, if $Z(\bar{x}\hspace{0.02cm},\sigma)$ is given by (\ref{eq:gaussnorm}),
\begin{equation} \label{eq:zhadamardk}
Z_{\scriptscriptstyle 0}(\sigma)  \leq\, Z(\bar{x}\hspace{0.02cm},\sigma) \,\leq Z_{\scriptscriptstyle c}(\sigma)
\end{equation}
where $Z_{\scriptscriptstyle 0}(\sigma) = \left(2\pi\sigma^2\right)^{\!\frac{n}{2}}$ and $Z_{\scriptscriptstyle c}(\sigma)$ is positive and given by (recall $n$ is the dimension of $M$)
\begin{equation} \label{eq:zcsigma}
  Z_{\scriptscriptstyle c}(\sigma) = \omega_{n-1}\hspace{0.02cm}\frac{\sigma}{(2c)^{n-1}}\hspace{0.03cm}\sum^{n-1}_{k=0}(-1)^k\hspace{0.03cm}\left(\!\!\begin{array}{c}n-1 \\ k \end{array}\!\!\right)\frac{\Phi\left((n-1-2k)\hspace{0.03cm}\sigma c\right)}{\mathstrut\Phi^\prime\left((n-1-2k)\hspace{0.03cm}\sigma c\right)}
\end{equation}
with $\omega_{n-1}$ the area of the unit sphere $S^{\hspace{0.02cm}n-1}$, and $\Phi$ the standard normal distribution function. 
\end{proposition}
\begin{proposition} \label{prop:zhomogeneous}
 If $M$ is a Riemannian homogeneous space, and $Z(\bar{x}\hspace{0.02cm},\sigma)$ is given by (\ref{eq:gaussnorm}), then $Z(\bar{x}\hspace{0.02cm},\sigma)$ does not depend on $\bar{x}$. In other words, $Z(\bar{x}\hspace{0.02cm},\sigma) = Z(\sigma)$.
\end{proposition}
If $M$ is a Hadamard manifold, and also a homogeneous space, then both Propositions \ref{prop:zhadamardk} and \ref{prop:zhomogeneous} apply to $M$. Indeed, if $M$ is a Riemannian homogeneous space, then its sectional curvatures lie within a bounded subset of the real line. Therefore, Proposition \ref{prop:zhadamardk} implies $Z(\bar{x}\hspace{0.02cm},\sigma)$ is finite for all $\bar{x} \in M$ and $\sigma > 0$. On the other hand, Proposition \ref{prop:zhomogeneous} implies that $Z(\bar{x}\hspace{0.02cm},\sigma) = Z(\sigma)$. \vfill\pagebreak

Thus, if $M$ is a Hadamard manifold, and also a homogeneous space, then (\ref{eq:pregaussdensity}), reduces to
\begin{equation} \label{eq:gaussdensity}
 P(dx|\bar{x},\sigma) \,=\, \left(Z(\sigma)\right)^{-1}\hspace{0.03cm}\exp\left[ -\frac{d^{\hspace{0.03cm}2}(x,\bar{x})}{2\sigma^2}\right]\hspace{0.04cm}\mathrm{vol}(dx)
\end{equation}
and yields a well-defined probability distribution $P(\bar{x}\hspace{0.02cm},\sigma)$ on $M$. This will be the main focus, throughout the following. \\[0.1cm]
\textbf{Proof of Proposition \ref{prop:zhadamardk}\,:} (\ref{eq:zhadamardk}) is a direct application of (\ref{eq:integralcomp}). Let $f(y) = f(y|\bar{x}\hspace{0.02cm},\sigma)$, and $\kappa_{\max} = 0$, $\kappa_{\min} = \kappa$. Also, since $M$ is a Hadamard manifold, note that $\min\lbrace\mathrm{inj}(\bar{x})\hspace{0.03cm},\pi\hspace{0.03cm}c^{\scriptscriptstyle -1}\rbrace = \infty$. Therefore, (\ref{eq:integralcomp}) (applied with $x = \bar{x}$), yields
$$
\omega_{n-1}\,\int^{\infty}_{0}\,\exp\left[-\frac{r^{\hspace{0.03cm}\scriptscriptstyle 2}}{\mathstrut 2\sigma^{\hspace{0.03cm}\scriptscriptstyle 2}}\right]\mathrm{sn}^{n-1}_{\scriptscriptstyle 0}(r)\hspace{0.02cm}dr \leq\, Z(\bar{x}\hspace{0.02cm},\sigma) \,\leq
\omega_{n-1}\,\int^{\infty}_{ 0}\,\exp\left[-\frac{r^{\hspace{0.03cm}\scriptscriptstyle 2}}{\mathstrut 2\sigma^{\hspace{0.03cm}\scriptscriptstyle 2}}\right]\mathrm{sn}^{n-1}_{\kappa}(r)\hspace{0.02cm}dr
$$
However, $\mathrm{sn}_{\scriptscriptstyle 0}(r) = r$ and $\mathrm{sn}_{\kappa}(r) = c^{\scriptscriptstyle -1}\hspace{0.02cm}\sinh(c\hspace{0.03cm} r)$. Therefore, the expression for $Z_{\scriptscriptstyle 0}(\sigma)$ follows easily. For $Z_{\scriptscriptstyle c}(\sigma)$, on the other hand, note that
$$
\int^{\infty}_{ 0}\,\exp\left[-\frac{r^{\hspace{0.03cm}\scriptscriptstyle 2}}{\mathstrut 2\sigma^{\hspace{0.03cm}\scriptscriptstyle 2}}\right]\mathrm{sn}^{n-1}_{\kappa}(r)\hspace{0.02cm}dr = \frac{1}{(2c)^{n-1}}\hspace{0.03cm}\int^{\infty}_{ 0}\,
\exp\left[-\frac{r^{\hspace{0.03cm}\scriptscriptstyle 2}}{\mathstrut 2\sigma^{\hspace{0.03cm}\scriptscriptstyle 2}}\right]\left( e^{c\hspace{0.03cm}r} - e^{-c\hspace{0.03cm}r\hspace{0.03cm}} \right)^{n-1}\hspace{0.03cm}dr
$$
Then, (\ref{eq:zcsigma}) follows by performing a binomial expansion, and using
$$
\int^{\infty}_{ 0}\,\exp\left[-\frac{r^{\hspace{0.03cm}\scriptscriptstyle 2}}{\mathstrut 2\sigma^{\hspace{0.03cm}\scriptscriptstyle 2}} +
(n-1-2k)\hspace{0.03cm}c\hspace{0.03cm}r \right]dr = \sigma\,\frac{\Phi\left((n-1-2k)\hspace{0.03cm}\sigma c\right)}{\mathstrut\Phi^\prime\left((n-1-2k)\hspace{0.03cm}\sigma c\right)}
$$
\textbf{Remark\,:} clearly, $Z_{\scriptscriptstyle 0}(\sigma)$ is the normalising factor of a Gaussian distribution, when $M$ is a Euclidean space, $M = \mathbb{R}^n\,$. On the other hand, $Z_{\scriptscriptstyle c}(\sigma)$ is the normalising factor of a Gaussian distribution, when $M$ is a hyperbolic space of dimension $n$, and constant negative curvature $\kappa = -c^{\hspace{0.02cm}\scriptscriptstyle 2}$. This will become clear in \ref{sec:z}, below. \\[0.1cm]
\noindent \textbf{Proof of Proposition \ref{prop:zhomogeneous}\,:} assume $M$ is a homogeneous space, and fix some point $o \in M$. There exists an isometry $g$ of $M$ such that $g\cdot\bar{x} = o$. In the integral (\ref{eq:gaussnorm}), introduce the new variable of integration $z = g \cdot x$. Since $g$ (being an isometry) preserves Riemannian volume,
$$
  Z(\bar{x}\hspace{0.02cm},\sigma) = \int_M\,   f(g^{\scriptscriptstyle -1}\cdot z|\bar{x}\hspace{0.02cm},\sigma)\,\mathrm{vol}(dz)
= \int_M\,   f( z|o\hspace{0.02cm},\sigma)\,\mathrm{vol}(dz) = Z(o\hspace{0.02cm},\sigma)
$$
where the second equality follows from (\ref{eq:gaussprofile}). Thus, $Z(\bar{x}\hspace{0.02cm},\sigma) = Z(o\hspace{0.02cm},\sigma)$ does not depend on $\bar{x}$.
\section{The normalising factor $Z(\sigma)$} \label{sec:z}
Assume now $M = G/K$ is a Riemannian symmetric space which belongs to the non-compact case, described in \ref{ssec:sspace}. In particular, $M$ is a Hadamard manifold, and also a homogeneous space. Thus, for each $\bar{x} \in M$ and $\sigma > 0$, there is a well-defined probability distribution $P(\bar{x}\hspace{0.02cm},\sigma)$ on $M$, given by (\ref{eq:gaussdensity}). Here, the normalising factor $Z(\sigma)$ can be expressed as a multiple integral, using the integral formula (\ref{eq:ssvolncka2}), of Proposition \ref{prop:ssvolncka}. Applying this proposition (with $o=\bar{x}$), it is enough to note
$$
f(\varphi(s\hspace{0.02cm},a)|\bar{x}\hspace{0.02cm},\sigma) = \exp\left[ -\frac{\Vert a \Vert^{2}_{\scriptscriptstyle B}}{2\sigma^2}\right]
$$ 
where $\Vert a \Vert^{2}_{\scriptscriptstyle B} = B(a,a)$, in terms of the $\mathrm{Ad}(G)$-invariant symmetric bilinear form $B$. Since this expression only depends on $a$, it is possible to integrate $s$ out of (\ref{eq:ssvolncka2}), to obtain
\begin{equation} \label{eq:ssz}
  Z(\sigma) \,=\,
 \frac{\omega(S)}{|W|}\,\int_{\mathfrak{a}}\exp\left[ -\frac{\Vert a \Vert^{2}_{\scriptscriptstyle B}}{2\sigma^2}\right]\prod_{\lambda \in \Delta_+}\left| \sinh\hspace{0.02cm}\lambda(a)\right|^{ m_\lambda}\hspace{0.03cm}da
\end{equation}
This formula expresses the normalising factor $Z(\sigma)$ as a multiple integral on the vector space $\mathfrak{a}$. \vfill\pagebreak

\noindent \textbf{Example 1\,:} the easiet instance of (\ref{eq:ssz}) arises when $M$ is a hyperbolic space of dimension $n$, and constant sectional curvature equal to $-1$. Then, $M$ has rank equal to $1$, so that $\mathfrak{a} = \mathbb{R}\hspace{0.03cm}\hat{a}$ for some unit vector $\hat{a} \in \mathfrak{a}$. Since the sectional curvature is equal to $-1$, there is only one positive root $\lambda$, say $\lambda(\hat{a}) = 1$, with multiplicity $m_\lambda= n-1$. In addition, there are two Weyl chambers, $C_+ = \lbrace t\hspace{0.02cm}\hat{a}\,; t > 0\rbrace$ and $C_- = \lbrace t\hspace{0.02cm}\hat{a}\,; t < 0\rbrace$. In other words, $|W| = 2$. Now, (\ref{eq:ssz}) reads
$$
Z(\sigma) \,=\,
 \frac{\omega_{n-1}}{2}\,\int^{+\infty}_{-\infty}\exp\left[-\frac{r^{\hspace{0.03cm}\scriptscriptstyle 2}}{\mathstrut 2\sigma^{\hspace{0.03cm}\scriptscriptstyle 2}}\right]\left| \sinh(r)\right|^{n-1}\hspace{0.03cm}dr =  \omega_{n-1}\,\int^{+\infty}_{ 0}\exp\left[-\frac{r^{\hspace{0.03cm}\scriptscriptstyle 2}}{\mathstrut 2\sigma^{\hspace{0.03cm}\scriptscriptstyle 2}}\right]\sinh^{n-1}(r)\hspace{0.03cm}dr
$$
In general, if all distances are divided by $c > 0$, the sectional curvature $-1$ is replaced by $-c^{\hspace{0.02cm}\scriptscriptstyle 2}$. Thus, when $M$ is a hyperbolic space of dimension $n$, and sectional curvature $-c^{\hspace{0.02cm}\scriptscriptstyle 2}$,  one has
$$
Z(\sigma) =  \omega_{n-1}\,\int^{+\infty}_{ 0}\exp\left[-\frac{r^{\hspace{0.03cm}\scriptscriptstyle 2}}{\mathstrut 2\sigma^{\hspace{0.03cm}\scriptscriptstyle 2}}\right](c^{\scriptscriptstyle -1}\hspace{0.02cm}\sinh(c\hspace{0.03cm}r))^{n-1}\hspace{0.03cm}dr
$$
This is exactly $Z_{\scriptscriptstyle c}(\sigma)$, expressed analytically in (\ref{eq:zcsigma}). \\[0.1cm]
\textbf{Example 2\,:} another example, also susceptible of analytic expression, is when $M$ is a cone of positive-definite matrices (covariance matrices), with real, complex, or quaternion coefficients. Then, $M = G/K$ with $G = \mathrm{GL}(N,\mathbb{K})$, where $\mathbb{K} = \mathbb{R}, \mathbb{C}$ or $\mathbb{H}$ (real numbers, complex numbers, or quaternions), and $K$ is a maximal compact subgroup of $G$, say $K = U(N), O(N)$ or $Sp(n)$. 
In each of these three cases, $\mathfrak{a}$ is the space of $N \times N$ real diagonal matrices, and the positive roots are the linear maps $\lambda(a) = a_{ii} - a_{jj}$ where $i < j$, each one having its multiplicity $m_\lambda = \beta$, ($\beta = 1,2$ or $4$, according to $\mathbb{K} = \mathbb{R}, \mathbb{C}$ or $\mathbb{H}$). In addition, $\Vert a \Vert^{2}_{\scriptscriptstyle B} = 4\hspace{0.02cm}\mathrm{tr}(a^{\scriptscriptstyle 2}) = 4\hspace{0.02cm}a^{{\scriptscriptstyle 2}}_{\scriptscriptstyle 11} + \ldots + 4\hspace{0.02cm}a^{{\scriptscriptstyle 2}}_{\scriptscriptstyle NN\,}$. 

The Weyl group $W$ is the groupe of permutation matrices in $K$, so $|W| = N!$. Finally, $S = K/T_{\scriptscriptstyle N}$ where $T_{\scriptscriptstyle N}$ is the subgroup of all matrices $t$ which are diagonal and belong to $K$. Replacing all of this into (\ref{eq:ssz}), it follows that
\footnote{Curious readers will want to compute $\omega_{\beta}(N)$. But, how\,? For example, $\omega_{\scriptscriptstyle 2}(N)$ can be found using the Weyl integral formula on U(N)~\cite{knapp}. This yields $\omega_{\scriptscriptstyle 2}(N) = \mathrm{vol}(U\!(N))/(2\pi)^N$, where the volume of the unitary group $U(N)$ is $\mathrm{vol}(U\!(N)) = (2\pi)^{(N^2+N)/2}/G(N)$, in terms of $G(N) = 1!\times2!\times\ldots\times(N-1)!$, which can be found, just by looking at the normalising factor of a Gaussian unitary ensemble.}
\begin{equation} \label{eq:covzbeta}
  Z(\sigma) = 
 \frac{\omega_{\beta}(N)}{N!}\,\int_{\mathfrak{a}}\,\prod^N_{i=1}\,\exp\left[ -\frac{2\hspace{0.02cm}a^{{\scriptscriptstyle 2}}_{\scriptscriptstyle ii}}{\sigma^2}\right]\prod_{i<j}\left| \sinh(a_{ii} - a_{jj})\right|^{ \beta}\hspace{0.03cm}da
\end{equation}
where $\omega_{\beta}(N)$ stands for $\omega(S)$, and $da = da_{\scriptscriptstyle 11}\ldots da_{\scriptscriptstyle NN\,}$. Passing to the variables $x_i = \exp(2a_{\scriptscriptstyle ii})$, 
$$
Z(\sigma) = 
\frac{\omega_{\beta}(N)}{\mathstrut 2^{\scriptscriptstyle N N_\beta}N!} \,\int_{\mathbb{R}^{\scriptscriptstyle N}_+}\,\prod^N_{i=1}\left[\rho(x_i\hspace{0.02cm},2\sigma^{\scriptscriptstyle 2})\hspace{0.03cm} x^{\!-N_\beta\,}_i\right]\hspace{0.03cm}|V(x)|^\beta\hspace{0.03cm}\prod^N_{i=1} \hspace{0.02cm}dx_i
$$
where $N_\beta = (\beta/2)(N-1) + 1$, $\rho(x,k) = \exp(-\log^2(x)/k)$ and $V(x) = \prod_{i<j} (x_j - x_i)$ is the Vandermonde determinant. Finally, using the elementary identity 
$$
\rho(x,k)\hspace{0.03cm} x^\alpha = \exp\left[\frac{k}{4}\hspace{0.03cm}\alpha^{\scriptscriptstyle 2}\right]\omega\left(e^{-\frac{k}{2}\hspace{0.02cm}\alpha}\hspace{0.02cm}x \hspace{0.02cm},k\right)
$$
it is immediately found that 
\begin{equation} \label{eq:covzvv}
Z(\sigma) = \frac{\omega_{\beta}(N)}{\mathstrut 2^{\scriptscriptstyle N N_\beta}N!}\times \exp\left[-N\hspace{0.02cm}N^2_\beta\hspace{0.04cm}(\sigma^{2}\!/2)\right]\times\int_{\mathbb{R}^{\scriptscriptstyle N}_+}\,\prod^N_{i=1}\rho(u_i\hspace{0.02cm},2\sigma^{\scriptscriptstyle 2})\,|V(u)|^\beta\,\prod^N_{i=1} \hspace{0.02cm}du_i
\end{equation}
For the case $\beta = 2$, the integral in (\ref{eq:covzvv}) will be expressed analytically in \ref{sec:orthpo}, below. The cases $\beta = 1, 4$ should be pursued using the techniques in~\cite{mehta} (Chapter 5). \\[0.1cm]
\textbf{Example 3\,:} for this last example, I am still unaware of any valid means of analytic expression. Let $M = \mathrm{D}_{\scriptscriptstyle N}$ be the Siegel domain~\cite{siegelspgeo}. This is the set of $N \times N$ symmetric complex matrices $z$,\hfill\linebreak such that $\mathrm{I}_N - z^\dagger z \geq 0$ (where the inequality is understood in the sense of the Loewner order).  Here, $M = G/K$, where $G\simeq \mathrm{Sp}(N,\mathbb{R})$ (real symplectic group) and $K \simeq U(N)$ (unitary group). Precisely, $G$ is the group of all $2N \times 2N$ complex matrices $g$, with $g^\mathrm{t}\hspace{0.02cm}\Omega\hspace{0.02cm}g = \Omega$ and $g^\dagger\hspace{0.02cm}\Gamma\hspace{0.02cm}g = \Gamma$, where $^\mathrm{t}$ denotes the transpose, and where $\Omega$ and $\Gamma$ are the matrices
$$
\Omega = \left(\begin{array}{cc}\! & \mathrm{I}_N \\[0.1cm] - \mathrm{I}_N &  \end{array}\!\right) \hspace{0.25cm};\hspace{0.25cm}
\Gamma = \left(\!\begin{array}{cc} \mathrm{I}_N  & \\[0.1cm]  & - \mathrm{I}_N \end{array}\!\right)
$$
In addition, $K$ is the group of all block-diagonal matrices $k = \mathrm{diag}(U,U^*)$ where $U \in U(N)$, and $^*$ denotes the conjugate. The action of $G$ on $M$ is given by the matrix analogue of Möbius transformations,
\begin{equation} \label{eq:siegel1}
  g\cdot z = (A\hspace{0.02cm}z +B)(C\hspace{0.02cm}z+D)^{\scriptscriptstyle -1} \hspace{0.5cm} g = \left(\!\begin{array}{cc} A & B \\[0.1cm] C  & D\end{array}\,\right)
\end{equation}
This action preserves the Siegel metric, which is defined by
\begin{equation} \label{eq:siegel2}
 \langle v,\!v\rangle_{\scriptscriptstyle z} = \Vert (\mathrm{I}_N - z\hspace{0.02cm}z^\dagger)^{\scriptscriptstyle -1}\hspace{0.03cm} v\hspace{0.03cm}\Vert^2_{\scriptscriptstyle B} \hspace{1cm} \Vert v \Vert^2_{\scriptscriptstyle B} = \frac{1}{2}\hspace{0.02cm}\mathrm{tr}(v\hspace{0.02cm}v^\dagger)
\end{equation}
where each tangent vector $v$ is identified with a symmetric complex matrix (with this metric, it is easy to see that geodesic symmetry at $0 \in \mathrm{D}_{\scriptscriptstyle N}$ is given by $s_{\scriptscriptstyle 0}(z) = -z$ for $z \in \mathrm{D}_{\scriptscriptstyle N}$). Now~\cite{terras2}, 
\begin{equation} \label{eq:siegel3}
  \mathfrak{a}  = \left \lbrace \left(\begin{array}{cc} & a \\[0.1cm] a & \end{array}\right)\,;\hspace{0.05cm} a = \mathrm{diag}(a_{\scriptscriptstyle 11},\ldots,
a_{\scriptscriptstyle NN})\right\rbrace
\end{equation}
The positive roots are $\lambda(a) = a_{ii} - a_{jj}$ for $i < j$, and $\lambda(a) = a_{ii} + a_{jj}$ for $i \leq j$, all  with $m_\lambda = 1$. The order of the Weyl group is $|W| = 2^NN!$, and $\omega(S) = \mathrm{vol}(U(N))/2^N$. Replacing into (\ref{eq:ssz}),
\begin{equation} \label{eq:siegelz}
   Z(\sigma) \,=\,
 \frac{\mathrm{vol}(U(N))}{2^{\scriptscriptstyle 2N}N!}\,\int_{\mathfrak{a}}\,\prod^N_{i=1}\,\exp\left[ -\frac{\hspace{0.02cm}a^{{\scriptscriptstyle 2}}_{\scriptscriptstyle ii}}{2\sigma^2}\right]\prod_{i<j} \sinh|a_{ii} - a_{jj}|\prod_{i\leq j} \sinh|a_{ii} + a_{jj}|\,da
\end{equation}
In~\cite{said2}, a special Monte Carlo method, for computing (\ref{eq:siegelz}) was indicated (this method is owed to Paolo Zanini). The idea is the following\,: if $\mathbf{a}$ is a random variable with normal distribution, of mean zero and covariance $\sigma^2\hspace{0.03cm}\mathrm{I}_N$ in $\mathbb{R}^N$, then $Z(\sigma)$ is given by the expectation,
\begin{equation} \label{eq:paolo}
Z(\sigma) \,=\,
 \frac{\mathrm{vol}(U(N))}{2^{\scriptscriptstyle 2N}N!}\,\mathbb{E}\left[ \prod_{i<j} \sinh|\mathbf{a}_{i} - \mathbf{a}_{j}|\prod_{i\leq j} \sinh|\mathbf{a}_{i} + \mathbf{a}_{j}|\right]
\end{equation}
For a given value of $\sigma$, this is easily approximated by an empirical average. However, it then remains  to guarantee that $Z(\sigma)$ is a well-behaved function of $\sigma$. Precisely (see Proposition \ref{prop:gausscgf} in \ref{sec:gaussmle}), if 
$\eta = (-2\sigma^{\scriptscriptstyle 2})^{\scriptscriptstyle -1}$, then $\psi(\eta) = \log\hspace{0.02cm}Z(\sigma)$ is a strictly convex function, from the half-line $(-\infty,0)$ onto $\mathbb{R}$. By approximating $Z(\sigma_{\scriptscriptstyle n})$ at certain nodes $\sigma_{\scriptscriptstyle n\,}$, and then performing a suitable spline interpolation, it becomes possible to guarantee this behavior of $\psi(\eta)$.

This Monte Carlo method applies, with very little modification, not only to the computation of (\ref{eq:siegelz}), but to the computation of the general formula (\ref{eq:ssz}). It has been used to produce tables of the function $Z(\sigma)$, for various Riemannian symmetric spaces $M$, of rank $N$ up to $30$, which have been successfully used, in numerical computation (recall the rank of $M$ is the dimension of  $\mathfrak{a}$). Unfortunately, this method breaks down, when $N$ is larger (roughly $\approx 50$). Either an analytic expression of $Z(\sigma)$ (see \ref{sec:orthpo}, below), or an asymptotic formula, for large $N$ (see \ref{sec:rmt}, below), are then needed. 

\section{MLE and maximum entropy} \label{sec:gaussmle}
Let $M$ be a Hadamard manifold, which is also a homogeneous space. Propositions \ref{prop:zhadamardk} and \ref{prop:zhomogeneous} then imply that, for any $\bar{x} \in M$ and $\sigma > 0$, there exists a well-defined probability distribution $P(\bar{x}\hspace{0.02cm},\sigma)$ on $M$, given by (\ref{eq:gaussdensity}). The family of distributions $P(\bar{x}\hspace{0.02cm},\sigma)$ fits the definition of Gaussian distributions, stated at the beginning of \ref{sec:rgd}. 
\begin{proposition} \label{prop:mlebarycentre}
 Let $P(\bar{x}\hspace{0.02cm},\sigma)$ be given by (\ref{eq:gaussdensity}), for $\bar{x} \in M$ and $\sigma > 0$. The maximum-likelihood estimate of the parameter $\bar{x}$, based on independent samples $(x_n\,;n=1,\ldots,N)$ from $P(\bar{x}\hspace{0.02cm},\sigma)$, is unique and equal to the empirical barycentre $\hat{x}_{\scriptscriptstyle N}$ of the samples $(x_n)$.
\end{proposition}
The proof of this proposition is immediate. From (\ref{eq:gaussdensity}), one has the log-likelihood function
\begin{equation} \label{eq:rgdll}
\ell(\bar{x}\hspace{0.02cm},\sigma) \,=\, - N\log Z(\sigma) - \frac{1}{2\sigma^2}\hspace{0.03cm}\sum^N_{n=1}\,d^{\hspace{0.03cm}2}(x_n\hspace{0.02cm},\bar{x})
\end{equation}
Since the first term does not depend on $\bar{x}$, one may maximise $\ell(\bar{x}\hspace{0.02cm},\sigma)$, first over $\bar{x}$ and then over $\sigma$. Clearly, maximising over $\bar{x}$ is equivalent to minimising the sum of squared distances $d^{\hspace{0.03cm}2}(x_n\hspace{0.02cm},\bar{x})$. This is just the least-squares problem (\ref{eq:empiricalfrechet}), whose solution is the empirical barycentre $\hat{x}_{\scriptscriptstyle N}\,$. Moreover, $\hat{x}_{\scriptscriptstyle N}$ is unique, since $M$ is a Hadamard manifold (as discussed in \ref{subsec:afsari})

Consider now maximum-likelihood estimation of $\sigma$. This is better carried out in terms of the natural parameter $\eta = (-2\sigma^{\scriptscriptstyle 2})^{\scriptscriptstyle -1}$, or in terms of the moment parameter $\delta = \psi^\prime(\eta)$, where $\psi(\eta) = \log\hspace{0.02cm}Z(\sigma)$ and the prime denotes the derivative.
\begin{proposition} \label{prop:gausscgf}
 The function $\psi(\eta)$, just defined, is a strictly convex function, which maps the half-line $(-\infty,0)$ onto $\mathbb{R}$. The maximum-likelihood estimates of the parameters $\eta$ and $\delta$ are
\begin{equation} \label{eq:mlsigma}
  \hat{\eta}_{\scriptscriptstyle N} = (\psi^\prime)^{\scriptscriptstyle -1}(\hat{\delta}_{\scriptscriptstyle N}) \hspace{0.5cm}\text{and}\hspace{0.5cm}
  \hat{\delta}_{\scriptscriptstyle N} = \frac{1}{N} \sum^N_{n=1}\,d^{\hspace{0.03cm}2}(x_n\hspace{0.02cm},\hat{x}_{\scriptscriptstyle N})
\end{equation}
where $(\psi^\prime)^{\scriptscriptstyle -1}$ denotes the reciprocal function. 
\end{proposition}
The proof of this proposition is given below. For now, note the maximum-entropy property of Gaussian distributions, stated in the following proposition.
\begin{proposition} \label{prop:maxent}
  The Gaussian distribution $P(\bar{x}\hspace{0.02cm},\sigma)$, given by (\ref{eq:gaussdensity}), is the unique distribution on $M$, having maximum Shannon entropy, among all distributions with given barycentre $\bar{x}$ and dispersion $\delta = \mathbb{E}_{\hspace{0.03cm}\scriptscriptstyle x \sim P}[d^{\hspace{0.03cm}\scriptscriptstyle 2}(x\hspace{0.02cm},\bar{x})]$. Its entropy is equal to $\psi^*(\delta)$ where $\psi^*$ is the Legendre transform of  $\psi$. 
\end{proposition}
\noindent \textbf{Proof of Proposition \ref{prop:gausscgf}\,:} denote $\mu$ the image of the distribution $P(\bar{x}\hspace{0.02cm},\sigma)$ under the mapping $x \mapsto d^{\hspace{0.03cm}2}(x\hspace{0.02cm},\bar{x})$. Then, $\psi(\eta)$ is the cumulant generating function of $\mu$,
\begin{equation}
  \psi(\eta) = \log\,\int^\infty_0\hspace{0.01cm}e^{\eta\hspace{0.02cm}s}\mu(ds)
\end{equation}
and is therefore strictly convex. Note from (\ref{eq:zhadamardk}) and (\ref{eq:zcsigma}) that $Z(\sigma) = 0$ when $\sigma = 0$ and $Z(\sigma)$ increases to $+\infty$ when $\sigma$ increases to $+\infty$. Recalling $\eta = (-2\sigma^{\scriptscriptstyle 2})^{\scriptscriptstyle -1}$ and $\psi(\eta) = \log\hspace{0.02cm}Z(\sigma)$, it becomes clear that $\psi$ is (in fact, strictly increasing, and) maps the half-line $(-\infty,0)$ onto $\mathbb{R}$. 

After maximisation with respect to $\bar{x}$, the log-likelihood function (\ref{eq:rgdll}) becomes,
\begin{equation} \label{eq:rdglegendre1}
\ell(\eta) \,=\, N\left\lbrace \eta\hspace{0.03cm}\hat{\delta}_{\scriptscriptstyle N} - \psi(\eta)\right\rbrace 
\end{equation}
which is a strictly concave function. Differentiating, and setting the derivative equal to $0$, directly yields the maximum-likelihood estimates  (\ref{eq:mlsigma}). \vfill\pagebreak
\noindent \textbf{Remark\,:} $\hat{\eta}_{\scriptscriptstyle N}$ in (\ref{eq:mlsigma}) is well-defined, since the range of $\psi^\prime$ is equal to $(0,\infty)$. Indeed, it is possible to use (\ref{eq:integralcomp}), as in the proof of (\ref{eq:zhadamardk}), to show that 
\begin{equation} \label{eq:psihadamardk}
  \psi^\prime_{\scriptscriptstyle 0}(\eta) \,\leq \psi^\prime(\eta) \leq\, \psi^\prime_c(\eta)
\end{equation}
where $\psi_{\scriptscriptstyle 0}(\eta) = \log Z_{\scriptscriptstyle 0}(\sigma)$, and $\psi_{\scriptscriptstyle c}(\eta) = \log Z_{\scriptscriptstyle c}(\sigma)$, with $\kappa = -c^{\hspace{0.02cm}\scriptscriptstyle 2}$ a lower bound on the sectional curvatures of $M$. Precisely, (\ref{eq:psihadamardk}) can be obtained by replacing $f(y) = d^{\hspace{0.03cm}\scriptscriptstyle 2}(y\hspace{0.02cm},\bar{x})\hspace{0.03cm}p(y|\bar{x}\hspace{0.02cm},\sigma)$ into (\ref{eq:integralcomp}), where $p(y|\bar{x}\hspace{0.02cm},\sigma)$ is the probability density function in (\ref{eq:gaussdensity}). Now, $\psi^\prime_{\scriptscriptstyle 0}(\eta) = n\hspace{0.02cm}\sigma^{\scriptscriptstyle 2}$, which increases to $+\infty$ when $\sigma$ increases to $+\infty$. On the other hand, by a straightforward application of the chain rule, it is seen that
\begin{equation} \label{eq:psicsigma}
 \psi^\prime_c(\eta) \,=\, \sigma^3\hspace{0.03cm}\frac{d}{d\sigma}\!\left(\log Z_{\scriptscriptstyle c}(\sigma)\right)
\end{equation}
which, from (\ref{eq:zcsigma}), is $=0$ when $\sigma = 0$. Now, it follows from (\ref{eq:psihadamardk}), $\psi^\prime$ maps the half-line $(-\infty,0)$ onto the half-line $(0,+\infty)$. \\[0.1cm]
\noindent \textbf{Proof of Proposition \ref{prop:maxent}\,:} let $Q(dx)$ be a probability distribution on $M$ with barycentre $\bar{x}$ and dispersion $\delta = \mathbb{E}_{\hspace{0.03cm}\scriptscriptstyle x \sim Q}[d^{\hspace{0.03cm}\scriptscriptstyle 2}(x\hspace{0.02cm},\bar{x})]$. Assume $Q(dx)$ has probablity density function $q(x)$, with respect to Riemannian volume. The Shannon entropy of $Q$ is given by
\begin{equation} \label{eq:shannon}
 S(q) \,=\, \int_M\,\log(q(x))\hspace{0.02cm}q(x)\hspace{0.02cm}\mathrm{vol}(dx)   
\end{equation}
Since $M$ is a homogeneous space, $S(q)$ does not depend on $\bar{x}$. Fixing some point $o \in M$, it is possible to assume, without loss of generality, that $\bar{x} = o$. Then, it is enough to maximise $S(q)$, subject to the constraints,
$$
\int_M\,q(x)\hspace{0.02cm}\mathrm{vol}(dx) \,=\,1 \hspace{0.4cm}\text{ and }\hspace{0.3cm} 
\int_M\,d^{\hspace{0.03cm}\scriptscriptstyle 2}(x\hspace{0.02cm},o)\hspace{0.03cm}q(x)\hspace{0.02cm}\mathrm{vol}(dx) \,=\, \delta
$$
Using the method of Lagrange multipliers, this leads to a stationary point 
\begin{equation} \label{eq:maxent1}
q(x) \,=\, \exp\left(\hspace{0.02cm}\eta\hspace{0.03cm}d^{\hspace{0.03cm}\scriptscriptstyle 2}(x\hspace{0.02cm},o) - \psi(\eta)\right) 
\end{equation}
where the Lagrange multiplier $\eta$ is finally given by $\eta = (\psi^\prime)^{\scriptscriptstyle -1}(\delta)$, in terms of the cumulant generating function,
$$
\psi(\eta) = \log\,\int_M\exp\left(\eta\hspace{0.03cm}d^{\hspace{0.03cm}\scriptscriptstyle 2}(x\hspace{0.02cm},o)\right)\mathrm{vol}(dx)
$$
Of course, $q(x)$ in (\ref{eq:maxent1}) is just $p(x|o\hspace{0.02cm},\sigma)$, once the parameter $\sigma > 0$ is defined by $\eta = (-2\sigma^{\scriptscriptstyle 2})^{\scriptscriptstyle -1}$. Since the Shannon entropy is strictly concave, this stationary point $q(x)$ is a unique maximum, over the (convex) set of probability density functions on $M$, which satisfy the above constraints. Its entropy is equal to
\begin{equation} \label{eq:maxent2}
S(q) \,=\,\int_M\left(\hspace{0.02cm}\eta\hspace{0.03cm}d^{\hspace{0.03cm}\scriptscriptstyle 2}(x\hspace{0.02cm},o) - \psi(\eta)\right)q(x)\hspace{0.02cm}\mathrm{vol}(dx) = \eta\hspace{0.02cm}\delta - \psi(\eta)
\end{equation}
To show that this is $\psi^*(\delta)$, as stated in the proposition, it is enough to show
\begin{equation} \label{eq:maxent3}
 S(q) \,=\, \sup_\eta\hspace{0.03cm}\lbrace \eta\hspace{0.02cm}\delta - \psi(\eta)\rbrace
\end{equation}
However, since $\psi$ is a strictly convex function, it is seen by differentiation that the $\sup$ is achieved when $\psi^\prime(\eta) = \delta$, exactly as in (\ref{eq:maxent1}). Accordingly, the right-hand side of (\ref{eq:maxent3}) is equal to $\eta\hspace{0.02cm}\delta - \psi(\eta)$, as in (\ref{eq:maxent2}). 

\section{Barycentre and covariance} \label{sec:gaussschur}

\subsection{The Riemannian barycentre}
Let $M$ be a Hadamard manifold, which is also a homogeneous space. Here, it is shown that the barycentre of the Gaussian distribution $P(\bar{x}\hspace{0.03cm},\sigma)$ on $M$, given by (\ref{eq:gaussdensity}), is equal to $\bar{x}$. 

First, it should be noted $P(\bar{x}\hspace{0.03cm},\sigma)$ does indeed have a well-defined Riemannian barycentre, since it has finite second-order moments. To see that this is true, it is enough to note that
$$
\int_M\,d^{\hspace{0.03cm}\scriptscriptstyle 2}(\bar{x}\hspace{0.02cm},x)\hspace{0.03cm}p(x|\bar{x}\hspace{0.02cm},\sigma)\hspace{0.03cm}\mathrm{vol}(dx) \,<\, \infty
$$
Ineded, this integral is just $\psi^\prime(\eta)$ in (\ref{eq:psihadamardk}). This means $\pi = P(\bar{x}\hspace{0.03cm},\sigma)$ satisfies (\ref{eq:secondordermoment}) for $y_o = \bar{x}$. 
\begin{proposition} \label{prop:gaussbarycentre}
Let $P(\bar{x}\hspace{0.03cm},\sigma)$ be given by (\ref{eq:gaussdensity}), for $\bar{x} \in M$ and $\sigma > 0$. The Riemannian barycentre of 
$P(\bar{x}\hspace{0.03cm},\sigma)$ is equal to $\bar{x}$.
\end{proposition}
\noindent \textbf{First proof\,:} the proof of this proposition relies on the fact that the variance function,
$$
\mathcal{E}(y) \,=\, \frac{1}{2}\hspace{0.03cm}
\int_M\,d^{\hspace{0.03cm}\scriptscriptstyle 2}(y\hspace{0.02cm},x)\hspace{0.03cm}p(x|\bar{x}\hspace{0.02cm},\sigma)\hspace{0.03cm}\mathrm{vol}(dx)
$$
is $1/2$-strongly convex. In particular, it has a unique stationary point, $\hat{x}$ with $\mathrm{grad}\,\mathcal{E}(\hat{x}) = 0$, which is also its unique global minimum, and (by definition) the Riemannian barycentre of $P(\bar{x}\hspace{0.03cm},\sigma)$. Now, let $f(\bar{x})$ be the function given by
$$
f(\bar{x}) \,=\,\int_M\,p(x|\bar{x}\hspace{0.02cm},\sigma)\hspace{0.03cm}\mathrm{vol}(dx) 
$$
Clearly, this is a constant function, equal to $1$ for all $\bar{x}$. On the other hand, its gradient may be written down, by differentiating under the integral, with respect to $\bar{x}$, using (\ref{eq:gradfx}) and (\ref{eq:gaussdensity}),
$$
\mathrm{grad}\,f(\bar{x}) \,=\, \sigma^{-2}\hspace{0.03cm}\int_M\,\mathrm{Exp}^{-1}_{\bar{x}}(x)\hspace{0.04cm}
p(x|\bar{x}\hspace{0.02cm},\sigma)\hspace{0.03cm}\mathrm{vol}(dx)
$$
Now, $\mathrm{grad}\,f(\bar{x})$ is identically zero. But, the right-hand side of the above expression is equal to $-\sigma^{-2}\hspace{0.03cm}\mathrm{grad}\,\mathcal{E}(\bar{x})$, by (\ref{eq:gradepsilonhadamard}). This shows that $\mathrm{grad}\,\mathcal{E}(\bar{x}) = 0$, and therefore $\bar{x}$ is the Riemannian barycentre of $P(\bar{x}\hspace{0.03cm},\sigma)$. \\[0.1cm]
\textbf{Second proof\,:} this proof works if $M$ is a Riemannnian symmetric space which belongs to the non-compact case. From (\ref{eq:gradepsilonhadamard}),
$$
\mathrm{grad}\,\mathcal{E}(\bar{x}) \,=\,-\hspace{0.03cm}\int_M\,\mathrm{Exp}^{-1}_{\bar{x}}(x)\hspace{0.04cm}
p(x|\bar{x}\hspace{0.02cm},\sigma)\hspace{0.03cm}\mathrm{vol}(dx)
$$
Let $s_{\bar{x}}$ be the geodesic symmetry at $\bar{x}$. From the definition of $s_{\bar{x}\,}$, $s_{\bar{x}}\cdot \mathrm{grad}\,\mathcal{E}(\bar{x}) = -\hspace{0.03cm}\mathrm{grad}\,\mathcal{E}(\bar{x})$. On the other hand, 
$$
s_{\bar{x}}\cdot \mathrm{grad}\,\mathcal{E}(\bar{x}) \,=\,
-\hspace{0.03cm}\int_M\left(s_{\bar{x}}\cdot\mathrm{Exp}^{-1}_{\bar{x}}(x)\right)
p(x|\bar{x}\hspace{0.02cm},\sigma)\hspace{0.03cm}\mathrm{vol}(dx)
$$
Since $s_{\bar{x}}$ is an isometry and fixes $\bar{x}$, it follows that 
$$
s_{\bar{x}}\cdot\mathrm{Exp}^{-1}_{\bar{x}}(x) =
\mathrm{Exp}^{-1}_{\bar{x}}(s_{\bar{x}}\cdot x) \text{ and } 
p(x|\bar{x}\hspace{0.02cm},\sigma) = 
p(s_{\bar{x}}\cdot x|\bar{x}\hspace{0.02cm},\sigma) 
$$
Therefore,
$$
s_{\bar{x}}\cdot \mathrm{grad}\,\mathcal{E}(\bar{x}) \,=\,
-\hspace{0.03cm}\int_M\,\mathrm{Exp}^{-1}_{\bar{x}}(s_{\bar{x}}\cdot x)\hspace{0.03cm}
p(s_{\bar{x}}\cdot x|\bar{x}\hspace{0.02cm},\sigma)\hspace{0.03cm}\mathrm{vol}(dx)
$$
and, introducing the variable of integration $z = s_{\bar{x}}\cdot x$, it follows that $s_{\bar{x}}\cdot \mathrm{grad}\,\mathcal{E}(\bar{x}) = 
\mathrm{grad}\,\mathcal{E}(\bar{x})$.
Now, it has been shown that $s_{\bar{x}}\cdot \mathrm{grad}\,\mathcal{E}(\bar{x}) = -\hspace{0.03cm}\mathrm{grad}\,\mathcal{E}(\bar{x})$ and that $s_{\bar{x}}\cdot \mathrm{grad}\,\mathcal{E}(\bar{x}) = 
\mathrm{grad}\,\mathcal{E}(\bar{x})$. Thus, $\mathrm{grad}\,\mathcal{E}(\bar{x}) = 0$ and one may conclude as in the first proof. 

%$d\hspace{0.02cm}g^{\scriptscriptstyle -1}\cdot u$ is denoted $g^{\scriptscriptstyle -1}\cdot u$.

%existence and uniqueness

\subsection{The covariance tensor}
The covariance form of the distribution $P(\bar{x}\hspace{0.03cm},\sigma)$ is the symmetric bilinear form $C_{\bar{x}}$ on $T_{\bar{x}}M$, 
\begin{equation} \label{eq:covariance}
  C_{\bar{x}\hspace{0.02cm}}(u\hspace{0.02cm},v) \,=\, \int_M\,\langle u\hspace{0.03cm},\mathrm{Exp}^{-1}_{\bar{x}}(x)\rangle\hspace{0.03cm}
                                                         \langle \mathrm{Exp}^{-1}_{\bar{x}}(x),v\rangle \,
                                                         p(x|\bar{x}\hspace{0.02cm},\sigma)\hspace{0.03cm}\mathrm{vol}(dx) \hspace{1cm}
u\,,v \in T_{\bar{x}}M
\end{equation}
With $\sigma > 0$ fixed, the map which assigns to $\bar{x} \in M$ the covariance form $C_{\bar{x}}$ is a (0,2)-tensor field on $M$, here called the covariance tensor of $P(\bar{x}\hspace{0.03cm},\sigma)$. In order to compute this tensor field, consider the following situation. 

Assume $M = G/K$ is a Riemannian symmetric space which belongs to the non-compact case. Here, $K = K_o\hspace{0.04cm}$, the stabiliser in $G$ of $o \in M$. For $k \in K$ and $u \in T_oM$, it is clear $k\cdot u \in T_oM$. This defines a representation of $K$ in the tangent space $T_oM$, called the isotropy representation. One says that $M$ is an irreducible symmetric space, if this isotropy representation is irreducible.

If $M$ is not irreducible, then it is a product of irreducible Riemannian symmetric spaces $M = M_{\scriptscriptstyle 1}\times\ldots\times M_{ s}$~\cite{helgason} (Proposition 5.5, Chapter VIII. This is the de Rham decomposition of $M$). 
Accordingly, for $x \in M$ and $u \in T_x M$, one may write $x = (x_{\scriptscriptstyle 1},\ldots,x_s)$ and $u = (u_{\scriptscriptstyle 1},\ldots,u_s)$, where $x_r \in M_r$ and $u_r \in T_{x_r}M_r\hspace{0.04cm}$. Now, looking back at  (\ref{eq:gaussdensity}), it may be seen that
\begin{equation} \label{eq:gdensitydrh}
p(x|\bar{x}\hspace{0.02cm},\sigma) \,=\, \prod^s_{r=1}p(x_r|\bar{x}_r\hspace{0.02cm},\sigma) \hspace{1cm} 
p(x_r|\bar{x}_r\hspace{0.02cm},\sigma) \,=\, 
\left(Z_{r}(\sigma)\right)^{-1}\hspace{0.03cm}\exp\left[ -\frac{d^{\hspace{0.03cm}2}(x_r\hspace{0.02cm},\bar{x}_r)}{2\sigma^2}\right]
\end{equation}
For the following proposition, let $\eta = (-2\sigma^{\scriptscriptstyle 2})^{\scriptscriptstyle -1}$ and $\psi_r(\eta) = \log Z_r(\sigma)$.
\begin{proposition} \label{prop:gausscovariance}
Assume that $M$ is a product of irreducible Riemannian symmetric spaces,  $M = M_{\scriptscriptstyle 1}\times\ldots\times M_{ s\hspace{0.04cm}}$. The covariance tensor $C$ in (\ref{eq:covariance}) is given by
\begin{equation} \label{eq:gausscovariance}
C_{\bar{x}}(u\hspace{0.02cm},u) \,=\, \sum^s_{r=1}\frac{\psi^\prime_r(\eta)}{\dim\hspace{0.03cm} M_r}\hspace{0.04cm} \Vert u_r \Vert^2_{\bar{x}_r}
\end{equation}
for $u \in T_{\bar{x}}M$ where $\bar{x} = (\bar{x}_{\scriptscriptstyle 1},\ldots,\bar{x}_s)$ and $u = (u_{\scriptscriptstyle 1},\ldots,u_s)$, with $\bar{x}_r \in M_r$ and $u_r \in T_{\bar{x}_r}M_r\hspace{0.04cm}$.  
\end{proposition}
\noindent \textbf{Example\,:} let $M = \mathrm{H}(N)$, so $M = \mathrm{GL}(N,\mathbb{C})/U(N)$, with $U(N)$ the stabiliser of $o = \mathrm{I}_N\,$.\hfill\linebreak The de Rham decomposition of $M$ is $M = M_{\scriptscriptstyle 1}\times M_{\scriptscriptstyle 2\hspace{0.04cm}}$, where $M_{\scriptscriptstyle 1} = \mathbb{R}$ and $M_{\scriptscriptstyle 2}$ is the submanifold whose elements are those $x \in M$ such that $\det(x) = 1$. Accordingly, each $\bar{x} \in M$ is identified with the couple $(\bar{x}_{\scriptscriptstyle 1}\hspace{0.02cm},\bar{x}_{\scriptscriptstyle 2})$, 
$$
\bar{x}_{\scriptscriptstyle 1} \,=\, \frac{1}{N}\log\det(\bar{x}) \hspace{0.5cm} 
\bar{x}_{\scriptscriptstyle 2} \,=\, (\det(\bar{x}))^{-{\scriptscriptstyle 1/N}}\hspace{0.03cm} \bar{x}
$$
and each $u \in T_{\bar{x}}M$ is written $u = u_{\scriptscriptstyle 1}\hspace{0.02cm}\bar{x} + u_{\scriptscriptstyle 2}$
$$
u_{\scriptscriptstyle 1} \,=\, \frac{1}{N}\hspace{0.03cm}\mathrm{tr}(\bar{x}^{\scriptscriptstyle -1}\hspace{0.02cm}u) \hspace{0.5cm}
u_{\scriptscriptstyle 2} \,=\, u - \frac{1}{N}\hspace{0.03cm}\mathrm{tr}(\bar{x}^{\scriptscriptstyle -1}\hspace{0.02cm}u)\hspace{0.04cm} \bar{x}
$$
These may be replaced into expression (\ref{eq:gausscovariance}), 
\begin{equation} \label{eq:gausshnprefim0}
C_{\bar{x}}(u\hspace{0.02cm},u) \,=\, \psi^\prime_{\scriptscriptstyle 1}(\eta)\hspace{0.03cm} u^2_{\scriptscriptstyle 1} \,+\, \frac{\psi^\prime_{\scriptscriptstyle 2}(\eta)}{N^{\scriptscriptstyle 2} - 1}\hspace{0.02cm}\Vert u_{\scriptscriptstyle 2}\Vert^2_{\bar{x}_{\scriptscriptstyle 2}}
\end{equation}
where $\psi_{\scriptscriptstyle 1}(\eta) = \log \left( 2\pi\hspace{0.03cm}\sigma^2\right)^{\!\frac{1}{2}}$, and $\psi_{\scriptscriptstyle 2}(\eta) = \log Z(\sigma) - \psi_{\scriptscriptstyle 1}(\eta)$ ($Z(\sigma)$ is given by (\ref{eq:sw_z}) in \ref{sec:orthpo}, below). After a direct calculation, this can be brought under the form
\begin{equation} \label{eq:gausshnprefim}
C_{\bar{x}}(u\hspace{0.02cm},u) \,=\, g_{\scriptscriptstyle 1}(\sigma)\hspace{0.03cm} \mathrm{tr}^2(\bar{x}^{\scriptscriptstyle -1}\hspace{0.02cm}u) +
g_{\scriptscriptstyle 2}(\sigma)\hspace{0.03cm} \mathrm{tr}(\bar{x}^{\scriptscriptstyle -1}\hspace{0.02cm}u)^2
\end{equation}
where $g_{\scriptscriptstyle 1}(\sigma)$ and $g_{\scriptscriptstyle 2}(\sigma) $ are certain functions of $\sigma$. \\[0.1cm]
\noindent \textbf{Remark\,:} as a corollary of Proposition \ref{prop:gausscovariance}, the covariance tensor $C$ is a $G$-invariant Riemannian metric on $M$. This is clear, for example, in the special case of (\ref{eq:gausshnprefim}), which coincides with the general expression of a $\mathrm{GL}(N,\mathbb{C})$-invariant metric. \vfill\pagebreak
\noindent \textbf{Proof of Proposition \ref{prop:gausscovariance}\,:} since $C_{\bar{x}}$ is bilinear 
\begin{equation} \label{eq:proofcovariance}
C_{\bar{x}}(u\hspace{0.02cm},u) \,=\, \sum^s_{r=1}\sum^s_{q=1}\, C_{\bar{x}}(u_r\hspace{0.02cm},u_q)
\end{equation}
It will be shown that
\begin{equation} \label{eq:proofcovariance1}
C_{\bar{x}}(u_r\hspace{0.02cm},u_q) = 0 \hspace{0.5cm} \text{for } r\neq q
\end{equation}
and, on the other hand, that
\begin{equation} \label{eq:proofcovariance2}
C_{\bar{x}}(u_r\hspace{0.02cm},u_r) = \frac{\psi^\prime_r(\eta)}{\dim\hspace{0.03cm} M_r}\hspace{0.04cm} \Vert u_r \Vert^2_{\bar{x}_r}
\end{equation}
Then, (\ref{eq:gausscovariance}) will follow immediately, by replacing (\ref{eq:proofcovariance1}) and (\ref{eq:proofcovariance2}) into (\ref{eq:proofcovariance}). \\[0.1cm]
\textbf{Proof of (\ref{eq:proofcovariance1})\,:} from (\ref{eq:covariance}),
\begin{equation} \label{eq:proofproofcovariance1}
C_{\bar{x}}(u_r\hspace{0.02cm},u_q) \,=\, \int_M\,\langle u_r\hspace{0.03cm},\mathrm{Exp}^{-1}_{\bar{x}}(x)\rangle\hspace{0.03cm}
                                                         \langle \mathrm{Exp}^{-1}_{\bar{x}}(x),u_q\rangle \,
                                                         p(x|\bar{x}\hspace{0.02cm},\sigma)\hspace{0.03cm}\mathrm{vol}(dx)
\end{equation}
However, since $M$ is given as a product Riemannian manifold,
\begin{equation} \label{eq:proofproofcovariance11}
\langle u_r\hspace{0.03cm},\mathrm{Exp}^{-1}_{\bar{x}}(x)\rangle \,=\,
\langle u_r\hspace{0.03cm},\mathrm{Exp}^{-1}_{\bar{x}_r}(x_r)\rangle \hspace{0.2cm}\text{and}\hspace{0.2cm}
\langle u_q\hspace{0.03cm},\mathrm{Exp}^{-1}_{\bar{x}}(x)\rangle \,=\,
\langle u_q\hspace{0.03cm},\mathrm{Exp}^{-1}_{\bar{x}_q}(x_q)\rangle
\end{equation}
Using (\ref{eq:gdensitydrh}) and (\ref{eq:proofproofcovariance11}), it follows from (\ref{eq:proofproofcovariance1}) that
$$
\begin{array}{rl}
C_{\bar{x}}(u_r\hspace{0.02cm},u_q) =& \int_{\scriptscriptstyle M_r}\langle u_r\hspace{0.03cm},\mathrm{Exp}^{-1}_{\bar{x}_r}(x_r)\rangle\hspace{0.03cm}
                                                         p(x_r|\bar{x}_r\hspace{0.02cm},\sigma)\hspace{0.03cm}\mathrm{vol}(dx_r)
\int_{\scriptscriptstyle M_q}\langle u_q\hspace{0.03cm},\mathrm{Exp}^{-1}_{\bar{x}_q}(x_q)\rangle\hspace{0.03cm}
                                                         p(x_q|\bar{x}_q\hspace{0.02cm},\sigma)\hspace{0.03cm}\mathrm{vol}(dx_q) \\[0.2cm]
=& \mathrm{grad}\,\mathcal{E}_r(\bar{x}_r)\,\mathrm{grad}\,\mathcal{E}_q(\bar{x}_q) \\[0.2cm]
=& 0
\end{array}
$$
where the second equality follows from (\ref{eq:gradepsilonhadamard}), applied to the variance functions
$$
\mathcal{E}_r(y) \,=\, \frac{1}{2}\hspace{0.03cm}
\int_{M_r}d^{\hspace{0.03cm}\scriptscriptstyle 2}(y\hspace{0.02cm},x_r)\hspace{0.03cm}p(x_r|\bar{x}_r\hspace{0.02cm},\sigma)\hspace{0.03cm}\mathrm{vol}(dx_r)
\hspace{0.1cm}\text{and}\hspace{0.1cm}
\mathcal{E}_q(y) \,=\, \frac{1}{2}\hspace{0.03cm}
\int_{M_q}d^{\hspace{0.03cm}\scriptscriptstyle 2}(y\hspace{0.02cm},x_q)\hspace{0.03cm}p(x_q|\bar{x}_q\hspace{0.02cm},\sigma)\hspace{0.03cm}\mathrm{vol}(dx_q)
$$
which, by Proposition \ref{prop:gaussbarycentre}, respectively have their global minima at $\bar{x}_r$ and $\bar{x}_q\hspace{0.04cm}$. \\[0.1cm]
\textbf{Proof of (\ref{eq:proofcovariance2})\,:} let $K_{\bar{x}}$ denote the stabiliser of $\bar{x}$ in $G$. For $k \in K_{\bar{x}}$ and $u_r \in T_{\bar{x}_r}M_r\hspace{0.04cm}$, note that $k\cdot u_r \in T_{\bar{x}_r}M_r\hspace{0.04cm}$. This defines an \textit{irreducible} representation of $K_{\bar{x}}$ in $T_{\bar{x}_r}M_r\hspace{0.04cm}$. The symmetric bilinear form $C_{\bar{x}}$ is invariant under this representation. Precisely, since any $k \in K_{\bar{x}}$ is an isometry which fixes $\bar{x}$, it follows from (\ref{eq:covariance}),
$$
\begin{array}{rll}
C_{\bar{x}}(k\cdot u_r\hspace{0.02cm},k\cdot u_r) = &
\int_{\scriptscriptstyle M_r}\,\langle k\cdot u_r\hspace{0.03cm},\mathrm{Exp}^{-1}_{\bar{x}_r}(x_r)\rangle^2\,                                                      
                                                         p(x_r|\bar{x}_r\hspace{0.02cm},\sigma)\hspace{0.03cm}\mathrm{vol}(dx_r) & \\[0.2cm]
=&
\int_{\scriptscriptstyle M_r}\,\langle u_r\hspace{0.03cm},\mathrm{Exp}^{-1}_{\bar{x}_r}(k^{\scriptscriptstyle -1}\cdot x_r)\rangle^2\,                                                      
                                                         p(x_r|\bar{x}_r\hspace{0.02cm},\sigma)\hspace{0.03cm}\mathrm{vol}(dx_r) & \\[0.2cm]
=& \int_{\scriptscriptstyle M_r}\,\langle u_r\hspace{0.03cm},\mathrm{Exp}^{-1}_{\bar{x}_r}(k^{\scriptscriptstyle -1}\cdot x_r)\rangle^2\,                                                      
                                                         p(k^{\scriptscriptstyle -1}\cdot x_r|\bar{x}_r\hspace{0.02cm},\sigma)\hspace{0.03cm}\mathrm{vol}(dx_r)  =&\!\!
C_{\bar{x}}(u_r\hspace{0.02cm},u_r)
\end{array}
$$
where the last equality follows by introducing the new variable of integration $z = k^{\scriptscriptstyle -1}\cdot x_r\hspace{0.04cm}$. Finally, from Schur's lemma~\cite{knapp}, $C_{\bar{x}}$ is a multiple of the metric,
$$
C_{\bar{x}}(u_r\hspace{0.02cm},u_r) \,=\, f(\eta)\,\hspace{0.04cm} \Vert u_r \Vert^2_{\bar{x}_r}
$$
where $f(\eta)$ may be found from $\mathrm{tr}(C_{\bar{x}})\,=\,(\dim\hspace{0.02cm}M_r)\hspace{0.02cm}f(\eta)$. To conclude, it is enough to note that
the trace may be evaluated by introducing an orthonormal basis of $T_{\bar{x}_r}M_r\hspace{0.04cm}$. It then follows that,
$$
\mathrm{tr}(C_{\bar{x}}) \,=\, 
\int_{M_r}\,\Vert \mathrm{Exp}^{-1}_{\bar{x}_r}(x_r)\Vert^2\hspace{0.03cm}p(x_r|\bar{x}_r\hspace{0.02cm},\sigma)\hspace{0.03cm}\mathrm{vol}(dx_r)
\,=\,\int_{M_r}\,d^{\hspace{0.03cm}\scriptscriptstyle 2}(\bar{x}_r\hspace{0.02cm},\hspace{0.03cm}x_r)\hspace{0.03cm}p(x_r|\bar{x}_r\hspace{0.02cm},\sigma)\hspace{0.03cm}\mathrm{vol}(dx_r)
$$  
which is equal to  $\psi^\prime_r(\eta)$, by the same argument as in the discussion before Proposition \ref{prop:gaussbarycentre}.

\vfill
\pagebreak

\section{An analytic formula for $Z(\sigma)$} \label{sec:orthpo}
Consider the special case where $M = \mathrm{H}(N)$, which corresponds to $\beta = 2$ in Example 2 of \ref{sec:z}. In this case, using the tools of random matrix theory (see~\cite{mehta}, Chapter 5), it is possible to provide an analytic formula for the normalising factor $Z(\sigma)$.
\begin{proposition} \label{prop:sw_z}
 When $M = \mathrm{H}(N)$, the normalising factor $Z(\sigma)$, given by (\ref{eq:covzvv}) with $\beta = 2$, admits of the following analytic formula
\begin{equation} \label{eq:sw_z}
  Z(\sigma) = 
\frac{\omega_{\scriptscriptstyle 2}(N)}{\mathstrut 2^{\scriptscriptstyle N^2}}\left( 2\pi\hspace{0.03cm}\sigma^2\right)^{\!\frac{N}{2}}\hspace{0.02cm} \exp\left[{\small\left(\frac{N^3 - N}{6}\right)}\sigma^{2}\right] \prod^{N-1}_{n=1}\left(1 - e^{-n\hspace{0.02cm}\sigma^2} \right)^{\!N-n}
\end{equation}
\end{proposition}
\noindent \textbf{Remark\,:} when $N = 2$, (\ref{eq:sw_z}) reduces to 
\begin{equation} \label{eq:sw_z_2}
  Z(\sigma) = \left(\frac{\pi\hspace{0.02cm}\sigma}{2}\right)^{\!2}\left( e^{\sigma^2} - 1\right)
\end{equation}
which can be checked, by directly calculating the integral (\ref{eq:covzvv}). \\[0.1cm]
\textbf{Proof of Proposition \ref{prop:sw_z}\,:} putting $\beta = 2$ in (\ref{eq:covzvv}), and noting that $N_{\scriptscriptstyle 2} = N$, it follows that
\begin{equation} \label{eq:sw_z_proof}
Z(\sigma) \,=\,
\frac{\omega_{\scriptscriptstyle 2}(N)}{\mathstrut 2^{\scriptscriptstyle N^2}N!}\hspace{0.02cm} \exp\left[-\frac{N^3}{2}\hspace{0.04cm}\sigma^{2}\right]\times I_{\scriptscriptstyle 2}
\end{equation}
where $I_{\scriptscriptstyle 2}$ is the integral
\begin{equation} \label{eq:i2sw}
I_{\scriptscriptstyle 2} = \int_{\mathbb{R}^{\scriptscriptstyle N}_+}\,\prod^N_{i=1}\rho(u_i\hspace{0.02cm},2\sigma^{\scriptscriptstyle 2})\,|V(u)|^2\,\prod^N_{i=1} \hspace{0.02cm}du_i
\end{equation}
This can be expressed using a well-known formula from random matrix theory~\cite{mehta} (Chapter 5, Page 79). Precisely, if $(p_{\hspace{0.02cm}n}\,; n = 0,1,\ldots)$ are orthonormal polynomials, with respect to the weight function $\rho(u\hspace{0.02cm},2\sigma^{\scriptscriptstyle 2})$ on $\mathbb{R}_+\hspace{0.04cm}$, then $I_{\scriptscriptstyle 2}$ is given by
\begin{equation} \label{eq:mehtasw}
I_{\scriptscriptstyle 2} \,=\, N!\hspace{0.03cm}\prod^{N-1}_{n=0} p^{-2}_{\hspace{0.02cm}nn}
\end{equation}
where $p_{\hspace{0.02cm}nn}$ is the leading coefficient in $p_{\hspace{0.02cm}n\hspace{0.04cm}}$. The required orthonormal polynomials $p_{\hspace{0.02cm}n}$ are given by $p_{\hspace{0.02cm}n} = (2\pi\sigma^2)^{-\frac{1}{4}}\hspace{0.03cm}s_{\hspace{0.02cm}n\hspace{0.04cm}}$, where $s_{\hspace{0.02cm}n}$ are the Stieltjes-Wigert polynomials~\cite{szegobook} (Page 33). Accordingly,
$$
p^{-2}_{\hspace{0.02cm}nn} \,=\, \left( 2\pi\hspace{0.03cm}\sigma^2\right)^{\!\frac{1}{2}}\hspace{0.02cm}\exp\left[\frac{(2n+1)^2}{2}\hspace{0.04cm}\sigma^2\right] \prod^{n}_{m=1}\left(1 - e^{-m\hspace{0.02cm}\sigma^2} \right)
$$
Then, working out the product (\ref{eq:mehtasw}), it easily follows
\begin{equation} \label{eq:i2final}
I_{\scriptscriptstyle 2} \,=\, N!\hspace{0.03cm}\left( 2\pi\hspace{0.03cm}\sigma^2\right)^{\!\frac{N}{2}}\hspace{0.02cm} \exp\left[{\small\left(\frac{4N^3-N}{6}\right)}\hspace{0.04cm}\sigma^{2}\right] \prod^{N-1}_{n=1}\left(1 - e^{-n\hspace{0.02cm}\sigma^2} \right)^{\!N-n}
\end{equation}
and (\ref{eq:sw_z}) may be obtained by replacing this into (\ref{eq:sw_z_proof}). \\[0.1cm]
\textbf{Remark\,:} the product appearing in (\ref{eq:i2final}) can be written as a product of $q$-Gamma functions. Letting $q = e^{-\sigma^2}$, and recalling the definition of the $q$-Gamma function~\cite{jackson}, it may be seen that
\begin{equation}
\prod^{N-1}_{n=1}\left(1 - e^{-n\hspace{0.02cm}\sigma^2} \right)^{\!N-n} \,=\, (1-q)^{\scriptscriptstyle (N^2-N)/2}\hspace{0.02cm}\prod^{N}_{n=2}\Gamma_{q\hspace{0.02cm}}(n) \hspace{1cm}\text{($\Gamma_q$ the $q$-Gamma function)}
\end{equation}
In other words, the product of $q$-Gamma functions plays, for the present problem, the same role that the product of classical Gamma functions (known as the Barnes function) plays, for the Gaussian unitary ensemble. 

\vfill
\pagebreak

\section{Large $N$ asymptotics} \label{sec:rmt}
Pursuing the development started in \ref{sec:orthpo}, it is possible to derive an asymptotic expression of $Z(\sigma)$, valid in the limit where $N$ goes to infinity, while the product $t = N\sigma^2$ remains constant. 
\begin{proposition} \label{prop:asymp_z}
 Let $Z(\sigma)$ be given by (\ref{eq:sw_z}). If $N \rightarrow \infty$, while $t = N\sigma^2$ remains constant, then the following equivalence holds,
\begin{equation} \label{eq:asymp_z}
  \frac{1}{N^{\scriptscriptstyle 2}}\hspace{0.03cm}\log Z(\sigma) \sim -\frac{1}{2}\log\left(\frac{2N}{\pi}\right) + \frac{3}{4} + \frac{t}{6} - \frac{\mathrm{Li}_{\scriptscriptstyle 3}(e^{-t}) - \zeta(3)}{t^{\scriptscriptstyle 2}}
\end{equation}
where $\mathrm{Li}_{\scriptscriptstyle 3}(x) = \sum^\infty_{k=1} x^k/k^{\scriptscriptstyle 3}$ for $|x| < 1$ (the trilogarithm), and $\zeta$ is the Riemann Zeta function.
\end{proposition}
The proposition follows by a direct calculation, once the following lemmas have been shown. 
\begin{lemma} \label{lemma_z_asymp_1}
  In the notation of (\ref{eq:sw_z}), if $N \rightarrow \infty$,
\begin{equation} \label{eq:lemma_z_asymp_1}
  \frac{1}{N^{\scriptscriptstyle 2}}\hspace{0.03cm}\log \omega_{\scriptscriptstyle 2}(N) \sim -  \frac{1}{2}\log\left(\frac{N}{2\pi}\right) + \frac{3}{4}
\end{equation}
\end{lemma}
\begin{lemma} \label{lemma_z_asymp_2}
 If $N \rightarrow \infty$, while $t = N\sigma^2$ remains constant, then
\begin{equation} \label{eq:lemma_z_asymp_2}
  \lim\hspace{0.04cm} \frac{1}{N^{\scriptscriptstyle 2}}\hspace{0.03cm}\log \prod^{N-1}_{n=1}\left(1 - e^{-n\hspace{0.02cm}\sigma^2} \right)^{\!N-n} \,=\, \int^{\scriptscriptstyle 1}_{\scriptscriptstyle 0}(1-x)\log\left(1 - e^{-tx}\right)\hspace{0.01cm}dx
\end{equation}
and this improper integral is equal to $-(\mathrm{Li}_{\scriptscriptstyle 3}(e^{-t}) - \zeta(3))/t^{\scriptscriptstyle 2}$.
\end{lemma}
\noindent \textbf{Proof of Lemma \ref{lemma_z_asymp_1}\,:} recall, from the footnote in \ref{sec:z}, that 
$$
\omega_{\scriptscriptstyle 2}(N) \,=\, (2\pi)^{(N^2-N)/2}/G(N) \hspace{0.4cm}\text{where } G(N) = 1!\times2!\times\ldots\times(N-1)!
$$
Then, from the asymptotic formula of the Barnes function~\cite{watson} (see Chapter XII, Exercice 49)
$$
\log \omega_{\scriptscriptstyle 2}(N) \,=\, \frac{N^2}{2}\hspace{0.03cm}\log(2\pi) - N^2\left[ \frac{1}{2}\log(N) -\frac{3}{4}\right] + o(N^2)
$$
which directly implies (\ref{eq:lemma_z_asymp_1}). \\[0.2cm]
\textbf{Proof of Lemma \ref{lemma_z_asymp_2}\,:} taking the logarithm of the product, the left-hand side of (\ref{eq:lemma_z_asymp_2}) reads
$$
\frac{1}{N^{\scriptscriptstyle 2}}\hspace{0.03cm}\log \prod^{N-1}_{n=1}\left(1 - e^{-n\hspace{0.02cm}\sigma^2} \right)^{\!N-n} \,=\, \frac{1}{N}\sum^{N-1}_{n=1}\left(1 - \frac{n}{N}\right)\log \left(1 - e^{-t\hspace{0.02cm}\frac{n}{N}}\right)
$$
which is a Riemann sum for the improper integral in the right-hand side. To evaluate this integral, one may resort to a symbolic computation software, or introduce the power series of the logarithm, under the integral,
$$
\int^{\scriptscriptstyle 1}_{\scriptscriptstyle 0}(1-x)\log\left(1 - e^{-tx}\right)\hspace{0.01cm}dx \,=\,
- \sum^\infty_{k=1}\frac{1}{k}\hspace{0.03cm}\int^{\scriptscriptstyle 1}_{\scriptscriptstyle 0}(1-x)\hspace{0.03cm}e^{-ktx}\hspace{0.01cm}dx
$$
and note that
$$
\int^{\scriptscriptstyle 1}_{\scriptscriptstyle 0}(1-x)\hspace{0.03cm}e^{-ktx}\hspace{0.01cm}dx = \frac{1 - e^{-ktx}}{(kt)^{\scriptscriptstyle 2}}
$$
in order to obtain $-(\mathrm{Li}_{\scriptscriptstyle 3}(e^{-t}) - \zeta(3))/t^{\scriptscriptstyle 2}$. \\[0.1cm]
\textbf{Remark\,:} from (\ref{eq:asymp_z}), it follows that $Z(\sigma) \rightarrow 0$ as $N \rightarrow \infty$, while $t = N\sigma^2$ remains constant. However, this is merely because $\omega_{\scriptscriptstyle 2}(N) \rightarrow 0$ as $N \rightarrow \infty$. Therefore, one should keep in mind,
\begin{equation} \label{eq:asymp_z_bis}
  \lim\hspace{0.04cm} \frac{1}{N^{\scriptscriptstyle 2}}\hspace{0.03cm}\log\left[\frac{Z(\sigma)}{\omega_{\scriptscriptstyle 2}(N)}\right] \,=\, -\frac{1}{2}\log(2) + \frac{3}{4} + \frac{t}{6} - \frac{\mathrm{Li}_{\scriptscriptstyle 3}(e^{-t}) - \zeta(3)}{t^{\scriptscriptstyle 2}}
\end{equation}
which may be thought of as the ``asymptotic cumulant generating function".

%$\log Z(\sigma) \sim - (N^{\scriptscriptstyle 2}/2)\hspace{0.02cm}\log(N)$, in the limit described in Proposition \ref{prop:asymp_z}, showing that $Z(\sigma)$ goes to zero. Essentially, this is due to $\omega_{\scriptscriptstyle 2}(N)$ going to zero.\hfill\linebreak All information relevant to the dependence of $Z(\sigma)$ on $\sigma$ is contained in the subsequent, constant terms of (\ref{eq:asymp_z}). It is these terms which should be thought of as the asymptotic cumulant generating function.

\section{The asymptotic distribution} \label{sec:rmtbis}
From the point of view of random matrix theory, a Gaussian distribution $P(\mathrm{I}_N\hspace{0.02cm},\sigma)$ on $M = \mathrm{H}(N)$ defines a unitary matrix ensemble. If $x$ is a random matrix, drawn from this ensemble, and $(x_i\,;i=1,\ldots, N)$ are its eigenvalues, which all belong to $(0,\infty)$, then the empirical distribution $\nu_{\scriptscriptstyle N\hspace{0.03cm}}$, which is given by (as usual, $\delta_{x_i}$ is the Dirac distribution at $x_i$)
\begin{equation} \label{eq:rr1}
\nu_{\scriptscriptstyle N}(B) = \mathbb{E}\left[\frac{1}{N}\sum^{\scriptscriptstyle N}_{i=1}\delta_{x_i}(B)\right]
\end{equation}
for measurable $B \subset (0,\infty)$, converges to an absolutely continuous distribution $\nu_{\scriptscriptstyle t\hspace{0.02cm}}$, when $N$ goes to infinity, while the product $t = N\sigma^2$ remains constant. 

%In random matrix theory, it is known that $\nu_{\scriptscriptstyle N}$ has a probability density function, equal to $\frac{1}{n}R^{\scriptscriptstyle \hspace{0.02cm}{(1)}}$, where $R^{\scriptscriptstyle \hspace{0.02cm}{(1)}}$ is the so-called one-point correlation function. 
\begin{proposition} \label{prop:rmtbis}
  Let $c = e^{-t}$ and $a(t) = c(1+\sqrt{1-c})^{\scriptscriptstyle -2}$ while $b(t) = c(1-\sqrt{1-c})^{\scriptscriptstyle -2}$. When $N$ goes to infinity, while the product $t = N\sigma^2$ remains constant, the empirical distribution $\nu_{\scriptscriptstyle N}$ converges weakly to the distribution $\nu_{\scriptscriptstyle t}$ with probability density function
\begin{equation} \label{eq:rmtbis}
 \frac{d\nu_{\scriptscriptstyle t}}{dx}(x) = \frac{1}{\pi\hspace{0.02cm}tx}\arctan\left(\frac{4\hspace{0.02cm}e^tx - (x+1)^2}{x+1}\right) \mathbf{1}_{[a(t),b(t)]}(x)
\end{equation}
where $\mathbf{1}_{[a(t),b(t)]}$ denotes the indicator function of the interval $[a(t),b(t)]$.
\end{proposition}
\noindent \textbf{Remark\,:} as one should expect, when $t = 0$ (so $\sigma^2 = 0$), $a(t) = b(t) = 1$.

The proof of Proposition \ref{prop:rmtbis} is a relatively direct application of a result in~\cite{asch} (Page 191). Recall the variables $u_i = e^{t}x_i$ which appear in (\ref{eq:covzvv}). Let $\tilde{\nu}_{\scriptscriptstyle N}$ be the empirical distribution of the $u_i$ (this is the same as (\ref{eq:rr1}), but with $u_i$ instead of $x_i$). By applying~\cite{mehta} (Chapter 5, Page 81),
\begin{equation} \label{eq:onepointcorr}
  \tilde{\nu}_{\scriptscriptstyle N}(B) = \frac{1}{N}\hspace{0.02cm}\int_B R^{\scriptscriptstyle \hspace{0.02cm}(1)}_{\scriptscriptstyle N}(u)\hspace{0.03cm}(du)
\end{equation}  
for measurable $B \subset (0,\infty)$, where the one-point correlation function $R^{\scriptscriptstyle \hspace{0.02cm}(1)}_{\scriptscriptstyle N}(u)$ is given by
\begin{equation} \label{eq:onepointcorrbis}
  R^{\scriptscriptstyle \hspace{0.02cm}(1)}_{\scriptscriptstyle N}(u) = \rho(u\hspace{0.02cm},2\sigma^{\scriptscriptstyle 2})\sum^{N-1}_{n=0}p^2_{\hspace{0.02cm} n}(u)
\end{equation}
in the notation of \ref{sec:orthpo} ($p_{\hspace{0.02cm}n}$ are orthonormal polynomials, with respect to the weight $\rho(u\hspace{0.02cm},2\sigma^{\scriptscriptstyle 2})$). According to~\cite{deift} (Page 133), $\tilde{\nu}_{\scriptscriptstyle N}$ given by (\ref{eq:onepointcorr}) converges weakly to the so-called equilibrium distribution $\tilde{\nu}_{\scriptscriptstyle t\hspace{0.02cm}}$, which minimises the electrostatic energy functional 
\begin{equation} \label{eq:electrostatic}
 E(\nu) = \frac{1}{t}\int^{\scriptscriptstyle \infty}_{\scriptscriptstyle 0}\frac{1}{2}\log^2(u)\nu(du) - \int^{\scriptscriptstyle \infty}_{\scriptscriptstyle 0}\int^{\scriptscriptstyle \infty}_{\scriptscriptstyle 0} \log|u-v|\nu(du)\nu(dv)
\end{equation}
over probability distributions $\nu$ on $(0,\infty)$. Also according to~\cite{deift} (Page 133), this equilibrium distribution is the asymptotic distribution of the zeros of the polynomial $p_{\hspace{0.02cm}\scriptscriptstyle N}$ (in the limit $N \rightarrow \infty$ while $N\sigma^2 = t$). Fortunately, $p_{\hspace{0.02cm}\scriptscriptstyle N}$ is just a constant multiple of the Stieltjes-Wigert polynomial $s_{\hspace{0.02cm}N}$~\cite{szegobook} (Page 33). Therefore, the required asymptotic distribution of zeros can be read from~\cite{asch} (Page 191). Finally, (\ref{eq:rmtbis}) follows by introducing the change of variables $x = e^{-t}u$. \\[0.1cm]
\textbf{Remark\,:} in~\cite{marino}, the equilibrium distribution $\tilde{\nu}_{\scriptscriptstyle t}$ is derived directly, by searching for stationary distributions of the energy functional (\ref{eq:electrostatic}). This leads to a singular integral equation, whose solution reduces to a Riemann-Hilbert problem. Astoundingly, the Gaussian distributions on $\mathrm{H}(N)$, as introduced in the present chapter, provide a matrix model for Chern-Simons quantum field theory (a detailed account is given in~\cite{marino}). 
%This asymptotic distribution $\nu_{\scriptscriptstyle t}$ has been computed both in mathematical literature (dealing with orthogonal polynomials~\cite{asch}) and in physical literature (dealing with Chern-Simons theory~\cite{marino}).

\vfill
\pagebreak

\section{Duality\,: the $\Theta$ distributions} \label{sec:THETA}
Recall the Riemannian symmetric space $M = \mathrm{H}(N)$ of \ref{sec:orthpo}. Its dual space is the unitary group $M^* = U(N)$.  Consider now a family of distributions on $M^*$, which will be called $\Theta$ distributions, and which display an interesting connection with Gaussian distributions on $M$, studied in \ref{sec:orthpo}. 
Recall Jacobi's $\vartheta$ function\footnote{To follow the original notation of Jacobi~\cite{watson}, this should be written $\vartheta(e^{i\phi}|q)$ where $q = e^{-\sigma^2}$. In other popular notations, this function is called $\vartheta_{\scriptscriptstyle 00}$ or $\vartheta_{\scriptscriptstyle 3\,}$.}, 
$$
\vartheta(e^{\scriptscriptstyle i\phi}|\sigma^{\scriptscriptstyle 2}) \,=\, \sum^{+\infty}_{m=-\infty} \exp(-m^2\hspace{0.02cm}\sigma^2 + 2m\hspace{0.03cm}i\phi)
$$
As a function of $\phi$, up to some minor modifications, this is just a wrapped normal distribution (in other words, the heat kernel of the unit circle),
$$
\frac{1}{2\pi}\hspace{0.03cm}
\vartheta\!\left(e^{\scriptscriptstyle i\phi}|{\scriptstyle \frac{\sigma^2}{2}}\right) \,=\,  \hspace{0.03cm}\sum^{\infty}_{m=-\infty} \exp\left[ - \frac{(2\phi - 2m\pi)^2}{2\sigma^2}\right]
$$
Each $x \in M^*$ can be written $x = k\cdot e^{i\theta}$ for some $k \in U(N)$ and $e^{i\theta} = \mathrm{diag}(e^{i\theta_i}\,;i=1,\ldots, N)$, where $k\cdot y = k\hspace{0.02cm}y\hspace{0.02cm}k^\dagger$, for $y \in M^*$. With this notation, define the following matrix $\vartheta$ function,
\begin{equation} \label{eq:THETAF}
  \Theta\left(x\middle|\sigma^2\right) \,=\, k\cdot \vartheta\!\left(e^{\scriptscriptstyle i\theta}|{\scriptstyle \frac{\sigma^2}{2}}\right)
\end{equation}
which is obtained from $x$ by applying Jacobi's $\vartheta$ function to each eigenvalue of $x$. Further, consider the positive function,
\begin{equation} \label{eq:THETAD}
 f_*(x|\bar{x}\hspace{0.02cm},\sigma) \,=\, \det\left[\left( 2\pi\hspace{0.03cm}\sigma^2\right)^{\!\frac{1}{2}}\hspace{0.03cm}\Theta\!\left(x\bar{x}^\dagger\middle|\sigma^2\right)\right]
\end{equation}
which is also equal to 
$$
\det\left[\left( 2\pi\hspace{0.03cm}\sigma^2\right)^{\!\frac{1}{2}}\hspace{0.03cm}\Theta\!\left(\bar{x}^\dagger x\middle|\sigma^2\right)\right]
$$ 
since the matrices $x\bar{x}^\dagger$ and $\bar{x}^\dagger x$ are similar. Then, let $Z_{\scriptscriptstyle M^*}(\sigma)$ denote the normalising constant
\begin{equation} \label{eq:zstar}
 Z_{\scriptscriptstyle M^*}(\sigma) = \int_{M^*}f_*(x|\bar{x}\hspace{0.02cm},\sigma)\,\mathrm{vol}(dx)
\end{equation}
which does not depend on $\bar{x}$, as can be seen, by introducing the new variable of integration $z = x\bar{x}^\dagger$, and using the invariance of $\mathrm{vol}(dx)$. (compare to the proof of Proposition \ref{prop:zhomogeneous}). 

Now, define a $\Theta$ distribution $\Theta(\bar{x},\sigma)$ as the probability distribution on $M^*$, whose probability density function, with respect to $\mathrm{vol}(dx)$, is given by 
\begin{equation} \label{eq:thetadensity}
  p_*(x|\bar{x}\hspace{0.02cm},\sigma) \,=\,\left(Z_{\scriptscriptstyle M^*}(\sigma)\right)^{-1}\hspace{0.03cm}f_*(x|\bar{x}\hspace{0.02cm},\sigma)
\end{equation}
\begin{proposition} \label{prop:thetadual}
  Let $Z_{\scriptscriptstyle M}(\sigma) = Z(\sigma)$, be given by (\ref{eq:sw_z}), and $Z_{\scriptscriptstyle M^*}(\sigma)$ be given by (\ref{eq:zstar}). Then, the following equality holds
\begin{equation} \label{eq:thetadual}
  \frac{Z_{\scriptscriptstyle M}(\sigma)}{Z_{\scriptscriptstyle M^*}(\sigma)} = \exp\left[{\small\left(\frac{N^3 - N}{6}\right)}\sigma^{2}\right]
\end{equation}
\end{proposition}
\noindent \textbf{Remark\,:} the Gaussian density (\ref{eq:gaussdensity}) on $M$, and the $\Theta$ distribution density (\ref{eq:thetadensity}) on $M^*$ are apparently unrelated. Therefore, it is interesting to note their normalising constants $Z_{\scriptscriptstyle M}(\sigma)$ and $Z_{\scriptscriptstyle M^*}(\sigma)$ scale together according to the simple relation (\ref{eq:thetadual}). The connection between the two distributions is due to the duality between the two spaces ($M$ and $M^*$). \vfill\pagebreak
%this provides yet another method of computing the normalising constant $Z(\sigma)$. Indeed, integrals with respect to the invariant Riemannian volume of a compact Lie group (such as $M^* = U(N)$) may be approximated using straightforward Monte Carlo methods~\cite{meckes} (see the ``Gauss-Gram-Schmidt" approach, Page 12). \\[0.1cm]
%\textbf{Remark\,:} of course, one is tempted to investigate how the relation (\ref{eq:thetadual}) extends to other pairs of dual Riemannian symmetric spaces, especially since the Monte Carlo method, just described, generalises quite easily to compact Riemannian symmetric spaces, beyond compact Lie groups. I hope to study this question, in as much depth as possible, after completing the present thesis.\vfill\pagebreak
\noindent \textbf{Proof of Proposition \ref{prop:thetadual}\,:} since $Z_{\scriptscriptstyle M^*}(\sigma)$ does not depend on $\bar{x}$, one may set $\bar{x} = o$ in (\ref{eq:zstar}), where $o = \mathrm{I}_N\hspace{0.03cm}$. Then,  $f_*(x|o\hspace{0.02cm},\sigma)$ is a class function, so (\ref{eq:zstar}) can be computed using (\ref{eq:integralunvmonde}). Note that $\omega(S_{\scriptscriptstyle N})$, which appears in (\ref{eq:integralunvmonde}), is equal to $\omega_{\scriptscriptstyle 2}(N)$, in the current notation. Therefore, 
\begin{equation} \label{eq:proofduality1}
 Z_{\scriptscriptstyle M^*}(\sigma) = 
\frac{\omega_{\scriptscriptstyle 2}(N)}{\mathstrut 2^{\scriptscriptstyle N^2} N!}\left( 2\pi\hspace{0.03cm}\sigma^2\right)^{\!\frac{N}{2}}\times I_{\scriptscriptstyle 2}
% \,\int_{[0\hspace{0.02cm},2\pi]^N}\,\prod^N_{i=1}\vartheta\left(e^{i\theta_i}\middle|\sigma^2\!/2\right)\hspace{0.03cm}|V(e^{i\theta})|^2\hspace{0.04cm}d\theta_{\scriptscriptstyle 1}\ldots\theta_{\scriptscriptstyle N}
\end{equation}
where $I_{\scriptscriptstyle 2}$ is the integral
\begin{equation} \label{eq:proofduality2}
I_{\scriptscriptstyle 2} \,=\,
\int_{[0\hspace{0.02cm},2\pi]^N}\,\prod^N_{i=1}\vartheta\!\left(e^{\scriptscriptstyle i\theta_i}|{\scriptstyle \frac{\sigma^2}{2}}\right)|V(e^{i\theta})|^2\hspace{0.04cm}d\theta_{\scriptscriptstyle 1}\ldots\theta_{\scriptscriptstyle N}
\end{equation}
which follows from the identity
$$
\det \Theta\!\left(x\middle|\sigma^2\right)
= \prod^N_{i=1}\vartheta\!\left(e^{\scriptscriptstyle i\theta_i}|{\scriptstyle \frac{\sigma^2}{2}}\right)
$$
Now, $I_{\scriptscriptstyle 2}$ can be expressed using~\cite{mehta} (Chapter 5, Page 79), as in the proof of Proposition \ref{prop:sw_z}. Precisely, if $(p_{\hspace{0.02cm}n}\,; n = 0,1,\ldots)$ are orthonormal trigonometric polynomials, with respect to the weight function $\vartheta\!\left(e^{\scriptscriptstyle i\theta}|{\scriptstyle \sigma^2\!/2}\right)$, on the unit circle, then $I_{\scriptscriptstyle 2}$ is given by (\ref{eq:mehtasw}), 
$$
I_{\scriptscriptstyle 2} \,=\, N!\hspace{0.03cm}\prod^{N-1}_{n=0} p^{-2}_{\hspace{0.02cm}nn}
$$
in terms of the leading coefficients $p_{\hspace{0.02cm}nn}$ of the polynomials $p_{\hspace{0.02cm}n}$ (these leading coefficients may always be chosen to be real). At present, the required orthonormal polynomials $p_{\hspace{0.02cm}n}$ are given by
\begin{equation} \label{eq:rogersz1}
  p_{\hspace{0.02cm}n}(z) \,=\, \left[q^{n}\!\prod^n_{m=1}( 1 - q^{m})^{-1}\right]^{\!\frac{1}{2}}r_{\hspace{0.02cm}n}(-q^{-\frac{1}{2}}z)
\end{equation}
where $q = e^{-\sigma^2}$ and $r_n(z)$ is the $n$-th Rogers-Szegö polynomial, which is monic~\cite{rogers}. Therefore, 
\begin{equation} \label{eq:rogerspnn}
   p^{-2}_{\hspace{0.02cm}nn} \,=\, \prod^{n}_{m=1}\left(1 - e^{-m\hspace{0.02cm}\sigma^2} \right)
\end{equation}
and, from (\ref{eq:proofduality2}), $I_{\scriptscriptstyle 2}$ is given by
\begin{equation} \label{eq:szegoi2}
I_{\scriptscriptstyle 2} \,=\, N! \prod^{N-1}_{n=1}\left(1 - e^{-n\hspace{0.02cm}\sigma^2} \right)^{\!N-n}
\end{equation}
which may be replaced into (\ref{eq:proofduality1}) to obtain
\begin{equation} \label{eq:zstarformula}
 Z_{\scriptscriptstyle M^*}(\sigma) = 
\frac{\omega_{\scriptscriptstyle 2}(N)}{\mathstrut 2^{\scriptscriptstyle N^2}}\left( 2\pi\hspace{0.03cm}\sigma^2\right)^{\!\frac{N}{2}}\hspace{0.02cm} \prod^{N-1}_{n=1}\left(1 - e^{-n\hspace{0.02cm}\sigma^2} \right)^{\!N-n}
\end{equation}
Finally, (\ref{eq:thetadual}) follows easily, by comparing (\ref{eq:zstarformula}) to (\ref{eq:sw_z}). \\[0.1cm]
\textbf{Remark\,:} the construction of the $\Theta$ distributions seems to indicate a general construction of ``dual distributions" on pairs of dual Riemannian symmetric spaces. Recalling the general notation of \ref{ssec:sspace}, it seems that Gaussian distributions arise from a classical Gaussian density profile on the maximal Abelian subspace $\mathfrak{a}$, while $\Theta$ distributions (``their duals") arise from wrapping this Gaussian density profile around the torus $\mathrm{Exp}_{o}(i\hspace{0.02cm}\mathfrak{a})$. \vfill\pagebreak

%\section{The FIM as a warped metric} \label{sec:warped}

\chapter{Bayesian inference and MCMC} \label{bayesian}

\minitoc
\vspace{0.1cm}
{\small
The present chapter is entirely made up of previously unpublished material. It continues the study of Gaussian distributions, from the previous chapter, in a new direction\,: Bayesian inference, and the Markov chain Monte Carlo (MCMC) techniques, useful in Bayesian inference. 
\begin{itemize}
 \item \ref{sec:mapvsmms} introduces two Bayesian estimators, the MAP and the MMS, for Gaussian distributions on a Riemannian symmetric space $M$. Proposition \ref{prop:mapmse} states these two estimators are equal, if the likelihood and prior densities are identical.
\item \ref{sec:mapmmsdistance} discusses a surprising experimental result\,: when $M$ is a space of constant negative curvature, numerical computation shows the MAP and the MMS are so close to each other that they appear to be equal, even if the likelihood and prior densities are different. 
\item \ref{sec:mc} states the original Proposition \ref{prop:gergodic}, which provides easy-to-verify sufficient conditions, for the geometric ergodicity of an isotropic Metropolis-Hastings Markov chain, in a Riemannian symmetric space which belongs to the non-compact case. This is then applied to the computation of the MMS, via the subsequent Proposition \ref{prop:bhatta1}.
\item \ref{sec:markovclt} discusses the Riemanian gradient descent method. Proposition \ref{prop:pgd} states this method has an exponential rate of convergence, when used to find the global mnimum of a strongly convex function, defined on a Hadamard manifold. Propositions \ref{prop:gergodic}, \ref{prop:bhatta1}, and \ref{prop:pgd} are the three essential ingredients of the recipe, used here for the computation of the MMS.
\item \ref{sec:mcmc_vollema} gives Lemma \ref{lemm:mcmc_volume}, to be used in the proof of Proposition \ref{prop:gergodic}. This lemma states that the logarithmic rate of growth of the volume density function is bounded at infinity, for a Riemannian symmetric space which belongs to the non-compact case.
\item \ref{sec:mcmcproof} is devoted to the proof of Proposition \ref{prop:gergodic}. The proof is a generalisation of the proof in~\cite{jarner}, carried out in the special case of Metropolis algorithms in a Euclidean space.
\end{itemize}
}

\vfill
\pagebreak

\section{MAP versus MMS} \label{sec:mapvsmms}
Let $M$ be a Riemannian symmetric space, which belongs to the non-compact case (see \ref{ssec:sspace}). Recall the Gaussian distribution $P(x\hspace{0.03cm},\sigma)$ on $M$ is given by its probability density function (\ref{eq:gaussdensity})
\begin{equation} \label{eq:likelihood}
p(y|x\hspace{0.02cm},\sigma)\,=\,  \left(Z(\sigma)\right)^{-1}\hspace{0.03cm}\exp\left[ -\frac{d^{\hspace{0.03cm}2}(y,x)}{2\sigma^2}\right]
\end{equation}
In \ref{sec:gaussmle}, it was seen that maximum-likelihood estimation of the parameter $x$, based on independent samples $(y_n\,;n=1,\ldots,N)$, amounts to computing the Riemannian barycentre of these samples. The one-sample maximum-likelihood estimate, given a single observation $y$, is therefore $\hat{x}_{\scriptscriptstyle ML} = y$. 

Instead of maximum-likelihood estimation, consider a Bayesian approach to estimating $x$, based on the observation $y$. To do so, assign to $x$ a prior density, which is also Gaussian,
\begin{equation} \label{eq:prior}
p(x|z\hspace{0.02cm},\tau)\,=\,\left(Z(\tau)\right)^{-1}\hspace{0.03cm}\exp\left[ -\frac{d^{\hspace{0.03cm}2}(x,z)}{2\tau^2}\right]
\end{equation}
Upon observation of $y$, Bayesian inference concerning $x$ is carried out, using the posterior density
\begin{equation} \label{eq:posterior}
\pi(x) \propto \exp\left[ -\frac{d^{\hspace{0.03cm}2}(y,x)}{2\sigma^2}-\frac{d^{\hspace{0.03cm}2}(x,z)}{2\tau^2}\right]
\end{equation}
where $\propto$ indicates a missing (unknown) normalising factor. 

In particular, the maximum \textit{a posteriori} estimator $\hat{x}_{\scriptscriptstyle MAP}$ of $x$ is equal to the mode of the posterior density $\pi(x)$. In other words, $\hat{x}_{\scriptscriptstyle MAP}$ minimises the weighted sum of squared distances $d^{\hspace{0.03cm}2}(y,x)/\sigma^2+d^{\hspace{0.03cm}2}(x,z)/\tau^2$. This is expressed in the following notation\footnote{If $p\hspace{0.03cm},q \in M$ and $c:[0,1]\rightarrow M$ is a geodesic curve with $c(0) = p$ and $c(1) = q$, then $p\,\#_{\scriptscriptstyle t}\, q = c(t)$, for $t \in [0,1]$. Therefore, $p\,\#_{\scriptscriptstyle t}\, q$ is a geodesic convex combination of $p$ and $q$, with respective weights $(1-t)$ and $t$.},
\begin{equation} \label{eq:mapformula}
  \hat{x}_{\scriptscriptstyle MAP} \,=\, z\, \#_{\scriptscriptstyle \rho}\, y \hspace{0.5cm} \text{where } \rho = \frac{\tau^2}{\sigma^2+\tau^2}
\end{equation}
Thus, $\hat{x}_{\scriptscriptstyle MAP}$ is a geodesic convex combination of the prior barycentre $z$ and the observation $y$, with respective weights $\sigma^2/(\sigma^2+\tau^2)$ and $\tau^2/(\sigma^2+\tau^2)$.

On the other hand, the minimum mean square error estimator $\hat{x}_{\scriptscriptstyle MMS}$ is the barycentre of the posterior density $\pi(x)$. That is, $\hat{x}_{\scriptscriptstyle MMS}$ is the global minimiser of 
\begin{equation} \label{eq:posteriorvariance}
  \mathcal{E}_{\pi}(y) \,=\,
\frac{1}{2}\hspace{0.03cm}
\int_M\,d^{\hspace{0.03cm}\scriptscriptstyle 2}(y\hspace{0.02cm},x)\hspace{0.03cm}\pi(x)\hspace{0.03cm}\mathrm{vol}(dx)
\end{equation}
whose existence and uniqueness are established in the remark below. While it is easy to compute $\hat{x}_{\scriptscriptstyle MAP}$ from (\ref{eq:mapformula}), it is much harder to find $\hat{x}_{\scriptscriptstyle MMS\,}$, as this requires minimising the integral (\ref{eq:posteriorvariance}), where the density $\pi(x)$ is known only up to normalisation. 

Still, there is one special case where these two estimators are equal.
\begin{proposition} \label{prop:mapmse}
  In the above notation, if $\sigma^2 = \tau^2$ (that is $\rho = 1/2$), then $\hat{x}_{\scriptscriptstyle MMS} = \hat{x}_{\scriptscriptstyle MAP\,}$.
\end{proposition}
This relies on the following (intuitively quite obvious) lemma. 
\begin{lemma}   \label{lemma:bcentre_isometry}
  Assume that $\pi$ is a probability distribution on $M$ with Riemannian barycentre $b$. If $g$ is an isometry of $M$ such that $g^*\pi = \pi$ ($g^*\pi$ denotes the image of the distribution $\pi$ under the mapping $g:M\rightarrow M$), then $g\cdot b = b$.
\end{lemma}
This lemma is proved by noting that, for any isometry $g$ of $M$, one has $\mathcal{E}_{g^*\pi} = \mathcal{E}_{\pi}\circ g^{\scriptscriptstyle -1}$. Accordingly, if $b$ is the Riemannian barycentre of $\pi$, $g\cdot b$ is the Riemannian barycentre of $g^*\pi$. \vfill\pagebreak

\noindent \textbf{Proof of Proposition \ref{prop:mapmse}\,:} in this case, 
$$
\pi(x) \propto \exp\left[ -\frac{d^{\hspace{0.03cm}2}(y,x)+d^{\hspace{0.03cm}2}(x,z)}{2\sigma^2}\right]
$$
On the other hand, $\hat{x}_{\scriptscriptstyle MAP} \,=\, z\, \#_{\scriptscriptstyle 1/2}\, y$ is the midpoint of the geodesic segment connecting $z$ to $y$ (note that $\rho = 1/2$). Let $s$ denote the geodesic symmetry at $\hat{x}_{\scriptscriptstyle MAP\,}$. Then, $s$ permutes $z$ and $y$, and therefore leaves invariant $\pi(x)$. Lemma \ref{lemma:bcentre_isometry} (applied with $g = s$) implies the Riemannian barycentre $\hat{x}_{\scriptscriptstyle MMS}$ of $\pi$ verifies $s\cdot \hat{x}_{\scriptscriptstyle MMS} = \hat{x}_{\scriptscriptstyle MMS\,}$. However, $\hat{x}_{\scriptscriptstyle MAP}$ is the unique fixed point of $s$. Therefore, $\hat{x}_{\scriptscriptstyle MMS} = \hat{x}_{\scriptscriptstyle MAP\,}$.\\[0.1cm]
\textbf{Remark\,:} to see that $\hat{x}_{\scriptscriptstyle MMS}$ is well-defined, it is enough to show the posterior density $\pi$ in (\ref{eq:posterior}) satisfies (\ref{eq:secondordermoment}). Indeed, this implies that $\pi$ has a well-defined Riemannian barycentre. 

Consider then the second-order moment in (\ref{eq:secondordermoment}), with $y_o = \hat{x}_{\scriptscriptstyle MAP\,}$. Specifically, this is 
\begin{equation} \label{eq:mmap}
m_{\scriptscriptstyle 2}(\hat{x}_{\scriptscriptstyle MAP}) = \int_M d^{\hspace{0.03cm}\scriptscriptstyle 2}(\hat{x}_{\scriptscriptstyle MAP}\hspace{0.02cm},x)\hspace{0.03cm}\pi(x)\hspace{0.03cm}\mathrm{vol}(dx)
\end{equation}
Rearrange (\ref{eq:posterior}) to obtain
\begin{equation} \label{eq:rearrangepi}
\pi(x) \propto \exp\left[ -h\left(\rho f_y(x)+(1-\rho)f_z(x)\right)\right] \hspace{1cm} \left(h = 1/\sigma^2 + 1/\tau^2\right)
\end{equation}
in the notation of \ref{sec:squareddistance}. Now, let $f(x) = \rho f_y(x)+(1-\rho)f_z(x)$. For $x \in M$, let $x = \mathrm{Exp}_{\hat{x}_{\scriptscriptstyle MAP}}(v)$, and recall the Taylor expansion (\ref{eq:taylor2}),

\begin{equation} \label{eq:rearrangetaylor}
  f\left(x\right) = f(\hat{x}_{\scriptscriptstyle MAP}) + \langle \mathrm{grad}\,f,v\rangle_{\scriptscriptstyle \hat{x}_{\scriptscriptstyle MAP}} + \frac{1}{2}\,\mathrm{Hess}\,f_{\scriptscriptstyle c(t^*)}(\dot{c},\dot{c})
\end{equation}
where $c(t^*)$ is a point along the geodesic $c(t) = \mathrm{Exp}_x(t\,v)$, corresponding to an instant $t^* \in (0,1)$. Note that $\mathrm{grad}\,f(\hat{x}_{\scriptscriptstyle MAP}) = 0$, as can be checked from (\ref{eq:mapformula}), and that, using (\ref{eq:hesfxcomp}),
$$
\mathrm{Hess}\,f(x) = \rho\hspace{0.03cm}\mathrm{Hess}\,f_y(x) + (1-\rho)\hspace{0.03cm}\mathrm{Hess}\,f_z(x)\,\geq\,g(y) 
$$ 
Replacing these into (\ref{eq:rearrangetaylor}), it follows that
$$
  f(x) \,\geq\, \rho\hspace{0.02cm}(1-\rho)\hspace{0.03cm}d^{\hspace{0.03cm}2}(z,y) + \frac{1}{2}\hspace{0.02cm}d^{\hspace{0.03cm}2}(\hat{x}_{\scriptscriptstyle MAP}\hspace{0.02cm},x)
$$
Then, if $C^{-1}_\pi$ is the missing normalising factor in (\ref{eq:rearrangepi}),
\begin{equation} \label{eq:rearrangepi1}
\pi(x) \leq C^{-1}_\pi\hspace{0.03cm}\exp\left[ -\frac{\rho}{\tau^{\scriptscriptstyle 2}}\hspace{0.03cm}d^{\hspace{0.03cm}2}(z,y)-\frac{h}{2}\hspace{0.03cm}d^{\hspace{0.03cm}2}(\hat{x}_{\scriptscriptstyle MAP}\hspace{0.02cm},x)\right]
\end{equation}
From (\ref{eq:mmap}) and (\ref{eq:rearrangepi1}),
\begin{equation} \label{eq:rearrangepi2}
m_{\scriptscriptstyle 2}(\hat{x}_{\scriptscriptstyle MAP}) \,\leq\, C^{-1}_\pi\exp\left[ -\frac{\rho}{\tau^{\scriptscriptstyle 2}}\hspace{0.03cm}d^{\hspace{0.03cm}2}(z,y)\right]\hspace{0.02cm}\int_Md^{\hspace{0.03cm}\scriptscriptstyle 2}(\hat{x}_{\scriptscriptstyle MAP}\hspace{0.02cm},x)\hspace{0.03cm}
\exp\left[-\frac{h}{2}\hspace{0.03cm}d^{\hspace{0.03cm}2}(\hat{x}_{\scriptscriptstyle MAP}\hspace{0.02cm},x)\right]
\hspace{0.03cm}\mathrm{vol}(dx)
\end{equation}
which is finite, as required in (\ref{eq:secondordermoment}). In fact, by a direct application of the integral formula (\ref{eq:ssvolncka2}), it is possible to show that
$$
\int_Md^{\hspace{0.03cm}\scriptscriptstyle 2}(\hat{x}_{\scriptscriptstyle MAP}\hspace{0.02cm},x)\hspace{0.03cm}
\exp\left[-\frac{h}{2}\hspace{0.03cm}d^{\hspace{0.03cm}2}(\hat{x}_{\scriptscriptstyle MAP}\hspace{0.02cm},x)\right]
\hspace{0.03cm}\mathrm{vol}(dx) =
h^{\scriptscriptstyle -3/2}\hspace{0.02cm}Z^\prime(h^{\scriptscriptstyle -1/2})
$$
where $Z(\sigma)$ was given in (\ref{eq:ssz}), and the prime denotes the derivative. Finally, replacing this into (\ref{eq:rearrangepi2}), it follows that
\begin{equation} \label{eq:prewasserstein}
m_{\scriptscriptstyle 2}(\hat{x}_{\scriptscriptstyle MAP}) \,\leq\, C^{-1}_\pi\exp\left[ (-\rho/\tau^{\scriptscriptstyle 2})\hspace{0.03cm}d^{\hspace{0.03cm}2}(z,y)\right]\hspace{0.02cm}h^{\scriptscriptstyle -3/2}\hspace{0.02cm}Z^\prime(h^{\scriptscriptstyle -1/2})
\end{equation}
%which will be used in the following paragraph. 
\section{Bounding the distance} \label{sec:mapmmsdistance}
Proposition \ref{prop:mapmse} states that $\hat{x}_{\scriptscriptstyle MMS} = \hat{x}_{\scriptscriptstyle MAP\,}$, if $\rho = 1/2$. When $M$ is a Euclidean space, it is famously known that $\hat{x}_{\scriptscriptstyle MMS} = \hat{x}_{\scriptscriptstyle MAP}$ for any value of $\rho$. In general, one expects these two estimators to be different from one another, if $\rho \neq 1/2$. 

However, when $M$ is a space of constant negative curvature, numerical experiments show that $\hat{x}_{\scriptscriptstyle MMS}$ and $\hat{x}_{\scriptscriptstyle MAP}$ lie surprisingly close to each other, and that they even appear to be equal. I am still unaware of any mathematical explanation of this phenomenon. 

It is possible to bound the distance between $\hat{x}_{\scriptscriptstyle MMS}$ and $\hat{x}_{\scriptscriptstyle MAP\,}$, using the so-called fundamental contraction property~\cite{sturm} (this is an immediate application of Jensen's inequality, as explained in the proof of Theorem 6.3 in~\cite{sturm}).
\begin{equation} \label{eq:contraction}
  d(\hat{x}_{\scriptscriptstyle MMS}\hspace{0.02cm},\hat{x}_{\scriptscriptstyle MAP}) \leq W(\pi,\delta_{\hat{x}_{\scriptscriptstyle MAP}})
\end{equation}
where $W$ denotes the Kantorovich ($L^{\scriptscriptstyle 1}$-Wasserstein) distance, and $\delta_{\hat{x}_{\scriptscriptstyle MAP}}$ denotes the Dirac probability distribution concentrated at $\hat{x}_{\scriptscriptstyle MAP\,}$. Now, the right-hand side of (\ref{eq:contraction}) is equal to the first-order moment 
\begin{equation} \label{eq:mmap1}
m_{\scriptscriptstyle 1}(\hat{x}_{\scriptscriptstyle MAP}) = \int_M d(\hat{x}_{\scriptscriptstyle MAP}\hspace{0.02cm},x)\hspace{0.03cm}\pi(x)\hspace{0.03cm}\mathrm{vol}(dx)
\end{equation}
Of course, the upper bound in (\ref{eq:contraction}) is not tight, since it is strictly positive, even when $\rho = 1/2$, as one may see from (\ref{eq:mmap1}). 

It will be shown below that a Metropolis-Hastings algorithm, with Gaussian proposals, can be used to generate (geometrically ergodic) samples $(x_n\,;n\geq 1)$ from the posterior density $\pi$.\hfill\linebreak Using these samples, it is possible to approximate (\ref{eq:mmap1}), by an empirical average, 
\begin{equation} \label{eq:mmapbis}
\bar{m}_{\scriptscriptstyle 1}(\hat{x}_{\scriptscriptstyle MAP}) = \frac{1}{N}\hspace{0.03cm}\sum^N_{n=1}d(\hat{x}_{\scriptscriptstyle MAP}\hspace{0.02cm},x_n)
\end{equation}
In addition, the samples $(x_n)$ can be used to compute a convergent approximation of $\hat{x}_{\scriptscriptstyle MMS\,}$. Precisely, the empirical barycentre $\bar{x}_{\scriptscriptstyle MMS}$ of the samples $(x_{\scriptscriptstyle 1},\ldots,x_{\scriptscriptstyle N})$ converges almost-surely to $\hat{x}_{\scriptscriptstyle MMS}$ (this is proved in \ref{ssec:recursiveb}).

Numerical experiments were conducted in the case when $M$ is a space of constant curvature, equal to $-1$, and of dimension $n$.  The following table was obtained for the values $\sigma^2 = \tau^2 = 0.1$, using samples $(x_{\scriptscriptstyle 1},\ldots,x_{\scriptscriptstyle N})$ where $N = 2\times10^5$.
$$
\begin{array}{rlllllllll}
\text{dimension }n & 2 & 3 & 4 & 5 & 6 & 7 & 8 & 9 & 10 \\[0.2cm]
\bar{m}_{\scriptscriptstyle 1}(\hat{x}_{\scriptscriptstyle MAP})     &         0.28    &  0.35  &  0.41   &  0.47       &   0.50       & 0.57          &  0.60        &        0.66& 0.70 \\[0.2cm]
  d(\bar{x}_{\scriptscriptstyle MMS}\hspace{0.02cm},\hat{x}_{\scriptscriptstyle MAP})     & 0.00  & 0.00   & 0.00   & 0.01  &  0.01  & 0.02  & 0.02   & 0.02 & 0.03
\end{array}
$$
and the following table for $\sigma^2 = 1$ and $\tau^2 = 0.5$, again using $N = 2\times10^5$.
$$
\begin{array}{rlllllllll}
\text{dimension }n & 2 & 3 & 4 & 5 & 6 & 7 & 8 & 9 & 10 \\[0.2cm]
\bar{m}_{\scriptscriptstyle 1}(\hat{x}_{\scriptscriptstyle MAP})    &         0.75    &  1.00  &  1.12   &  1.44       &   1.73       &  1.97         & 2.15     &   2.54     &2.91 \\[0.2cm]
  d(\bar{x}_{\scriptscriptstyle MMS}\hspace{0.02cm},\hat{x}_{\scriptscriptstyle MAP})      & 0.00   &  0.00  & 0.03   & 0.02  & 0.02  &  0.03  &  0.04  &  0.03 &  0.12
\end{array}
$$
The first table, as expected, confirms Proposition \ref{prop:mapmse}. The second table, more surprisingly, shows that $\hat{x}_{\scriptscriptstyle MMS}$ and $\hat{x}_{\scriptscriptstyle MAP}$ can be quite close to each other, even when $\rho \neq 1/2$. 

Other values of $\sigma^2$ and $\tau^2$ lead to similar orders of magnitude for $\bar{m}_{\scriptscriptstyle 1}(\hat{x}_{\scriptscriptstyle MAP})$ and $  d(\bar{x}_{\scriptscriptstyle MMS}\hspace{0.02cm},\hat{x}_{\scriptscriptstyle MAP})$. While $\bar{m}_{\scriptscriptstyle 1}(\hat{x}_{\scriptscriptstyle MAP})$  increases with the dimension $n$, $d(\bar{x}_{\scriptscriptstyle MMS}\hspace{0.02cm},\hat{x}_{\scriptscriptstyle MAP})$ does not appear sensitive to increasing dimension. 

Based on these experimental results, one may be tempted to conjecture that $\hat{x}_{\scriptscriptstyle MMS} = \hat{x}_{\scriptscriptstyle MAP\,}$, even when $\rho \neq 1/2$. Naturally, numerical experiments do not equate to a mathematical proof.

\section{Computing the MMS} \label{sec:mc}

%samples from the posterior density $\pi$, generated using the Metropolis-Hastings
\subsection{Metropolis-Hastings algorithm} \label{ssec:mh}
A crucial step, in Bayesian inference, is sampling from the posterior density. Here, this is $\pi(x)$, given by (\ref{eq:posterior}). 
Since $\pi(x)$ is known only up to normalisation, a suitable sampling method is afforded by the Metropolis-Hastings algorithm. This algorithm generates a Markov chain $(x_n\,;n\geq 1)$, with transition kernel~\cite{roberts}
\begin{equation} \label{eq:hmP}
  Pf(x) \,=\, \int_M\,\alpha(x\hspace{0.02cm},y)\hspace{0.02cm}q(x\hspace{0.02cm},y)f(y)\hspace{0.02cm}\mathrm{vol}(dy) + \rho(x)\hspace{0.02cm}f(x)
\end{equation}
for any bounded measurable function $f:M\rightarrow \mathbb{R}$, where $\alpha(x\hspace{0.02cm},y)$ is the probability of accepting a transition from $x$ to $dy$, and $\rho(x)$ is the probability of staying at $x$, and where $q(x,y)$ is the proposed transition density
\begin{equation} \label{eq:hmTr}
  q(x\hspace{0.02cm},y) \geq 0 \hspace{0.2cm}\text{and}\hspace{0.2cm} \int_M\,q(x\hspace{0.02cm},y)\hspace{0.02cm}\mathrm{vol}(dy) = 1 \hspace{0.5cm} \text{for } x\in M
\end{equation}
In the following, $(x_n)$ will always be an isotropic Metropolis-Hastings chain, in the sense that $q(x\hspace{0.02cm},y) = q(d(x\hspace{0.02cm},y))$, so $q(x\hspace{0.02cm},y)$ only depends on the distance $d(x\hspace{0.02cm},y)$. In this case, the acceptance probability $\alpha(x\hspace{0.02cm},y)$ is given by $\alpha(x\hspace{0.02cm},y) = \mathrm{min}\left\lbrace1\hspace{0.02cm},\pi(y)/\pi(x) \right\rbrace$.

The aim of the Metropolis-Hastings algorithm is to produce a Markov chain $(x_n)$ which is geometrically ergodic. Geometric ergodicity means the distribution $\pi_n$ of $x_n$ converges to $\pi$, with a geometric rate, in the sense that there exist $\beta \in (0,1)$ and $R(x_{\scriptscriptstyle 1}) \in (0,\infty)$, as well as a function $V:M\rightarrow \mathbb{R}$, such that (in the following, $\pi(dx) = \pi(x)\hspace{0.02cm}\mathrm{vol}(dx)$)
\begin{equation} \label{eq:VV}
  V(x) \geq \max\left\lbrace 1\hspace{0.02cm},d^{\hspace{0.03cm} 2}(x,x^*)\right\rbrace \text{ for some } x^* \in M
\end{equation}
\begin{equation} \label{eq:gergodic}
\left| \int_M\, f(x)\hspace{0.02cm}(\pi_n(dx) - \pi(dx)) \right| \,\leq\,R(x_{\scriptscriptstyle 1})\hspace{0.02cm}\beta^n \hspace{0.9cm}
\end{equation}
for any function $f:M\rightarrow \mathbb{R}$ with $|f|\leq V$. If the chain $(x_n)$ is geometrically ergodic, then it satisfies the strong law of large numbers~\cite{tweedie} 
\begin{equation} \label{eq:mcmclln}
\frac{1}{N}\hspace{0.03cm}\sum^N_{n=1}f(x_n) \longrightarrow \int_M\,f(x)\hspace{0.02cm}\pi(dx) \text{ \;\;(almost-surely)}
\end{equation}
as well as a corresponding central limit theorem (see Theorem 17.0.1, in~\cite{tweedie}). Then, in practice, the Metropolis-Hastings algorithm generates samples $(x_n)$ from the posterior density $\pi(x)$. 

In \ref{sec:mcmcproof}, the following general statement will be proved, concerning the geometric ergodicity of isotropic Metropolis-Hastings chains. 
\begin{proposition} \label{prop:gergodic}
Let $M$ be a Riemannian symmetric space, which belongs to the non-compact case. Assume $(x_n\,;n\geq 1)$ is a Markov chain in $M$, with transition kernel given by (\ref{eq:hmP}), with proposed transition density $q(x\hspace{0.02cm},y) = q(d(x\hspace{0.02cm},y))$, and with strictly positive invariant density $\pi$. 

The chain $(x_n)$ satisfies (\ref{eq:VV}) and (\ref{eq:gergodic}), if the following assumptions hold,\\[0.1cm]
(a1) there exists $x^* \in M$, such that $r(x) = d(x^*,x)$ and $\ell(x) = \log\hspace{0.02cm}\pi(x)$ satisfy
$$
\limsup_{r(x)\rightarrow \infty}\hspace{0.02cm}\frac{\langle\mathrm{grad}\,r,\mathrm{grad}\,\ell\rangle_{\scriptscriptstyle x}}{r(x)}\,<\,0
$$%which V
(a2) if $n(x) = \left.\mathrm{grad}\,\ell(x)\middle/\Vert \mathrm{grad}\,\ell(x)\Vert\right.$, then $n(x)$ satisfies
$$
\limsup_{r(x)\rightarrow \infty}\hspace{0.02cm}\langle\mathrm{grad}\,r,n\rangle_{\scriptscriptstyle x}\,<\,0
$$
(a3) there exist $\delta_{\scriptscriptstyle q} > 0$ and $\varepsilon_{\scriptscriptstyle q} > 0$ such that $d(x\hspace{0.02cm},y) <\delta_{\scriptscriptstyle q}$ implies $q(x\hspace{0.02cm},y) > \varepsilon_{\scriptscriptstyle q\,}$
\end{proposition}
\noindent \textbf{Remark\,:} the posterior density $\pi$ in (\ref{eq:posterior}) verifies Assumptions (a1) and (a2). To see this, let $x^* = z$, and note from (\ref{eq:gradr}) and (\ref{eq:gradfx}) that
$$
\mathrm{grad}\,\ell(x) \,=\, -\frac{1}{\tau^2}\hspace{0.02cm}r(x)\hspace{0.02cm}\mathrm{grad}\,r(x) - \frac{1}{\sigma^2}\hspace{0.02cm}\mathrm{grad}\,f_y(x)
$$
Then, taking the scalar product with $\mathrm{grad}\,r$,
\begin{equation} \label{eq:proofconditiona1}
\langle\mathrm{grad}\,r,\mathrm{grad}\,\ell\rangle_{x} = -\frac{1}{\tau^2}\hspace{0.02cm}r(x) - 
\frac{1}{\sigma^2}\hspace{0.02cm}\langle\mathrm{grad}\,r,\mathrm{grad}\,f_y\rangle_{x}
\end{equation}
since $\mathrm{grad}\,r(x)$ is a unit vector, for all $x \in M$. Now, $\mathrm{grad}\,f_y(x) = -\mathrm{Exp}^{-1}_x(y)$, by (\ref{eq:gradfx}). But, 
since $r(x)$ is a convex function of $x$,
$$
\langle \mathrm{grad}\,r,\mathrm{Exp}^{-1}_x(y)\rangle \,\leq\, r(y) - r(x)
$$
for any $y \in M$. Thus, the right-hand side of (\ref{eq:proofconditiona1}) is strictly negative, as soon as $r(x) > r(y)$, and Assumption (a1) is indeed verified. That Assumption (a2) is also verified can be proved by a similar reasoning. \\[0.1cm]
\textbf{Remark\,:} on the other hand, Assumption (a3) holds, if the proposed transition density $q(x\hspace{0.02cm},y)$ is a Gaussian density, $q(x\hspace{0.02cm},y) = p(y|x,\tau^{\scriptscriptstyle 2}_{\scriptscriptstyle q})$. \\[0.1cm]
\indent With this choice of $q(x\hspace{0.02cm},y)$, all the assumptions of Proposition \ref{prop:gergodic} are verified, for the posterior density $\pi$ in (\ref{eq:posterior}).  Proposition \ref{prop:gergodic} therefore implies that the Metropolis-Hastings algorithm generates geometrically ergodic samples $(x_n\,;n\geq 1)$, from this posterior density.

\subsection{The empirical barycentre} \label{ssec:recursiveb}
Let $(x_n\,;n\geq 1)$ be a Metropolis-Hastings Markov chain in $M$, with its transition kernel (\ref{eq:hmP}), and invariant density $\pi$. Assume the chain $(x_n)$ is geometrically ergodic, so it satisfies the strong law of large numbers (\ref{eq:mcmclln}). 

Then, let $\bar{x}_{\scriptscriptstyle N}$ denote the empirical barycentre of the first $N$ samples $(x_{\scriptscriptstyle 1},\ldots,x_{\scriptscriptstyle N})$. This is the unique global minimum of the variance function
\begin{equation} \label{eq:emprecursive}
\mathcal{E}_{\scriptscriptstyle N}(y) \,=\,\frac{1}{2N}\sum^N_{n=1}d^{\hspace{0.03cm}2}(y\hspace{0.02cm},x_n)
\end{equation}
Assuming it is well-defined, let $\hat{x}$ denote the Riemannian barycentre of the invariant density $\pi$. It turns out that $\bar{x}_{\scriptscriptstyle N}$ converges almost-surely to $\hat{x}$. 
\begin{proposition} \label{prop:bhatta1}
 Let $(x_n)$ be any Markov chain in a Hadamard manifold $M$, with invariant distribution $\pi$. Denote
$\bar{x}_{\scriptscriptstyle N}$ the empirical barycentre of $(x_{\scriptscriptstyle 1},\ldots,x_{\scriptscriptstyle N})$, and $\hat{x}$ the Riemannian barycentre of $\pi$ (assuming it is well-defined). If $(x_n)$ satisfies the strong law of large numbers (\ref{eq:mcmclln}), then $\bar{x}_{\scriptscriptstyle N}$ converges to $\hat{x}$, almost-surely. 
\end{proposition}
According to the remarks after Proposition \ref{prop:gergodic}, the Metropolis-Hastings Markov chain $(x_n)$, whose invariant density is the posterior density $\pi(x)$, given by (\ref{eq:posterior}), is geometrically ergodic. Therefore, by Proposition \ref{prop:bhatta1}, the empirical barycentre $\bar{x}_{\scriptscriptstyle MMS\,}$, of the samples $(x_{\scriptscriptstyle 1},\ldots,x_{\scriptscriptstyle N})$, converges almost-surely to the minimum mean square error estimator $\hat{x}_{\scriptscriptstyle MMS}$ (since this is just the barycentre of the posterior density $\pi$). This provides a practical means of approximating $\hat{x}_{\scriptscriptstyle MMS\,}$. Indeed, $\bar{x}_{\scriptscriptstyle MMS}$ can be computed using the Riemannian gradient descent method (this method is discussed in \ref{sec:markovclt}, below).

The proof of Proposition \ref{prop:bhatta1} is nearly a word-for-word repetition of the proof in~\cite{bhatta1} (that of Theorem 2.3). \\[0.1cm]
\textbf{Proof of Proposition \ref{prop:bhatta1}\,:} denote $\mathcal{E}_{\pi}$ the variance function of the invariant distribution $\pi$,
$$
\mathcal{E}_{\pi}(y) = \frac{1}{2}\hspace{0.03cm}
\int_M\,d^{\hspace{0.03cm}\scriptscriptstyle 2}(y\hspace{0.02cm},x)\hspace{0.03cm}\pi(dx)
$$
First, for any compact $K \subset M$, it will be proved that
\begin{equation} \label{eq:proofbhatta1}
  \sup_{y\in K}\hspace{0.03cm}\left| \mathcal{E}_{\scriptscriptstyle N}(y) - \mathcal{E}_{\pi}(y)\right| \longrightarrow 0\text{ \;\;(almost-surely)}
\end{equation}
To do so, let $\delta > 0$ and let $\lbrace w_j\,;j=1,\ldots,J\rbrace$ be a $\delta$-net in $K$ (for any $y \in  K$, there exists $w_j$ such that $d(w_j\hspace{0.02cm},y) < \delta$). By the strong law of large numbers (\ref{eq:mcmclln}), 
\begin{equation} \label{eq:proofbhatta11}
  \max_{j=1,\ldots, J}\hspace{0.03cm}\left| \mathcal{E}_{\scriptscriptstyle N}(w_j) - \mathcal{E}_{\pi}(w_j)\right| \longrightarrow 0\text{ \;\;(almost-surely)}
\end{equation}
Using the elementary identity
$$
\left| d^{\hspace{0.03cm}\scriptscriptstyle 2}(y\hspace{0.02cm},x_n) - d^{\hspace{0.03cm}\scriptscriptstyle 2}(w\hspace{0.02cm},x_n)\right| \leq \left(d(y\hspace{0.02cm},x_n) + d(w\hspace{0.02cm},x_n)\right)\left|d(y\hspace{0.02cm},x_n) - d(w\hspace{0.02cm},x_n)\right|
$$
it follows by the triangle inequality that
\begin{equation} \label{eq:squarelipschitz}
\left| d^{\hspace{0.03cm}\scriptscriptstyle 2}(y\hspace{0.02cm},x_n) - d^{\hspace{0.03cm}\scriptscriptstyle 2}(w\hspace{0.02cm},x_n)\right| \leq \left(d(y\hspace{0.02cm},x_n) + d(w\hspace{0.02cm},x_n)\right)d(w\hspace{0.02cm},y)
\end{equation}
From (\ref{eq:squarelipschitz}), it is possible to show that, for $y$ and $w$ in $K$,
\begin{equation} \label{eq:proofbhatta12}
\left| \mathcal{E}_{\scriptscriptstyle N}(y) - \mathcal{E}_{\scriptscriptstyle N}(w)\right| \leq \sup_{z\in K}\left( \frac{1}{N}\sum^N_{n=1}d(z\hspace{0.02cm},x_n)\right)d(w\hspace{0.02cm},y)
\end{equation}
However, by the strong law of large numbers (\ref{eq:mcmclln}), if $y_o \in K$ and $N$ is sufficiently large,
$$   
\frac{1}{N}\sum^N_{n=1}d(z\hspace{0.02cm},x_n) \leq 1 +
\int_M\,d(y_o\hspace{0.02cm},x)\hspace{0.03cm}\pi(dx) + \mathrm{diam}\,K\text{ \;\;(almost-surely)}
$$
Calling this quantity $A$, it follows that for $N$ sufficiently large  (note that this is the same $N$, for all $y$ and $w$ in $K$),
\begin{equation} \label{eq:proofbhatta13}
\left| \mathcal{E}_{\scriptscriptstyle N}(y) - \mathcal{E}_{\scriptscriptstyle N}(w)\right| \leq A\hspace{0.03cm}d(w\hspace{0.02cm},y)\text{ \;\;(almost-surely)}
\end{equation}
From (\ref{eq:squarelipschitz}), it is also possible to show that, for $y$ and $w$ in $K$,
\begin{equation} \label{eq:proofbhatta14}
\left| \mathcal{E}_{\pi}(y) - \mathcal{E}_{\pi}(w)\right| \leq A\hspace{0.03cm}d(w\hspace{0.02cm},y)
\end{equation}
Now, if $y \in K$, let $w(y) \in \lbrace  w_j\rbrace$ be such that $d(w(y),y) < \delta$. Then, for $y$ in $K$,
$$
\left| \mathcal{E}_{\scriptscriptstyle N}(y) - \mathcal{E}_{\pi}(y)\right| \leq
\left| \mathcal{E}_{\scriptscriptstyle N}(y) - \mathcal{E}_{\scriptscriptstyle N}(w(y))\right|+
\left| \mathcal{E}_{\scriptscriptstyle N}(w(y)) - \mathcal{E}_{\pi}(w(y))\right|+
\left| \mathcal{E}_{\pi}(w(y)) - \mathcal{E}_{\pi}(y)\right|
$$
By (\ref{eq:proofbhatta13}) and (\ref{eq:proofbhatta14}), if $N$ is sufficiently large, it follows that
$$
\left| \mathcal{E}_{\scriptscriptstyle N}(y) - \mathcal{E}_{\pi}(y)\right| \leq 2A\delta + \max_{j=1,\ldots, J}\hspace{0.03cm}\left| \mathcal{E}_{\scriptscriptstyle N}(w_j) - \mathcal{E}_{\pi}(w_j)\right|
$$
and (\ref{eq:proofbhatta1}) follows from (\ref{eq:proofbhatta11}), since $\delta > 0$ is arbitrary. 

Second, for $N$ sufficiently large, and for any $C > 0$, it will be proved that there exists a compact $K \subset M$, such that
\begin{equation} \label{eq:proofbhatta2}
  y\notin K \;\;\Longrightarrow\;\; \mathcal{E}_{\scriptscriptstyle N}(y) > C \text{ \;\;(almost-surely)}
\end{equation}
To do so, note from (\ref{eq:emprecursive}), by the triangle inequality 
$$
\mathcal{E}_{\scriptscriptstyle N}(y) \geq \frac{1}{2N}\sum^N_{n=1}(d(y\hspace{0.02cm},\hat{x})-d(\hat{x}\hspace{0.02cm},x_n))^2 
                                                     \, \geq \,\frac{1}{2}d^{\hspace{0.03cm}2}(y\hspace{0.02cm},\hat{x}) - \left(\frac{1}{N}\sum^N_{n=1}d(\hat{x}\hspace{0.02cm},x_n) \right)d(y\hspace{0.02cm},\hat{x})
$$
However, by the strong law of large numbers (\ref{eq:mcmclln}), if $N$ is sufficiently large
$$
\frac{1}{N}\sum^N_{n=1}d(\hat{x}\hspace{0.02cm},x_n) \leq 1 + \int_M\,d(\hat{x}\hspace{0.02cm},x)\hspace{0.03cm}\pi(dx)
$$
Calling this quantity $B$, it follows that for $N$ sufficiently large, 
\begin{equation} \label{eq:proofbhatta21}
\mathcal{E}_{\scriptscriptstyle N}(y) \geq
\frac{1}{2}d^{\hspace{0.03cm}2}(y\hspace{0.02cm},\hat{x}) - B\hspace{0.02cm}d(y\hspace{0.02cm},\hat{x})
\end{equation}
and this directly yields (\ref{eq:proofbhatta2}), since closed and bounded sets are compact (as a consequence of the Hopf-Rinow theorem~\cite{chavel}).

Now, to complete the proof, note the following. By (\ref{eq:proofbhatta2}), for $N$ sufficiently large, there exists a compact $K \subset M$, such that 
$\mathcal{E}_{\scriptscriptstyle N}(y) > \mathcal{E}_{\pi}(\hat{x})+1$ almost-surely, whenever $y \notin K$. That is,
\begin{equation} \label{eq:proofbhatta3}
  \inf_{y\notin K}\hspace{0.02cm}\mathcal{E}_{\scriptscriptstyle N}(y) > \mathcal{E}_{\pi}(\hat{x})+1\text{ \;\;(almost-surely)}
\end{equation}
Moreover, one may always assume that $K$ is a neighborhood of $\hat{x}$. Then, if $B(\hat{x},\epsilon) \subset K$, it follows from (\ref{eq:proofbhatta1}) that, for $N$ sufficiently large, 
$$
 \inf_{y \in B(\hat{x},\epsilon)}\hspace{0.02cm}\mathcal{E}_{\scriptscriptstyle N}(y) <
 \inf_{y \in B(\hat{x},\epsilon)}\hspace{0.02cm}\mathcal{E}_{\pi}(y) + \frac{\epsilon^2}{4}\text{ \;\;(almost-surely)}
$$
or, since $\hat{x}$ is the unique global minimum of $\mathcal{E}_{\pi}(y)$,
\begin{equation} \label{eq:proofbhatta4}
 \inf_{y \in B(\hat{x},\epsilon)}\hspace{0.02cm}\mathcal{E}_{\scriptscriptstyle N}(y) <
 \mathcal{E}_{\pi}(\hat{x}) + \frac{\epsilon^2}{4}\text{ \;\;(almost-surely)}
\end{equation}
But, also by (\ref{eq:proofbhatta1}), for $N$ sufficiently large,
$$
 \inf_{y \in K- B(\hat{x},\epsilon)}\hspace{0.02cm}\mathcal{E}_{\scriptscriptstyle N}(y) >
 \inf_{y \in K- B(\hat{x},\epsilon)}\hspace{0.02cm}\mathcal{E}_{\pi}(y) - \frac{\epsilon^2}{4}\text{ \;\;(almost-surely)}
$$
However, since the variance function $\mathcal{E}_{\pi}$ is $1/2$-strongly convex, with its global minimum at $\hat{x}$,
$$
\mathcal{E}_{\pi}(y) \geq \mathcal{E}_{\pi}(\hat{x}) + \frac{1}{2}d^{\hspace{0.03cm}2}(y\hspace{0.02cm},\hat{x})
$$
and this implies
\begin{equation} \label{eq:proofbhatta5}
 \inf_{y \in K- B(\hat{x},\epsilon)}\hspace{0.02cm}\mathcal{E}_{\scriptscriptstyle N}(y) >
\mathcal{E}_{\pi}(\hat{x}) + \frac{\epsilon^2}{4}\text{ \;\;(almost-surely)}
\end{equation}
Finally, (\ref{eq:proofbhatta3}), (\ref{eq:proofbhatta4}) and (\ref{eq:proofbhatta5}) show that, for $N$ sufficiently large
$$
\inf_{y\in M}\hspace{0.02cm}\mathcal{E}_{\scriptscriptstyle N}(y) = 
\inf_{y\in B(\hat{x},\epsilon)}\hspace{0.02cm}\mathcal{E}_{\scriptscriptstyle N}(y)\text{ \;\;(almost-surely)}
$$
Since $\mathcal{E}_{\scriptscriptstyle N}$ has a unique global minimum $\bar{x}_{\scriptscriptstyle N\,}$, it follows  that $\bar{x}_{\scriptscriptstyle N}$ belongs to the closure of $B(\hat{x},\epsilon)$, almost-surely, when $N$ is sufficiently large. The proof is now complete, since $\epsilon$ is arbitrary. 

%%%metric space
%%%verify all \mathcal{E}
\vfill
\pagebreak

\section{Riemannian gradient descent} \label{sec:markovclt}
Since the minimum mean square error estimator $\hat{x}_{\scriptscriptstyle MMS}$ could not be computed directly, it was approximated by $\bar{x}_{\scriptscriptstyle MMS\hspace{0.03cm}}$, the global minimum of the variance function $\mathcal{E}_{\scriptscriptstyle N\hspace{0.03cm}}$, defined as in (\ref{eq:emprecursive}). This function $\mathcal{E}_{\scriptscriptstyle N}$ being $1/2$-strongly convex, its global minimum can be computed using the Riemannian gradient descent method, which even guarantees an exponential rate of convergence.

This method is here studied from a general point of view. The aim is to minimise a function $f:M\rightarrow \mathbb{R}$, where $M$ is a Hadamard manifold, with sectional curvatures in the interval $[-c^{\hspace{0.02cm}\scriptscriptstyle 2},0]$, and $f$ is an $(\alpha/2)$-strongly convex function.

Recall from (\ref{eq:strongconv}) in \ref{ssec:sqd}, this means $f$ is $(\alpha/2)$-strongly convex along any geodesic in $M$. In particular, for $x\hspace{0.02cm}, y \in M$,
\begin{equation} \label{eq:strongconvtan}
  f(y) - f(x) \geq \langle\mathrm{Exp}^{-1}_x(y),\mathrm{grad}\,f(x)\rangle_x\,+\,(\alpha/2)\hspace{0.02cm}d^{\hspace{0.03cm}2}(x,y)
\end{equation}
This implies that $f$ has compact sublevel sets. Indeed, let $x^*$ be the global minimum of $f$, so $\mathrm{grad}\,f(x^*) = 0$. Putting $x = x^*$ and $y = x$ in (\ref{eq:strongconvtan}), it follows that
\begin{equation} \label{eq:properconv}
  f(x) - f(x^*) \geq (\alpha/2)\hspace{0.02cm}d^{\hspace{0.03cm}2}(x^*\!,x)
\end{equation}
Accordingly, if $S(y)$ is the sublevel set of $y$, then $S(y)$ is contained in the closed ball $\bar{B}(x^*\!,R_y)$, where $R_y = (2/\alpha)(f(y) - f(x^*))$. Therefore, $S(y)$ is compact, since it is closed and bounded~\cite{chavel}.

The Riemannian gradient descent method is based on the iterative scheme
\begin{equation} \label{eq:pgd}
  x^{t+1} = \mathrm{Exp}_{x^t}(-\mu\hspace{0.02cm}\mathrm{grad}\,f(x^t))
\end{equation}
where $\mu$ is a positive step-size, $\mu \leq 1$. If this is chosen sufficiently small, then the iterates $x^t$ remain within the sublevel set $S(x^{\scriptscriptstyle 0})$.
In fact, let $\bar{B}_{\scriptscriptstyle 0} = \bar{B}(x^*\!,R_{x^{\scriptscriptstyle 0}})$ and  
$\bar{B}^\prime_{\scriptscriptstyle 0} = \bar{B}(x^*\!,R_{x^{\scriptscriptstyle 0}} + G)$, where $G$ denotes the supremum of the norm of $\mathrm{grad}\,f(x)$, taken over $x \in \bar{B}_{\scriptscriptstyle 0}$.  Then, let $H^\prime_{\scriptscriptstyle 0}$ denote the supremum of the operator norm of 
$\mathrm{Hess}\,f(x)$, taken over $x \in \bar{B}^\prime_{\scriptscriptstyle 0}$.
\begin{lemma} \label{lem:sublevelgrad}
For the Riemannian gradient descent method (\ref{eq:pgd}), if $\mu \leq 2/\!H^\prime_{\scriptscriptstyle 0\hspace{0.02cm}}$, then the iterates $x^t$ remain within the sublevel set $S(x^{\scriptscriptstyle 0})$. 
\end{lemma} 
Once it has been ensured that the iterates $x^t$ remain within $S(x^{\scriptscriptstyle 0})$, it is even possible to choose $\mu$ in such a way that these iterates achieve an exponential rate of convergence towards $x^*$. This relies on the fact that $x^*$ is a ``strongly attractive" critical point of the vector field $\mathrm{grad}\,f$. Precisely, putting $y = x^*$ in (\ref{eq:strongconvtan}), it follows that
\begin{equation}\label{eq:strongtangattract}
  \langle\mathrm{Exp}^{-1}_x(x^*),\mathrm{grad}\,f(x)\rangle_x \leq -\,(\alpha/2)\hspace{0.02cm}d^{\hspace{0.03cm}2}(x,x^*)+ (f(x^*) - f(x))
\end{equation}
Now, let $C_{\scriptscriptstyle 0} = c\hspace{0.02cm}R_{x^{\scriptscriptstyle 0}}\coth(c\hspace{0.02cm}R_{x^{\scriptscriptstyle 0}})$.
\begin{proposition} \label{prop:pgd}
Let $\bar{H}^\prime_{\scriptscriptstyle 0}  =\max\lbrace H^\prime_{\scriptscriptstyle 0\hspace{0.02cm}},1\rbrace$. If $\mu \leq 1/\!(\bar{H}^\prime_{\scriptscriptstyle 0}C^{\phantom{\prime}}_{\scriptscriptstyle 0})$ (this implies $\mu \leq 2/\!H^\prime_{\scriptscriptstyle 0}$) and $\mu \leq 1/\alpha$, 
\begin{equation} \label{eq:proppgd}
d^{\hspace{0.03cm} 2}(x^t,x^*) \leq (1- \mu\alpha)^t\hspace{0.02cm}d^{\hspace{0.03cm} 2}(x^{\scriptscriptstyle 0},x^*)
\end{equation}
\end{proposition}
The proof of Proposition \ref{prop:pgd} will employ the following lemma.
\begin{lemma}\label{lem:pgd}
Let $\bar{H}^\prime_{\scriptscriptstyle 0}  =\max\lbrace H^\prime_{\scriptscriptstyle 0\hspace{0.02cm}},1\rbrace$.  For any $x \in \bar{B}_{\scriptscriptstyle 0\hspace{0.02cm}}$,
\begin{equation} \label{eq:lempgd}
   \Vert\mathrm{grad}\,f\Vert^2_x \leq 2\bar{H}^\prime_{\scriptscriptstyle 0}(f(x) - f(x^*))
\end{equation}
\end{lemma}
\noindent \textbf{Remark\,:} the rate of convergence predicted by (\ref{eq:proppgd}) is exponential, but depends on the initial guess $x^{\scriptscriptstyle 0}$, through the constants $\bar{H}^\prime_{\scriptscriptstyle 0}$ and $C^{\phantom{\prime}}_{\scriptscriptstyle 0\hspace{0.02cm}}$. This rate can become arbitrarily bad, if $x^{\scriptscriptstyle 0}$ is chosen sufficiently far from $x^*$, since both $\bar{H}^\prime_{\scriptscriptstyle 0}$ and $C^{\phantom{\prime}}_{\scriptscriptstyle 0\hspace{0.02cm}}$ may then become arbitrarily large. By contrast, if $M$ is a Euclidean space (that is, in the limit $c = 0$), $C_{\scriptscriptstyle 0} = 1$, is a constant. \\[0.1cm]
\textbf{Remark\,:} I have never met with a function $f:M \rightarrow \mathbb{R}$ ($M$ a non-Euclidean Hadamard manifold), which is strongly convex, and also has a bounded Hessian. I do not even know whether it is possible or not to construct such a function. \vfill\pagebreak
\noindent \textbf{Proof of Lemma \ref{lem:sublevelgrad}\,:} let $c:[0,1]\rightarrow M$ be the geodesic curve with $c(0) = x^t$ and $c(1) = x^{t+1}$. From (\ref{eq:pgd}), $\dot{c}(0) = -\mu\hspace{0.02cm}\mathrm{grad}\,f(x^t)$. Then, by the Taylor expansion (\ref{eq:taylor2}),
\begin{equation} \label{eq:proofsublevel1}
f(x^{t+1}) = f(x^t) - \mu\hspace{0.02cm}\Vert \mathrm{grad}\,f\Vert^2_{x^t}
+ \frac{1}{2}\,\mathrm{Hess}\,f_{\scriptscriptstyle c(u)}(\dot{c},\dot{c})
\end{equation}
for some $u \in (0,1)$. Assume that $x^t$ belongs to $S(x^{\scriptscriptstyle 0}) \subset \bar{B}_{\scriptscriptstyle 0\hspace{0.02cm}}$. Then, by the triangle inequality,
$$
d(x^*,c(u)) \leq d(x^*,x^t) + d(x^t,c(u)) \leq R_{x^{\scriptscriptstyle 0}} + \mu\hspace{0.02cm}G
$$
where the second inequality follows from the definition of $G$, because $d(x^t,c(u)) = u\Vert \dot{c}(0)\Vert$. Since $\mu \leq 1$, it follows that $d(x^*,c(u))  \leq R_{x^{\scriptscriptstyle 0}} + G$. Therefore, $c(u) \in \bar{B}^\prime_{\scriptscriptstyle 0\hspace{0.02cm}}$. Then, from the definition of $H^\prime_{\scriptscriptstyle 0}$,
$$
\mathrm{Hess}\,f_{\scriptscriptstyle c(u)}(\dot{c},\dot{c}) \leq H^\prime_{\scriptscriptstyle 0}\hspace{0.02cm}\Vert \dot{c}\Vert^2_{\scriptscriptstyle c(u)}
= H^\prime_{\scriptscriptstyle 0}\hspace{0.02cm}\mu^2\Vert \mathrm{grad}\,f\Vert^2_{x^t}
$$
Replacing this into (\ref{eq:proofsublevel1}), 
\begin{equation} \label{eq:proofsublevel2}
f(x^{t+1}) \leq f(x^t) - \mu(1-\mu\hspace{0.02cm}(H^\prime_{\scriptscriptstyle 0}/2))\hspace{0.02cm}\Vert \mathrm{grad}\,f\Vert^2_{x^t}
\end{equation}
Clearly, then, taking $\mu \leq 2/\!H^\prime_{\scriptscriptstyle 0\hspace{0.03cm}}$, it follows that $f(x^{t+1}) \leq f(x^t)$ so that $x^{t+1}$ belongs to $S(x^{\scriptscriptstyle 0})$. The lemma is proved by induction. \\[0.1cm]
\textbf{Proof of Proposition \ref{prop:pgd}\,:} let $c:[0,1]\rightarrow M$ be the geodesic with $c(0) = x^t$ and $c(1) = x^{t+1}$.\hfill\linebreak Note from (\ref{eq:pgd}) that $\dot{c}(0) = -\mu\hspace{0.02cm}\mathrm{grad}\,f(x^t)$. Let $W(x) = d^{\hspace{0.03cm} 2}(x,x^*)/2$, and write down its Taylor expansion (\ref{eq:taylor2}), 
\begin{equation} \label{eq:proofpgd1}
W(x^{t+1}) = W(x^t) - \mu\hspace{0.02cm}\langle\mathrm{grad}W,\mathrm{grad}\,f\rangle_{x^t} + 
 \frac{1}{2}\,\mathrm{Hess}\,W_{\scriptscriptstyle c(u)}(\dot{c},\dot{c})
\end{equation}
for some $u \in (0,1)$. Note that $\mathrm{grad}\,W$ and $\mathrm{Hess}\,W$ are given by (\ref{eq:gradfx}) and (\ref{eq:hesfxcomp}), and also that $x^t$ and $x^{t+1}$ belong to $S(x^{\scriptscriptstyle 0}) \subset \bar{B}_{\scriptscriptstyle 0}\hspace{0.02cm}$, by Lemma \ref{lem:sublevelgrad}, since $\mu \leq 2/\!H^\prime_{\scriptscriptstyle 0}$. Since $S(x^{\scriptscriptstyle 0})$ is a convex set (recall the definition from \ref{ssec:sqd}), $c(u)$ also belongs to $S(x^{\scriptscriptstyle 0}) \subset \bar{B}_{\scriptscriptstyle 0}\hspace{0.02cm}$. By the definition of $C_{\scriptscriptstyle 0\hspace{0.02cm}}$,
$$
\mathrm{Hess}\,W_{\scriptscriptstyle c(u)}(\dot{c},\dot{c}) \leq 
 C_{\scriptscriptstyle 0}\hspace{0.02cm}\Vert \dot{c}\Vert^2_{\scriptscriptstyle c(u)}
= C_{\scriptscriptstyle 0}\hspace{0.02cm}\mu^2\Vert \mathrm{grad}\,f\Vert^2_{x^t}
$$
Replacing into (\ref{eq:proofpgd1}), one now has
\begin{equation} \label{eq:proofpgd2}
W(x^{t+1}) \leq W(x^t) + \mu\hspace{0.02cm}\langle\mathrm{Exp}^{-1}_{x^t}(x^*),\mathrm{grad}\,f\rangle_{x^t} + 
(C_{\scriptscriptstyle 0}/2)\hspace{0.02cm}\mu^2\Vert \mathrm{grad}\,f\Vert^2_{x^t}
\end{equation}
Therefore, by (\ref{eq:strongtangattract}) and (\ref{eq:lempgd}),
\begin{equation} \label{eq:proofpgd3}
W(x^{t+1}) \leq W(x^t)(1-\mu\alpha) + \mu(1-\mu(\bar{H}^\prime_{\scriptscriptstyle 0}C^{\phantom{\prime}}_{\scriptscriptstyle 0}))(f(x^*) - f(x))
\end{equation}
If $\mu \leq 1/\!(\bar{H}^\prime_{\scriptscriptstyle 0}C^{\phantom{\prime}}_{\scriptscriptstyle 0})$, then (\ref{eq:proofpgd3}) implies $W(x^{t+1}) \leq (1-\mu\alpha)W(x^t)$, because $f(x^*) - f(x) \leq 0$.\hfill\linebreak The proposition easily follows by induction, since $1-\mu\alpha \geq 0$. \\[0.1cm]
\textbf{Proof of Lemma \ref{lem:pgd}\,:} let $c$ denote the geodesic with $c(0) = x$ and $\dot{c}(0) = (-1/\bar{H}^\prime_{\scriptscriptstyle 0})\hspace{0.02cm}\mathrm{grad}\,f(x)$. By the same arguments as in the proof of Lemma \ref{lem:sublevelgrad}, one has that $c(u) \in \bar{B}^\prime_{\scriptscriptstyle 0}$ for all $u \in [0,1]$. Therefore, letting $y = c(1)$ and writing down the Taylor expansion (\ref{eq:taylor2}),
$$
f(y) - f(x) \leq  (-1/\bar{H}^\prime_{\scriptscriptstyle 0})\Vert\mathrm{grad}\,f\Vert^2_x + (\bar{H}^\prime_{\scriptscriptstyle 0}/2)\Vert (1/\bar{H}^\prime_{\scriptscriptstyle 0})\hspace{0.02cm}\mathrm{grad}\,f(x)\Vert^2_x = (-1/2\bar{H}^\prime_{\scriptscriptstyle 0})\Vert\mathrm{grad}\,f\Vert^2_x
$$
Multiplying this inequality by $-2\bar{H}^\prime_{\scriptscriptstyle 0\hspace{0.02cm}}$, 
$$
2\bar{H}^\prime_{\scriptscriptstyle 0}(f(x) - f(y)) \geq \Vert \mathrm{grad}\,f\Vert^2_x
$$
Now, (\ref{eq:lempgd}) obtains by noting that $f(x) - f(x^*) \geq f(x) - f(y)$.

%%% correct for the last lemma
%%%only H prime
%%%max it with 1

%let $x \in C$ and $y = \mathrm{Exp}_x((-1/L_f)\hspace{0.02cm}\mathrm{grad}\,f(x))$. Note that, since $L_f \geq 1$,
%$$
%d(x,y) = (1/L_f)\Vert \mathrm{grad}\,f(x)\Vert \leq G
%$$
%so, by the triangle inequality, 
%$$
%d(x^*,y) \leq d(x^*,x) + G \leq \mathrm{diam}\,C + G
%$$
%and, therefore, $y \in C^\prime$. Let $c:[0,1]\rightarrow M$ be the geodesic curve with $c(0) = x$ and $c(1) = y$. Since $C^\prime$ is convex, with $c(0) \in C^\prime$ and $c(1) \in C^\prime$, it follows that $c(t) \in C^\prime$ for all $t \in [0,1]$.
%
%Write down the Taylor expansion (\ref{eq:taylor2}) of $f$ along the geodesic curve $c$. From (\ref{eq:LVLf}),
%$$
%f(c(1)) - f(c(0)) \leq (-1/L_f)\Vert\mathrm{grad}\,f(x)\Vert^2_x\,+\, (L_f/2)\Vert (-1/L_f)\mathrm{grad}\,f(x)\Vert^2_x
%$$
%In other words,
%$$
%f(y) - f(x) \leq -\frac{1}{2L_f}\Vert\mathrm{grad}\,f(x)\Vert^2
%$$
%Finally, since $f(x^*) \leq f(y)$, for any $y \in M$,
%$$
%f(x^*) - f(x) \leq -\frac{1}{2L_f}\Vert\mathrm{grad}\,f(x)\Vert^2
%$$
%which is the same as (\ref{eq:lempgd}). 

\section{A volume growth lemma} \label{sec:mcmc_vollema}
Lemma \ref{lemm:mcmc_volume} will be used in the proof of Proposition \ref{prop:gergodic}, to be carried out in  \ref{sec:mcmcproof}. This lemma is of a purely geometric content, and is therefore considered separately, beforehand. 

Let $M$ be a Riemannian symmetric space, which belongs to the non-compact case (see \ref{ssec:sspace}). Then, in particular, $M$ is a Hadamard manifold. 

Fix $x^* \in M$, and let $(r,\theta)$ be geodesic spherical coordinates, with origin at $x^*$. Any $z \in M$, other than $x^*$, is uniquely determined by its coordinates $(r,\theta)$, and will be written $z(r,\theta)$. 

Recall the volume density function $\det(\mathcal{A}(r,\theta))$, from the integral formula (\ref{eq:integralsphericalhadamard}). This will be denoted $\lambda(r,\theta) =\det(\mathcal{A}(r,\theta))$. 

Essentially, the following lemma states the logarithmic rate of growth of the volume density function $\lambda(r,\theta)$ is bounded at infinity.
\begin{lemma} \label{lemm:mcmc_volume}
Let $M$ be a Riemannian symmetric space, which belongs to the non-compact case. Fix $x^* \in M$ and denote $r(x) = d(x^*,x)$ for $x \in M$. Then, for any $R > 0$,
\begin{equation} \label{eq:volmcmc}
\limsup_{r(x) \rightarrow \infty}\hspace{0.04cm}\frac{\sup_{\scriptscriptstyle z(r,\theta) \in B(x,R)}\hspace{0.04cm} \lambda(r,\theta)}
{\inf_{\scriptscriptstyle z(r,\theta) \in B(x,R)}\hspace{0.04cm} \lambda(r,\theta)}\,<\,\infty
\end{equation}
\end{lemma}
The proof of this lemma proceeds in the following way. Identify the unit sphere in $T_{x^*}M$ with $S^{n-1}$, and consider for $\theta \in S^{n-1}$ the self-adjoint curvature operator $R_\theta:T_{x^*M} \rightarrow T_{x^*}M$, given by
$$
R_\theta(v) = -R(\theta,v)\hspace{0.02cm} \theta \hspace{0.2cm};\hspace{0.2cm} v \in T_{x^*}M
$$
Recall that the Riemann curvature tensor is parallel (because $M$ is a symmetric space). Then, from (\ref{eq:jacobioperator}) and the definition of $\mathcal{A}(r,\theta)$, it follows that $\mathcal{A}(r,\theta)$ solves the Jacobi equation
\begin{equation} \label{eq:jacobiss}
  \mathcal{A}^{\prime\prime} - R_{\theta}\hspace{0.02cm}\mathcal{A} = 0 \hspace{1cm} \mathcal{A}(0) = 0 \,,\, \mathcal{A}^\prime(0) = \mathrm{Id}_{x^*}
\end{equation}
where the prime denotes differentiation with respect to $r$. At present, all the eigenvalues of $R_\theta$ are positive. If $c^{\scriptscriptstyle 2}(\theta)$ runs through these eigenvalues, then it follows from (\ref{eq:jacobiss}) that
\begin{equation} \label{eq:lambdamc1}
  \lambda(r,\theta) \,=\, \prod_{c(\theta)}\left( \frac{\sinh(c(\theta)\hspace{0.02cm}r)}{c(\theta)}\right)^{\!\!m_{c(\theta)}}
\end{equation}
where $m_{c(\theta)}$ denotes the multiplicity of the eigenvalue $c^{\scriptscriptstyle 2}(\theta)$ of $R_\theta\hspace{0.02cm}$.

It is possible to express (\ref{eq:lambdamc1}) in a different form. Let $M = G/K$ where $K$ is the stabiliser in $G$ of $x^*$. Let $\mathfrak{g}$ and $\mathfrak{k}$ be the Lie algebras of $G$ and $K$, and $\mathfrak{g} = \mathfrak{k} + \mathfrak{p}$ the corresponding Cartan decomposition. Let $\mathfrak{a}$ be a maximal Abelian subspace of $\mathfrak{p}$, and recall that it is always possible to write $r\theta = \mathrm{Ad}(k)\,a$ for some $k \in K$ and $a \in \mathfrak{a}$ (see Lemma 6.3, Chapter V, in~\cite{helgason}). In this notation, $r = \Vert a \Vert_{x^*}$ and $c(\theta) = \lambda(a)/ \Vert a \Vert_{x^*}\hspace{0.02cm}$, where $\lambda$ is a positive roots of $\mathfrak{g}$ with respect to $\mathfrak{a}$, with multiplicity $m_\lambda = m_{c(\theta)}$\hspace{0.02cm} (see Lemma 2.9, Chapter VII, in~\cite{helgason}). Replacing into (\ref{eq:lambdamc1}) gives
\begin{equation} \label{eq:lambdamc2}
 \lambda(r,\theta) \,=\, \prod_{\lambda \in \Delta_+}\left(  \frac{\sinh(\lambda(a))}{\lambda(a)/ \Vert a \Vert}\right)^{\!\!m_\lambda}
\end{equation}
Here, if the right-hand side is denoted by $f(a)$, then it is elementary that $\log\hspace{0.02cm}f(a)$ is a Lipschitz function, on the complement of any bounded subset of $\mathfrak{a}$ which contains the zero element of $\mathfrak{a}$.

Returning to (\ref{eq:volmcmc}), let the supremum in the numerator be achieved at $(r_{\max\hspace{0.02cm}},\theta_{\max})$ and the infimum in the denominator be achieved at $(r_{\min\hspace{0.02cm}},\theta_{\min})$. Let $(k_{\max},a_{\max})$ and $(k_{\min},a_{\min})$ be corresponding values of $k$ and $a$. Note that for $z(r,\theta) \in B(x,R)$, by the triangle inequality, $r \geq r(x) - R$. But, since $r = \Vert a \Vert_{x^*\hspace{0.02cm}}$, this also means $\Vert a \Vert_{x^*\hspace{0.02cm}}  \geq r(x) - R$. 

Therefore, if $r(x) > R$ then,  as stated above, $\log\hspace{0.02cm}f(a)$ is a Lipschitz function, on the set of  $a$ such that $\Vert a \Vert_{x^*\hspace{0.02cm}}  \geq r(x) - R$. If $\mathrm{C}$ is the corresponding Lipschitz constant,
\begin{equation} \label{eq:prooflemvolmc1}
\frac{\sup_{\scriptscriptstyle z(r,\theta) \in B(x,R)}\hspace{0.04cm} \lambda(r,\theta)}
{\inf_{\scriptscriptstyle z(r,\theta) \in B(x,R)}\hspace{0.04cm} \lambda(r,\theta)}\,\leq\, \exp[\mathrm{C}\hspace{0.02cm}\Vert a_{\max} - a_{\min}\Vert_{x^*}]
\end{equation}
\indent Now, (\ref{eq:volmcmc}) will follow by showing that $\Vert a_{\max} - a_{\min}\Vert_{x^*} < 2R$ wherever $r(x) > R$. 

To do so, let $z_{\max} = z(r_{\max\hspace{0.02cm}},\theta_{\max})$ and $z_{\min} = z(r_{\min\hspace{0.02cm}},\theta_{\min})$, and note $d(z_{\max\hspace{0.02cm}},z_{\min}) \leq 2R$.\hfill\linebreak If $c: [0,1]\rightarrow M$ is a geodesic curve with $c(0) = z_{\min}$ and $c(1) = z_{\max\hspace{0.02cm}}$, then
\begin{equation} \label{eq:prooflemvolmc2}
\int^{\scriptscriptstyle 1}_{\scriptscriptstyle 0}\,\Vert \dot{c}(t)\Vert_{\scriptscriptstyle c(t)}\hspace{0.03cm}dt = d(z_{\max\hspace{0.02cm}},z_{\min}) \leq 2R
\end{equation}
On the other hand, if $c(t) = c(r(t)\hspace{0.02cm},\theta(t))$, then it is possible to write $r(t)\hspace{0.02cm}\theta(t) = \mathrm{Ad}(k(t))\,a(t)$, where $k(t)$ and $a(t)$ are differentiable curves in $K$ and $\mathfrak{a}$. It will be shown below that this implies
\begin{equation} \label{eq:prooflemvolmc3}
\Vert \dot{c}(t)\Vert^2_{\scriptscriptstyle c(t)} = \Vert \dot{a}(t)\Vert^2_{\scriptscriptstyle x^*} + \sum_{\lambda \in \Delta_+} \sinh^2(\lambda(a(t))\hspace{0.02cm}\Vert \dot{k}_{\lambda}(t)\Vert^2_{\scriptscriptstyle x^*}
\end{equation}
where $\dot{k}_{\lambda}(t)$ is defined following (\ref{eq:ncmetric3}), below. Finally, from (\ref{eq:prooflemvolmc2}) and (\ref{eq:prooflemvolmc3}), it follows that
$$
\Vert a_{\max} - a_{\min}\Vert_{x^*} \,\leq\, \int^{\scriptscriptstyle 1}_{\scriptscriptstyle 0}\Vert \dot{a}(t)\Vert_{\scriptscriptstyle x^*}\hspace{0.03cm}dt
\,\leq\, \int^{\scriptscriptstyle 1}_{\scriptscriptstyle 0}\Vert \dot{c}(t)\Vert_{\scriptscriptstyle c(t)}\hspace{0.03cm}dt \leq 2R
$$
Replacing into (\ref{eq:prooflemvolmc1}), this yields
$$
\frac{\sup_{\scriptscriptstyle z(r,\theta) \in B(x,R)}\hspace{0.04cm} \lambda(r,\theta)}
{\inf_{\scriptscriptstyle z(r,\theta) \in B(x,R)}\hspace{0.04cm} \lambda(r,\theta)}\,\leq\, \exp(2\mathrm{C}R)
$$
for all $x$ such that $r(x) > R$. However, this immediately implies (\ref{eq:volmcmc}). \\[0.1cm]
\textbf{Proof of (\ref{eq:prooflemvolmc3})\,:} in the notation of \ref{ssec:sspace}, $c(t) = \varphi(s(t)\hspace{0.02cm},a(t))$, where $s(t)$ is the representative of $k(t)$ in the quotient $K/K_{\mathfrak{a}\hspace{0.03cm}}$. Recall that $\varphi(s\hspace{0.02cm},a) = \mathrm{Exp}_o(\beta(s\hspace{0.02cm},a))$ where $\beta(s,a) = \mathrm{Ad}(s)\,a$\hfill\linebreak (the dependence on $t$ is now suppressed). Then, by differentiating with respect to $t$, 
$$
\dot{\beta}(s,a) \,=\, \mathrm{Ad}(s)\left(\dot{a} \,+\, [\dot{s},a]\right) 
$$
Further, by replacing from (\ref{eq:dexpss}),
$$
\dot{c} \,=\, \exp(r\theta)\cdot\mathrm{sh}(R_{r\theta})( \dot{\beta}(s,a))
$$
However, $\mathrm{Ad}(s)$ preserves norms, and $\mathrm{Ad}(s^{\scriptscriptstyle -1})\circ R_{r\theta}\circ \mathrm{Ad}(s) = R_{a\hspace{0.02cm}}$, as in \ref{eq:raeigen}). Therefore,
\begin{equation} \label{eq:ncmetric1}
\Vert \dot{c}\Vert^2_{\scriptscriptstyle c} \,=\, \left\Vert\mathrm{sh}(R_{a})\left(\dot{a} \,+\, [\dot{s},a]\right)\right\Vert^2_{\scriptscriptstyle x^*}
\end{equation}
and from the definition of $\mathrm{sh}(R_{a})$,
\begin{equation} \label{eq:ncmetric2}
\mathrm{sh}(R_{a}) \,=\, \Pi_{\mathfrak{a}} + \sum_{\lambda \in \Delta_+}\frac{\sinh(\lambda(a))}{\lambda(a)}\,\Pi_{\lambda}
\end{equation}
Now, one has the orthogonal decomposition $\dot{s} = \sum_{\lambda \in \Delta_+} (\xi_\lambda + d\hspace{0.02cm}\theta(\xi_\lambda))$ where $[a,\xi_\lambda] = \lambda(a)\hspace{0.03cm}\xi_\lambda$ and $d\hspace{0.02cm}\theta$ was introduced before (\ref{eq:sscommute}) (see Lemma 3.6, Chapter VI, in~\cite{helgason}). In turn, this yields the orthogonal decomposition
\begin{equation} \label{eq:ncmetric3}
[a,\dot{s}] \,=\,\sum_{\lambda \in \Delta_+} \lambda(a)\hspace{0.03cm}(\xi_\lambda - d\hspace{0.02cm}\theta(\xi_\lambda))
\end{equation}
Letting $\dot{s}_\lambda = (\xi_\lambda - d\hspace{0.02cm}\theta(\xi_\lambda))$, it follows from (\ref{eq:ncmetric1}) and (\ref{eq:ncmetric2}) that
$$
\Vert \dot{c}\Vert^2_{\scriptscriptstyle c} \,=\, \Vert \dot{a}\Vert^2_{\scriptscriptstyle x^*} + \sum_{\lambda \in \Delta_+} \sinh^2(\lambda(a))\hspace{0.02cm}\Vert \dot{s}_{\lambda}\Vert^2_{\scriptscriptstyle x^*}
$$
This is the same as (\ref{eq:prooflemvolmc3}), once $\dot{k}$ is identified with its representative $\dot{s}$. 

\section{Proof of geometric ergodicity} \label{sec:mcmcproof}
The proof of Proposition \ref{prop:gergodic} relies on the so-called geometric drift condition. This condition requires that there exist a function $V:M \rightarrow \mathbb{R}$ such that
\begin{equation} \label{eq:VVbis}
  V(x) \geq \max\left\lbrace 1\hspace{0.02cm},d^{\hspace{0.03cm} 2}(x,x^*)\right\rbrace \text{ for some } x^* \in M
\end{equation}
\begin{equation} \label{eq:gdrift}
PV(x) \leq \,\lambda\hspace{0.02cm}V(x) + b\hspace{0.02cm}\mathbf{1}_{\scriptscriptstyle C}(x) \hspace{3cm}
\end{equation}
for some $\lambda \in (0,1)$ and $b \in (0,\infty)$, and where $C$ is a small set for $P$ (for the definition, see~\cite{tweedie}).
If the geometric drift condition (\ref{eq:gdrift}) is verified, then the geometric  ergodicity condition (\ref{eq:gergodic}) holds~\cite{tweedie}.

The proof is a generalisation of the proof carried out in the special case where $M$ is a Euclidean space, in~\cite{jarner}. The idea is to use Assumptions (a1)--(a3) to show that the following two conditions hold,
\begin{equation} \label{eq:ergod1}
\limsup_{r(x)\rightarrow \infty}\hspace{0.03cm}\frac{PV(x)}{V(x)}\,<\,1
\end{equation}
\begin{equation} \label{eq:ergod2}
\phantom{li}\, \sup_{x \in M}\hspace{0.03cm} \frac{PV(x)}{V(x)}\,<\,\infty
\end{equation}
where $r(x) = d(x^*,x)$, and $V(x) = a\hspace{0.03cm}\pi^{\scriptscriptstyle -\frac{1}{2}}(x)$ with $a$ chosen so $V(x) \geq 1$ for all $x \in M$. However, under Assumption (a3), these  two conditions are shown to imply (\ref{eq:gdrift}).
\begin{lemma} \label{lem:mcmc1}
 Let $(x_n)$ be a Markov chain in $M$, with transition kernel (\ref{eq:hmP}), with proposed transition density $q(x\hspace{0.02cm},y) = q(d(x\hspace{0.02cm},y))$, and continuous, strictly positive invariant density $\pi$. Moreover, assume the proposed transition density satisfies Assumption (a3). If Conditions (\ref{eq:ergod1}) and (\ref{eq:ergod2})  are verified, then the geometric drift condition (\ref{eq:gdrift}) holds.
\end{lemma}
On the other hand, (\ref{eq:VVbis}) (which is just the same as (\ref{eq:VV})) is a straightforward result of Assumption (a1), which implies the existence of strictly positive $\mu, R$ and $\pi_{\scriptscriptstyle R}$ such that
\begin{equation} \label{eq:preVV}
  r(x) \geq R \;\Longrightarrow\; \pi(x) \leq \pi_{\scriptscriptstyle R}\hspace{0.02cm}\exp\left(-\mu\hspace{0.02cm}r^2(x)\right)
\end{equation}
Then, to obtain (\ref{eq:VVbis}), it is enough to chose $a = \max\left\lbrace 1,R^{\hspace{0.03cm}\scriptscriptstyle 2},\pi^{\scriptscriptstyle 1/2}_{\scriptscriptstyle R},2\hspace{0.02cm}\mu^{\scriptscriptstyle -1}\right\rbrace\hspace{0.02cm}$. \\[0.1cm]
\textbf{Proof of Lemma \ref{lem:mcmc1}\,:} the proof is almost identical to the proofs for random-walk Metropolis chains in Euclidean space~\cite{jarner}\cite{rroberts}. The main point is that Assumption (a3) implies that every non-empty bounded subset of $M$ is a small set for the transition kernel $P$ in (\ref{eq:hmP}). With this in mind, the geometric drift condition (\ref{eq:gdrift}) follows almost directly from the two conditions (\ref{eq:ergod1}) and (\ref{eq:ergod2}). Indeed, (\ref{eq:ergod1}) implies that there exist $\lambda \in (0,1)$ and $R \in (0,\infty)$ such that
$$
r(x) \geq R \;\Longrightarrow\; PV(x) \leq \lambda\hspace{0.02cm}V(x)
$$
That is, (\ref{eq:gdrift}) is verified on $M - C$, where $C$ is the open ball $B(x^*,R)$. In addition, by (\ref{eq:ergod2}),
$$
b = \left[ \sup_{x \in B(x^*,R)} V(x)\right]\,\left[\hspace{0.04cm} \sup_{x \in M}\hspace{0.03cm} \frac{PV(x)}{V(x)}\right]\,<\,\infty
$$ 
Therefore, (\ref{eq:gdrift}) is also verified on $C$, since for $x \in C$,
$$
PV(x) \leq b \leq \lambda\hspace{0.02cm}V(x) + b
$$
Thus, (\ref{eq:gdrift}) is verified throughout $M$. It remains to note that $C$ is a small set, since it is bounded. \vfill\pagebreak

Now, the aim is to establish the two conditions (\ref{eq:ergod1}) and (\ref{eq:ergod2}). These will follow from Propositions \ref{prop:proofgergodic1} and \ref{prop:proofgergodic2}, below. Consider the proposed transition kernel
\begin{equation}\label{eq:propkernelQ} 
Qf(x) \,=\, \int_M\,q(x\hspace{0.02cm},y)f(y)\hspace{0.02cm}\mathrm{vol}(dy)
\end{equation}
for any bounded measurable function $f:M\rightarrow \mathbb{R}$. If $f$ is the indicator function of a measurable set $A$, then it is usual to write $Qf(x) = Q(x,A)$. For $x \in M$, consider its acceptance region
$$
A(x) = \left\lbrace y \in M: \pi(y) \geq \pi(x)\right\rbrace
$$
\begin{proposition} \label{prop:proofgergodic1}
 Under the assumptions of Proposition \ref{prop:gergodic}, the following limit holds
\begin{equation} \label{eq:limQ}
  \liminf_{r(x) \rightarrow \infty} Q(x\hspace{0.02cm},A(x))\,> \, 0
\end{equation}
\end{proposition}
\begin{proposition} \label{prop:proofgergodic2}
 Under the assumptions of Proposition \ref{prop:gergodic}, if (\ref{eq:limQ}) holds, then the two conditions (\ref{eq:ergod1}) and (\ref{eq:ergod2}) are verified, where $V(x) = a\hspace{0.03cm}\pi^{\scriptscriptstyle -\frac{1}{2}}(x)$ with $a$ chosen so $V(x) \geq 1$ for all $x \in M$.
\end{proposition}
The proof of these two propositions will use the following fact, concerning the contour manifolds of the probability density function $\pi(x)$. For $x \in M$, the contour manifold of $x$ is the set $C_{x}$ of all $y \in M$ such that $\pi(y) = \pi(x)$. This is a hypersurface in $M$, whenever $\pi(x)$ is a regular value of  $\pi$ (by the ``regular level set theorem"~\cite{lee}). \\[0.1cm]
\indent \textbf{fact\,:} if $r(x)$ is sufficiently large, then $C_{x}$ can be parameterised by the unit sphere in $T_{x^*}M$. Precisely, it is possible to write 
\begin{equation} \label{eq:contourm}
  C_{x} = \left\lbrace \mathrm{Exp}_{x^*}\left( c(v)\hspace{0.02cm}v\right)\,; v \in S_{x^*}M\right\rbrace
\end{equation}
where $c$ is a positive continuous function on $S_{x^*}M$,  the set of unit vectors $v$ in $T_{x^*}M$. Moreover, $A(x)$ is exactly the region inside of $C_{x\hspace{0.03cm}}$. Precisely, $y \in A(x)$ if and only if $y = \mathrm{Exp}_{x^*}(c\hspace{0.02cm}v)$ where $v \in S_{x^*}M$ and $c \leq c(v)$. \\[0.1cm]
\textbf{Proof of Proposition \ref{prop:proofgergodic1}\,:} by Assumption (a2), there exist $\delta > 0$ and $R > 0$ such that
\begin{equation} \label{eq:proofproofgergodic11}
 r(y) \geq R\;\Longrightarrow\; \langle\mathrm{grad}\,r,n\rangle_{y}\,<\,-\delta
\end{equation}
Let $-c^{\scriptscriptstyle 2}$ be a lower bound on the sectional curvatures of $M$, and $\Lambda$ be a positive number with
\begin{equation} \label{eq:Lambda}
 (\dim\hspace{0.03cm}M)^{\scriptscriptstyle \frac{1}{2}}\hspace{0.03cm}\Lambda \,\leq\, \frac{\delta}{2c}\hspace{0.02cm}\tanh(c\hspace{0.02cm}R)
\end{equation}
Now, for any $x \in M$ with $r(x) \geq R + \Lambda$, consider the set
$$
\Omega(x) \,=\,\left\lbrace \mathrm{Exp}_{x}(-a\hspace{0.02cm}u)\,; a \in (0,\Lambda)\,,u \in S_xM\,, \Vert \mathrm{grad}\,r(x) - u\Vert_{x} \hspace{0.02cm}\leq \frac{\delta}{2}\right\rbrace
$$
Let $y = \mathrm{Exp}_{x}(-a\hspace{0.02cm}u)$ be a point in $\Omega(x)$, and $\gamma(t)$ the unit-speed geodesic with $\gamma(0) = x$ and $\gamma(a) = y$. It is first proved that
\begin{equation} \label{eq:proofproofgergodic12}
  \langle\dot{\gamma}\hspace{0.02cm},n\rangle_{\gamma(t)}\, > 0 \hspace{0.3cm} \text{for } t \in (0,a)
\end{equation}
Indeed, the left-hand side of (\ref{eq:proofproofgergodic12}) may be written
$$
  \langle\dot{\gamma}\hspace{0.02cm},n\rangle_{\gamma(t)} =   -\,\langle\mathrm{grad}\,r\hspace{0.02cm},n\rangle_{\gamma(t)}\hspace{0.03cm}
+\langle\dot{\gamma}+\mathrm{grad}\,r\hspace{0.02cm},n\rangle_{\gamma(t)}
$$
Then, if $\Pi^t_{\scriptscriptstyle 0}$ denotes the parallel transport along $\gamma$ from $\gamma(0) = x$ to $\gamma(t)$,
\begin{equation} \label{eq:proofproofgergodic13}
  \langle\dot{\gamma}\hspace{0.02cm},n\rangle_{\gamma(t)} =   -\,\langle\mathrm{grad}\,r\hspace{0.02cm},n\rangle_{\gamma(t)}
+\langle\Pi^t_{\scriptscriptstyle 0}(\mathrm{grad}\,r(x) - u)\hspace{0.02cm},n\rangle_{\gamma(t)}
+\langle\mathrm{grad}\,r - \Pi^t_{\scriptscriptstyle 0}(\mathrm{grad}\,r(x))\hspace{0.02cm},n\rangle_{\gamma(t)}
\end{equation}
which may be checked by adding together the three terms, and noting that $\dot{\gamma}(t) = \Pi^t_{\scriptscriptstyle 0}(-u)$, since $\gamma$ is a geodesic with $\dot{\gamma}(0) = -u$. But, by the triangle inequality
$$
r(\gamma(t)) \,\geq\, r(x) - d(x\hspace{0.02cm},\gamma(t))\,>\, (R+\Lambda) - \Lambda = R
$$
since $d(x^*,x) = r(x) \geq  R + \Lambda$ and $d(x\hspace{0.02cm},\gamma(t)) \leq a \leq \Lambda$. Thus, it follows from (\ref{eq:proofproofgergodic11})
\begin{equation} \label{eq:proofproofgergodic14}
 - \langle\mathrm{grad}\,r,n\rangle_{\gamma(t)}\,>\,\delta
\end{equation}
Moreover, since the parallel transport $\Pi^t_{\scriptscriptstyle 0}$ preserves norms, and since by definition of $\Omega(x)$, $\Vert\mathrm{grad}\,r(x) - u\Vert_{x} \leq \delta/2$, it follows from the Cauchy-Schwarz inequality
\begin{equation} \label{eq:proofproofgergodic15}
\langle\Pi^t_{\scriptscriptstyle 0}(\mathrm{grad}\,r(x) - u)\hspace{0.02cm},n\rangle_{\gamma(t)} \,\geq - \Vert \Pi^t_{\scriptscriptstyle 0}(\mathrm{grad}\,r(x) -u)\Vert_{x}= -\Vert\mathrm{grad}\,r(x) - u\Vert_{x} \geq -\delta/2
\end{equation}
On the other hand, let $(e_{\scriptscriptstyle i}\,;1,\ldots,n)$ be a parallel orthonormal base, along the geodesic $\gamma$. Then,
$$
\langle \mathrm{grad}\,r - \Pi^t_{\scriptscriptstyle 0}(\mathrm{grad}\,r(x))\hspace{0.02cm},e_{\scriptscriptstyle i}\rangle_{\gamma(t)} \,=\,
\int^t_{\scriptscriptstyle 0}\left\langle\mathrm{Hess}\,r\cdot\dot{\gamma}\hspace{0.02cm},e_{\scriptscriptstyle i}\right\rangle_{\gamma(s)}ds
$$
But, according to (\ref{eq:hescomp}) from Theorem \ref{th:comp},
$$
\int^t_{\scriptscriptstyle 0}\left\langle\mathrm{Hess}\,r\cdot\dot{\gamma}\hspace{0.02cm},e_{\scriptscriptstyle i}\right\rangle_{\gamma(s)}ds\, \leq \int^t_{\scriptscriptstyle 0}c\hspace{0.02cm}\coth\left(c\hspace{0.02cm} r(\gamma(s))\right)ds\,\leq \Lambda\hspace{0.03cm}c\hspace{0.02cm}\coth\left(c\hspace{0.02cm} R\right)
$$
Thus, using (\ref{eq:Lambda}), it follows by the Cauchy-Schwarz inequality 
\begin{equation} \label{eq:proofproofgergodic16}
\langle \mathrm{grad}\,r - \Pi^t_{\scriptscriptstyle 0}(\mathrm{grad}\,r(x))\hspace{0.02cm},n\rangle_{\gamma(t)}\hspace{0.03cm}
\geq -\delta/2
\end{equation}
Finally, by adding (\ref{eq:proofproofgergodic14}) to (\ref{eq:proofproofgergodic15}) and (\ref{eq:proofproofgergodic16}), it follows from (\ref{eq:proofproofgergodic13})
$$
 \langle\dot{\gamma}\hspace{0.02cm},n\rangle_{\gamma(t)}\, > \delta - \delta/2 - \delta/2 = 0
$$
which is the same as (\ref{eq:proofproofgergodic12}). Moving on, from (\ref{eq:proofproofgergodic12}), it is possible to prove that
\begin{equation} \label{eq:proofproofgergodic17}
  \Omega(x) \subset A(x)
\end{equation}
for all $x$ such that $r(x) \geq R + \Lambda$, where $A(x)$ is the acceptance region of $x$, defined after (\ref{eq:propkernelQ}). 

To prove (\ref{eq:proofproofgergodic17}), consider $y \in \Omega(x)$ and $\gamma(t)$ as before, with $\gamma(0) = x$ and $\gamma(a) = y$. Now, assume that $y \in C_{x\hspace{0.03cm}}$, the contour manifold of $x$, defined in (\ref{eq:contourm}). Then, $\pi(\gamma(0)) = \pi(\gamma(a))$, so that, by the mean-value theorem, there exists $t \in (0,a)$ such that
$$
\frac{d}{dt}\pi(\gamma(t)) = \langle\dot{\gamma}(t),\mathrm{grad}\,\pi\rangle_{\gamma(t)} = 0
$$
But, from the definition of $n(x)$, this implies
$$
\langle \dot{\gamma}(t),n\rangle_{\gamma(t)}=\,\Vert\mathrm{grad}\,\pi(x)\Vert^{\scriptscriptstyle -1}\hspace{0.03cm}\langle\dot{\gamma}(t),\mathrm{grad}\,\pi\rangle_{\gamma(t)} = 0
$$
in contradiction with (\ref{eq:proofproofgergodic12}). Thus, the assumption that $y \in C_{x}$ cannot hold. Since $y \in \Omega(x)$ is arbitrary, this means that
\begin{equation} \label{eq:proofproofgergodic18}
 \Omega(x)\,\cap\, C_{x}\,=\, \varnothing
\end{equation}
However, note that $y_* = \mathrm{Exp}_{x}(-a\hspace{0.02cm}\mathrm{grad}\,r(x))$ belongs to $\Omega(x)$, as can be seen from the definition of $\Omega(x)$. Also, since $r(y_*) = r(x) - a$, it follows that $y_*$ is inside of $C_{x\hspace{0.03cm}}$. Therefore, $y_* \in A(x)$, and the intersection of $\Omega(x)$ and $A(x)$ is non-empty. Finally, it is enouh to note that the set $\Omega(x)$ is connected, since it is the image under $\mathrm{Exp}_{x}$ of a connected set. This implies that, if the intersection of $\Omega(x)$ and $R(x)$, the complement of $A(x)$, were non-empty, then $\Omega(x)$ would also intersect $C_{x\hspace{0.03cm}}$. Clearly, this would be in contradiction with (\ref{eq:proofproofgergodic18}). \vfill\pagebreak
Using (\ref{eq:proofproofgergodic17}), it is now possible to prove (\ref{eq:limQ}).  Indeed, for $x$ such that $r(x) \geq R+\Lambda$, it follows from (\ref{eq:proofproofgergodic17}) that 
\begin{equation} \label{eq:proofproofgergodic19}
  Q(x\hspace{0.02cm},A(x)) \geq Q(x\hspace{0.02cm},\Omega(x)) = \int_{\Omega(x)}\,q(x\hspace{0.02cm},y)\hspace{0.02cm}\mathrm{vol}(dy)
\end{equation}
where the last equality follows from (\ref{eq:propkernelQ}). However, by Assumption (a3),
\begin{equation} \label{eq:proofproofgergodic20}
\int_{\Omega(x)}\,q(x\hspace{0.02cm},y)\hspace{0.02cm}\mathrm{vol}(dy) \,\geq\,\varepsilon_{\scriptscriptstyle q}\times\mathrm{vol}\left(\Omega(x)\cap B(x,\delta_{\scriptscriptstyle q})\right)
\end{equation}
Now, to prove (\ref{eq:limQ}), it only remains to show that
\begin{equation} \label{eq:proofproofgergodic21}
 \mathrm{vol}\left(\Omega(x)\cap B(x,\delta_{\scriptscriptstyle q})\right) \,\geq \mathrm{c}> 0
\end{equation}
where the constant $\mathrm{c}$ does not depend on $x$. Indeed, it is then clear from (\ref{eq:proofproofgergodic19}) and (\ref{eq:proofproofgergodic20}) that
$$
\liminf_{r(x) \rightarrow \infty} Q(x\hspace{0.02cm},A(x))\,> \varepsilon_{\scriptscriptstyle q}\times\mathrm{c}> 0
$$
To obtain (\ref{eq:proofproofgergodic21}), let $(r,\theta)$ be geodesic spherical coordinates, with origin at $x$. Using the integral formula (\ref{eq:integralsphericalhadamard}), after noting $\lambda(r,\theta) =\det(\mathcal{A}(r,\theta))$, it follows
\begin{equation} \label{eq:proofproofgergodic22}
\mathrm{vol}\left(\Omega(x)\cap B(x,\delta_{\scriptscriptstyle q})\right) \,=\, \int^{\tau}_{\scriptscriptstyle 0}\!\!\!\int_{\scriptscriptstyle S^{n-1}}\mathbf{1}\lbrace\Vert \mathrm{grad}\,r(x) - u(\theta)\Vert_{x} \leq \delta/2 \rbrace\hspace{0.04cm}\lambda(r,\theta)\hspace{0.03cm}dr\hspace{0.02cm}\omega_{n-1}(d\theta)
\end{equation}
where $\tau = \min\lbrace\Lambda,\delta_{\scriptscriptstyle q}\rbrace$, and the map $\theta \mapsto u(\theta)$ identifies $S^{n-1}$ with $S_xM$. Here, by (\ref{eq:volcomp}) from Theoreom \ref{th:compvol}, $\lambda(r,\theta) \geq r^{n-1}$. Therefore, (\ref{eq:proofproofgergodic22}) implies
$$
\mathrm{vol}\left(\Omega(x)\cap B(x,\delta_{\scriptscriptstyle q})\right) \geq\left(\tau^{\scriptscriptstyle n}/n\right)\times\omega_{n-1}\!\left(\lbrace\Vert \mathrm{grad}\,r(x) - u(\theta)\Vert_{x} \leq \delta/2 \rbrace\right)
$$
However, since the area measure $\omega$ is invariant by rotation, the area
$$
\omega_{n-1}\!\left(\lbrace\Vert \mathrm{grad}\,r(x) - u(\theta)\Vert_{x} \leq \delta/2 \rbrace\right) = \varsigma
$$
does not depend on $x$. Precisely, $\varsigma$ is equal to the area of a spherical cap, with angle equal to $2\hspace{0.02cm}\mathrm{acos}(1-\delta^{\scriptscriptstyle 2}/8)$.
Finally, (\ref{eq:proofproofgergodic21}) is immediately obtained, by letting $\mathrm{c} = \left(\tau^{\scriptscriptstyle d}/d\right)\times \varsigma$. \\[0.1cm]
\textbf{Proof of Proposition \ref{prop:proofgergodic2}\,:} let $V(x) = a\hspace{0.03cm}\pi^{\scriptscriptstyle -\frac{1}{2}}(x)$, as in the proposition. Recall the transition kernel $P$ is given by (\ref{eq:hmP}), which implies
$$
\rho(x) = \int_M(1-\alpha(x\hspace{0.02cm},y))\hspace{0.02cm}q(x\hspace{0.02cm},y)\hspace{0.02cm}\mathrm{vol}(dy)
$$
since the right-hand side of (\ref{eq:hmP}) should integrate to $1$ when $f(x)$ is the constant function $f(x) = 1$. But, since $\alpha(x\hspace{0.02cm},y) = \mathrm{min}\left\lbrace1\hspace{0.02cm},\pi(y)/\pi(x) \right\rbrace$, it follows that $1 - \alpha(x\hspace{0.02cm},y) = 0$ when $y \in A(x)$, the acceptance region of $x$, defined after (\ref{eq:propkernelQ}). Thus,
$$
\rho(x) = \int_{R(x)}\left[1 - \frac{\pi(y)}{\pi(x)}\right]q(x\hspace{0.02cm},y)\hspace{0.02cm}\mathrm{vol}(dy)
$$
where $R(x)$, the complement of $A(x)$, is the rejection region of $x$. With this expression of $\rho(x)$, putting $f(x) = V(x)$ in (\ref{eq:hmP}), it follows by a direct calculation that $PV(x)/V(x)$ is equal to 
\begin{equation} \label{eq:proofproofgergodic21}
\int_{A(x)}q(x\hspace{0.02cm},y)\left[\frac{\pi(x)}{\pi(y)}\right]^{\! \frac{1}{2}\hspace{0.02cm}}\mathrm{vol}(dy) + 
\int_{R(x)}q(x\hspace{0.02cm},y)\left( 1 - \left[\frac{\pi(y)}{\pi(x)}\right] + \left[\frac{\pi(y)}{\pi(x)}\right]^{\! \frac{1}{2}}\right)\mathrm{vol}(dy)
\end{equation}
Here, all the ratios are less than or equal to $1$, so that (\ref{eq:hmTr}) immediately implies (\ref{eq:ergod2}). 

In order to prove (\ref{eq:ergod1}), it is enough to prove that
\begin{equation}\label{eq:proofproofgergodic22}
 \lim_{r(x)\rightarrow\infty}\,\int_{A(x)}q(x\hspace{0.02cm},y)\left[\frac{\pi(x)}{\pi(y)}\right]^{\! \frac{1}{2}\hspace{0.02cm}}\mathrm{vol}(dy)\,=\,0\hspace{2.3cm}
\end{equation}
\begin{equation}\label{eq:proofproofgergodic23}
 \lim_{r(x)\rightarrow\infty}\,\int_{R(x)}q(x\hspace{0.02cm},y)\left(\left[\frac{\pi(y)}{\pi(x)}\right]^{\! \frac{1}{2}}- \left[\frac{\pi(y)}{\pi(x)}\right] \right)\mathrm{vol}(dy) \,=\,0
\end{equation}
Indeed, if these two limits are replaced in (\ref{eq:proofproofgergodic21}), it will follow that
$$
\limsup_{r(x)\rightarrow\infty}\, \frac{PV(x)}{V(x)} = \limsup_{r(x)\rightarrow\infty}\,Q(x,R(x))
                                                                      = \limsup_{r(x)\rightarrow\infty}\, 1 - Q(x,A(x)) < 1
$$
where the inequality is obtained using (\ref{eq:limQ}). However, this is the same as  (\ref{eq:ergod1}). Thus, to complete the proof, it is enough to prove (\ref{eq:proofproofgergodic22}) and (\ref{eq:proofproofgergodic23}). The proofs of (\ref{eq:proofproofgergodic22}) and (\ref{eq:proofproofgergodic23}) being very similar, only the proof of (\ref{eq:proofproofgergodic22}) is presented. \\[0.1cm]
\textbf{Proof of (\ref{eq:proofproofgergodic22})\,:} this is divided into three steps. First, it is proved that
\begin{equation} \label{eq:LLL1}
  \lim_{L\rightarrow \infty}\,\int_{\scriptscriptstyle A(x) - B(x,L)}q(x\hspace{0.02cm},y)\hspace{0.02cm}(\alpha(y,x))^{\scriptscriptstyle \frac{1}{2}}\hspace{0.03cm}\mathrm{vol}(dy) = 0 \hspace{0.5cm} \text{uniformly in $x$}
\end{equation}
where $\alpha(y,x) = \pi(x)/\pi(y)$. To prove (\ref{eq:LLL1}) note that $\alpha(y,x) \leq 1$ for $y \in A(x)$, and that $A(x) - B(x,L) \subset M - B(x,L)$. It follows that, for any $x \in M$,
\begin{equation} \label{eq:LLL2}
 \int_{\scriptscriptstyle A(x) - B(x,L)}q(x\hspace{0.02cm},y)\hspace{0.02cm}(\alpha(y,x))^{\scriptscriptstyle \frac{1}{2}}\hspace{0.03cm}\mathrm{vol}(dy) \,\leq\,
\int_{\scriptscriptstyle M - B(x,L)}q(x\hspace{0.02cm},y)\hspace{0.03cm}\mathrm{vol}(dy)
\end{equation}
Since $M$ is a symmetric space, there exists an isometry $g$ of $M$ such that $g\cdot x^* = x$. Since $g$ preserves Riemannian volume,
$$
\int_{\scriptscriptstyle M - B(x,L)}q(x\hspace{0.02cm},y)\hspace{0.03cm}\mathrm{vol}(dy) = 
\int_{\scriptscriptstyle M - B(x^*,L)}q(x\hspace{0.02cm},g\cdot y)\hspace{0.03cm}\mathrm{vol}(dy)
$$
But, $q(x\hspace{0.02cm},y) = q(d(x\hspace{0.02cm},y))$ depends only on the Riemannian distance $d(x\hspace{0.02cm},y)$. This implies that $q(x\hspace{0.02cm},g\cdot y) = q(x^*,y)$, since $g$ is an isometry. Thus,
$$
\int_{\scriptscriptstyle M - B(x,L)}q(x\hspace{0.02cm},y)\hspace{0.03cm}\mathrm{vol}(dy) = 
\int_{\scriptscriptstyle M - B(x^*,L)}q(x^*,y)\hspace{0.03cm}\mathrm{vol}(dy)
$$
Here, the right-hand side does not depend on $x$, and tends to zero as $L \rightarrow \infty$, as can be seen by putting $x = x^*$ in (\ref{eq:hmTr}). Now (\ref{eq:LLL1}) follows directly from (\ref{eq:LLL2}).

Second, assume that $r(x)$ is so large that the level set $C_x$ verifies (\ref{eq:contourm}) and $A(x)$ is equal to the region inside $C_{x\hspace{0.03cm}}$. It is then proved that, for any $L > 0$, 
\begin{equation} \label{eq:KKK1}
\lim_{r(x)\rightarrow\infty}\,\int_{\scriptscriptstyle A(x) \cap B(x,L) - C_x(\varepsilon)}
q(x\hspace{0.02cm},y)\hspace{0.02cm}(\alpha(y,x))^{\scriptscriptstyle \frac{1}{2}}\hspace{0.03cm}\mathrm{vol}(dy) = 0
\end{equation}
where $C_x(\varepsilon)$ is the tubular neighborhood of $C_x$ given by
$$
C_x(\varepsilon) \,=\,\left\lbrace \mathrm{Exp}_y\left(s\hspace{0.02cm}\mathrm{grad}\,r(y)\right)\,;y \in C_x\,,|s|<\varepsilon\right\rbrace
$$
Because of (\ref{eq:hmTr}), to prove (\ref{eq:KKK1}) it is enough to prove that
\begin{equation} \label{eq:KKK2}
\lim_{r(x)\rightarrow\infty}\,\alpha(y,x) = 0 \hspace{0.5cm} \text{uniformly in $y \in A(x) \cap B(x,L) - C_x(\varepsilon)$}
\end{equation}
However, this follows by Assumption (a1). Indeed, this assumption guarantees the existence of some strictly positive $\mu\hspace{0.02cm},R$ and $\pi_{\scriptscriptstyle R\hspace{0.04cm}}$, as in (\ref{eq:preVV}). Then, take $r(x) \geq R + \varepsilon$ and note that,\hfill\linebreak by (\ref{eq:preVV}), for $y$ as in (\ref{eq:KKK2}), if $r(y) \leq R$,
\begin{equation} \label{eq:KKK3}
\alpha(y,x)  \leq \frac{\pi_{\scriptscriptstyle R}\hspace{0.02cm}\exp\left(-\mu\hspace{0.02cm}r^2(x)\right)}{\pi(y)} \leq
\frac{\pi_{\scriptscriptstyle R}\hspace{0.02cm}\exp\left(-\mu\hspace{0.02cm}r^2(x)\right)}{\min_{r(y)\leq R}\pi(y)}
\end{equation}
where the right-hand side converges to zero as $r(x) \rightarrow \infty$, uniformly in $y$. On the other hand, if $r(y) > R$, let $c$ be the unit-speed geodesic connecting $x^*$ to $y$. Since $y \in A(x)$ (so $y$ lies inside $C_x$) there exists some $r \geq r(y)$ such that $c(r) \in C_{x\hspace{0.03cm}}$. Moreover, since $y \notin C_x(\varepsilon)$, it follows that $r > r(y) + \varepsilon$. Then, it is possible to show, by Assumption (a1),
$$
\alpha(y,x) = \frac{\pi(c(r))}{\pi(c(r(y)))} \leq \,\exp[-\mu\left( r^{\scriptscriptstyle 2} - r^{\scriptscriptstyle 2}(y)\right)]
$$
By a direct calculation, this implies
\begin{equation} \label{eq:KKK4}
  \alpha(y,x) \leq\exp[-\mu\left( 2\hspace{0.03cm}\varepsilon r - \varepsilon^{\scriptscriptstyle 2}\right)] \leq
                        \exp[-\mu\left( 2\hspace{0.03cm}\varepsilon r(w) - \varepsilon^{\scriptscriptstyle 2}\right)] 
\end{equation}
where $w \in C_x$ is such that $r(w)$ is the minimum of $r(w^\prime)$, taken over all $w^\prime \in C_{x\hspace{0.03cm}}$. Note that the right-hand side of (\ref{eq:KKK4}) does not depend on $y$. Moreover, $\pi(w)$ tends to zero as $r(x) \rightarrow \infty$, since $\pi(w) = \pi(x)$, and $\pi(x)$ tends to zero as $r(x) \rightarrow \infty$. Therefore, because $\pi(w)$ is positive, it follows that $r(w) \rightarrow \infty$ as $r(x)\rightarrow \infty$. But, this implies the right-hand side of (\ref{eq:KKK4}) converges to zero as $r(x) \rightarrow \infty$, uniformly in $y$. Now, (\ref{eq:KKK2}) follows from (\ref{eq:KKK4}).

The third, and final, step is to show that, for any $L>0$,
\begin{equation} \label{eq:III1}
  \lim_{\varepsilon\rightarrow 0}\limsup_{r(x)\rightarrow \infty}\,\int_{\scriptscriptstyle A(x) \cap B(x,L) \cap C_x(\varepsilon)}
q(x\hspace{0.02cm},y)\hspace{0.02cm}(\alpha(y,x))^{\scriptscriptstyle \frac{1}{2}}\hspace{0.03cm}\mathrm{vol}(dy) = 0
\end{equation}
For brevity, the proof is carried out under the assumption that $q(x\hspace{0.02cm},y)$ is bounded, uniformly in $x$ and $y$. If this assumption holds, then (\ref{eq:III1}) follows immediately by showing
\begin{equation} \label{eq:III2}
  \lim_{\varepsilon\rightarrow 0}\limsup_{r(x)\rightarrow \infty}\,\mathrm{vol}\left( B(x,L) \cap C_x(\varepsilon)\right) = 0
\end{equation}
To show (\ref{eq:III2}), let $\theta \mapsto v(\theta)$ identify the Euclidean unit sphere $S^{n-1}$ with $S_{x^*}M$, and consider the following sets
$$
\begin{array}{rl}
  T(x) =& \lbrace \theta \in S^{n-1}: \mathrm{Exp}_{x^*}(rv(\theta)) \in B(x,L)\;\;\text{for some $r\geq 0$}\rbrace \\[0.2cm]
  S(x) =& \lbrace \mathrm{Exp}_{x^*}(rv(\theta))\,; \theta \in T(x)\;\text{and}\;|r-r(x)|\leq L\rbrace
\end{array}
$$
Using the triangle inequality, it is possible to show that
\begin{equation}\label{eq:III3}
  B(x,L) \subset S(x) \subset B(x,3L)
\end{equation}
To estimate the volume in (\ref{eq:III2}), let $(r,\theta)$ be geodesic spherical coordinates, with origin at $x^*$. The first inclusion in (\ref{eq:III3}) implies $\mathrm{vol}\left(B(x,L)\cap C_x(\varepsilon)\right) \leq \mathrm{vol}\left(S(x)\cap C_x(\varepsilon)\right)$, and this yields
$$
\mathrm{vol}\left(B(x,L)\cap C_x(\varepsilon)\right) \,\leq\, \int^{\scriptscriptstyle r(x)+L}_{\scriptscriptstyle r(x)-L}\!\!\int_{\scriptscriptstyle T(x)}\mathbf{1}_{C_x(\varepsilon)}\left( \mathrm{Exp}_{x^*}(rv(\theta))\right)\hspace{0.03cm}\lambda(r,\theta)\hspace{0.02cm}dr\hspace{0.02cm}\omega_{n-1}(d\theta)
$$
in the notation of (\ref{eq:integralsphericalhadamard}), where $\lambda(r,\theta) = \det(\mathcal{A}(r,\theta))$. Bounding the last integral from above,
\begin{equation} \label{eq:III4}
\mathrm{vol}\left(B(x,L)\cap C_x(\varepsilon)\right) \,\leq\,2\varepsilon\hspace{0.02cm}\omega_{n-1}(T(x))\hspace{0.02cm}\sup_{z(r,\theta)\in B(x,3L)}\hspace{0.04cm} \lambda(r,\theta)
\end{equation}
where $z(r,\theta) =  \mathrm{Exp}_{x^*}(rv(\theta))$.  Similarly, the second inclusion in (\ref{eq:III3}) implies
$$
\mathrm{vol}\left( B(x,3L)\right)\,\geq\,\mathrm{vol}(S(x))= 
\int^{\scriptscriptstyle r(x)+L}_{\scriptscriptstyle r(x)-L}\!\!\int_{\scriptscriptstyle T(x)}\lambda(r,\theta)\hspace{0.02cm}dr\hspace{0.02cm}\omega_{n-1}(d\theta)
$$
and bounding the last integral from below gives
\begin{equation} \label{eq:III5}
\mathrm{vol}\left( B(x,3L)\right)\,\geq\,2L
\hspace{0.02cm}\omega_{n-1}(T(x))\hspace{0.02cm}\inf_{z(r,\theta)\in B(x,3L)}\hspace{0.04cm} \lambda(r,\theta)
\end{equation}
From (\ref{eq:III4}) and (\ref{eq:III5}), it follows that
\begin{equation} \label{eq:III6}
  \mathrm{vol}\left(B(x,L)\cap C_x(\varepsilon)\right) \,\leq\,
  \left(\varepsilon\middle/L\right)\mathrm{vol}\left(B(x,3L)\right)\frac{\sup_{z(r,\theta)\in B(x,3L)}\hspace{0.04cm} \lambda(r,\theta)}{\inf_{z(r,\theta)\in B(x,3L)}\hspace{0.04cm} \lambda(r,\theta)}
\end{equation}
However, by the volume growth lemma \ref{lemm:mcmc_volume}, from \ref{sec:mcmc_vollema},
$$
\limsup_{r(x) \rightarrow \infty}\hspace{0.04cm}\frac{\sup_{z(r,\theta)\in B(x,3L)}\hspace{0.04cm} \lambda(r,\theta)}{\inf_{z(r,\theta)\in B(x,3L)}\hspace{0.04cm} \lambda(r,\theta)}= \mathrm{R}\,<\,\infty
$$
Replacing into (\ref{eq:III6}), and noting that, since $M$ is a symmetric space, 
$$\mathrm{vol}\left(B(x,3L)\right) = \mathrm{vol}\left(B(x^*,3L)\right)
$$,
it follows that
$$
\limsup_{r(x) \rightarrow \infty}\hspace{0.04cm}\mathrm{vol}\left(B(x,L)\cap C_x(\varepsilon)\right) \leq
\left(\varepsilon\middle/L\right)\mathrm{vol}\left(B(x^*,3L)\right)\mathrm{R}
$$
This immediately implies (\ref{eq:III2}), and therefore (\ref{eq:III1}). \\[0.1cm]
\textbf{Conclusion\,:} finally, (\ref{eq:proofproofgergodic22}) can be obtained by combining (\ref{eq:LLL1}), (\ref{eq:KKK1}) and (\ref{eq:III1}). Precisely, the integral under the limit in (\ref{eq:proofproofgergodic22}) can be decomposed into the sum of three integrals
$$
\left(\int_{\scriptscriptstyle A(x) - B(x,L)}+
\int_{\scriptscriptstyle A(x) \cap B(x,L) - C_x(\varepsilon)}+
\int_{\scriptscriptstyle A(x) \cap B(x,L) \cap C_x(\varepsilon)\hspace{0.02cm}}
\right)q(x\hspace{0.02cm},y)\hspace{0.02cm}(\alpha(y,x))^{\scriptscriptstyle \frac{1}{2}}\hspace{0.03cm}\mathrm{vol}(dy)
$$
By (\ref{eq:LLL1}), for any $\Delta > 0$, it is possible to choose $L$ to make the first integral less than $\Delta/3$, irrespective of $x$ and $\varepsilon$. By (\ref{eq:III1}), it is possible to choose $\varepsilon$ to make the third integral less than $\Delta/3$, for all $x$ with sufficiently large $r(x)$. With $L$ and $\varepsilon$ chosen in this way, (\ref{eq:KKK1}) implies the second integral is less than $\Delta/3$, if $r(x)$ is sufficiently large. Then, the sum of the three integrals is  less than $\Delta$, and (\ref{eq:proofproofgergodic22}) follows, because $\Delta$ is arbitrary.

%%CHANGE NAME OF GEODESIC
%%dimension is n
%%subscripts

\chapter{Stochastic approximation} \label{stocha}

\minitoc
\vspace{0.1cm}

{\small
The present chapter is based on~\cite{colt}\cite{aistats}. It aims to give a general treatment, under realistic assumptions, of two problems related to stochastic approximation on Riemannian manifolds. \\[0.1cm]
The first problem is to estimate the rate of convergence of a stochastic approximation scheme, to the set of critical points (\textit{i.e.} zeros) of its mean field. 
\begin{itemize}
  \item \ref{sec:critapprox} introduces the concept of an approximate critical point of the mean field.
  \item \ref{sec:expscheme} and \ref{sec:retscheme} provide non-asymptotic upper bounds, for the number of iterations necessary, for a stochastic approximation scheme to find an approximate critical point of its mean field (specifically, exponential schemes are considered in \ref{sec:expscheme}, and retraction schemes in \ref{sec:retscheme}).
 \item  \ref{sec:example_mix} and \ref{sec:example_pca} apply the results of \ref{sec:expscheme} and \ref{sec:retscheme} to two examples\,: estimation of a mixture of Gaussian densities, and principal component analysis (PCA). 
\end{itemize}
The second problem is to derive a central limit theorem, describing the asymptotic behavior of constant-step-size exponential schemes, defined on Hadamard manifolds.
\begin{itemize}
\item \ref{sec:clt1} states the central limit theorem\,: under realistic assumptions, a constant-step-size exponential scheme defines a geometrically ergodic Markov chain. As the step-size goes to zero, a re-scaled version of this Markov chain has the same asymptotic behavior as a linear diffusion process, with a multivariate normal invariant distribution. 
\item \ref{sec:cltproofs} details the proof of this central limit theorem.
\end{itemize}
As a follow-up on the second problem, one final example is studied, as part of the present chapter.
\begin{itemize}
\item \ref{sec:ar1} introduces the Riemannian AR(1) model\,: a Markov chain $(x_t\,;t=0,1,\ldots)$ with values in a Hadamard manifold $M$, where each $x_{t+1}$ is a geodesic convex combination of the old $x_t$ and of a new input $y_{\hspace{0.02cm}t+1}$, with respective weights $1-\mu$ (for $x_t$) and $\mu$ (for $y_{\hspace{0.02cm}t+1}$), for some $\mu \in (0,1)$.\hfill\linebreak If $(y_t\,;t=1,2,\ldots)$ are independent samples from a probability distribution $P$ on $M$, then the Markov chain $(x_t)$ is geometrically ergodic, and its invariant distribution concentrates at the Riemannian barycentre of $P$, as $\mu$ goes to zero. 
%sample is a geodesic convex combination of the previous sample and of geodesic convex
%applies the central limit theorem to the study of a Riemannian auto-regressive model ...
\end{itemize}
}
\vfill
\pagebreak

%%pflug

%\section{Convergence to critical points} \label{sec:critconv}

%\begin{equation} \label{eq:sageneral}
%  x_{t+1} = \mathrm{Ret}_{x^t}\left(\mu_{\scriptscriptstyle t+1}\hspace{0.03cm}X_{y_{\scriptscriptstyle \hspace{0.02cm}t+1}}(x_t)\right)
%\end{equation}

\section{Approximate critical points} \label{sec:critapprox}
Here, the main object of study will be a stochastic approximation scheme, on a Riemannian manifold $M$. Given some initial value $x_{\scriptscriptstyle 0} \in M$, and independent observations $(y_t\,;t=1,2,\ldots)$, drawn from a probability distribution $P$ on a measurable space $Y$, this computes a sequence of iterates $(x_t\,;t = 1,2,\ldots)$, according to the update rule
\begin{equation} \label{eq:retscheme}
  x_{t+1} = \mathrm{Ret}_{x_t}\!\left(\mu_{\scriptscriptstyle t+1}\hspace{0.02cm}X_{y_{\scriptscriptstyle\hspace{0.02cm}t+1}}(x_t)\right)
\end{equation}
where $\mathrm{Ret}:TM \rightarrow M$ is a retraction, $(\mu_{\scriptscriptstyle t}\,;t=1,2,\ldots)$ is a sequence of (positive) step-sizes, and the map $X:Y\times M \rightarrow TM$ is such that $X(y,x) = X_y(x)$ always belongs to $T_xM$.

One says that $X:Y\times M \rightarrow TM$ is a random vector field. The corresponding mean vector field $X:M \rightarrow TM$ is given by
\begin{equation} \label{eq:meanfield}
X(x) = \int_{Y}\,X_y(x)\hspace{0.03cm}P(dy)
\end{equation}
which means that the noise vector field, given by $e_y(x) = X_y(x) - X(x)$, has zero expectation. In the following, it will be assumed the variance of this noise vector field is not too large, 
\begin{equation} \label{eq:variancecontrol}
   \int_{Y}\,\Vert e_y(x)\Vert^2_x\hspace{0.04cm}P(dy) \leq \sigma^2_{\scriptscriptstyle 0} + \sigma^2_{\scriptscriptstyle 1}\hspace{0.02cm}\Vert X\Vert^2_x
\end{equation}
for some constants $\sigma^2_{\scriptscriptstyle 0}\hspace{0.02cm},\sigma^2_{\scriptscriptstyle 1}$. 

The scheme (\ref{eq:retscheme}) is often used to search for zeros (critical points) of the mean vector field $X$\footnote{Zeros of vector fields are also called ``singular points", and ``stationary points". The term ``critical points" seems more in line with the context of stochastic approximation, where the mean vector field is often a gradient vector field, so the scheme (\ref{eq:retscheme}) is a stochastic gradient scheme, used to solve some optimisation problem.}.\hfill\linebreak After $t$ iterations, this scheme will have generated the iterates $(x_s\,;s = 1,\ldots,t)$. One may randomly sample these, by looking at $x_{\tau_t}$ where $\tau_t$ follows a discrete probability distribution 
\begin{equation} \label{eq:tautt}
  \mathbb{P}(\tau_t = s)  = \frac{\mu_{s+1}}{\sum^t_{s=1}\mu_{s+1}} \hspace{1cm} s = 1,\ldots, t
\end{equation}
Then, the scheme is said to have found an approximate critical point (precisely, an $\epsilon$-critical point, for some suitable accuracy $\epsilon > 0$) in expectation, if $\mathbb{E}\Vert X(x_{\tau_t})\Vert^2 \leq \epsilon$. 

For example, note that if $\mu_{\scriptscriptstyle t} = \mu$ is a constant, so (\ref{eq:retscheme}) is a constant-step-size scheme,
$$
\mathbb{E}\left[\Vert X(x_{\tau_t})\Vert^2_{x_{\tau_t}}\right] = \frac{1}{t}\sum^t_{s=1} \mathbb{E}\left[\Vert X(x_{s})\Vert^2_{x_{s}}\right]
$$
is just the average, over the first $t$ iterates, of the expected norm of the mean field. 

In order to study the stochastic approximation scheme (\ref{eq:retscheme}), it is helpful to introduce a Lyapunov function $V:M\rightarrow \mathbb{R}$. This is a positive function, which is continuously differentiable, and has $\ell$-Lipschitz gradient, in the sense of (\ref{eq:lipschitzgrad}). It is moreover assumed to satisfy
\begin{equation} \label{eq:lyapunov}
 c\hspace{0.02cm}\Vert X\Vert^2_x \leq - \langle \mathrm{grad}\,V,X\rangle_x %\hspace{0.5cm}\text{and}\hspace{0.5cm} \Vert \mathrm{grad}\,V\Vert_x \leq \bar{c}\hspace{0.02cm}\Vert X\Vert_x
\end{equation}
for some constant $c > 0$. \\[0.1cm]
\textbf{Example 1\,:} let $M = S^n \subset \mathbb{R}^{n+1}$, the unit sphere of dimension $n$. If $x^*$ is some critical point of the mean field $X$, then one may choose $V(x) = 1-\cos d(x\hspace{0.02cm},x^*)$, where $d(x\hspace{0.02cm},x^*)$ denotes the Riemannian distance between $x$ and $x^*$. In this case, $V$ is positive and has $1$-Lipschitz gradient. \\[0.1cm]
\textbf{Example 2\,:} let $M$ be a Hadamard manifold, with sectional curvatures bounded below by $\kappa_{\min} = -c^{\hspace{0.02cm}\scriptscriptstyle 2}$. If $x^*$ is some critical point of the mean field $X$, then one may choose $V(x) = V_{x^*}(x)$, for some $\delta >0$, as in (\ref{eq:hdist}). From Proposition \ref{prop:hdist} and Lemma \ref{lem:lipschitzlemma}, $V$ is positive and has $(1+\delta\hspace{0.02cm}c)$-Lipschitz gradient. \vfill\pagebreak
Lemma \ref{lem:lipschitztaylor} will be essential to all further analysis of the scheme (\ref{eq:retscheme}). The proof of this lemma, being somewhat elementary, is not given in detail. 
\begin{lemma} \label{lem:lipschitztaylor}
 If $V:M\rightarrow \mathbb{R}$ has $\ell$-Lipschitz gradient, then
\begin{equation} \label{eq:liptaylor}
  \left|V(\mathrm{Exp}_x(v)) - V(x) -  \langle\mathrm{grad}\,V,v\rangle_{x}\right|\leq   (\ell/2)\hspace{0.02cm}\Vert v\Vert^2_x
\end{equation}
for any $x \in M$ and $v \in T_xM$.
\end{lemma}
\noindent \textbf{Sketch of proof\,:} consider the geodesic $c:[0,1]\rightarrow M$, given by $c(t) = \mathrm{Exp}_x(t\hspace{0.02cm}v)$. Then, let $V(t) = V(c(t))$ and note that $V^\prime(t) = \langle\mathrm{grad}\,V,\dot{c}\rangle_{c(t)\hspace{0.02cm}}$. Let $\Pi^{\scriptscriptstyle 0}_{t}$ denote parallel transport along $c$, from $c(t)$ to $c(0)$. Since this preserves scalar products, and $\dot{c}$ is parallel,
\begin{equation} \label{eq:vprime1}
V^\prime(t) = \langle\mathrm{grad}\,V,\dot{c}\rangle_{c(0)} + \langle \Pi^{\scriptscriptstyle 0}_{t}\left(\mathrm{grad}\,V_{c(t)}\right) - \mathrm{grad}\,V_{c(0)},\dot{c}\rangle_{c(0)}
\end{equation}
Then,  using (\ref{eq:lipschitzgrad}), it may be shown that
\begin{equation} \label{eq:vprime2}
\left| \langle \Pi^{\scriptscriptstyle 0}_{t}\left(\mathrm{grad}\,V_{c(t)}\right) - \mathrm{grad}\,V_{c(0)},\dot{c}\rangle_{c(0)}\right| \leq
\ell\hspace{0.02cm}t(L(c))^2
\end{equation}
Since $c(0) = x$ and $\dot{c}(0) = v$, (\ref{eq:liptaylor}) follows by replacing (\ref{eq:vprime2}) into (\ref{eq:vprime1}), and integrating over $t$. 

\section{Exponential schemes} \label{sec:expscheme}
Consider now the case where $\mathrm{Ret} = \mathrm{Exp}$, in (\ref{eq:retscheme}). That is,
\begin{equation} \label{eq:expscheme}
  x_{t+1} = \mathrm{Exp}_{x_t}\!\left(\mu_{\scriptscriptstyle t+1}\hspace{0.02cm}X_{y_{\scriptscriptstyle\hspace{0.02cm}t+1}}(x_t)\right)
\end{equation}
For this exponential scheme, Proposition \ref{prop:nasymp_exp} provides a non-asymptotic bound on $\mathbb{E}\Vert X(x_{\tau_t})\Vert^2$, where
$\tau_t$ was defined in (\ref{eq:tautt}). This proposition uses the notation
\begin{equation} \label{eq:barmup}
%  \lbrace \mu\hspace{0.02cm}t\rbrace_{\scriptscriptstyle -1} = \frac{1}{\sum^t_{s=1}\mu^{\scriptscriptstyle \phantom{p+1}}_s}\hspace{0.3cm};\hspace{0.2cm}
\lbrace\mu^{\scriptscriptstyle p}\rbrace_{t} = \frac{\sum^t_{s=1} \mu^{\scriptscriptstyle p+1}_{s+1}}{\sum^t_{s=1}\mu^{\scriptscriptstyle \phantom{p+1}}_{s+1}} %\hspace{0.5cm} \text{for $p \geq -1$} 
\end{equation}
which is motivated by the fact that if $\mu_{\scriptscriptstyle t} = \mu$ is a constant, so (\ref{eq:expscheme}) is a constant-step-size scheme,\hfill\linebreak then $\lbrace\mu^{\scriptscriptstyle p}\rbrace_{t} = \mu^{\scriptscriptstyle p}$. In this spirit, $\lbrace\mu^{\scriptscriptstyle 1}\rbrace_{t}$ will be written $\lbrace\mu^{\scriptscriptstyle 1}\rbrace_{ t} = \lbrace\mu\rbrace_{ t}$, throughout the following.
\begin{proposition} \label{prop:nasymp_exp}
 Consider the exponential scheme (\ref{eq:expscheme}), with mean vector field (\ref{eq:meanfield}), and where the noise variance satisfies (\ref{eq:variancecontrol}). Assume that there exists a positive function $V:M\rightarrow \mathbb{R}$, with $\ell$-Lipschitz gradient, which verifies (\ref{eq:lyapunov}). If $\mu_t \leq c\hspace{0.02cm}(2\ell(1+\sigma^{\scriptscriptstyle 2}_{\scriptscriptstyle 1}))^{\scriptscriptstyle -1}$ for all $t$, then
\begin{equation} \label{eq:nasymp_exp}
   \mathbb{E}\left[\Vert X(x_{\tau_t})\Vert^2_{x_{\tau_t}}\right] \leq (2/\!\hspace{0.02cm}c)\!\,\left[\left(V(x_{\scriptscriptstyle 0})\middle/t\right)\!\hspace{0.02cm}\lbrace \mu^{\scriptscriptstyle -1}\rbrace_{t} + (\ell\hspace{0.02cm}\sigma^2_{\scriptscriptstyle 0})\hspace{0.03cm}\lbrace\mu\rbrace_{t}\hspace{0.02cm} \right]
\end{equation}
\end{proposition}
\noindent \textbf{Remark\,:} the simplest application of this proposition is to a stochastic gradient scheme, whith mean field $X(x) = -\mathrm{grad}\,f(x)$ for a cost function $f:M\rightarrow \mathbb{R}$. If $f$ is positive (or just bounded below), and has $\ell_f$-Lipschitz gradient, then $V = f$ can be introduced, as a Lyapunov function, since (\ref{eq:lyapunov}) then holds with $c = 1$. In the case of a constant-step-size scheme, with $\mu \leq (2\ell_f(1+\sigma^{\scriptscriptstyle 2}_{\scriptscriptstyle 1}))^{\scriptscriptstyle -1}$, it follows from (\ref{eq:nasymp_exp}) that
\begin{equation} \label{eq:nasymp_exp_grad}
\frac{1}{2t}\sum^t_{s=1} \mathbb{E}\left[\Vert \mathrm{grad}\,f(x_{s})\Vert^2_{x_{s}}\right] \leq \left(f(x_{\scriptscriptstyle 0})\middle/t\hspace{0.02cm}\mu\right) + (\ell_f\hspace{0.02cm}\sigma^2_{\scriptscriptstyle 0})\hspace{0.03cm}\mu
\end{equation}
In particular, if $t$ is sufficiently large, then one must have $\mathbb{E}\Vert \mathrm{grad}\,f(x_{s})\Vert^2 \leq 3(\ell_f\hspace{0.02cm}\sigma^2_{\scriptscriptstyle 0})\hspace{0.03cm}\mu$, for at least one $s$ in the range $s = 1,\ldots, t$. \\[0.1cm]
\textbf{Remark\,:} Proposition \ref{prop:nasymp_exp} provides an estimate of the rate of convergence of a stochastic approximation scheme, to the set of critical points of its mean field, which is applicable even when this set of critical points is complicated. This is especially helpful for stochastic gradient schemes, with a cost function that has many global minima (see the the example in \ref{sec:example_mix}). Proposition \ref{prop:nasymp_ret} will extend Proposition \ref{prop:nasymp_exp}, from exponential schemes, to retraction schemes. \\[0.1cm]
\textbf{Proof of Proposition \ref{prop:nasymp_exp}\,:} for each $s = 0,1,\ldots,$ it follows from Lemma \ref{lem:lipschitztaylor} that
$$
V(x_{s+1}) - V(x_s) \leq \mu_{\scriptscriptstyle s+1}\hspace{0.02cm}\langle\mathrm{grad}\,V,X_{y_{\scriptscriptstyle\hspace{0.02cm}s+1}}\rangle_{x_s} + \mu^2_{\scriptscriptstyle s+1}(\ell/2)\Vert X_{y_{\scriptscriptstyle\hspace{0.02cm}s+1}}\Vert^2_{x_s}
%\langle\mathrm{grad}\,V,\dot{c}\rangle_{c(0)} + (\ell/2)\hspace{0.02cm}(L(c))^2
$$
Then, since $X_{y_{\scriptscriptstyle\hspace{0.02cm}s+1}}(x_s) = X(x_s) + e_{y_{\scriptscriptstyle\hspace{0.02cm}s+1}}(x_s)$, 
\begin{equation} \label{eq:proofexpscheme1}
V(x_{s+1}) - V(x_s) \leq \mu_{\scriptscriptstyle s+1}\hspace{0.02cm}\langle\mathrm{grad}\,V,X_{y_{\scriptscriptstyle\hspace{0.02cm}s+1}}\rangle_{x_s} + \mu^2_{\scriptscriptstyle s+1}\ell\left(\Vert X \Vert^2_{x_s} + \Vert e_{y_{\scriptscriptstyle\hspace{0.02cm}s+1}}\Vert^2_{x_s}\right)
%\langle\mathrm{grad}\,V,\dot{c}\rangle_{c(0)} + (\ell/2)\hspace{0.02cm}(L(c))^2
\end{equation}
Let $\mathcal{Y}_s$ be the $\sigma$-algebra generated by $y_{\scriptscriptstyle 1},\ldots, y_s\hspace{0.03cm}$. Taking conditional expectations in (\ref{eq:proofexpscheme1}), it follows from (\ref{eq:meanfield}) that
$$
-\mu_{\scriptscriptstyle s+1}\hspace{0.02cm}\langle\mathrm{grad}\,V,X\rangle_{x_s} \leq \mathbb{E}\left[V(x_{s}) - V(x_{s+1})\middle|\mathcal{Y}_s\right] + \mu^2_{\scriptscriptstyle s+1}\ell\left(\Vert X \Vert^2_{x_s} + \mathbb{E}\left[\Vert e_{y_{\scriptscriptstyle\hspace{0.02cm}s+1}}\Vert^2_{x_s}\middle|\mathcal{Y}_s\right]\right)
$$
Then, from (\ref{eq:variancecontrol}),
$$
-\mu_{\scriptscriptstyle s+1}\hspace{0.02cm}\langle\mathrm{grad}\,V,X\rangle_{x_s} \leq \mathbb{E}\left[V(x_{s}) - V(x_{s+1})\middle|\mathcal{Y}_s\right] + \mu^2_{\scriptscriptstyle s+1}\ell\left(\sigma^2_{\scriptscriptstyle 0} + (1+\sigma^2_{\scriptscriptstyle 1})\Vert X \Vert^2_{x_s}\right)
$$
Therefore, using (\ref{eq:lyapunov}), and rearranging terms,
\begin{equation} \label{eq:proofexpscheme2}
(c - \ell(1+\sigma^2_{\scriptscriptstyle 1})\mu_{\scriptscriptstyle s+1})\hspace{0.03cm}\mu_{\scriptscriptstyle s+1}\hspace{0.02cm}\Vert X \Vert^2_{x_s} \leq \mathbb{E}\left[V(x_{s}) - V(x_{s+1})\middle|\mathcal{Y}_s\right] + (\ell\hspace{0.02cm}\sigma^2_{\scriptscriptstyle 0})\hspace{0.03cm}\mu^2_{\scriptscriptstyle s+1}
\end{equation}
If $\mu_{s+1} \leq c\hspace{0.02cm}(2\ell(1+\sigma^{\scriptscriptstyle 2}_{\scriptscriptstyle 1}))^{\scriptscriptstyle -1}$, this becomes
$$
(c/\!\hspace{0.04cm}2)\hspace{0.03cm}\mu_{\scriptscriptstyle s+1}\hspace{0.02cm}\Vert X \Vert^2_{x_s} \leq \mathbb{E}\left[V(x_{s}) - V(x_{s+1})\middle|\mathcal{Y}_s\right] + (\ell\hspace{0.02cm}\sigma^2_{\scriptscriptstyle 0})\hspace{0.03cm}\mu^2_{\scriptscriptstyle s+1}
$$
Finally, (\ref{eq:nasymp_exp}) follows by summing over $s = 1,\ldots, t$ and dividing by $\sum^t_{s=1}\mu_{s+1} = t/\!\hspace{0.02cm}\lbrace \mu^{\scriptscriptstyle -1}\rbrace_t\,$. 
\section{Retraction schemes} \label{sec:retscheme}
Consider now the case where $\mathrm{Ret}$ in (\ref{eq:retscheme}) is a regular retraction, in the sense of \ref{sec:retractions}. Then (\ref{eq:retscheme}) can be written under an exponential form,
\begin{equation} \label{eq:retexpscheme}
x_{t+1} = \mathrm{Exp}_{x_t}\!\left(\Phi_{x_t}\!\left(\mu_{\scriptscriptstyle t+1}\hspace{0.02cm}X_{y_{\scriptscriptstyle\hspace{0.02cm}t+1}}(x_t)\right)\right)
\end{equation}
where the maps $\Phi_x :T_xM\rightarrow T_xM$ were defined in (\ref{eq:PHI}). This new exponential form is useful, since it renders possible the application of Lemma \ref{lem:lipschitztaylor}, as in the proof of Proposition \ref{prop:nasymp_exp}.

In addition to being regular, the retraction $\mathrm{Ret}$ is assumed to verify
\begin{equation} \label{eq:retract_asump}
   \Vert\Phi_x(v)\Vert_x \leq \Vert v\Vert_x \hspace{0.3cm}\text{and}\hspace{0.2cm}
   \Vert\Phi_x(v) - v \Vert_x \leq \delta\hspace{0.02cm}\Vert v\Vert^3_x
\end{equation}
for all $x \in M$ and $v\in T_x M$, where $\delta > 0$ is a constant. This assumption holds true for the retractions studied in \ref{sec:retractions} and \ref{sec:grassret}, as may be verified, using elementary properties of the arctangent function.

It will also be assumed that the random vector field $X_y(x)$ has bounded third-order moments, 
\begin{equation} \label{eq:thirdordermoments}
  \int_{Y}\,\Vert X_y(x)\Vert^a_x\hspace{0.04cm}P(dy) \leq \tau_{\scriptscriptstyle a} \hspace{0.2cm}; a = 2,3
\end{equation}
for some constants $\tau_{\scriptscriptstyle 2}\hspace{0.02cm},\tau{\scriptscriptstyle 3} > 0$. This implies that it is possible to replace $\sigma^2_{\scriptscriptstyle 1} = 0$ in (\ref{eq:variancecontrol}). 

The following Proposition \ref{prop:nasymp_ret} is obtained by applying Lemma \ref{lem:lipschitztaylor} to the exponential form (\ref{eq:retexpscheme}) of the retraction scheme (\ref{eq:retscheme}), and taking advantage of the assumptions (\ref{eq:retract_asump}) and (\ref{eq:thirdordermoments}). 
\begin{proposition} \label{prop:nasymp_ret}
 Consider the retraction scheme (\ref{eq:retscheme}), where $\mathrm{Ret}$ is a regular retraction, which satisfies (\ref{eq:retract_asump}). Assume that (\ref{eq:thirdordermoments}) holds, so it is possible replace $\sigma^2_{\scriptscriptstyle 1} = 0$ in (\ref{eq:variancecontrol}). Assume also that there exists a positive function, with bounded and $\ell$-Lipschitz gradient, which verifies (\ref{eq:lyapunov}).\hfill\linebreak If $\mu_t \leq (c/2\ell)$ for all $t$, then
\begin{equation} \label{eq:nasymp_ret}
   \mathbb{E}\left[\Vert X(x_{\tau_t})\Vert^2_{x_{\tau_t}}\right] \leq (2/\!\hspace{0.02cm}c)\!\,\left[\left(V(x_{\scriptscriptstyle 0})\middle/t\right)\!\hspace{0.02cm}\lbrace \mu^{\scriptscriptstyle -1}\rbrace_{t} + (\ell\hspace{0.02cm}\sigma^2_{\scriptscriptstyle 0})\hspace{0.03cm}\lbrace\mu\rbrace_{t}+(\delta\tau_{\scriptscriptstyle 3}\hspace{0.02cm}\Vert V\Vert_{\scriptscriptstyle 1,\infty})\hspace{0.03cm}\lbrace\mu^{\scriptscriptstyle 2}\rbrace_{t}\hspace{0.02cm} \right]
\end{equation}
where $\Vert V\Vert_{\scriptscriptstyle 1,\infty} = \sup_{x \in M}\Vert \mathrm{grad}\,V\Vert_x\hspace{0.02cm}$.
\end{proposition}
\indent The assumptions of Proposition \ref{prop:nasymp_ret} (namely, that $X_y(x)$ has bounded third order moments, and that $\mathrm{grad}\,V(x)$ is uniformly bounded), can seem a bit too strong. In fact, these assumptions are quite natural, in several applications, where the underlying Riemannian manifold $M$ is compact. One such application, to the PCA problem, is presented in \ref{sec:example_pca}. \\[0.1cm]
\noindent \textbf{Remark\,:} the first two terms on the right-hand side of (\ref{eq:nasymp_ret}) are the same as on the right-hand side of (\ref{eq:nasymp_exp}). Thus, replacing the Riemannian exponential $\mathrm{Exp}$ by a regular retraction $\mathrm{Ret}$ has the effect of introducing a second-order term (\textit{i.e.} a constant multiple of $\lbrace\mu^{\scriptscriptstyle 2}\rbrace_{t}$) into (\ref{eq:nasymp_ret}). This additional term vanishes, in the limit where $\delta$ goes to zero. \\[0.1cm]
\textbf{Proof of Proposition \ref{prop:nasymp_ret}\,:} for $s = 0,1,\ldots,$ it follows by applying Lemma \ref{lem:lipschitztaylor} to (\ref{eq:retexpscheme}) that
\begin{equation} \label{eq:proof_nasympret1}
V(x_{s+1}) - V(x_s) \leq \left\langle \mathrm{grad}\,V,\Phi_{x_s}\!\left(\mu_{\scriptscriptstyle s+1}\hspace{0.02cm}X_{y_{\scriptscriptstyle\hspace{0.02cm}s+1}}\right)\right\rangle_{x_s} + (\ell/2)\left\Vert 
\Phi_{x_s}\!\left(\mu_{\scriptscriptstyle s+1}\hspace{0.02cm}X_{y_{\scriptscriptstyle\hspace{0.02cm}s+1}}\right)\right\Vert^2_{x_s}
\end{equation}
Here, the right-hand side may also be written
$$
\mu_{\scriptscriptstyle s+1}\hspace{0.02cm}\left\langle \mathrm{grad}\,V,X_{y_{\scriptscriptstyle\hspace{0.02cm}s+1}}\right\rangle_{x_s} + (\ell/2)\left\Vert 
\Phi_{x_s}\!\left(\mu_{\scriptscriptstyle s+1}\hspace{0.02cm}X_{y_{\scriptscriptstyle\hspace{0.02cm}s+1}}\right)\right\Vert^2_{x_s}
+
\left\langle \mathrm{grad}\,V,\Phi_{x_s}\!\left(\mu_{\scriptscriptstyle s+1}\hspace{0.02cm}X_{y_{\scriptscriptstyle\hspace{0.02cm}s+1}}\right) - \mu_{\scriptscriptstyle s+1}\hspace{0.02cm}X_{y_{\scriptscriptstyle\hspace{0.02cm}s+1}}\right\rangle_{x_s}
$$
However, by (\ref{eq:retract_asump}),
\begin{equation} \label{eq:proof_nasympret2}
\left\Vert 
\Phi_{x_s}\!\left(\mu_{\scriptscriptstyle s+1}\hspace{0.02cm}X_{y_{\scriptscriptstyle\hspace{0.02cm}s+1}}\right)\right\Vert^2_{x_s} \leq 
\mu^2_{s+1}\hspace{0.02cm}
\Vert X_{y_{\scriptscriptstyle\hspace{0.02cm}s+1}}\Vert^2_{x_s}
\end{equation}
and, in addition,
\begin{equation} \label{eq:proof_nasympret3}
\left\Vert 
\Phi_{x_s}\!\left(\mu_{\scriptscriptstyle s+1}\hspace{0.02cm}X_{y_{\scriptscriptstyle\hspace{0.02cm}s+1}}\right) - \mu_{\scriptscriptstyle s+1}\hspace{0.02cm}X_{y_{\scriptscriptstyle\hspace{0.02cm}s+1}}\right\Vert_{x_s} \leq 
\delta\hspace{0.02cm}\mu^3_{s+1}\hspace{0.02cm}
\Vert X_{y_{\scriptscriptstyle\hspace{0.02cm}s+1}}\Vert^3_{x_s}
\end{equation}
Replacing (\ref{eq:proof_nasympret2}) and (\ref{eq:proof_nasympret3}) into (\ref{eq:proof_nasympret1}), and using the Cauchy-Schwarz inequality,
$$
V(x_{s+1}) - V(x_s) \leq \mu_{\scriptscriptstyle s+1}\hspace{0.02cm}\langle \mathrm{grad}\,V,X_{y_{\scriptscriptstyle\hspace{0.02cm}s+1}}\rangle_{x_s} + \mu^2_{s+1}(\ell/2)\hspace{0.02cm}
\Vert X_{y_{\scriptscriptstyle\hspace{0.02cm}s+1}}\Vert^2_{x_s} + 
(\delta\Vert V\Vert_{\scriptscriptstyle 1,\infty})\hspace{0.02cm}\mu^3_{s+1}\hspace{0.02cm}
\Vert X_{y_{\scriptscriptstyle\hspace{0.02cm}s+1}}\Vert^3_{x_s}
$$
Now, it is possible to proceed as in the proof of Proposition \ref{prop:nasymp_exp}. Taking conditional expectations with respect to $\mathcal{Y}_s\hspace{0.03cm}$, and using (\ref{eq:meanfield}) and (\ref{eq:variancecontrol}),
$$
- \mu_{\scriptscriptstyle s+1}\hspace{0.02cm}\langle \mathrm{grad}\,V,X\rangle_{x_s} - 
\mu^2_{s+1}\ell\hspace{0.02cm}
\Vert X_{y_{\scriptscriptstyle\hspace{0.02cm}s+1}}\Vert^2_{x_s} \leq 
-\Delta V_s +  
(\ell\hspace{0.02cm}\sigma^2_{\scriptscriptstyle 0})\hspace{0.03cm}\mu^2_{s+1} + 
(\delta\Vert V\Vert_{\scriptscriptstyle 1,\infty})\hspace{0.02cm}\mu^3_{s+1}\hspace{0.02cm}
\mathbb{E}\left[\Vert X_{y_{\scriptscriptstyle\hspace{0.02cm}s+1}}\Vert^3_{x_s}\middle|\mathcal{Y}_s\right]
$$
where $\Delta V_s = \mathbb{E}\left[V(x_{s+1}) - V(x_{s})\middle|\mathcal{Y}_s\right]$.
Then, using (\ref{eq:lyapunov}), it follows that
$$
(c - \ell\mu_{\scriptscriptstyle s+1})\hspace{0.03cm}\mu_{\scriptscriptstyle s+1}\hspace{0.02cm}\Vert X \Vert^2_{x_s} \leq
-\Delta V_s +  
(\ell\hspace{0.02cm}\sigma^2_{\scriptscriptstyle 0})\hspace{0.03cm}\mu^2_{s+1} + 
(\delta\Vert V\Vert_{\scriptscriptstyle 1,\infty})\hspace{0.02cm}\mu^3_{s+1}\hspace{0.02cm}
\mathbb{E}\left[\Vert X_{y_{\scriptscriptstyle\hspace{0.02cm}s+1}}\Vert^3_{x_s}\middle|\mathcal{Y}_s\right]
$$
so, inserting (\ref{eq:thirdordermoments}), one obtains the inequality
$$
(c - \ell\mu_{\scriptscriptstyle s+1})\hspace{0.03cm}\mu_{\scriptscriptstyle s+1}\hspace{0.02cm}\Vert X \Vert^2_{x_s} \leq
-\Delta V_s +  
(\ell\hspace{0.02cm}\sigma^2_{\scriptscriptstyle 0})\hspace{0.03cm}\mu^2_{s+1} + 
(\delta\tau_{\scriptscriptstyle 3}\hspace{0.02cm}\Vert V\Vert_{\scriptscriptstyle 1,\infty})\hspace{0.03cm}\mu^3_{s+1}
$$
Here, if $\mu_t \leq (c/2\ell)$ for all $t$, then
$$
(c/\!\hspace{0.04cm}2)\hspace{0.03cm}\mu_{\scriptscriptstyle s+1}\hspace{0.02cm}\Vert X \Vert^2_{x_s} \leq
-\Delta V_s +  
(\ell\hspace{0.02cm}\sigma^2_{\scriptscriptstyle 0})\hspace{0.03cm}\mu^2_{s+1} + 
(\delta\tau_{\scriptscriptstyle 3}\hspace{0.02cm}\Vert V\Vert_{\scriptscriptstyle 1,\infty})\hspace{0.03cm}\mu^3_{s+1}
$$
Finally, (\ref{eq:nasymp_ret}) follows by summing over $s = 1,\ldots, t$ and dividing by $\sum^t_{s=1}\mu_{s+1} = t/\!\hspace{0.02cm}\lbrace \mu^{\scriptscriptstyle -1}\rbrace_t\,$. 

\vfill
\pagebreak
\section{Example\,: mixture estimation} \label{sec:example_mix}
Let $M$ be a Riemannian symmetric space, which belongs to the non-compact case, (see \ref{ssec:sspace}). Consider a probability density $m$ on $M$, which is a mixture of Gaussian densities  (of the kind defined in \ref{sec:rgd}),
\begin{equation} \label{eq:mixture}
  m(y|x) = \frac{1}{K}\sum^K_{\kappa = 1}p(y|x_{\kappa}) \hspace{0.5cm} \text{where }\, p(y|x_{\kappa}) = (Z(1))^{-1}\exp\left[ -\frac{d^{\hspace{0.03cm}2}(y,x_\kappa)}{2}\right]
\end{equation}
where $K$ is the number of mixture components, and the normalising factor $Z(1)$ is given by (\ref{eq:ssz}).\hfill\linebreak
The parameters $x = (x_\kappa\,;\kappa=1,\ldots,K)$ are to be estimated, by fitting the mixture density (\ref{eq:mixture}) to data $y_{\scriptscriptstyle 1},\ldots, y_{\scriptscriptstyle N}\hspace{0.03cm}$. Then, maximum-likelihood estimation amounts to minimising the negative log-likelihood function
\begin{equation} \label{eq:neglh}
  f(x) = -\log\hspace{0.02cm}Z(1)-\frac{1}{N}\sum^N_{n=1}\log\hspace{0.02cm} m(y_n|x)
\end{equation}
where the first term,  $-\log\hspace{0.02cm}Z(1)$, has been added to ensure that $f(x)$ is positive. The function $f$ is defined on the product Riemannian manifold, $M^{\scriptscriptstyle K} = M\times\ldots\times M$. Its gradient is then $\mathrm{grad}\,f = (\mathrm{grad}_{\kappa}\,f\,;\kappa = 1,\ldots, K)$, where $\mathrm{grad}_{\kappa}\,f$ denotes the gradient with respect to $x_{\kappa}\hspace{0.03cm}$.
\begin{lemma} \label{lem:mixgrad}
  For the negative log-likelihood function (\ref{eq:neglh}), 
\begin{equation} \label{eq:mixgrad}
  \mathrm{grad}_{\kappa}\,f(x) = -\frac{1}{N}\sum^N_{n=1}\omega_\kappa(y_n)\hspace{0.03cm}\mathrm{Exp}^{-1}_{x_\kappa}(y_n)
\end{equation}
where $\omega_\kappa(y) \propto p(y|x_{\kappa})$ are positive weights, which add up to $1$.
\end{lemma}
Let $(y_t\,;t=1,2,\ldots)$ be chosen at random among the data $y_{\scriptscriptstyle 1},\ldots, y_{\scriptscriptstyle N}\hspace{0.03cm}$. By Lemma \ref{lem:mixgrad},
\begin{equation} \label{eq:mixstochgrad}
  x^{t+1}_\kappa = \mathrm{Exp}_{x^t_\kappa}\!\left(\mu\hspace{0.03cm}X_\kappa(y_{\scriptscriptstyle\hspace{0.02cm}t+1},x^t_\kappa)\right) \hspace{0.5cm} \text{where }\,X_\kappa(y_{\scriptscriptstyle\hspace{0.02cm}t+1},x^t_\kappa) = \omega_\kappa(y_{\scriptscriptstyle\hspace{0.02cm}t+1})\hspace{0.03cm}\mathrm{Exp}^{-1}_{x^t_\kappa}(y_{\scriptscriptstyle\hspace{0.02cm}t+1})
\end{equation}
is a constant-step-size stochastic gradient scheme, for the cost function $f$. Here, the step-size $\mu$ is assumed to be less than  $1$ (in comparison to (\ref{eq:expscheme}), $t$ and $t+1$ have been written as superscripts, rather than subscripts, in order to accommodate the appearance of $\kappa$). 

Now, let $C$ be a compact and convex subset of $M$, which contains all of the data points $y_{n}\hspace{0.03cm}$, as well as all of the initial values $x^{\scriptscriptstyle 0}_\kappa$ (since $M$ is a Hadamard manifold, one may take $C$ to be any sufficiently large closed geodesic ball). The diameter of $C$ will be denoted $\mathrm{D}_{\scriptscriptstyle C}\hspace{0.03cm}$. From (\ref{eq:mixstochgrad}),
$$
x^{t+1}_\kappa = x^t_\kappa\; \#_{\scriptscriptstyle \rho_{\hspace{0.01cm}t+1}}\,y_{\scriptscriptstyle\hspace{0.02cm}t+1} \hspace{0.5cm} \text{where } \rho_{\hspace{0.02cm}t+1} = 
\mu\hspace{0.03cm}\omega_\kappa(y_{\scriptscriptstyle\hspace{0.02cm}t+1})
$$
in the notation of (\ref{eq:mapformula}), from \ref{sec:mapvsmms}. Accordingly, the iterates $x^t_\kappa$ remain within $C$, for all $t$ and $\kappa$.\hfill\linebreak
Because $C$ is compact and convex, it then becomes possible to derive the following result, by repeating, with very minor changes, the arguments leading to (\ref{eq:nasymp_exp_grad}). 
\begin{proposition} \label{prop:mixgrad}
 For the stochastic gradient scheme (\ref{eq:mixstochgrad}), let $C$ be a compact and convex subset of $M$, which contains all of the data points $y_{n}\hspace{0.03cm}$, as well as all of the initial values $x^{\scriptscriptstyle 0}_\kappa$.\hfill\linebreak If $\mu \leq (1/2\ell_{\scriptscriptstyle C})$, where $\ell_{\scriptscriptstyle C}$ denotes the supremum of the operator norm of $\mathrm{Hess}\,f(x)$, taken over $x = (x_{\kappa}\,;\kappa = 1,\ldots, K) \in C^{\scriptscriptstyle K}$, then for all $t = 1,2,\ldots$,
\begin{equation} \label{eq:mixnasymp_exp_grad}
\frac{1}{2t}\sum^t_{s=1} \mathbb{E}\left[\Vert \mathrm{grad}\,f(x^{s})\Vert^2_{x^{s}}\right] \leq \left(f_{\scriptscriptstyle C}\middle/t\hspace{0.02cm}\mu\right) + (\ell_{\scriptscriptstyle C}\hspace{0.02cm}\sigma^2_{\scriptscriptstyle 0})\hspace{0.03cm}\mu
\end{equation}
Here, $f_{\scriptscriptstyle C} = \sup_{x \in C^{\scriptscriptstyle K}} f(x)$ and $\sigma_{\scriptscriptstyle 0} = \sup_{x \in C^{\scriptscriptstyle K}}\Vert \mathrm{grad}\,f\Vert_x$  (explicit bounds on $f_{\scriptscriptstyle C\hspace{0.02cm}}$, $\sigma_{\scriptscriptstyle 0}$ and $\ell_{\scriptscriptstyle C}$ are given in the remark below).
\end{proposition}
\noindent \textbf{Remark\,:} tedious, but straightforward, calculations provide the upper bounds
\begin{equation} \label{eq:Cconstants1}
f_{\scriptscriptstyle C} \leq \frac{\mathrm{D}^2_{\scriptscriptstyle C}}{2} \hspace{0.27cm};\hspace{0.23cm} \sigma_{\scriptscriptstyle 0} \leq K\hspace{0.02cm}\mathrm{D}_{\scriptscriptstyle C}
\end{equation}
\begin{equation} \label{eq:Cconstants2}
\ell_{\scriptscriptstyle C} \leq (1+c\hspace{0.02cm}\mathrm{D}_{\scriptscriptstyle C}) + (1 + Z(1)\exp(\mathrm{D}_{\scriptscriptstyle C}/2))\mathrm{D}^2_{\scriptscriptstyle C}
\end{equation}
where $c$ is such that the sectional curvatures of $M$ lie within $[-c^{\hspace{0.02cm}\scriptscriptstyle 2},0]$. \\[0.2cm]
\textbf{Proof of Lemma \ref{lem:mixgrad}\,:} taking the gradient of (\ref{eq:neglh}), it is clear that
\begin{equation} \label{eq:proofmixgrad0}
  \mathrm{grad}_{\kappa}\,f(x) = -\frac{1}{N}\sum^N_{n=1}   \mathrm{grad}_{\kappa}\log m(y_n|x)
\end{equation}
Now, $\mathrm{grad}_{\kappa}\log m(y|x)$ can be computed as follows. If $\lambda$ is a random variable, independent from $y$, with $\mathbb{P}(\lambda = \kappa) = K^{\scriptscriptstyle -1}$, for $\kappa = 1,\ldots, K$, then $\mathbb{P}(\lambda = \kappa|y) = \omega_\kappa(y)$, with $\omega_\kappa(y)$ as in (\ref{eq:mixgrad}). Therefore, using Bayes rule,
$$
  \frac{p(\lambda\hspace{0.02cm},y)}{m(y|x)} = \sum^K_{\nu = 1}\mathbf{1}\lbrace\lambda = \nu\rbrace\hspace{0.03cm}\omega_\nu(y)
$$
where $p(\lambda\hspace{0.02cm},y)$ is the joint distribution of the couple $(\lambda\hspace{0.02cm},y)$. Taking logarithms,
\begin{equation} \label{eq:proofmixgrad1}
  \log p(\lambda\hspace{0.02cm},y) - \log m(y|x) = \sum^K_{\nu = 1}\mathbf{1}\lbrace\lambda = \nu\rbrace\hspace{0.03cm}\log \omega_\nu(y)
\end{equation}
If $\mathbb{E}_y$ denotes conditional expectation with respect to $y$,
$$
\mathbb{E}_y\left[\mathrm{grad}_{\kappa}\sum^K_{\nu = 1}\mathbf{1}\lbrace\lambda = \nu\rbrace\hspace{0.03cm}\log \omega_\nu(y)\right] = \sum^K_{\nu = 1}\omega_\nu(y)\hspace{0.03cm}\mathrm{grad}_{\kappa}\log \omega_\nu(y) = 0
$$
where the second equality follows since the conditional probabilities $\omega_\nu(y)$ always add up to $1$. By replacing this into (\ref{eq:proofmixgrad1}),
\begin{equation} \label{eq:proofmixgrad2}
\mathrm{grad}_{\kappa}\log m(y|x) = \mathbb{E}_y\left[ \mathrm{grad}_{\kappa}  \log p(\lambda\hspace{0.02cm},y)\right]
\end{equation}
But, since $\lambda$ and $y$ are independent, the joint distribution $p(\lambda\hspace{0.02cm},y)$ reads
$$
  p(\lambda\hspace{0.02cm},y) = \frac{1}{K}\sum^K_{\nu = 1}\mathbf{1}\lbrace\lambda = \nu\rbrace\hspace{0.03cm} p(y|x_{\nu})
$$
Therefore, taking logarithms, 
$$
 \log p(\lambda\hspace{0.02cm},y) = -\log(K) + \sum^K_{\nu = 1}\mathbf{1}\lbrace\lambda = \nu\rbrace\hspace{0.03cm} \log p(y|x_{\nu})
$$
This immediately yields,
\begin{equation} \label{eq:proofmixgrad4}
\mathbb{E}_y\left[ \mathrm{grad}_{\kappa}  \log p(\lambda\hspace{0.02cm},y)\right] = 
\omega_\kappa(y)\hspace{0.03cm} \mathrm{grad}_{\kappa}\log p(y|x_{\kappa}) = 
\omega_\kappa(y)\hspace{0.03cm}\mathrm{Exp}^{-1}_{x_\kappa}(y)
\end{equation}
where the second equality follows from (\ref{eq:gradfx}), and from the definition of $p(y|x_{\kappa})$ in (\ref{eq:mixture}). Finally, replacing (\ref{eq:proofmixgrad4}) into (\ref{eq:proofmixgrad2}), 
\begin{equation} \label{eq:proofmixgrad5}
\mathrm{grad}_{\kappa}\log m(y|x) = \omega_\kappa(y)\hspace{0.03cm}\mathrm{Exp}^{-1}_{x_\kappa}(y)
\end{equation}
so that (\ref{eq:mixgrad}) follows by plugging (\ref{eq:proofmixgrad5}) into (\ref{eq:proofmixgrad0}). \vfill\pagebreak

\section{Example\,: the PCA problem} \label{sec:example_pca}
Here, the notation will be the same as in \ref{sec:pca} and \ref{sec:grassret}. The aim is to apply Proposition \ref{prop:nasymp_ret}, to a constant-step-size stochastic gradient scheme, for the objective function (\ref{eq:pcaf}),
\begin{equation} \label{eq:pcafbis}
  f(x) = \mathrm{tr}\left(x\Delta\right) \hspace{1cm} x \in \mathrm{Gr}_{\scriptscriptstyle \mathbb{R}}(p\,,q)
\end{equation}
where $\Delta$ is the covariance matrix of a zero-mean random vector $y$, with values in $\mathbb{R}^d$ ($d= p+q)$. It is assumed that $y$ has finite moments of order $6$. 

The gradient of the objective function $f$ was given by (\ref{eq:pcapx}) and (\ref{eq:grasspx}). These can be written
\begin{equation} \label{eq:pcapxbis}
\mathrm{grad}\,f(x) = g\cdot \tilde{\omega}(x) \hspace{0.5cm} \text{where } \tilde{\omega}(x) = \mathrm{P}_o(g^\dagger\cdot \Delta)
\end{equation}
Now, let $b \in \mathrm{St}_{\scriptscriptstyle \mathbb{R}}(p\,,q)$ be such that $x = [b]$. By the discussion before (\ref{eq:retgrassman1}), choosing $g = (b,b^{\scriptscriptstyle \perp})$, it follows that $\mathrm{grad}\,f(x) = [X(b)]$, where $X(b) =  b^{\scriptscriptstyle \perp}\omega(b)\hspace{0.02cm}$. From (\ref{eq:grassto}) and (\ref{eq:grasspx}), it is clear that $\omega(b) = (b^{\scriptscriptstyle \perp})^\dagger\Delta\hspace{0.03cm} b$. Therefore, using the fact that $x = bb^\dagger$ (this is the definition of $[b]$),
\begin{equation} \label{eq:pcameanfield}
X(b) = (\mathrm{I}_d - x)\Delta\hspace{0.03cm} b
\end{equation}
In terms of the random vector $y$, $X(b)$ is the expectation of $X_y(b)$, where
\begin{equation} \label{eq:pcaxy}
  X_y(b) =  (\mathrm{I}_d - x)(yy^\dagger)\hspace{0.03cm} b
\end{equation}
Let $X_y(x) = [X_y(b)]$, and note that the expectation of $X_y(x)$ is equal to $\mathrm{grad}\,f(x)$ (by linearity). 
Then, consider the constant-step-size stochastic gradient scheme
\begin{equation} \label{eq:pcascheme1}
    x_{t+1} = \mathrm{Ret}_{x_t}\!\left(\mu\hspace{0.02cm}X_{y_{\scriptscriptstyle\hspace{0.02cm}t+1}}(x_t)\right)
\end{equation}
where $(y_t\,;t = 1,2,\ldots)$ are independent copies of $y$. If the retraction $\mathrm{Ret}$ is given as in (\ref{eq:retgrassman}), this becomes
\begin{equation} \label{eq:pcascheme}
  x_{t+1} = \mathrm{Span}(b_{t+1}) \hspace{0.25cm};\hspace{0.21cm} b_{t+1} = b_t +  \mu\hspace{0.03cm}(\mathrm{I}_d - b^{\phantom{\dagger}}_{\scriptscriptstyle\hspace{0.02cm}t}b^\dagger_{\scriptscriptstyle\hspace{0.02cm}t})(y^{\phantom{\dagger}}_{\scriptscriptstyle\hspace{0.02cm}t+1}y^\dagger_{\scriptscriptstyle\hspace{0.02cm}t+1})\hspace{0.03cm} b_t
\end{equation}
Proposition \ref{prop:nasymp_ret}, applied to this scheme, yields the following bound.
\begin{proposition} \label{prop:pca_nasymp}
Consider the constant-step-size scheme (\ref{eq:pcascheme1})-(\ref{eq:pcascheme}). For all $t = 1,2,\ldots,$
\begin{equation} \label{eq:pca_nasymp}
\frac{1}{2t}\sum^t_{s=1} \mathbb{E}\left[\Vert \mathrm{grad}\,f(x_{s})\Vert^2_{x_{s}}\right] \leq\left(p\Vert\Delta\Vert_{\scriptscriptstyle \mathrm{op}})\middle/t\hspace{0.02cm}\mu\right)\!\hspace{0.02cm} + (4\Vert\Delta\Vert_{\scriptscriptstyle \mathrm{op}}\hspace{0.02cm}m^{4}_y )\hspace{0.03cm}\mu+(\sqrt{8}\Vert \Delta\Vert_{\scriptscriptstyle F}\hspace{0.02cm}m^{6}_y)\hspace{0.03cm}\mu^{\scriptscriptstyle 2}
\end{equation}
where $\Vert\Delta\Vert_{\scriptscriptstyle \mathrm{op}}$ and $\Vert\Delta\Vert_{\scriptscriptstyle F}$ denote the operator norm and Frobenius norm of the matrix $\Delta$, while $m^{4}_y$ and $m^{6}_y$ denote the fourth-order and sixth-order moments of the random vector $y$.
\end{proposition}
Proposition \ref{prop:pca_nasymp} follows directly from Proposition \ref{prop:nasymp_ret}, by introducing $V(x) = p\Vert\Delta\Vert_{\scriptscriptstyle \mathrm{op}}-f(x)$, which satisfies $0 \leq V(x) \leq p\Vert\Delta\Vert_{\scriptscriptstyle \mathrm{op}\hspace{0.02cm}}$. Since $-\mathrm{grad}\,V(x) =\mathrm{grad}\,f(x)$, (\ref{eq:lyapunov}) now holds with $c = 1$. 

The function $V$ has $2\Vert\Delta\Vert_{\scriptscriptstyle \mathrm{op}}$-Lipschitz gradient, as will be shown in the remark below, and the norm of its gradient can be computed from (\ref{eq:pcapxbis}),
\begin{equation} \label{eq:pcav10}
\Vert\mathrm{grad}\,V\Vert_x = \Vert\mathrm{grad}\,f\Vert_x = 
\Vert\tilde{\omega}_x\Vert_o \leq \Vert g^\dagger\cdot \Delta \Vert_{\scriptscriptstyle F}
\end{equation}
where the inequality follows from (\ref{eq:grasssp}), since $\mathrm{P}_o$ is an orthogonal projection. But, since $g$ is orthogonal, (\ref{eq:pcav10}) implies that $\Vert\mathrm{grad}\,V\Vert_x$ is bounded by $\Vert \Delta\Vert_{\scriptscriptstyle F\hspace{0.02cm}}$, uniformly in $x$.

Thus, to obtain (\ref{eq:pca_nasymp}), it is possible to replace into (\ref{eq:nasymp_ret}), $V(x_{\scriptscriptstyle 0}) \leq p\Vert\Delta\Vert_{\scriptscriptstyle \mathrm{op}\hspace{0.02cm}}$, $\ell = 2\Vert\Delta\Vert_{\scriptscriptstyle \mathrm{op}\hspace{0.02cm}}$, and $\Vert V\Vert_{\scriptscriptstyle 1,\infty} = \Vert \Delta\Vert_{\scriptscriptstyle F\hspace{0.02cm}}$. 

For $\sigma^2_{\scriptscriptstyle 0}$ and $\tau_{\scriptscriptstyle 3\hspace{0.03cm}}$, recall that $X_y(x) = [X_y(b)]$, where $X_y(b) =  b^{\scriptscriptstyle \perp}\omega_y(b)$, with $\omega_y(b) = 
(b^{\scriptscriptstyle \perp})^\dagger(yy^\dagger)\hspace{0.03cm} b$. However, this implies $X_y(x) = g\cdot \tilde{\omega}_y(b)$, ($\tilde{\omega}_y(b)$ is obtained from $\omega_y(b)$, according to (\ref{eq:grassto})). Thus, $\Vert X_y \Vert_x = \Vert \tilde{\omega}_y(b)\Vert_o = \sqrt{2}\Vert \omega_y(b)\Vert_{\scriptscriptstyle F\hspace{0.02cm}}$. By evaluating the Frobenius norm,
$$
\Vert \omega_y(b)\Vert^2_{\scriptscriptstyle F} = \mathrm{tr}\left((\mathrm{I}_d - x)(yy^\dagger)x(yy^\dagger)\right) \leq
\Vert(\mathrm{I}_d - x)(yy^\dagger)\Vert_{\scriptscriptstyle F}\hspace{0.02cm}
\Vert (yy^\dagger)x\Vert_{\scriptscriptstyle F} \leq \Vert yy^\dagger\Vert^2_{\scriptscriptstyle F}
$$
where the first inequality follows from the Cauchy-Schwarz inequality, and the second inequality follows because $x$ and $\mathrm{I}_d - x$ are orthogonal projectors. Since $\Vert yy^\dagger\Vert_{\scriptscriptstyle F}= \Vert y \Vert^2$ (the squared Euclidean norm of $y$), this implies $\Vert X_y \Vert_x \leq \sqrt{2}\Vert y \Vert^2$. Therefore, it is possible to set $\sigma^2_{\scriptscriptstyle 0} =2m^{4}_y$ and $\tau_{\scriptscriptstyle 3} = \sqrt{8}m^{6}_y\hspace{0.03cm}$.  

Finally, the constant $\delta$ in (\ref{eq:retract_asump}) can be taken equal to $1$. Indeed, if $\mathrm{Ret}_x$ is given by (\ref{eq:retgrassman}) and $\Phi_x$ is given by (\ref{eq:phigrass}), then for $v \in T_x\mathrm{Gr}_{\scriptscriptstyle \mathbb{R}}(p\,,q)$, where $v = g\cdot \tilde{\omega}$ and $\omega$ has s.v.d. $\omega = ras^\dagger$, $\Phi_x(v) = g\cdot \tilde{\varphi}$ where $\varphi$ has s.v.d. $\varphi= r \arctan(a) s^\dagger$. Therefore,
\begin{equation} \label{eq:proofdeltaa1}
   \Vert\Phi_x(v) - v \Vert_x = \Vert g\cdot\tilde{\varphi} - g\cdot\tilde{\omega}\Vert_x = \Vert \tilde{\varphi} - \tilde{\omega}\Vert_o
%\leq \delta\hspace{0.02cm}\Vert v\Vert^3_x
\end{equation}
If $k$ is given by (\ref{eq:proofgrass3}), then
$$
 \Vert \tilde{\varphi} - \tilde{\omega}\Vert_o = 
 \Vert k\cdot \arctan(\tilde{a}) - k \cdot \tilde{a}\Vert_o = \Vert \arctan(\tilde{a}) - \tilde{a}\Vert_o
$$
By an elementary property of the $\arctan$ function, $\Vert \arctan(\tilde{a}) - \tilde{a}\Vert_o \leq \Vert \tilde{a}\Vert^3_o\hspace{0.03cm}$. Therefore,
\begin{equation} \label{eq:proofdeltaa2}
 \Vert \tilde{\varphi} - \tilde{\omega}\Vert_o \leq \Vert \tilde{a}\Vert^3_o
\end{equation}
Replacing (\ref{eq:proofdeltaa2}) into (\ref{eq:proofdeltaa1}), and noting that $\Vert \tilde{a} \Vert_o = \Vert v \Vert_x\hspace{0.03cm}$, it follows that
$\Vert\Phi_x(v) - v \Vert_x \leq \Vert v \Vert^3_x\hspace{0.03cm}$. This is the second inequality in (\ref{eq:retract_asump}), with $\delta = 1$. The first inequality in (\ref{eq:retract_asump}) is obtained by an analogous reasoning, using once more the properties of the  $\arctan$ function. \\[0.1cm]
\textbf{Remark\,:} it was claimed that the function $V$ has $2\Vert\Delta\Vert_{\scriptscriptstyle \mathrm{op}}$-Lipschitz gradient (this means $\mathrm{grad}\,V$ satisfies (\ref{eq:lipschitzgrad}), with $\ell =  2\Vert\Delta\Vert_{\scriptscriptstyle \mathrm{op}})$. To prove this claim, let 
$c(t)$ be a geodesic, with $c(0) = x$ and $\dot{c}(0) = v$. In the notation of (\ref{eq:grasslift}), $
c(t) = \exp(t\, \hat{\omega}_{\scriptscriptstyle v})\cdot x$. From~\cite{helgason} (Theorem 3.3, Chapter IV),
\begin{equation} \label{eq:grass_parallel}
  \Pi^{\scriptscriptstyle 1}_{\scriptscriptstyle 0}\left(\mathrm{grad}\,V(x)\right) = \exp(\hat{\omega}_{\scriptscriptstyle v})\cdot \mathrm{grad}\,V(x)
\end{equation}
But, $\mathrm{grad}\,V(x) = - \mathrm{grad}\,f(x)$, which is given by (\ref{eq:pcapxbis}). Therefore,
\begin{equation} \label{eq:grass_parallel1}
 \Pi^{\scriptscriptstyle 1}_{\scriptscriptstyle 0}\left(\mathrm{grad}\,V(x)\right) = (\exp(\hat{\omega}_{\scriptscriptstyle v})g)\cdot \mathrm{P}_o(g^\dagger\cdot \Delta)
\end{equation}
On the other hand, letting $y = c(1)$, one has $y = (\exp(\hat{\omega}_{\scriptscriptstyle v})g)\cdot o$. Thus, from (\ref{eq:pcapxbis}),
\begin{equation} \label{eq:grass_parallel2}
 \mathrm{grad}\,V(y) = 
(\exp(\hat{\omega}_{\scriptscriptstyle v})g)\cdot \mathrm{P}_o((\exp(\hat{\omega}_{\scriptscriptstyle v})g)^\dagger\cdot \Delta)
\end{equation}
From (\ref{eq:grass_parallel1}) and (\ref{eq:grass_parallel2}),
$$
\Vert \mathrm{grad}\,V(y)  -  \Pi^{\scriptscriptstyle 1}_{\scriptscriptstyle 0}\left(\mathrm{grad}\,V(x)\right)\Vert_y = 
\Vert \mathrm{P}_o((\exp(\hat{\omega}_{\scriptscriptstyle v})g)^\dagger\cdot \Delta) - \mathrm{P}_o(g^\dagger\cdot \Delta) \Vert_o \leq
\Vert (\exp(\hat{\omega}_{\scriptscriptstyle v})g)^\dagger\cdot \Delta - g^\dagger\cdot \Delta \Vert_{\scriptscriptstyle F}
$$
where the inequality holds since $\mathrm{P}_o$ is an orthogonal projection. Using the fact that
$$
(\exp(\hat{\omega}_{\scriptscriptstyle v})g)^\dagger\cdot \Delta - g^\dagger\cdot \Delta = \int^1_0 \left(\Delta(t)\hspace{0.02cm}\hat{\omega}_{\scriptscriptstyle v} - \hat{\omega}_{\scriptscriptstyle v}\hspace{0.02cm}\Delta(t)\right) dt
$$
where $\Delta(t) = (\exp(t\,\hat{\omega}_{\scriptscriptstyle v})g)^\dagger\cdot \Delta$, so $\Vert \Delta(t)\Vert_{\scriptscriptstyle op} = \Vert \Delta \Vert_{\scriptscriptstyle op}\hspace{0.02cm}$, it follows that
$$
\Vert \mathrm{grad}\,V(y)  -  \Pi^{\scriptscriptstyle 1}_{\scriptscriptstyle 0}\left(\mathrm{grad}\,V(x)\right)\Vert_y \leq 
\int^1_0 \Vert\Delta(t)\hspace{0.02cm}\hat{\omega}_{\scriptscriptstyle v} - \hat{\omega}_{\scriptscriptstyle v}\hspace{0.02cm}\Delta(t)\Vert_{\scriptscriptstyle F}\hspace{0.04cm} dt \leq 2\Vert \Delta\Vert_{\scriptscriptstyle op}\hspace{0.02cm}\Vert \hat{\omega}_{\scriptscriptstyle v}\Vert_{\scriptscriptstyle F}
$$
By the remark at the end of \ref{sec:pca}, the right-hand side is $2\Vert \Delta\Vert_{\scriptscriptstyle op}\hspace{0.02cm}\Vert v\Vert_{x\hspace{0.03cm}}$. In other words,
$$
\Vert \mathrm{grad}\,V(y)  -  \Pi^{\scriptscriptstyle 1}_{\scriptscriptstyle 0}\left(\mathrm{grad}\,V(x)\right)\Vert_y \leq 
2\Vert \Delta\Vert_{\scriptscriptstyle op}\hspace{0.02cm}L(c)
$$
This is equivalent to the required form of (\ref{eq:lipschitzgrad}), as can be seen by applying $\Pi^{\scriptscriptstyle 0}_{\scriptscriptstyle 1}$ under the norm.
%\left\Vert \Pi^{\scriptscriptstyle 0}_{\scriptscriptstyle 1}\left(\mathrm{grad}\,V_{c(1)}\right) - \mathrm{grad}\,V_{c(0)}\hspace{0.03cm}\right\Vert_{c(0)}\,\leq\,\ell\hspace{0.02cm}L(c)

% let $X(c) = \mathrm{grad}\,V(x)$. If 

%it is contractive
%

%\delta from the retraction

\section{A central limit theorem (CLT)} \label{sec:clt1}
Here, the aim will be to derive a central limit theorem, describing the asymptotic behavior of certain constant-step-size exponential schemes, defined on Hadamard manifolds. This is a generalisation of the central limit theorem, which holds in Euclidean space, found in~\cite{pflug}.

%Let $M$ be a Hadamard manifold, all of whose sectional curvatures lie within the interval $[-c^{\hspace{0.02cm}\scriptscriptstyle 2},0]$. Consider the constant-step-size exponential scheme,

\subsection{Geometric ergodicity}
Consider the constant-step-size exponential scheme, defined on a Hadamard manifold $M$,
\begin{equation} \label{eq:cssexp}
   x_{t+1} = \mathrm{Exp}_{x_t}\!\left(\mu\hspace{0.02cm}X_{y_{\scriptscriptstyle\hspace{0.02cm}t+1}}(x_t)\right)
\end{equation}
Since the observations $(y_t\,;t=1,2,\ldots)$ are independent and identically distributed, it follows that $(x_t\,;t=0,1,\ldots)$ is a time-homogeneous Markov chain with values in $M$. The following assumptions ensure that $(x_t)$ is geometrically ergodic, and therefore has a unique invariant distribution $\pi_\mu$.\\[0.1cm]
\noindent \textbf{e1.} the noise vector $e_y(x)$ satisfies (\ref{eq:variancecontrol}). \\[0.1cm]
\noindent \textbf{e2.} $e_y(x)$ is $P$-almost-surely a continuous function of $x$, and the distribution of $e_y(x)$ has strictly positive density, with respect to the Lebesgue measure on $T_xM$. \\[0.1cm]
\noindent \textbf{v1.} there exists a positive function $V:M\rightarrow \mathbb{R}$, with compact sublevel sets, and $\ell$-Lipschitz gradient, which satisfies (\ref{eq:lyapunov}). \\[0.1cm]
\noindent \textbf{v2.} there exist $x^* \in M$ and $\lambda > 0$, such that
\begin{equation} \label{eq:attractive}
\langle \mathrm{grad}\,V,X\rangle_x \leq -\lambda\hspace{0.02cm}V(x)  \hspace{0.5cm} \text{for } x \neq x^*
\end{equation}
\noindent \textbf{v3.} $V(x) = 0$ if and only if $x = x^*$. 
\begin{proposition} \label{prop:harris}
 Consider the constant-step-size scheme (\ref{eq:cssexp}), on a Hadamard manifold $M$.\hfill\linebreak Assume that \textbf{e1}, \textbf{e2}, \textbf{v1}, \textbf{v2} hold. If $\mu \leq c\hspace{0.02cm}(2\ell(1+\sigma^{\scriptscriptstyle 2}_{\scriptscriptstyle 1}))^{\scriptscriptstyle -1}$, then the Markov chain $(x_t)$ is geometrically ergodic, with a unique invariant distribution $\pi_\mu$.
\end{proposition} 
As the step-size $\mu$ goes to zero, the invariant distribution $\pi_\mu$ concentrates on the point $x^*$.
\begin{proposition} \label{prop:todirac}
  Under the same conditions as in Proposition \ref{prop:harris}, if \textbf{v3} holds, then $\pi_{\mu}\,\Rightarrow\,\delta_{x^*}$ as $\mu \rightarrow 0$ (here, $\Rightarrow$ denotes weak convergence of probability measures). 
\end{proposition}

%%%geometric drift !!

\subsection{A diffusion limit} \label{subsec:functional}
Consider now the re-scaled sequence $(u_t\,;t=0,1,\ldots)$, with values in $T_{x^*}M$, 
\begin{equation} \label{eq:re-scaledu}
  u_t = \psi_\mu(x_t) \hspace{0.5cm} \text{where } \psi_\mu(x) = \mu^{-\frac{1}{2}}\mathrm{Exp}^{-1}_{x^*}(x)
\end{equation}
This is the image of $(x_t\,;t=0,1,\ldots)$, under the diffeomorphism $\psi_\mu:M\rightarrow T_{x^*}M$. It is therefore a time-homogeneous Markov chain with values in $T_{x^*}M$. The trasition kernels of $(x_t)$ and $(u_t)$ will be denoted $Q_\mu$ and $\tilde{Q}_\mu\hspace{0.03cm}$, respectively. Note that
\begin{equation} \label{eq:qtoqtilde}
 \tilde{Q}_\mu\phi(u) = Q_\mu(\phi\circ \psi_\mu)(\psi^{\scriptscriptstyle -1}(u))
\end{equation}
for any measurable function $\phi:T_xM \rightarrow \mathbb{R}$.

The following assumptions ensure that, as $\mu$ goes to zero, $(u_t\,;t = 0,1,\ldots)$ behave like samples, taken at evenly spaced times $\tau_t = t\mu$, from a linear diffusion process $(U_\tau\,;\tau \geq 0)$. \\[0.1cm]
\textbf{d1.} the (2,0)-tensor field $\Sigma$, defined by
\begin{equation} \label{eq:fieldsigma}
\Sigma(x) = \int_{Y}\, e_y(x) \otimes e_y(x)\hspace{0.03cm}P(dy) \hspace{0.5cm} \text{for } x\in M
\end{equation}
is continuous on $M$. \vfill\pagebreak
\noindent \textbf{d2.} there exists a linear map $A:T_{x^*}M\rightarrow T_{x^*}M$, such that for $x \in M$,
\begin{equation} \label{eq:fieldA}
  X(x) = \Pi^{\scriptscriptstyle 1}_{\scriptscriptstyle 0}\left(A\left(\mathrm{Exp}^{-1}_{x^*}(x)\right) + R(x)\right)
\end{equation}
where $\Pi^{\scriptscriptstyle 1}_{\scriptscriptstyle 0}$ denotes parallel transport along the unique geodesic $c:[0,1] \rightarrow M$, connecting $x^*$ to $x$,\hfill\linebreak and $\Vert R(x) \Vert_{x^*} = o(d(x\hspace{0.02cm},x^*))$. \\[0.1cm]
\indent Now, let $(U_\tau\,:\tau \geq 0)$, be the linear diffusion process, with generator,
\begin{equation} \label{eq:generator}
 \mathcal{L}\phi(u) = A^i_ju^j\hspace{0.02cm}\frac{\partial\phi}{\partial u^i}(u) + \frac{1}{2} \Sigma^{ij}_*\hspace{0.02cm}\frac{\partial^2\phi}{\partial u^i\partial u^j}(u)
\end{equation}
where $(A^i_j)$ and $(\Sigma^{ij}_*)$ are the matrices which represent the linear map $A$ and the tensor $\Sigma(x^*)$, in a basis of normal coordinates centred at $x^*$. 
\begin{proposition} \label{prop:functional}
Consider the constant-step-size scheme (\ref{eq:cssexp}), on a Hadamard manifold $M$.\hfill\linebreak Let $(u_t\,:t=0,1,\ldots)$ be given by (\ref{eq:re-scaledu}), and assume that \textbf{e1}, \textbf{d1}, \textbf{d2} hold. For any compactly-supported, smooth function $\phi:T_xM \rightarrow \mathbb{R}$,
\begin{equation} \label{eq:functional}
\mu^{-1}\left[ \tilde{Q}_\mu\phi(u) - \phi(u)\right] =  \mathcal{L}\phi(u) \,+\, \varepsilon_{\mu}(u) 
\end{equation}
where $\varepsilon_{\mu}(u) \rightarrow 0$ as $\mu \rightarrow 0$, uniformly on compact subsets of $T_{x^*}M$. 
\end{proposition}
\noindent \textbf{Remark\,:} this proposition implies a \textit{functional central limit theorem}, by application of~\cite{kallenberg} (Theorem 19.28). This says that the process $(U^\mu_\tau\,;\tau \geq 0)$, equal to $u_t$ for $t\mu \leq \tau < (t+1)\mu$, converges in distribution to the linear diffusion $U$, with generator (\ref{eq:generator}), in Skorokhod space. This functional central limit theorem can be used to study the asymptotic behavior of $(x_t)$, near a general critical point, which satisfies \textbf{d2}. Sadly, I have not yet had time to develop this idea. 

\subsection{The stable case}
A central limit theorem can be derived, in the case where $x^*$ is a stable critical point of the mean field $X$, in the following sense. \\[0.1cm]
\textbf{t1.} the linear map $A$ in (\ref{eq:fieldA}) has its spectrum contained in the open left half-plane.\\[0.1cm] 
In this case, the generator $\mathcal{L}$ in (\ref{eq:generator}) admits of a unique invariant distribution, which is multivariate normal with mean zero and covariance matrix $V$, the solution of the Lyapunov equation $AV + V\!\hspace{0.01cm}A^\dagger = \Sigma_*$~\cite{duflo} ($A = (A^i_j)$ and $\Sigma_* = (\Sigma^{ij}_*)$). This will be denoted $\mathrm{N}(0,V)$.

Under the conditions of Proposition \ref{prop:harris}, the Markov chain $(x_t)$ has a unique invariant distribution $\pi_\mu$. Then, the same holds for the Markov chain $(u_t)$, which will have a unique invariant distribution $\tilde{\pi}_\mu$. This is $\tilde{\pi}_\mu(A) = \pi_\mu(\mathrm{Exp}_{x^*}(\mu^{1/2}A))$, for any measurable $A \subset T_{x^*}M$.

The following assumptions will be essential for the central limit theorem, which is stated in Proposition \ref{prop:clt1}. Assumption \textbf{t2} ensures tightness of the family $(\tilde{\pi}_\mu\,;\mu \leq c\hspace{0.02cm}(2\ell(1+\sigma^{\scriptscriptstyle 2}_{\scriptscriptstyle 1}))^{\scriptscriptstyle -1})$. \\[0.1cm]
\textbf{t2.} for each $r > 0$ there exists $v(r) > 0$ such that $V(x) \geq v(r)$ if $d(x\hspace{0.02cm},x^*) > r$. Moreover, $v(r) \rightarrow \infty$ as $r \rightarrow \infty$ and $\left. a\middle/v(a^{\scriptscriptstyle 1/2}r)\right.$ is a non-descreasing function of $a > 0$, for any $r$. \\[0.1cm]
\textbf{t3.} there exists $\alpha > 0$ such that 
\begin{equation} \label{eq:extravariancecontrol}
  \int_{Y}\,\Vert e_y(x)\Vert^{2+\alpha}_x\hspace{0.04cm}P(dy) \leq \tilde{\sigma}^2_{\scriptscriptstyle 0} + \tilde{\sigma}^2_{\scriptscriptstyle 1}\hspace{0.02cm}V(x)
\end{equation}
for some constants $\tilde{\sigma}^2_{\scriptscriptstyle 0}\hspace{0.02cm},\tilde{\sigma}^2_{\scriptscriptstyle 1}$. 
\begin{proposition} \label{prop:clt1}
 Under the conditions of Propositions \ref{prop:todirac} and \ref{prop:functional}, if \textbf{t1}, \textbf{t2}, \textbf{t3} hold, then $\tilde{\pi}_\mu \Rightarrow \mathrm{N}(0,V)$ as $\mu \rightarrow 0$. 
\end{proposition}

\section{Proof of the CLT} \label{sec:cltproofs}

\subsection{Proof of Proposition \ref{prop:harris}}
The proof relies on the following two lemmas, which will be proved below.
\begin{lemma} \label{lem:feller}
 Assume that \textbf{e2} holds. Then, the Markov chain $(x_t)$ is Feller, and $|\mathrm{vol}|$-irreducible and aperiodic (where $|\mathrm{vol}|$ denotes the Riemannian volume measure on $M$). 
\end{lemma}
\begin{lemma} \label{lem:drift}
 Assume that \textbf{e1}, \textbf{v1}, \textbf{v2} hold. If $\mu \leq c\hspace{0.02cm}(2\ell(1+\sigma^{\scriptscriptstyle 2}_{\scriptscriptstyle 1}))^{\scriptscriptstyle -1}$, then
\begin{equation} \label{eq:driftzero}
 Q_\mu V(x) \leq (1-\lambda\mu/2)V(x) + (\ell\hspace{0.02cm}\sigma^2_{\scriptscriptstyle 0})\hspace{0.03cm}\mu^2
\end{equation}
for all $x \in M$.
\end{lemma}
Admitting these lemmas, the fact that the chain $(x_t)$ is geometrically ergodic follows from~\cite{tweedie} (Theorem 16.0.1). Specifically, let $\tilde{V}(x) = V(x) + 1$. Then, by (\ref{eq:driftzero}),
$$
Q_\mu \tilde{V}(x) \leq (1-\lambda\mu/2)\tilde{V}(x) + b \hspace{0.5cm} \text{where } b = (\ell\hspace{0.02cm}\sigma^2_{\scriptscriptstyle 0})\hspace{0.03cm}\mu^2 + \lambda\mu/2
$$
Let $C = \lbrace x:V(x) \leq 4b/(\lambda\mu)\rbrace$. Clearly, 
\begin{equation} \label{eq:drift1}
Q_\mu \tilde{V}(x) \leq (1-\lambda\mu/4)\tilde{V}(x) + b\hspace{0.02cm}\mathbf{1}_{\scriptscriptstyle C}(x)
\end{equation}
By \textbf{v1}, $\tilde{V}$ has compact sublevel sets, so $C$ is a compact subset of $M$. Therefore, since $(x_t)$ is Feller, $C$ is a small set for $Q_\mu$~\cite{tweedie} (Theorem 6.0.1). Then,  (\ref{eq:drift1}) is a geometric drift condition towards the small set $C$. This is equivalent to $(x_t)$ being geometrically ergodic. \\[0.1cm]
\textbf{Proof of Lemma \ref{lem:feller}\,:} let $f:M\rightarrow \mathbb{R}$ be a bounded continuous function. By a slight abuse of notation, let $y$ denote a random variable with distribution $P$. From (\ref{eq:cssexp}),
\begin{equation} \label{eq:qmuf}
Q_\mu f(x) = \mathbb{E}\left[f\!\left(\mathrm{Exp}_x\!\left(\mu\hspace{0.02cm}X(x) + \mu\hspace{0.02cm}e_y(x)\right)\right)\right]
\end{equation}
By \textbf{e2}, $e_y(x)$ is $P$-almost-surely a continuous function of $x$. By the dominated convergence theorem, $Q_\mu f(x)$ is a bounded continuous function of $x$. In other words, the transition kernel $Q_\mu$ is a Feller kernel, so the chain $(x_t)$ is Feller. 

To show that $(x_t)$ is $|\mathrm{vol}|$-irreducible and aperiodic, it is enough to show that $Q_\mu(x,B) > 0$ whenever $\mathrm{vol}(B) > 0$, where $Q_\mu(x,B) = Q_\mu\mathrm{1}_{\scriptscriptstyle B}(x)$. By \textbf{e2}, if $w = e_y(x)$ then the distribution of $w$ is of the form $\gamma(w)\hspace{0.03cm}dw$, where $\gamma(w) > 0$, and $dw$ denotes the Lebesgue measure on $T_xM$. Therefore, from (\ref{eq:cssexp}), 
$$
Q_\mu(x,B) = \int_{T_xM} \mathbf{1}_{\scriptscriptstyle B}\!\left(\mathrm{Exp}_x\!\left(\mu\hspace{0.02cm}X(x) + \mu\hspace{0.02cm} w\right)\right)\,\gamma(w)\hspace{0.03cm}dw
$$
Since $M$ is a Hadamard manifold, $\mathrm{Exp}_x$ is a diffeomorphism of $T_xM$ onto $M$. Accordingly, 
$$
Q_\mu(x,B) = (1/\mu)^{n}\hspace{0.05cm}\int_{B} \gamma\!\left((1/\mu)\mathrm{Exp}^{-1}_x(z) - X(x)\right)\!\left| J_x(z)\right|^{-1}\mathrm{vol}(dz)
$$
where $n$ is the dimension of $M$, and $\mathrm{Exp}^*_x(\mathrm{vol})(dw) = (\left|J_x\right|\circ \mathrm{Exp}_x)(w)dw$, so that $|J_x(z)| > 0$.
Now, if $\mathrm{vol}(B) > 0$, it is clear that $Q_\mu(x,B) > 0$. \\[0.1cm]
\textbf{Proof of Lemma \ref{lem:drift}\,:} for any $x_{\scriptscriptstyle 0} \in M$, it follows from (\ref{eq:cssexp}) that
\begin{equation} \label{eq:qmuV}
Q_\mu V(x_{\scriptscriptstyle 0}) = \mathbb{E}\left[V(x_{\scriptscriptstyle 1})\right] \hspace{0.5cm} \text{where } x_{\scriptscriptstyle 1} = 
\mathrm{Exp}_{x_{\scriptscriptstyle 0}}\!\left(\mu\hspace{0.02cm}X(x_{\scriptscriptstyle 0}) + \mu\hspace{0.02cm}e_y(x_{\scriptscriptstyle 0})\right)
\end{equation}
where $y$ denotes a random variable with distribution $P$ (this is the same abuse of notation made in (\ref{eq:qmuf})). However, using Lemma \ref{lem:lipschitztaylor}, it is possible to write, as in  (\ref{eq:proofexpscheme1}), \vfill\pagebreak
$$
V(x_{\scriptscriptstyle 1})  \leq V(x_{\scriptscriptstyle 0})  + \mu\hspace{0.02cm}\langle\mathrm{grad}\,V,X_y\rangle_{x_{\scriptscriptstyle 0}} + \mu^2\ell\left(\Vert X \Vert^2_{x_{\scriptscriptstyle 0}} + \Vert e_{y}\Vert^2_{x_{\scriptscriptstyle 0}}\right)
$$
By \textbf{e1} and (\ref{eq:qmuV}), it follows after taking expectations,
\begin{equation} \label{eq:prooflemdrifto1}
Q_\mu V(x_{\scriptscriptstyle 0}) \leq V(x_{\scriptscriptstyle 0}) + \mu\hspace{0.02cm}\langle\mathrm{grad}\,V,X\rangle_{x_{\scriptscriptstyle 0}} +
\mu^2\ell(1+\sigma^2_{\scriptscriptstyle 1})\Vert X \Vert^2_{x_{\scriptscriptstyle 0}}
 + (\ell\hspace{0.02cm}\sigma^2_{\scriptscriptstyle 0})\hspace{0.03cm}\mu^2
\end{equation} 
By \textbf{v1}, $V$ satisfies (\ref{eq:lyapunov}), so that (\ref{eq:prooflemdrifto1}) implies
\begin{equation} \label{eq:prooflemdrifto2}
Q_\mu V(x_{\scriptscriptstyle 0}) \leq V(x_{\scriptscriptstyle 0}) + \mu\hspace{0.02cm}(1 - \mu\hspace{0.02cm}
\ell(1+\sigma^2_{\scriptscriptstyle 1})/c)
\hspace{0.02cm}\langle\mathrm{grad}\,V,X\rangle_{x_{\scriptscriptstyle 0}} +
 + (\ell\hspace{0.02cm}\sigma^2_{\scriptscriptstyle 0})\hspace{0.03cm}\mu^2
\end{equation} 
Since $\langle\mathrm{grad}\,V,X\rangle$ is negative, if $\mu \leq c\hspace{0.02cm}(2\ell(1+\sigma^{\scriptscriptstyle 2}_{\scriptscriptstyle 1}))^{\scriptscriptstyle -1}$, then (\ref{eq:prooflemdrifto2}) implies
$$
Q_\mu V(x_{\scriptscriptstyle 0}) \leq V(x_{\scriptscriptstyle 0}) + (\mu/2)
\hspace{0.02cm}\langle\mathrm{grad}\,V,X\rangle_{x_{\scriptscriptstyle 0}} 
 + (\ell\hspace{0.02cm}\sigma^2_{\scriptscriptstyle 0})\hspace{0.03cm}\mu^2
$$
Finally, by \textbf{v2}, this yields
\begin{equation} \label{eq:prooflemdrift3}
Q_\mu V(x_{\scriptscriptstyle 0}) \leq (1-\lambda\mu/2)V(x_{\scriptscriptstyle 0}) + (\ell\hspace{0.02cm}\sigma^2_{\scriptscriptstyle 0})\hspace{0.03cm}\mu^2
\end{equation}
which is the same as (\ref{eq:driftzero}), since $x_{\scriptscriptstyle 0}$ is arbitrary.

\subsection{Proof of Proposition \ref{prop:todirac}}
Proposition \ref{prop:gergodic} implies that the chain $(x_t)$ has a unique invariant distribution, here denoted $\pi_\mu\hspace{0.02cm}$, for any $\mu \leq c\hspace{0.02cm}(2\ell(1+\sigma^{\scriptscriptstyle 2}_{\scriptscriptstyle 1}))^{\scriptscriptstyle -1}$. Integrating both sides of (\ref{eq:driftzero}) with respect to $\pi_\mu\hspace{0.02cm}$, it follows that
$$
\int_M Q_\mu V(x)\hspace{0.02cm}\pi_\mu(dx) \leq (1-\lambda\mu/2)\int_M V(x)\hspace{0.02cm}\pi_\mu(dx) + (\ell\hspace{0.02cm}\sigma^2_{\scriptscriptstyle 0})\hspace{0.03cm}\mu^2
$$
Since $\pi_\mu$ is an invariant distribution of the transition kernel $Q_\mu\hspace{0.02cm}$, this means
$$
\int_M V(x)\hspace{0.02cm}\pi_\mu(dx) \leq (1-\lambda\mu/2)\int_M V(x)\hspace{0.02cm}\pi_\mu(dx) + (\ell\hspace{0.02cm}\sigma^2_{\scriptscriptstyle 0})\hspace{0.03cm}\mu^2
$$
In other words,
\begin{equation} \label{eq:Vmoment}
 \int_M V(x)\hspace{0.02cm}\pi_\mu(dx) \leq 2(\ell\hspace{0.02cm}\sigma^2_{\scriptscriptstyle 0}/\lambda)\hspace{0.02cm}\mu
\end{equation}
so, by Markov's inequality, 
\begin{equation} \label{eq:Vmarkov}
  \pi_{\mu}(V > v) \leq 2(\ell\hspace{0.02cm}\sigma^2_{\scriptscriptstyle 0}/\lambda)\hspace{0.02cm}\frac{\mu}{v} \hspace{0.5cm} \text{for all $v > 0$}
\end{equation}
By \textbf{v1}, $V$ has compact sublevel sets, so (\ref{eq:Vmarkov}) implies the family $(\pi_\mu\,;\mu \leq c\hspace{0.02cm}(2\ell(1+\sigma^{\scriptscriptstyle 2}_{\scriptscriptstyle 1}))^{\scriptscriptstyle -1})$ is tight. If $\pi_*$ is a limit point of this family at $\mu = 0$, then by the Portmaneau theorem, $\pi_*(V > v) = 0$ for all $v > 0$. In other words, $\pi_*(V = 0) = 1$. By \textbf{v3}, this is equivalent to $\pi_*(\lbrace x^*\rbrace) = 1$, or to $\pi_* = \delta_{x^*}\hspace{0.02cm}$. 

\subsection{Proof of Proposition \ref{prop:functional}}
The proof exploits the relation (\ref{eq:qtoqtilde}), between the transition kernels $Q_\mu$ and $\tilde{Q}_\mu\hspace{0.02cm}$, using the following Lemmas \ref{lem:qtaylor}, \ref{lem:normalderivatives}, and \ref{lem:xnormals}. 

In Lemma \ref{lem:qtaylor}, $[H\!:\!T]$ will denote the contraction of a (0,2)-tensor $H$ with a (2,0)-tensor $T$.\hfill\linebreak This is  $[H\!:\!T] = H_{ij}T^{ij}$, in any local coordinates. Moreover, if $f :M\rightarrow \mathbb{R}$ is compactly-supported and smooth, $A_f\hspace{0.02cm},B_f$ denote positive numbers such that
\begin{equation} \label{eq:AfBf}
  |\mathrm{Hess}\,f_x(w,w)| \leq A_f\Vert w \Vert^2_x \hspace{0.25cm}\text{and}\hspace{0.2cm}
 |\nabla\mathrm{Hess}\,f_x(w,w,w)| \leq B_f\Vert w \Vert^3_x
\end{equation}
for any $x \in M$ and $w \in T_xM$, where $\nabla\mathrm{Hess}\,f$ is the covariant derivative of the Hessian of $f$ (respectively, $\mathrm{Hess}\,f$ and $\nabla\mathrm{Hess}\,f$ are (0,2)- and (0,3)-tensor fields).

\begin{lemma} \label{lem:qtaylor}
  For any compactly-supported, smooth $f :M\rightarrow \mathbb{R}$,
\begin{equation} \label{eq:qtaylor}
  Q_\mu f(x) = f(x) + \mu\hspace{0.02cm}Xf(x) + \frac{\mu^2}{2}[\mathrm{Hess}\,f\!:\!\Sigma + X \otimes X]_x + \mu^2\hspace{0.03cm}\mathcal{R}_x(f,\mu)
\end{equation}
where the remainder term $\mathcal{R}_x(f,\mu)$ satisfies
$$
|\mathcal{R}_x(f,\mu)| 
 \leq 2A_f\hspace{0.03cm}\mathbb{E}\left[\mathbf{1}\lbrace \Vert e_y\Vert_x> K\rbrace\Vert e_y\Vert^2_x\right] + \hspace{8cm}
$$
\begin{equation} \label{eq:qremainder}
\phantom{abcd}6B_f\hspace{0.03cm}
\mu\hspace{0.02cm}(2\Vert X \Vert^2_x + 2K^2)\hspace{0.03cm}\mathbb{E}\left[\mathbf{1}\lbrace \Vert e_y\Vert_x> K\rbrace\Vert e_y\Vert_x\right]+2B_f\hspace{0.03cm}\mu\hspace{0.02cm}(4\Vert X \Vert^3_x + 4K^3) 
\end{equation}
for any (arbitrarily chosen) $K > 0$.
\end{lemma}
Given normal coordinates $(x^i\,;i=1,\ldots, n)$ on $M$, with origin at $x^*$, recall the coordinate vector fields $\partial_i = \left.\partial\middle/\partial x^i\right.$. Any function $\phi:T_{x^*}M\rightarrow \mathbb{R}$ may be identified with a function of\hfill\linebreak $n$ variables, $\phi(u) = \phi(u^{\scriptscriptstyle 1},\ldots,u^n)$, for $u \in T_{x^*}M$ where $u = u^i\partial_i(x^*)$. 
\begin{lemma} \label{lem:normalderivatives}
Let $(x^i\,;i=1,\ldots, n)$ be normal coordinates on $M$ with origin at $x^*$. For any smooth function $\phi:T_{x^*}M\rightarrow \mathbb{R}$, if $\psi_\mu$ is given by (\ref{eq:re-scaledu}), then 
\begin{equation} \label{eq:normalderivatives}
  \partial_i(\phi\circ \psi_\mu)(\psi^{\scriptscriptstyle -1}(u)) = \mu^{-\frac{1}{2}}\frac{\partial\phi}{\partial u^i}(u) \hspace{0.36cm};\hspace{0.3cm}
  \partial_{ij}(\phi\circ \psi_\mu)(\psi^{\scriptscriptstyle -1}(u)) = \mu^{-1}\frac{\partial^2\phi}{\partial u^i\partial u^j}(u)
\end{equation}
\end{lemma}
\begin{lemma} \label{lem:xnormals}
 Let $X^i(x)$ denote the components of the mean field $X$, with respect to the normal coordinates $(x^i\,;i=1,\ldots,n)$. If \textbf{d2} holds, then 
\begin{equation} \label{eq:xnormals}
  X^i(\psi^{\scriptscriptstyle -1}(u)) = \mu^{\frac{1}{2}}\hspace{0.02cm}A^i_ju^j + R^i(\mu^{\frac{1}{2}}u)
\end{equation}
where $|R^i(u)| = o(\Vert u \Vert_{x^*})$. 
\end{lemma}
Lemmas \ref{lem:qtaylor}, \ref{lem:normalderivatives}, and \ref{lem:xnormals} will be proved below. Accepting them to be true, recall (\ref{eq:qtoqtilde})
$$
 \tilde{Q}_\mu\phi(u) = Q_\mu(\phi\circ \psi_\mu)(\psi^{\scriptscriptstyle -1}(u))
$$
Replacing (\ref{eq:qtaylor}) into the right-hand side gives
$$
\mu^{-1}\left[ \tilde{Q}_\mu\phi(u) - \phi(u)\right] = \hspace{7.6cm}
$$
\begin{equation} \label{eq:proofunctional}
X(\phi\circ \psi_\mu)(\psi^{\scriptscriptstyle -1}(u)) + \frac{\mu}{2}\hspace{0.02cm}[\mathrm{Hess}\,(\phi\circ \psi_\mu)\!:\!T]_{\psi^{-1}(u)} + \mu\hspace{0.02cm}\mathcal{R}_{\psi^{-1}(u)}(\phi\circ \psi_\mu\hspace{0.02cm},\mu)
\end{equation}
where $T = \Sigma + X \otimes X$. However, working in normal coordinates,
$$
X(\phi\circ \psi_\mu)(\psi^{\scriptscriptstyle -1}(u)) = 
  X^i(\psi^{\scriptscriptstyle -1}(u)) \partial_i(\phi\circ \psi_\mu)(\psi^{\scriptscriptstyle -1}(u)) 
$$
so that, by (\ref{eq:normalderivatives}) and (\ref{eq:xnormals}),
$$
X(\phi\circ \psi_\mu)(\psi^{\scriptscriptstyle -1}(u)) = \left\lbrace A^i_ju^j 
 + 
\mu^{-\frac{1}{2}} R^i(\mu^{\frac{1}{2}}u)
\right\rbrace
\frac{\partial\phi}{\partial u^i}(u)
$$
Since $\phi$ is compactly-supported, this can be written
\begin{equation} \label{eq:proofunctional1}
X(\phi\circ \psi_\mu)(\psi^{\scriptscriptstyle -1}(u)) = A^i_ju^j\hspace{0.02cm}\frac{\partial\phi}{\partial u^i}(u) + \varepsilon^{\scriptscriptstyle 1}_\mu(u)
\end{equation}
where $\varepsilon^{\scriptscriptstyle 1}_\mu(u) \rightarrow 0$, uniformly on $T_{x^*}M$, as  $\mu \rightarrow 0$. For the second term in (\ref{eq:proofunctional}), using (\ref{eq:hesscoordinates}),
$$
[\mathrm{Hess}\,(\phi\circ \psi_\mu)\!:\!T]_{\psi^{-1}(u)} = T^{ij}(\psi^{-1}(u))\left[   \partial_{ij}(\phi\circ \psi_\mu)(\psi^{\scriptscriptstyle -1}(u)) - \Gamma^k_{ij}(\psi^{-1}(u))\hspace{0.02cm}\partial_k(\phi\circ \psi_\mu)(\psi^{\scriptscriptstyle -1}(u))\right]
$$
so that, by (\ref{eq:normalderivatives}),
$$
[\mathrm{Hess}\,(\phi\circ \psi_\mu)\!:\!T]_{\psi^{-1}(u)} = \mu^{-1}\hspace{0.03cm}T^{ij}(\psi^{-1}(u))\left[\frac{\partial^2\phi}{\partial u^i\partial u^j}(u) - \mu^{\frac{1}{2}}\Gamma^k_{ij}(\psi^{-1}(u))\frac{\partial\phi}{\partial u^k}(u)\right]
$$
where $(\Gamma^i_{jk})$ denote the Christoffel symbols. Since $\phi$ is compactly-supported, this can be written
$$
[\mathrm{Hess}\,(\phi\circ \psi_\mu)\!:\!T]_{\psi^{-1}(u)} = \mu^{-1}\hspace{0.03cm}T^{ij}_*\hspace{0.02cm}\frac{\partial^2\phi}{\partial u^i\partial u^j}(u) + \varepsilon^{\scriptscriptstyle 2}_\mu(u)
%- \mu^{\frac{1}{2}}\Gamma^k_{ij}(\psi^{-1}(u))\frac{\partial\phi}{\partial u^i}(u)\right]
$$
where $(T^{ij}_*)$ is the matrix which represents the tensor $T(x^*)$ in normal coordinates, and where $\varepsilon^{\scriptscriptstyle 2}_\mu(u) \rightarrow 0$, uniformly on $T_{x^*}M$, as  $\mu \rightarrow 0$. Since (clearly, from (\ref{eq:fieldA})), $T(x^*) = \Sigma(x^*)$, it follows
\begin{equation} \label{eq:proofunctional2}
[\mathrm{Hess}\,(\phi\circ \psi_\mu)\!:\!T]_{\psi^{-1}(u)} = \mu^{-1}\hspace{0.03cm}\Sigma^{ij}_*\hspace{0.02cm}\frac{\partial^2\phi}{\partial u^i\partial u^j}(u) + \varepsilon^{\scriptscriptstyle 2}_\mu(u)
\end{equation}
Then, replacing (\ref{eq:proofunctional1}) and (\ref{eq:proofunctional2}) into (\ref{eq:proofunctional}), and recalling the definition of $\mathcal{L}$ from (\ref{eq:generator}),
\begin{equation} \label{eq:proofunctional11}
\mu^{-1}\left[ \tilde{Q}_\mu\phi(u) - \phi(u)\right] =
\mathcal{L}\phi(u) \,+\, \varepsilon^{\scriptscriptstyle 1}_\mu(u) + \varepsilon^{\scriptscriptstyle 2}_\mu(u) + 
\mu\hspace{0.02cm}\mathcal{R}_{\psi^{-1}(u)}(\phi\circ \psi_\mu\hspace{0.02cm},\mu)
\end{equation}
To conclude, let $\varepsilon_\mu(u) = \varepsilon^{\scriptscriptstyle 1}_\mu(u) + \varepsilon^{\scriptscriptstyle 2}_\mu(u) + 
\mu\hspace{0.02cm}\mathcal{R}_{\psi^{-1}(u)}(\phi\circ \psi_\mu\hspace{0.02cm},\mu)$, and recall that $\varepsilon^{\scriptscriptstyle 1}_\mu(u)$ and
$\varepsilon^{\scriptscriptstyle 2}_\mu(u)$ converge to zero, uniformly on $T_{x^*}M$. Moreover, using (\ref{eq:variancecontrol}) and (\ref{eq:qremainder}), it is straightforward that $\mathcal{R}_{\psi^{-1}(u)}(\phi\circ \psi_\mu\hspace{0.02cm},\mu)$ is bounded on compact subsets of $T_{x^*}M$, (by an upper bound which is independent of $\mu$). Therefore, $\varepsilon_{\mu}(u) \rightarrow 0$ as $\mu \rightarrow 0$, uniformly on compact subsets of $T_{x^*}M$. \\[0.1cm]
\textbf{Proof of Lemma \ref{lem:qtaylor}\,:} the proof will rely on the following variant of Taylor expansion\hfill\linebreak (compare to~\cite{pflug}, Section 2). Let $f:M \rightarrow \mathbb{R}$ be a compactly-supported, smooth function. If $A_f\hspace{0.02cm},B_f$ are given by (\ref{eq:AfBf}), 
$x \in M$ and $\xi\hspace{0.02cm},\eta \in T_x M$, then
$$
f(\mathrm{Exp}_x(\xi + \eta)) = f(x) + (\xi + \eta)f + \frac{1}{2}\hspace{0.02cm}[\mathrm{Hess}\,f\!:\!(\xi + \eta)\otimes(\xi+\eta)] + \mathcal{R}_f(x) \hspace{0.62cm}
$$
\begin{equation} \label{eq:pflug}
\text{where }\left|\mathcal{R}_f(x)\right| \leq  2A_f\hspace{0.02cm}\Vert \eta\Vert^2_x + 6B_f\hspace{0.02cm}\Vert \xi\Vert^2_x\Vert \eta\Vert^{\phantom{2}}_x +  2B_f\hspace{0.02cm}\Vert \xi\Vert^3_x \hspace{2.965cm}
\end{equation}
To apply (\ref{eq:pflug}), recall (\ref{eq:qmuf}),
$$
Q_\mu f(x) = \mathbb{E}\left[f\!\left(\mathrm{Exp}_x\!\left(\mu\hspace{0.02cm}X(x) + \mu\hspace{0.02cm}e_y(x)\right)\right)\right]
$$
and let $\xi = \mu\hspace{0.02cm}X(x) + \mu\hspace{0.02cm}\mathbf{1}\lbrace \Vert e_y\Vert_x \leq K\rbrace\hspace{0.02cm}e_y(x)$, $\eta = \mu\hspace{0.02cm}\mathbf{1}\lbrace \Vert e_y\Vert_x > K\rbrace\hspace{0.02cm}e_y(x)$. Taking the expectation of the Taylor expansion in (\ref{eq:pflug}) and using (\ref{eq:meanfield}) and (\ref{eq:fieldsigma}), it follows that, as in (\ref{eq:qtaylor}),
$$
Q_\mu f(x) = f(x) + \mu\hspace{0.02cm}Xf(x) + \frac{\mu^2}{2}[\mathrm{Hess}\,f\!:\!\Sigma + X \otimes X]_x + \mu^2\hspace{0.03cm}\mathcal{R}_x(f,\mu)
$$
where $|\mathcal{R}_x(f,\mu)|$ is less than or equal to 
%$$
% \leq 2A_f\hspace{0.03cm}\mathbb{E}\left[\mathbf{1}\lbrace \Vert e_y\Vert_x> K\rbrace\Vert e_y\Vert^2_x\right] + \hspace{8cm}
%$$
%\begin{equation} \label{eq:qremainder}
%\phantom{abcd}6B_f\hspace{0.03cm}
%\mu\hspace{0.02cm}(2\Vert X \Vert^2_x + 2K^2)\hspace{0.03cm}\mathbb{E}\left[\mathbf{1}\lbrace \Vert e_y\Vert_x> K\rbrace\Vert e_y\Vert_x\right]+2B_f\hspace{0.03cm}\mu\hspace{0.02cm}(4\Vert X \Vert^3_x + 4K^3) 
$$
\begin{array}{l}
2A_f\hspace{0.03cm}\mathbb{E}\left[\mathbf{1}\lbrace \Vert e_y\Vert_x> K\rbrace\Vert e_y\Vert^2_x\right] + \\[0.15cm]
6B_f\hspace{0.03cm}
\mu\hspace{0.02cm}\mathbb{E}\left[\Vert \hspace{0.02cm}X(x) + \mathbf{1}\lbrace \Vert e_y\Vert_x \leq K\rbrace\hspace{0.02cm}e_y(x)\Vert^2_x \times \Vert \mathbf{1}\lbrace \Vert e_y\Vert_x > K\rbrace\hspace{0.02cm}e_y(x)\Vert_x\right] + \\[0.15cm]
2B_f\hspace{0.03cm}\mu\mathbb{E}\left[\Vert \hspace{0.02cm}X(x) + \mathbf{1}\lbrace \Vert e_y\Vert_x \leq K\rbrace\hspace{0.02cm}e_y(x)\Vert^3_x \right]
%\Vert \xi\Vert^2_x\Vert \eta\Vert^{\phantom{2}}_x + \Vert \eta\Vert^2_x]
\end{array}
$$
Then, to obtain (\ref{eq:qremainder}), it is enough to note 
$$
\begin{array}{l}
\Vert \hspace{0.02cm}X(x) + \mathbf{1}\lbrace \Vert e_y\Vert_x \leq K\rbrace\hspace{0.02cm}e_y(x)\Vert^3_x \leq 4\Vert X \Vert^3_x + 4K^3 \\[0.2cm]
\Vert \hspace{0.02cm}X(x) + \mathbf{1}\lbrace \Vert e_y\Vert_x \leq K\rbrace\hspace{0.02cm}e_y(x)\Vert^2_x \leq 2\Vert X \Vert^2_x + 2K^2
\end{array}
$$
which follow from the elementary inequalities $(a+b)^3 \leq 4a^3 + 4b^3$ and $(a+b)^2 \leq 2a^2 + 2b^2$. \vfill\pagebreak
\noindent \textbf{Proof of (\ref{eq:pflug})\,:} if $f:M\rightarrow \mathbb{R}$ is smooth and compactly-supported, for $x \in M$ and $\zeta \in T_xM$, 
one has from the second- and third-order Taylor expansions of $f$ at $x$, that 
$$
f(\mathrm{Exp}_x(\zeta)) = f(x) + \zeta f + \frac{1}{2}\hspace{0.02cm}[\mathrm{Hess}\,f\!:\!\zeta\otimes \zeta] + \mathcal{R}_f(x)
$$
where, simultaneously, $|\mathcal{R}_f(x)| \leq A_f\Vert \zeta\Vert^2_x$ and $|\mathcal{R}_f(x)| \leq B_f\Vert \zeta\Vert^3_x\hspace{0.2cm}$. If $\zeta = \xi + \eta$, then
$$
\begin{array}{rl}
|\mathcal{R}_f(x)| \leq 2A_f\Vert \eta\Vert^2_x & \text{if $\Vert \eta\Vert_x \geq \Vert \xi\Vert_x$}\\[0.12cm]
|\mathcal{R}_f(x)| \leq 2B_f\Vert \xi\Vert^3_x  + 6B_f\Vert \xi\Vert^2_x\Vert \eta\Vert_x& \text{if $\Vert \eta\Vert_x < \Vert \xi\Vert_x$}
\end{array}
$$
and (\ref{eq:pflug}) is obtained by adding up these two cases. \\[0.2cm] 
\noindent \textbf{Proof of Lemma \ref{lem:normalderivatives}\,:} let $f:M\rightarrow \mathbb{R}$ be a smooth function. From the definition of coordinate vector fields~\cite{lee} (Page 49),
\begin{equation} \label{eq:coofields}
\partial_if(x) = (f\circ \mathrm{Exp}_{x^*})^\prime(\mathrm{Exp}^{-1}_{x^*}(x))(\partial_i(x^*))
\end{equation}
where the prime denotes the Fr\'echet derivative. To obtain (\ref{eq:normalderivatives}), set $f = \phi\circ \psi_\mu$ and $x = \psi^{\scriptscriptstyle -1}(u)$, so that $f\circ \mathrm{Exp}_{x^*}(w) = \phi(\mu^{-\frac{1}{2}}w)$ (for $w \in T_{x^*}M$) and $\mathrm{Exp}^{-1}_{x^*}(x) = \mu^{\frac{1}{2}}u$. Then, in particular, $(f\circ \mathrm{Exp}_{x^*})^\prime = \mu^{-\frac{1}{2}}\phi^\prime$. Replacing into (\ref{eq:coofields}), it follows that
$$
\partial_i(\phi\circ \psi_\mu)(\psi^{\scriptscriptstyle -1}(u)) = \mu^{-\frac{1}{2}}\phi^\prime(u)(\partial_i(x^*))
$$
Now, if $\phi$ is identified with a function of $n$ variables, $\phi(u) = \phi(u^{\scriptscriptstyle 1},\ldots,u^n)$ where $u = u^i\partial_i(x^*)$, then $\phi^\prime(u)(\partial_i(x^*)) = \partial\phi(u)/\partial u^i$. This yields the first identity in (\ref{eq:normalderivatives}). The second identity follows from the first by repeated application. \\[0.1cm]
\noindent \textbf{Proof of Lemma \ref{lem:xnormals}\,:} assume that \textbf{d2} holds. Using the same notation as in (\ref{eq:fieldA}), consider the Taylor expansion of the coordinate vector fields $\partial_i$ (see \cite{chavel}, Page 90), 
$$
\partial_i(x) = \Pi^{\scriptscriptstyle 1}_{\scriptscriptstyle 0}\left(\partial_i(x^*) + \nabla\partial_i(x^*)\left(\mathrm{Exp}^{-1}_{x^*}(x)\right) + o(d(x\hspace{0.02cm},x^*))\right)
$$
where $\nabla \partial_i(x^*):T_{x^*}M\rightarrow T_{x^*}M$ is the covariant derivative of $\partial_i$ at $x^*$. From (\ref{eq:christoffel}) and (\ref{eq:normal0}), it is clear that $\nabla \partial_i(x^*) = 0$, and therefore
\begin{equation} \label{eq:coordinateA}
\partial_i(x) = \Pi^{\scriptscriptstyle 1}_{\scriptscriptstyle 0}\left(\partial_i(x^*) + o(d(x\hspace{0.02cm},x^*))\right)
\end{equation}
Take the scalar product of (\ref{eq:fieldA}) and (\ref{eq:coordinateA}). Since parallel transport preserves scalar products,
$$
  \langle X\hspace{0.02cm},\partial_i\rangle_x = \langle A\left(\mathrm{Exp}^{-1}_{x^*}(x)\right),\partial_i(x^*)\rangle_{x^*} + o(d(x\hspace{0.02cm},x^*))
$$
However, $\mathrm{Exp}^{-1}_{x^*}(x) = x^i\partial_i(x^*)$, and $\partial_i(x^*)$ form an orthonormal basis of $T_{x^*}M$. Therefore,
\begin{equation} \label{eq:proofxnormals1}
  \langle X\hspace{0.02cm},\partial_i\rangle_x = A^i_jx^j + o(d(x\hspace{0.02cm},x^*))
\end{equation}
where $A(\partial_i(x^*)) = A^k_i\partial^{\phantom{k}}_k(x^*)$. Finally, note that (in normal coordinates), the metric coefficients satisfy 
$$
g_{ij}(x) = \delta_{ij} + o(d(x\hspace{0.02cm},x^*))
$$ 
Using these to express the scalar product in (\ref{eq:proofxnormals1}), it can be seen that
\begin{equation} \label{eq:proofnormals2}
X^i(x) = A^i_jx^j + o(d(x\hspace{0.02cm},x^*))
\end{equation}
Thus, (\ref{eq:xnormals}) follows by putting $x = \psi^{-1}_\mu(u)$ in (\ref{eq:proofnormals2}). Then, $x^j = \mu^{\frac{1}{2}}u^j$ and $d(x\hspace{0.02cm},x^*) = \mu^{\frac{1}{2}}\Vert u \Vert_{x^*\hspace{0.02cm}}$. 

\subsection{Proof of Proposition \ref{prop:clt1}}
To begin, it is helpful to establish tightness of the family $(\tilde{\pi}_\mu\,;\mu \leq c\hspace{0.02cm}(2\ell(1+\sigma^{\scriptscriptstyle 2}_{\scriptscriptstyle 1}))^{\scriptscriptstyle -1})$. 
\begin{lemma} \label{lem:tildetight}
Assume that \textbf{e1}, \textbf{e2}, \textbf{v1}, \textbf{v2}, \textbf{t2} hold. Then, the family of probability distributions $(\tilde{\pi}_\mu\,;\mu \leq c\hspace{0.02cm}(2\ell(1+\sigma^{\scriptscriptstyle 2}_{\scriptscriptstyle 1}))^{\scriptscriptstyle -1})$ is tight.
\end{lemma}
Accepting this lemma, let $\tilde{\pi}_*$ be some limit point of the family $(\tilde{\pi}_\mu\,;\mu \leq c\hspace{0.02cm}(2\ell(1+\sigma^{\scriptscriptstyle 2}_{\scriptscriptstyle 1}))^{\scriptscriptstyle -1})$\hfill\linebreak at $\mu = 0$. By integrating both sides of (\ref{eq:functional}) with respect to 
$\tilde{\pi}_\mu\hspace{0.02cm}$, and recalling that $\tilde{\pi}_\mu$ is an invariant distribution of $\tilde{Q}_\mu$ (so the integral of the left-hand side is zero), it follows that
\begin{equation} \label{eq:proofclt11}
  \int_{T_{x^*}M}\mathcal{L}\phi(u)\hspace{0.03cm}\tilde{\pi}_\mu(du) = - \int_{T_{x^*}M} \varepsilon_{\mu}(u)\hspace{0.03cm}\tilde{\pi}_\mu(du)
\end{equation}
where $\varepsilon_\mu(u) = \varepsilon^{\scriptscriptstyle 1}_\mu(u) + \varepsilon^{\scriptscriptstyle 2}_\mu(u) + \mu\hspace{0.02cm}\mathcal{R}_{\psi^{-1}(u)}(\phi\circ \psi_\mu\hspace{0.02cm},\mu)$, in the notation of (\ref{eq:proofunctional}), (\ref{eq:proofunctional1}) and (\ref{eq:proofunctional2}), from the proof of Proposition \ref{prop:functional}. Since both $\varepsilon^{\scriptscriptstyle 1}_\mu(u)$ and
$\varepsilon^{\scriptscriptstyle 2}_\mu(u)$ converge to zero as $\mu \rightarrow 0$, uniformly on $T_{x^*}M$, it follows from (\ref{eq:proofclt11}) that,
$$
\left|\int_{T_{x^*}M}\mathcal{L}\phi(u)\hspace{0.03cm}\tilde{\pi}_*(du)\right| \leq \limsup_{\mu \rightarrow 0} \int_{T_{x^*}M}\,
\mu\left|\mathcal{R}_{\psi^{-1}(u)}(\phi\circ \psi_\mu\hspace{0.02cm},\mu)\right|\tilde{\pi}_\mu(du)
$$
Since $\tilde{\pi}_\mu$ is the image of $\pi_\mu$ under $\psi_\mu\hspace{0.02cm}$, this is the same as
\begin{equation} \label{eq:proofclt12}
\left|\int_{T_{x^*}M}\mathcal{L}\phi(u)\hspace{0.03cm}\tilde{\pi}_*(du)\right| \leq \limsup_{\mu \rightarrow 0} \int_{M}\,
\mu\left|\mathcal{R}_{x}(\phi\circ \psi_\mu\hspace{0.02cm},\mu)\right|{\pi}_\mu(dx)
\end{equation}
To bound the right-hand side, put $f = \phi\circ \psi_\mu$ in (\ref{eq:qremainder}). If $\bar{f} = \phi \circ \mathrm{Exp}_{x^*}$, then $\bar{f}$ is compactly-supported and smooth. Moreover, applying the chain rule, it follows from (\ref{eq:AfBf}) that $A_f = \mu^{-1}\!A_{\bar{f}}$ and $B_f = 
\mu^{-\frac{3}{2}}B_{\bar{f}\hspace{0.03cm}}$. Therefore, by (\ref{eq:qremainder}),
$$
\left|\mathcal{R}_{x}(\phi\circ \psi_\mu\hspace{0.02cm},\mu)\right|
\leq 2\mu^{-1}\!A_{\bar{f}}\hspace{0.04cm}\mathbb{E}\left[\mathbf{1}\lbrace \Vert e_y\Vert_x> K\rbrace\Vert e_y\Vert^2_x\right] + \hspace{8cm}
$$
\begin{equation} \label{eq:qremainderters}
\phantom{abcd}6\mu^{-\frac{1}{2}}B_{\bar{f}}\hspace{0.03cm}
(2\Vert X \Vert^2_x + 2K^2)\hspace{0.03cm}\mathbb{E}\left[\mathbf{1}\lbrace \Vert e_y\Vert_x> K\rbrace\Vert e_y\Vert_x\right]+2\mu^{-\frac{1}{2}}B_{\bar{f}}\hspace{0.03cm}(4\Vert X \Vert^3_x + 4K^3) 
\end{equation}
Now, since \textbf{t3} holds, it follows from (\ref{eq:extravariancecontrol}) that
\begin{equation} \label{eq:proofclt13}
\mathbb{E}\left[\mathbf{1}\lbrace \Vert e_y\Vert_x> K\rbrace\Vert e_y\Vert^2_x\right] \leq K^{-\alpha}\hspace{0.02cm}(\tilde{\sigma}^2_{\scriptscriptstyle 0} + \tilde{\sigma}^2_{\scriptscriptstyle 1}\hspace{0.02cm}V(x))
\end{equation}
Moreover, by (\ref{eq:variancecontrol}) (assuming that $K > 1$),
\begin{equation} \label{eq:proofclt14}
(2\Vert X \Vert^2_x + 2K^2)\hspace{0.03cm}\mathbb{E}\left[\mathbf{1}\lbrace \Vert e_y\Vert_x> K\rbrace\Vert e_y\Vert_x\right] \leq
(2\Vert X \Vert^2_x + 2K^2)(\sigma^2_{\scriptscriptstyle 0} + \sigma^2_{\scriptscriptstyle 1}\hspace{0.02cm}\Vert X\Vert^2_x)
\end{equation}
Then, it follows from (\ref{eq:qremainderters}), (\ref{eq:proofclt13}) and (\ref{eq:proofclt14}) that
$$
\begin{array}{l}
%\limsup_{\mu \rightarrow 0} 
\mu\left|\mathcal{R}_{x}(\phi\circ \psi_\mu\hspace{0.02cm},\mu)\right| \leq \\[0.2cm]
2A_{\bar{f}}\hspace{0.04cm}K^{-\alpha}\hspace{0.02cm}(\tilde{\sigma}^2_{\scriptscriptstyle 0} + \tilde{\sigma}^2_{\scriptscriptstyle 1}\hspace{0.02cm}V(x)) \,+ 6\mu^{\frac{1}{2}}B_{\bar{f}}(2\Vert X \Vert^2_{x} + 2K^2)(\sigma^2_{\scriptscriptstyle 0} + \sigma^2_{\scriptscriptstyle 1}\hspace{0.02cm}\Vert X\Vert^2_{x}) +
2\mu^{\frac{1}{2}}B_{\bar{f}}\hspace{0.03cm}(4\Vert X \Vert^3_{x} + 4K^3) 
\end{array}
$$
Integrate this inequality with respect to $\pi_\mu\hspace{0.02cm}$, and recall from Proposition \ref{prop:todirac} that $\pi_\mu$ converges weakly to $\delta_{x^*}$ as $\mu \rightarrow 0$. It follows that,
$$
\limsup_{\mu \rightarrow 0} 
\int_{M}\,
\mu\left|\mathcal{R}_{x}(\phi\circ \psi_\mu\hspace{0.02cm},\mu)\right|{\pi}_\mu(dx) \leq 
2A_{\bar{f}}\hspace{0.04cm}K^{-\alpha}\hspace{0.02cm}(\tilde{\sigma}^2_{\scriptscriptstyle 0}
+ \tilde{\sigma}^2_{\scriptscriptstyle 1}\hspace{0.02cm}V(x^*))
$$
However, since $K$ can be chosen arbitrarily large, and $\alpha > 0$, the limit superior is equal to zero, and (\ref{eq:proofclt12}) becomes
$$
\int_{T_{x^*}M}\mathcal{L}\phi(u)\hspace{0.03cm}\tilde{\pi}_*(du) = 0
$$
This means that $\tilde{\pi}_*$ is an invariant distribution of the generator $\mathcal{L}$, and therefore $\tilde{\pi}_* = \mathrm{N}(0,V)$, as required. \vfill\pagebreak
%\[0.1cm]
\noindent \textbf{Proof of Lemma \ref{lem:tildetight}\,:} by Proposition \ref{prop:harris}, \textbf{e1}, \textbf{e2}, \textbf{v1}, \textbf{v2} ensure that the chain $(x_t)$ has a unique invariant distribution $\pi_\mu\hspace{0.02cm}$, whenever $\mu \leq c\hspace{0.02cm}(2\ell(1+\sigma^{\scriptscriptstyle 2}_{\scriptscriptstyle 1}))^{\scriptscriptstyle -1}$. Then, the chain $(u_t)$ has a unique invariant distribution $\tilde{\pi}_\mu\hspace{0.02cm}$. According to (\ref{eq:re-scaledu}), this is $\tilde{\pi}_\mu(A) = \pi_\mu(\mathrm{Exp}_{x^*}(\mu^{1/2}A))$,\hfill\linebreak for any measurable $A \subset T_{x^*}M$. The same \textbf{e1}, \textbf{e2}, \textbf{v1}, \textbf{v2} also imply (\ref{eq:Vmoment}), in the proof of Proposition (\ref{prop:todirac}). Now, for $u \in T_{x^*} M$, let $x = \mathrm{Exp}_{x^*}(\mu^{1/2}u)$, and note that $\Vert u \Vert_{x^*} > r$ if and only if $d(x\hspace{0.02cm},x^*) > \mu^{\scriptscriptstyle 1/2}r$. It then follows from Assumption \textbf{t2} that
$$
\tilde{\pi}_\mu(\Vert u \Vert_{x^*} > r) \leq \pi_\mu(V > v(\mu^{\scriptscriptstyle 1/2}r))
$$
so, using Markov's inequality and (\ref{eq:Vmoment}),
$$
\tilde{\pi}_\mu(\Vert u \Vert_{x^*} > r) \leq 2(\ell\hspace{0.02cm}\sigma^2_{\scriptscriptstyle 0}/\lambda)\left(\mu\middle/v(\mu^{\scriptscriptstyle 1/2}r)\right)
$$
To conclude, let $\bar{\mu} = c\hspace{0.02cm}(2\ell(1+\sigma^{\scriptscriptstyle 2}_{\scriptscriptstyle 1}))^{\scriptscriptstyle -1}$. By \textbf{t2}, 
$\left.\mu\middle/v(\mu^{\scriptscriptstyle 1/2}r)\right. \leq \left.\bar{\mu}\middle/v(\bar{\mu}^{\scriptscriptstyle 1/2}r)\right.$. Therefore,
$$
\tilde{\pi}_\mu(\Vert u \Vert_{x^*} > r) \leq 2(\ell\hspace{0.02cm}\sigma^2_{\scriptscriptstyle 0}/\lambda)\left(\bar{\mu}\middle/v(\bar{\mu}^{\scriptscriptstyle 1/2}r)\right)
$$
However (again by \textbf{t2}), the right-hand side of this inequality is independent of $\mu$, and goes to zero as $r \rightarrow \infty$. This is equivalent to the required tightness.

\section{Riemannian AR(1)} \label{sec:ar1}
Let $M$ be a Hadamard manifold, and $P$ a probability distribution on $M$, which has a strictly positive probability density, with respect to Riemannian volume.   Then, let $(y_t\,;t = 1,2,\ldots)$ be independent samples from $P$, and $x_{\scriptscriptstyle 0}$ be a point in $M$. Consider the stochastic update rule, starting from $x_{\scriptscriptstyle 0}$ at $t = 0$,
\begin{equation} \label{eq:ar1}
  x_{t+1} = x_t\, \#_{\scriptscriptstyle \mu}\, y_{\scriptscriptstyle\hspace{0.02cm}t+1} \hspace{0.5cm}\text{where $\mu \in (0,1)$}
\end{equation}
where the notation is that of (\ref{eq:mapformula}) from \ref{sec:mapvsmms}. This will be called a Riemannian AR(1) model, since each new $x_{t+1}$ is a geodesic convex combination of the old $x_t$ and of the new sample $y_{\scriptscriptstyle\hspace{0.02cm}t+1\hspace{0.02cm}}$. If $M$ is a Euclidean space, $M = \mathbb{R}^n$, then (\ref{eq:ar1}) reads $x_{t+1} = (1-\mu)\hspace{0.02cm}x_t + \mu\hspace{0.02cm}y_{\scriptscriptstyle\hspace{0.02cm}t+1\hspace{0.02cm}}$, which is a first-order auto-regressive model (whence the name AR(1)).

The update rule (\ref{eq:ar1}) may be viewed as a constant-step-size exponential scheme, of the form (\ref{eq:cssexp}). Specifically, (\ref{eq:ar1}) is equivalent to
\begin{equation} \label{eq:arscheme}
     x_{t+1} = \mathrm{Exp}_{x_t}\!\left(\mu\hspace{0.02cm}X_{y_{\scriptscriptstyle\hspace{0.02cm}t+1}}(x_t)\right) \hspace{0.5cm} \text{where } X_y(x) = \mathrm{Exp}^{-1}_x(y)
\end{equation}
which defines a time-homogeneous Markov chain $(x_t)$ with values in $M$. 

One is tempted to apply the results of \ref{sec:clt1} (\textit{e.g.} on geometric ergodicity), directly to the scheme (\ref{eq:arscheme}). However, some of the assumptions in \ref{sec:clt1} (especially \textbf{e1}), turn out to be quite unnatural. Fortunately, it is possible to proceed along a different path, which only requires the existence of second-order moments. Specifically, it is merely required that
\begin{equation} \label{eq:arsecondorder}
  \mathcal{E}(x) = \frac{1}{2}\hspace{0.03cm}\int_M\,d^{\hspace{0.03cm} 2}(x\hspace{0.02cm},y)\hspace{0.03cm}P(dy) < \,\infty
\end{equation}
for some (and therefore all) $x \in M$. As discussed in \ref{subsec:frechhad}, (\ref{eq:arsecondorder}) guarantees existence and uniqueness of the Riemannian barycentre $x^*$ of $P$. This is enough for the following proposition.
\begin{proposition} \label{prop:arergodic}
 Consider the Riemannian AR(1) model (\ref{eq:ar1}) (or (\ref{eq:arscheme})), on a Hadamard manifold $M$. If (\ref{eq:arsecondorder}) is verified, then 
the Markov chain $(x_t)$ is geometrically ergodic, with a unique invariant distribution $\pi_\mu$. Moreover, $\pi_{\mu}\,\Rightarrow\,\delta_{x^*}$ as $\mu \rightarrow 0$. 
%\hfill\linebreak Assume that \textbf{e1}, \textbf{e2}, \textbf{v1}, \textbf{v2} hold. If $\mu \leq c\hspace{0.02cm}(2\ell(1+\sigma^{\scriptscriptstyle 2}_{\scriptscriptstyle 1}))^{\scriptscriptstyle -1}$, then the Markov chain $(x_t)$ is geometrically ergodic, with a unique invariant distribution $\pi_\mu$.
\end{proposition}
The proof of this proposition begins like that of Proposition \ref{prop:harris}, by noting that the Markov chain $(x_t)$ is Feller and $|\mathrm{vol}|$-irreducible and aperiodic. Indeed, since $X_y(x)$ is given by (\ref{eq:arscheme}), and since $P$ has a strictly positive probability density with respect to $|\mathrm{vol}|$, it follows that Assumption \textbf{e2} holds. Therefore, it is possible to argue exactly as in the proof of Lemma \ref{lem:feller}. \vfill\pagebreak

Now, let $V(x) = d^{\hspace{0.02cm} 2}(x^*,x)/2$. To prove that the chain $(x_t)$ is geometrically ergodic, it is enough to obtain the inequality 
\begin{equation} \label{eq:ardrift}
 Q_\mu V(x) \leq (1-\mu)^2\hspace{0.02cm}V(x) + \mu^2\hspace{0.02cm}\mathcal{E}(x^*)
\end{equation}
which is similar to (\ref{eq:driftzero}) of Lemma \ref{lem:drift}. This can then be used, exactly as in the proof of Proposition \ref{prop:harris}, based on~\cite{tweedie} (Theorem 16.0.1).\\[0.1cm]
\textbf{Proof of (\ref{eq:ardrift})\,:} for any $x \in M$, note from (\ref{eq:ar1}) that
\begin{equation} \label{eq:finalproof1}
Q_\mu V(x) = \mathbb{E}\left[V(x\, \#_{\scriptscriptstyle \mu}\, y)\right]
\end{equation}
where $y$ denotes a random variable with distribution $P$. Recall from \ref{ssec:sqd} that $V(x)$ is $1/2$-strongly convex. Therefore, by (\ref{eq:strongconv}),
$$
V(x\, \#_{\scriptscriptstyle \mu}\, y) \leq (1-\mu)V(x) + \mu V(y) - \mu(1-\mu)d^{\hspace{0.03cm} 2}(x\hspace{0.02cm},y)/2
$$
taking expectations, this yields,
\begin{equation} \label{eq:finalproof2}
Q_\mu V(x) \leq 
(1-\mu)V(x) + \mu\hspace{0.02cm}\mathcal{E}(x^*) - \mu(1-\mu)\mathcal{E}(x)
\end{equation}
after using the fact that $\mathcal{E}(x) = \mathbb{E}[d^{\hspace{0.03cm} 2}(x\hspace{0.02cm},y)/2]$ for any $x \in M$, which is clear from (\ref{eq:arsecondorder}). Now, recall from \ref{subsec:frechhad} that $\mathcal{E}$ is $1/2$-strongly convex. Therefore, by (\ref{eq:properconv})
\begin{equation} \label{eq:finalproof3}
\mathcal{E}(x) \geq \mathcal{E}(x^*) + d^{\hspace{0.02cm} 2}(x^*,x)/2
\end{equation}
where the right-hand side is just $\mathcal{E}(x^*) + V(x)$. Thus, replacing (\ref{eq:finalproof3}) into (\ref{eq:finalproof2}), one has
$$
Q_\mu V(x) \leq 
(1-\mu)V(x) + \mu\hspace{0.02cm}\mathcal{E}(x^*) - \mu(1-\mu) V(x)-  \mu(1-\mu)\mathcal{E}(x^*) 
$$
which immediately yields (\ref{eq:ardrift}). \\[0.1cm]
Geometric ergodicity ensures the chain $(x_t)$ has a unique invariant distribution $\pi_\mu$. To prove that $\pi_{\mu}\,\Rightarrow\,\delta_{x^*}$ as $\mu \rightarrow 0$, it is possible to argue as in the proof of Proposition \ref{prop:todirac}. Precisely, integrating both sides of (\ref{eq:ardrift}) with respect to $\pi_\mu\hspace{0.02cm}$, it follows that
$$
 \int_M Q_\mu V(x)\hspace{0.02cm}\pi_\mu(dx) \leq (1-\mu)^2\hspace{0.02cm}\int_M V(x)\hspace{0.02cm}\pi_\mu(dx) + \mu^2\hspace{0.02cm}\mathcal{E}(x^*)
%\int_M Q_\mu V(x)\hspace{0.02cm}\pi_\mu(dx) \leq (1-\lambda\mu/2)\int_M V(x)\hspace{0.02cm}\pi_\mu(dx) + (\ell\hspace{0.02cm}\sigma^2_{\scriptscriptstyle 0})\hspace{0.03cm}\mu^2
$$
Since $\pi_\mu$ is an invariant distribution of the transition kernel $Q_\mu\hspace{0.02cm}$, this means
$$
 \int_M V(x)\hspace{0.02cm}\pi_\mu(dx)  \leq (1-\mu)^2\hspace{0.02cm}\int_M V(x)\hspace{0.02cm}\pi_\mu(dx) + \mu^2\hspace{0.02cm}\mathcal{E}(x^*)
%\int_M V(x)\hspace{0.02cm}\pi_\mu(dx) \leq (1-\lambda\mu/2)\int_M V(x)\hspace{0.02cm}\pi_\mu(dx) + (\ell\hspace{0.02cm}\sigma^2_{\scriptscriptstyle 0})\hspace{0.03cm}\mu^2
$$
In other words,
\begin{equation} \label{eq:arVmoment}
 \int_M V(x)\hspace{0.02cm}\pi_\mu(dx) \leq \mathcal{E}(x^*)\mu/(2-\mu)
\end{equation}
Since $\mu/(2-\mu) \leq \mu$ for $\mu \leq 1$, (\ref{eq:arVmoment}) can be used like (\ref{eq:Vmoment}) in the proof of Proposition \ref{prop:todirac}. In this way, $\pi_{\mu}\,\Rightarrow\,\delta_{x^*}$ follows by noting that $V$, defined by $V(x)= d^{\hspace{0.02cm} 2}(x^*,x)/2$, has compact sublevel sets, by the Hopf-Rinow theorem and the fact that $M$ is complete (see~\cite{chavel}), and that $V(x) = 0$ if and only if $x = x^*$. \\[0.1cm]
\textbf{Remark\,:} thanks to Proposition \ref{prop:arergodic}, it is now possible to prove that a central limit theorem, identical to Proposition \ref{prop:clt1}, holds for the Riemannian AR(1) model (\ref{eq:ar1}). This only requires the additional condition (\ref{eq:extravariancecontrol}).

%\section{The CLT\,: Part II} \label{sec:clt2}

\chapter{Open problems} \label{chapopen}
While working on this thesis, there are several problems which I found very interesting, and important, but could not solve, or even attack in a meaningful way. I would therefore like to close the thesis with a list of these problems, in the hope that they will attract the attention of more people (not just myself). \\[0.1cm]
\textbf{In Chapter \ref{barycentre}\,:} the conclusion of Lemma \ref{lem:hddim} only holds for compact Riemannian symmetric spaces which are simply connected. Therefore, the subsequent Propositions \ref{prop:differentiability} and \ref{prop:unique} are restricted to this simply connected case. The problem is to describe, at least partially, what happens for compact Riemannian symmetric spaces which are not simply connected. It would be particularly interesting to give counterexamples to either one of Propositions \ref{prop:differentiability} and \ref{prop:unique}, in the non-simply-connected case. \\[0.1cm]
\textbf{In Chapter \ref{gaussian}\,:} formula (\ref{eq:ssz}) gives the normalising factor $Z(\sigma)$ for a Gaussian distribution on a symmetric space $M$, which belongs to the non-compact case. When $M$ is the space of $N \times N$ Hermitian positive-definite matrices, (\ref{eq:sw_z}) and (\ref{eq:asymp_z}) provide a closed form expression (valid for any $N$), and an asymptotic expression (valid for large $N$), of the multiple integral in (\ref{eq:ssz}).
The problem is to find similar expressions of this integral, for other symmetric spaces. It should be easiest to first deal with the spaces of $N \times N$ real or quaternion positive-definite matrices, and then move on to other spaces, such as the Siegel domaine (Example 3, in \ref{sec:z}). \hfill\linebreak
\textbf{In Chapter \ref{bayesian}\,:} the problem is to prove or disprove the conjecture, mentioned in \ref{sec:mapmmsdistance}. Namely, that the MAP and MMS Bayesian estimators are equal, for Gaussian distributions on a space of constant negative curvature. \\[0.1cm]
\textbf{In Chapter \ref{bayesian}\,:} as mentioned in \ref{sec:markovclt}, I have never met with a function $f:M \rightarrow \mathbb{R}$ ($M$ a non-Euclidean Hadamard manifold), which is strongly convex, and also has a bounded Hessian. The problem is to construct a function $f$ with these properties, or to show this is not possible. Another problem, which is quite important for convex optimisation, is to show that a function $f:M \rightarrow \mathbb{R}$, which is convex and has a bounded Hessian, has co-coercive gradient, in the sense of~\cite{nesterov} (Theorem 2.1.5, property (2.1.11)).
\\[0.1cm]
\textbf{In Chapter \ref{stocha}\,:} to state in clear and general terms the functional central limit theorem, which follows from Proposition \ref{prop:functional}, and also to derive a similar functional central limit theorem for decreasing-step-size schemes. These can be used in studying the behavior of stochastic approximation schemes, in the presence of unstable critical points (only stable critical points were considered, in the above). \\[0.1cm]
\textbf{In Chapter \ref{stocha}\,:} to generalise Proposition \ref{prop:clt1} to the case where $M$ is not a Hadamard manifold. I believe that, in this case, the asymptotic form of the invariant distribution will no longer be multivariate normal (roughly, the scheme can always ``jump across" the cut locus of a stable critical point).%\\[0.1cm]
%\textbf{In Chapter \ref{stocha}\,:}

%In , an analytic expression for $Z(\sigma)$ was obtained, and \ref{sec:rmt},

%nor simply connected

%normalisation

%the conjecture

%co-coercive

%functional clt
%barycentre of markov

\bibliographystyle{IEEEtran}    
\bibliography{adtsalem}

% Generated by IEEEtran.bst, version: 1.14 (2015/08/26)
\begin{thebibliography}{10}
\providecommand{\url}[1]{#1}
\csname url@samestyle\endcsname
\providecommand{\newblock}{\relax}
\providecommand{\bibinfo}[2]{#2}
\providecommand{\BIBentrySTDinterwordspacing}{\spaceskip=0pt\relax}
\providecommand{\BIBentryALTinterwordstretchfactor}{4}
\providecommand{\BIBentryALTinterwordspacing}{\spaceskip=\fontdimen2\font plus
\BIBentryALTinterwordstretchfactor\fontdimen3\font minus
  \fontdimen4\font\relax}
\providecommand{\BIBforeignlanguage}[2]{{%
\expandafter\ifx\csname l@#1\endcsname\relax
\typeout{** WARNING: IEEEtran.bst: No hyphenation pattern has been}%
\typeout{** loaded for the language `#1'. Using the pattern for}%
\typeout{** the default language instead.}%
\else
\language=\csname l@#1\endcsname
\fi
#2}}
\providecommand{\BIBdecl}{\relax}
\BIBdecl

\bibitem{warped}
S.~Said, L.~Bombrun, and Y.~Berthoumieu, ``Warped {Riemannian} metrics for
  location-scale models,'' in \emph{Geometric structures of information},
  F.~Nielsen, Ed.\hskip 1em plus 0.5em minus 0.4em\relax Springer Switzerland,
  2019.

\bibitem{colt}
A.~Durmus, P.~Jimenez, E.~Moulines, S.~Said, and H.~T. Wai, ``Convergence
  analysis of {Riemannian} stochastic approximation schemes,''
  \emph{\texttt{arXiv:2005.13284}}, 2020.

\bibitem{aistats}
A.~Durmus, P.~Jimenez, E.~Moulines, and S.~Said, ``On {Riemannian} stochastic
  approximation schemes with fixed step-size (under review),'' in
  \emph{Artificial Intelligence and Statistics}, 2021.

\bibitem{tierz}
L.~Santilli and M.~Tierz, ``Riemannian gaussian distributions, random matrix
  ensembles and diffusion kernels,'' \emph{\texttt{arXiv:2011.13680}}, 2020.

\bibitem{sanslarmes}
P.~A. Meyer, ``Géométrie stochastique sans larmes,'' \emph{Séminaire de
  probabilités (Strasbourg)}, vol.~15, pp. 44--102, 1981.

\bibitem{eschenburg}
J.~H. Eschenburg, \emph{Lecture notes on symmetric spaces (course material
  available online)}, \texttt{uni-augsburg.de/eschenbu/symspace.pdf}, 1997.

\bibitem{absil}
P.~A. Absil, R.~Mahony, and R.~Sepulchre, \emph{Optimization algorithms on
  matrix manifolds}.\hskip 1em plus 0.5em minus 0.4em\relax Princeton
  University Press, 2008.

\bibitem{absilmalick}
P.~A. Absil and J.~Malick, ``Projection-like retractions on matrix manifolds,''
  \emph{SIAM Journal on Optimization}, vol.~22, pp. 135--158, 2012.

\bibitem{sakai}
T.~Sakai, ``On cut loci of compact symmetric spaces,'' \emph{Hokkaido
  Mathematical Journal.}, vol.~6, pp. 136--161, 1977.

\bibitem{helgason}
S.~Helgason, \emph{Differential geometry and symmetric spaces}.\hskip 1em plus
  0.5em minus 0.4em\relax New York and London: Academic Press, 1962.

\bibitem{chavel}
I.~Chavel, \emph{Riemannian geometry, a modern introduction}.\hskip 1em plus
  0.5em minus 0.4em\relax Cambridge University Press, 2006.

\bibitem{distributional}
C.~Mantegazza, G.~Mascellani, and G.~Uraltsev, ``On the distributional hessian
  of the distance function,'' \emph{\texttt{arXiv:1303.1421}}, 2013.

\bibitem{huber}
P.~J. Huber and E.~M. Ronchetti, \emph{Robust statistics (2nd edition)}.\hskip
  1em plus 0.5em minus 0.4em\relax Wiley-Blackwell, 2009.

\bibitem{marsden}
J.~E. Marsden, T.~Ratiu, and R.~Abraham, \emph{Manifolds, tensor analysis, and
  applications}.\hskip 1em plus 0.5em minus 0.4em\relax Springer-Verlag, 2001.

\bibitem{besseclosed}
A.~L. Besse, \emph{Manifolds all of whose geodesics are closed}.\hskip 1em plus
  0.5em minus 0.4em\relax New York: Springer-Verlag, 1978.

\bibitem{bogachev}
V.~I. Bogachev, \emph{Measure Theory, Volume I}.\hskip 1em plus 0.5em minus
  0.4em\relax Springer-Verlag, 2007.

\bibitem{mehta}
M.~L. Mehta, \emph{Random matrices (3rd edition)}.\hskip 1em plus 0.5em minus
  0.4em\relax Elsevier Ltd., 2004.

\bibitem{meckes}
E.~S. Meckes, \emph{The random matrix theory of the classical compact
  groups}.\hskip 1em plus 0.5em minus 0.4em\relax Cambridge University Press,
  2019.

\bibitem{nomizu2}
S.~Kobayashi and K.~Nomizu, \emph{Foundations of differential geometry, Volume
  II}.\hskip 1em plus 0.5em minus 0.4em\relax Interscience Publishers, 1969.

\bibitem{higham}
H.~J. Higham, \emph{Functions of matrices\,: theory and computation}.\hskip 1em
  plus 0.5em minus 0.4em\relax SIAM Publications, 2008.

\bibitem{pennec2006}
X.~Pennec, P.~Fillard, and N.~Ayache, ``A {Riemannian} framework for tensor
  computing,'' \emph{International Journal of Computer Vision}, vol.~66, no.~1,
  pp. 41--66, 2006.

\bibitem{salemgibbs}
S.~Said and J.~H. Manton, ``Riemannian barycentres of {Gibbs} distributions\,:
  new results on concentration and convexity,'' \emph{Information Geometry
  (under review)}, 2020.

\bibitem{frechet}
M.~Fréchet, ``Les éléments aléatoires de nature quelconque dans un espace
  distancié,'' \emph{Annales de l'I.H.P.}, vol.~10, no.~4, pp. 215--210, 1948.

\bibitem{bhatta1}
R.~Bhattacharya and V.~Patrangenaru, ``Large sample theory of instrinsic and
  extrinsic sample means on manifolds {I},'' \emph{The annals of statistics},
  vol.~31, no.~1, pp. 1--29, 2003.

\bibitem{bhatta2}
------, ``Large sample theory of instrinsic and extrinsic sample means on
  manifolds {II},'' \emph{The annals of statistics}, vol.~33, no.~3, pp.
  1225--1259, 2005.

\bibitem{kendall}
W.~S. Kendall, ``Probability, convexity, and harmonic maps with small image
  {I}\,: uniqueness and fine existence,'' \emph{Proceedings of the London
  Mathematical Society}, vol.~61, no.~2, pp. 371--406, 1990.

\bibitem{afsari}
B.~Afsari, ``Riemannian {$L^{\scriptscriptstyle p}$} center of mass\,:
  existence, uniqueness and convexity,'' \emph{Proceedings of the American
  Mathematical Society}, vol. 139, no.~2, pp. 655--673, 2010.

\bibitem{propeller}
W.~S. Kendall, ``The propeller\,: a counterexample to a conjectured criterion
  for the existence of certain convex functions,'' \emph{Journal of the London
  Mathematical Society}, vol.~36, no.~2, pp. 364--374, 1992.

\bibitem{sturm}
K.~T. Sturm, ``Probability measures on metric spaces of nonpositive
  curvature,'' \emph{Contemporary mathematics}, vol. 338, pp. 1--34, 2003.

\bibitem{arnaumic}
M.~Arnaudon and L.~Miclo, ``Means in complete manifolds\,: completeness and
  approximation,'' \emph{ESAIM\,:Probability and Statistics}, vol.~18, pp.
  185--206, 2014.

\bibitem{wong}
R.~Wong, \emph{Asymptotic approximation of integrals}, Society of Industrial
  and Applied Mathematics, 2001.

\bibitem{hausdorff}
D.~Schleicher, \emph{Hausdorff dimension, its properties and its surprises
  (available online)}, \texttt{arXiv:0505099}, 2007.

\bibitem{watson}
E.~T. Whittaker and G.~N. Watson, \emph{A course of modern analysis (4th
  edition)}, Cambridge University Press, 1950.

\bibitem{petersen}
P.~Petersen, \emph{Riemannian geometry (2nd edition)}, Springer Science, 2006.

\bibitem{kantorovich}
L.~P. Kantorovich and G.~P. Akilov, \emph{Functional analysis (2nd edition)},
  Pergamon Press, 1982.

\bibitem{said2}
S.~Said, H.~Hajri, L.~Bombrun, and B.~C. Vemuri, ``Gaussian distributions on
  {Riemannian} symmetric spaces\,: statistical learning with structured
  covariance matrices,'' \emph{IEEE Transactions on Information Theory},
  vol.~64, no.~2, pp. 752--772, 2018.

\bibitem{historynormal}
S.~Stahl, ``The evolution of the normal distribution,'' \emph{Mathematics
  Magazine}, vol.~79, no.~2, pp. 96--113, 2006.

\bibitem{borel}
E.~Borel, \emph{Introduction géométrique à quelques théories physiques},
  Gauthier-Villars, 1914.

\bibitem{perrin}
J.~Perrin, ``Mouvement brownien et molécules,'' \emph{Journal de physique
  théorique et appliquée}, vol.~9, no.~1, pp. 5--39, 1910.

\bibitem{knapp}
A.~W. Knapp, \emph{Lie groups, beyond an introduction (2nd edition)}.\hskip 1em
  plus 0.5em minus 0.4em\relax Birkhauser, 2002.

\bibitem{siegelspgeo}
C.~L. Siegel, ``Symplectic geometry,'' \emph{American Journal of Mathematics},
  vol.~65, no.~1, pp. 1--86, 1943.

\bibitem{terras2}
A.~Terras, \emph{Harmonic analysis on symmetric spaces and applications, Vol.
  II}.\hskip 1em plus 0.5em minus 0.4em\relax Springer-Verlag, 1988.

\bibitem{szegobook}
G.~Szegö, \emph{Orthogonal Polynomials (1st edition)}.\hskip 1em plus 0.5em
  minus 0.4em\relax American Mathematical Society, 1939.

\bibitem{jackson}
B.~C. Berndt, \emph{What is a $q$-series? (appeared in \emph{Ramanujan
  rediscovered}, available online)},
  \texttt{faculty.math.illinois.edu/berndt/articles/q.pdf}, 2012.

\bibitem{asch}
A.~B.~J. Kuijlaars and W.~{Van Assche}, ``The asymptotic zero distribution of
  orthogonal polynomials with varying recurrence coefficients,'' \emph{Journal
  of Approximation theory}, vol.~99, pp. 167--197, 1999.

\bibitem{deift}
P.~Deift, \emph{Orthogonal polynomials and random matrices\,: a
  {Riemann}-{Hilber} approach}, American Mathematical Sociery, 1998.

\bibitem{marino}
M.~Mariño, \emph{Chern-Simons theory, matrix models, and topological strings},
  Oxford University Press, 2005.

\bibitem{rogers}
G.~E. Andrews, \emph{The theory of partitions}, Addison-Wesley Publishing
  Company, 1976.

\bibitem{jarner}
S.~F. Jarner and E.~Hansen, ``Geometric ergodicity of {Metropolis}
  algorithms,'' \emph{Stochastic Processes and Applications}, vol.~58, pp.
  341--361, 1998.

\bibitem{roberts}
R.~O. Roberts and J.~S. Rosenthal, ``General state-space {Markov} chains and
  {MCMC} algorithms,'' \emph{Probability Surveys}, vol.~1, pp. 20--71, 2004.

\bibitem{tweedie}
S.~Meyn and R.~L. Tweedie, \emph{Markov chains and stochastic stability},
  Cambrdge University Press, 2008.

\bibitem{rroberts}
G.~O. Roberts and R.~L. Tweedie, ``Geometric ergodicity and central limit
  theorems for multidimensional {Metropolis} and {Hastings} algorithms,''
  \emph{Biometrika}, vol.~82, no.~1, pp. 95--110, 1996.

\bibitem{lee}
J.~M. Lee, \emph{Introduction to smooth manifolds (2nd edition)}, Springer
  Science, 2012.

\bibitem{pflug}
G.~C. Pflug, ``Stochastic minimisation with constant step-size\,: asymptotic
  laws,'' \emph{SIAM Journal on Control and Optimization}, vol.~24, no.~4, pp.
  655--666, 1986.

\bibitem{kallenberg}
O.~Kallenberg, \emph{Foundations of modern probability (2nd edition)},
  Springer-Verlag, 2002.

\bibitem{duflo}
M.~Duflo, \emph{Algorithmes stochastiques}, Springer-Verlag, 1996.

\bibitem{nesterov}
Y.~Nesterov, \emph{Lectures on convex optimization}.\hskip 1em plus 0.5em minus
  0.4em\relax Springer Switzerland, 2018.

\end{thebibliography}

\end{document}